\documentclass[11pt,reqno,a4paper]{amsart}
\usepackage[T1]{fontenc}
\usepackage{lmodern}
\usepackage[margin=28mm]{geometry}
\usepackage{amsmath,amssymb,mathtools,microtype}
\usepackage[hidelinks]{hyperref}
\usepackage{tikz}
\usetikzlibrary{arrows.meta}
\usepackage{graphicx,placeins,needspace}
\usepackage{mathrsfs}
\newtheorem{theorem}{Theorem}[section]
\newtheorem{proposition}[theorem]{Proposition}
\newtheorem{lemma}[theorem]{Lemma}
\newtheorem{corollary}[theorem]{Corollary}
\theoremstyle{definition}
\newtheorem{example}{Example}
\theoremstyle{remark}
\newtheorem*{remark}{Remark}
\numberwithin{equation}{section}
\newcommand{\C}{\mathbb C}
\newcommand{\Q}{\mathbb Q}
\newcommand{\R}{\mathbb R}
\newcommand{\Z}{\mathbb Z}
\newcommand{\dd}{\,\mathrm d}
\DeclareMathOperator{\conv}{conv}
\DeclareMathOperator{\Rea}{Re}
\DeclareMathOperator{\Ima}{Im}
\title[Global zero distribution and pluripolar cones]{Global zero distribution and pluripolar exceptional cones for complex linear differential equations}
\author{Xing-Yu Li}
\date{}
\makeatletter
\let\manuscriptSubjectClass\@subjclass
\let\manuscriptKeywords\@keywords
\let\@subjclass\@empty
\let\@keywords\@empty
\let\manuscriptSetAuthors\@setauthors
\renewcommand{\@setauthors}{%
  \manuscriptSetAuthors
  \begingroup
  \centering\footnotesize\normalfont
  Department of Mathematics, The University of Manchester\\
  Oxford Road, Manchester M13 9PL, United Kingdom\\
  \href{mailto:xingyu.li@manchester.ac.uk}{xingyu.li@manchester.ac.uk}\quad
  \href{https://xyli9596.github.io/}{https://xyli9596.github.io/}\par
  \endgroup
}
\makeatother
\hypersetup{pdftitle={Global zero distribution and pluripolar exceptional cones for complex linear differential equations},pdfauthor={Xing-Yu Li},pdfsubject={Global zero distribution and pluripolar exceptional initial values for complex linear differential equations with exponential-sum coefficients}}
\begin{document}
\raggedbottom

\begin{abstract}
We study the zeros of solutions of complex linear differential equations
of arbitrary order $n\ge2$ with exponential-sum coefficients.
Using the initial-value space \(\mathbb C^n\) to parametrize the solutions, we show that solutions
with initial vectors outside a complete pluripolar cone of Lebesgue
measure zero share an explicit common leading term for the global zero
count and a common joint limit in radius and angle. The rescaled radial
distance from the disk boundary has an exponential limiting distribution,
while the rescaled angular deviations from each main ray are Gaussian.
Inspired by logarithmic tracts in the Eremenko--Lyubich class $\mathcal B$,
we pull exterior sectors back conformally to countably many curvilinear
strips extending in finitely many directions and obtain zero-counting
asymptotics in each strip, although the coefficients may lie outside this class.
Moreover, the asymptotic formula for the global zero count and the joint
limit in radius and angle persist under arbitrary entire perturbations
of strictly smaller weighted exponential type.
\end{abstract}
\maketitle

\begingroup
\par
\smallskip
\noindent\textit{2020 Mathematics Subject Classification.}
\manuscriptSubjectClass.\par
\endgroup

\section{Introduction}\label{sec:introduction}

Linear differential equations with exponential-sum coefficients
include Mathieu's equation \cite[\S\S28.2 and 28.33]{DLMF} and
one-dimensional Schr\"odinger equations with potentials of the form
$V(z)=c e^{\lambda z}$ \cite{BerensteinHulsey}. The zeros of their solutions describe
oscillation in the complex plane. For a fixed equation, we study
how the global zero distribution of its solutions varies with the
initial values. We consider linear differential equations on $\C$
\begin{equation}\label{eq:reduced}
 f^{(n)}+\sum_{j=0}^{n-2}A_j(z)f^{(j)}=0,\qquad
 A_j(z)=\sum_{\lambda\in\Gamma_j}P_{j,\lambda}(z)e^{\lambda z},
\end{equation}
where $n\ge2$, each $\Gamma_j\subset\C$ is finite, and the
multipliers $P_{j,\lambda}$ are nonzero polynomials. At least one
coefficient has a nonzero frequency. Every solution is entire.
The $f^{(n-1)}$ term can be removed by a zero-free entire factor,
which preserves all zeros; see \cite[p.~131]{Laine}.
For background on exponential polynomials and complex differential
equations, see \cite{HITW}.

The initial values can change the zero set substantially. For
example, consider the Schr\"odinger equation \cite[Example~26]{CGHRseparation}
\begin{equation}\label{eq:model-original}
 f''-\frac{1+e^{2z}}4f=0,\qquad z\in\C,
\end{equation}
with potential $V(z)=e^{2z}/4$ and energy $E=-1/4$.
It has the independent zero-free solutions
$f_\pm(z)=e^{-z/2\pm e^z/2}$. Their difference
$f_0=f_+-f_-$ has infinitely many zeros.
Its zeros satisfy $e^z=2\pi i k$,
$k\in\Z\setminus\{0\}$, which lie on horizontal rows.
On a row at height $y=O(\sqrt r)$, the disk extends to
$x=\sqrt{r^2-y^2}=r-y^2/(2r)+o(1)$. The row therefore contributes
asymptotically $e^r e^{-y^2/(2r)}/(2\pi)$ zeros.
There are of order $\sqrt r$ rows with comparable
contributions near the positive real direction, so summing over
the rows gives
\[
 n(r,1/f_0)\sim\frac{\sqrt r\,e^r}{\sqrt2\,\pi^{3/2}}.
\]
Here and below, $n(r,1/f)$ counts zeros in $|z|\le r$, with
multiplicities.

\begin{samepage}
\medskip
\noindent\textbf{Main conclusion.}
There are constants $a>0$, $\beta_*\ge0$ and $C>0$ such that every
solution $f$ of~\eqref{eq:reduced} whose initial values lie outside
a small exceptional set satisfies
\[
 n(r,1/f)\sim C r^{\beta_*+1/2}e^{ar}.
\]
The exceptional set need not be a linear subspace; we therefore introduce the more general notion of a complex cone to describe its geometry.
\end{samepage}

\subsection*{Initial values and the exceptional cone}
\enlargethispage{\baselineskip}
Fix $z_0\in\C$. For $\mathbf c=(c_1,\ldots,c_n)\in\C^n$, let $f_{\mathbf c}$ solve~\eqref{eq:reduced} with
\[
 f_{\mathbf c}^{(k)}(z_0)=c_{k+1},\qquad k=0,\ldots,n-1.
\]
Choose fundamental solutions $f_1,\ldots,f_n$ whose initial
vectors are the standard basis. Then
$f_{\mathbf c}=\sum_{j=1}^n c_jf_j$.
For a complex linear subspace $V$, let $\dim V$ denote its complex dimension.

In a real vector space, a \emph{real cone} with vertex at the
origin is a set containing $0$ and, with each nonzero point $\mathbf x$,
the whole ray $\{t\mathbf x:t\ge0\}$.
A \emph{complex cone} with vertex at the origin is a set
$K\subset\C^n$ containing $0$ and, with each nonzero point
$\mathbf c\in K$, the whole complex line
\[
 \C \mathbf c=\{\lambda \mathbf c:\lambda\in\C\}.
\]
Every complex linear subspace is a complex cone. Conversely,
a complex cone is a complex linear subspace if and only if it
is also closed under addition: $\mathbf c,\mathbf d\in K$ implies $\mathbf c+\mathbf d\in K$.
A \emph{complex linear hyperplane} in $\C^n$ has the form
$H=\{\mathbf c:\ell(\mathbf c)=0\}$, where $\ell$ is a nonzero complex linear
form. Then $\dim H=n-1$. In $\C^2$, the complex
lines through the origin are exactly the complex linear hyperplanes.

In the second-order example above, write $f=Af_++Bf_-$.
The coefficient vectors of the zero-free solutions, together with
the origin, form
\[
 E=\{AB=0\}=\{A=0\}\cup\{B=0\}\subset\C^2.
\]
Thus $E$ consists of exactly two complex lines, each also a
hyperplane. It is a cone but not a linear subspace:
$(1,0),(0,1)\in E$, whereas $(1,1)\notin E$.
The coordinates $(A,B)$ and the initial values are related by
an invertible complex linear map.
Identifying nonzero coefficient vectors that differ by a nonzero
scalar gives the complex projective line $\mathbb P^1(\C)$,
or Riemann sphere. With $t=A/B$ for $B\ne0$ and $t=\infty$
for $B=0$, the two lines in $E$ correspond to $0$ and $\infty$.
The remaining directions form
$\mathbb P^1(\C)\setminus\{0,\infty\}\cong\C\setminus\{0\}$.
More generally, the directions of a complex cone
$K\subset\C^n$ form a subset of $\mathbb P^{n-1}(\C)$.
Here $E$ is closed, nowhere dense and of Lebesgue measure zero,
while its complement $E^c=\{AB\ne0\}$ is open, dense and of
full measure. The set $E^c\cup\{0\}$ is also a complex cone,
but not a linear subspace. 

More generally, the complement of
any complex cone becomes another complex cone after adjoining
the origin.
The cone condition alone allows finitely, countably or
uncountably many complex lines. Even a countable union of
complex lines can be dense while having Lebesgue measure zero.

Pluripotential theory gives a way to describe small exceptional
sets in the initial-value space.
\phantomsection\label{par:pluripotential-definitions}
A function $\Psi:\C^n\to[-\infty,\infty)$ is
\emph{plurisubharmonic} if it is upper semicontinuous, is not
identically $-\infty$, and, for every $\mathbf a,\mathbf b\in\C^n$ with $\mathbf b\ne0$,
the function
\[
 v_{\mathbf a,\mathbf b}(\lambda)=\Psi(\mathbf a+\lambda \mathbf b),\qquad \lambda\in\C,
\]
is subharmonic or identically $-\infty$.
We write $\Psi\in\operatorname{PSH}(\C^n)$.
A set $S\subset\C^n$ is \emph{pluripolar} if there is such a
function $\Psi$ with
\[
 S\subset\{\mathbf c\in\C^n:\Psi(\mathbf c)=-\infty\}.
\]
It is \emph{complete pluripolar} if
$S=\{\mathbf c\in\C^n:\Psi(\mathbf c)=-\infty\}$ for some such $\Psi$;
see \cite[\S2.9 and Chapter~4]{Klimek}.
These sets have Lebesgue measure zero in $\C^n$.
The converse fails. In $\C$, where pluripolar means polar,
the unit circle has zero area but is not polar, whereas its countable
dense subset $\{e^{2\pi iq}:q\in\Q\}$ is polar; see
\cite[Corollary~4.7.7]{Klimek}.
The set $E$ above is complete pluripolar, since
$\log|AB|$ is plurisubharmonic and $E=\{\log|AB|=-\infty\}$.
Not every pluripolar set is complete pluripolar. For example,
$K=\{(w,w):|w|\le1\}\subset\C^2$ is pluripolar because the
plurisubharmonic function $(A,B)\mapsto\log|A-B|$ equals
$-\infty$ on $K$.
Any $\Psi\in\operatorname{PSH}(\C^2)$ equal to $-\infty$ on $K$
must be $-\infty$ on the whole line $A=B$, by restriction to
this line and \cite[Corollary~2.4.6]{Klimek}.

For equation~\eqref{eq:reduced}, we construct a function
$\Psi\in\operatorname{PSH}(\C^n)$, $\Psi\not\equiv-\infty$.
For every initial vector $\mathbf c$ with $\Psi(\mathbf c)>-\infty$, the
corresponding solution $f_{\mathbf c}$ satisfies the asymptotic formula
in the main conclusion, with the same constants $a$, $\beta_*$
and $C$.
The set $\mathcal E=\{\mathbf c\in\C^n:\Psi(\mathbf c)=-\infty\}$ is a
complete pluripolar cone.
Plurisubharmonic functions are used to describe exceptional
parameters in equidistribution and counting problems
\cite{RS,DinhVu}.
To identify the common asymptotic formula for our equation,
we also need to locate the zeros and control their combined
contribution in large disks.

\subsection*{Local zero counts and a global upper bound}

In the Schr\"odinger example~\eqref{eq:model-original}, we have
$a=1$ and may take $s=\varphi(z)=e^z$. An exterior sector
$|s|>R_0$, $\alpha_1<\arg s<\alpha_2$, of angular width less
than $\pi$, has the logarithmic pullbacks
\[
 \{z=x+iy:x>\log R_0,\quad
       \alpha_1+2\pi\ell<y<\alpha_2+2\pi\ell\},
 \qquad \ell\in\Z.
\]
These are parallel horizontal half-strips; compare
\cite[\S5.1]{LWW}.

A transcendental entire map belongs to the Eremenko--Lyubich
class $\mathcal B$ if its critical and finite asymptotic values
form a bounded set. The preimage of the exterior of a sufficiently
large disk then consists of at most countably many unbounded
logarithmic tracts, each mapped as a universal covering onto that
exterior. Cutting the exterior along a radial ray gives an argument
interval of width $2\pi$ and countably many conformal inverse
branches in each tract, with unbounded images~\cite[\S2, p.~993]{EL}.
Exponential sums need not belong to $\mathcal B$;
Sixsmith~\cite[Lemma~8.2]{Sixsmith} gives examples without
logarithmic tracts.

For the general coefficients in~\eqref{eq:reduced}, no
$\mathcal B$ assumption is made. Their frequencies select finitely
many directions; compare the discussion of convex hulls of
frequencies and Stokes rays in~\cite[\S\S2.2 and~3.1]{Li}.
Near each selected direction, we introduce a change of variable
$s=\varphi(z)$ for the whole equation, with countably many inverse
branches over suitable exterior sectors. Their images are
curvilinear strips whose boundaries are asymptotic to parallel
rays in that direction. On each fixed branch, the transformed
coefficients tend to constants as $|s|\to\infty$, uniformly on
closed subsectors. On a conformal branch, a solution and its
expression in the new coordinate have corresponding zeros with
unchanged multiplicities. Their zero counts therefore agree on
corresponding domains. These facts give local zero-counting
asymptotics for every nonzero solution in each strip.

Using the coefficient frequency sets $\Gamma_j$
in~\eqref{eq:reduced}, define the weighted frequency set
$\mathcal F=\{\lambda/(n-j):\lambda\in\Gamma_j,\ 0\le j\le n-2\}$
and put
\begin{equation}\label{eq:weighted}
 a=\max_{\omega\in\mathcal F}|\omega|>0,\qquad
 \Omega_a=\{\omega\in\mathcal F:|\omega|=a\}.
\end{equation}
For $\omega\in\mathcal F\setminus\{0\}$, set
\begin{equation}\label{eq:beta}
 \beta(\omega)=\max_{\substack{0\le j\le n-2\\(n-j)\omega\in\Gamma_j}}
 \frac{\deg P_{j,(n-j)\omega}}{n-j}.
\end{equation}
Fix $\omega\in\Omega_a$, write $\beta=\beta(\omega)$, and choose
an angular interval $J$ of length less than $\pi$. For sufficiently
large $R_0$, $\varphi$ has inverse branches on
$S_J=\{s:|s|>R_0,\ \arg s\in J\}$. Fix one inverse branch and let
$T$ be its image. Then $\varphi:T\to S_J$ is conformal, with
$\varphi^{-1}(s)=\omega^{-1}\psi(s)$. The transformed solution is
$v(s)=f(\omega^{-1}\psi(s))$ on $S_J$. Put
\[
 U(\alpha_1,\alpha_2)=
 \varphi^{-1}\{s\in S_J:\alpha_1<\arg s<\alpha_2\}.
\]
The strips extend in direction $-\arg\omega$, with
$\varphi(z)\sim(\omega z)^\beta e^{\omega z}$ on each fixed strip.

For $E\subset\C$, write
$n_E(r,1/f)=n(E\cap\{|z|\le r\},1/f)$, with multiplicities.
A nondecreasing step function on $J$ means a finite sum
$M(\alpha)=\sum_{\theta\le\alpha}m_\theta$, where
$\theta\in J$ and $m_\theta>0$. It is right-continuous and
zero near the left endpoint; the empty sum is allowed.

\Needspace{20\baselineskip}
\begin{samepage}
\begin{theorem}\label{thm:intro}
 Let $a,\beta,\varphi,J$ and $U(\alpha_1,\alpha_2)$ be as above.
 For every nonzero solution $f$ of~\eqref{eq:reduced}, there is a
 nondecreasing step function $M_f$ on $J$ such that
 \begin{equation}\label{eq:intro-count}
  n_{U(\alpha_1,\alpha_2)}(r,1/f)
  =r^\beta e^{ar}\bigl(M_f(\alpha_2)-M_f(\alpha_1)\bigr)
   +o(r^\beta e^{ar})
 \end{equation}
 whenever $\alpha_1<\alpha_2$ lie in $J$ and neither endpoint
 is a jump point of $M_f$. The possible jump directions lie
 in a finite set determined by the coefficients of the equation and by $\omega$.

 If $M_f$ is constant on a neighbourhood of
 $[\alpha_1,\alpha_2]\subset J$, then
 \[
  n_{U(\alpha_1,\alpha_2)}(r,1/f)=o(r^\beta e^{ar}),
 \]
 even if either or both angular boundary curves are included.
\end{theorem}
\end{samepage}

See \hyperref[subsec:proof-intro]{\S\ref*{subsec:proof-intro}} for the proof
and \hyperref[subsec:nonmaximal-frequencies]{\S\ref*{subsec:nonmaximal-frequencies}}
for a supporting-line criterion allowing $0<|\omega|<a$.

If $M_f$ has a positive jump, the zeros of $f$ have infinite
exponent of convergence; see~\cite[p.~9]{Levin} for the definition.

The direction of the strip $T$ is $-\arg\omega$, independently of
$J$. The length of $J$ determines its asymptotic width. 
The next
theorem makes suitable choices and gives the same positive density
outside two proper linear subspaces of initial values.

\begin{theorem}\label{thm:generic-density}
 Fix $\omega\in\Omega_a$ as in~\eqref{eq:weighted},
 and fix $z_0\in\C$. There are an angular interval $J$,
 numbers $\alpha_1<\alpha_2$ in $J$, a strip
 $U(\alpha_1,\alpha_2)$ as above, and two complex linear
 subspaces $V_-,V_+\subsetneq\C^n$ with $\dim V_\pm\le n-1$
 such that every solution $f_{\mathbf c}$ with $\mathbf c\notin V_-\cup V_+$ satisfies
 \begin{equation}\label{eq:generic-density}
  n_{U(\alpha_1,\alpha_2)}(r,1/f_{\mathbf c})
   \sim c_\omega r^\beta e^{ar},\qquad c_\omega>0.
 \end{equation}
 The sector, subspaces and constant are independent of the
 chosen solution.
\end{theorem}

The proof is given in \hyperref[subsec:proof-generic-density]{Section~\ref*{subsec:proof-generic-density}}.

For each $\omega\in\Omega_a$, write $U_\omega$ and $V_{\omega,\pm}$ for
the strip and subspaces in Theorem~\ref{thm:generic-density}.
These subspaces are complex cones. In~\eqref{eq:model-original}, they
are the hyperplanes $A=0$ and $B=0$ in the coordinates $f=Af_++Bf_-$.
Their union is a complete pluripolar cone and fails to be a linear subspace.
In Example~\ref{ex:fourth-order-subspaces}, both subspaces have
dimension $2$ in $\C^4$, so each has dimension strictly between that
of a line and that of a hyperplane.
Since $\Omega_a$ is finite, every $\mathbf c$ outside the finite union
$\bigcup_{\omega\in\Omega_a}(V_{\omega,-}\cup V_{\omega,+})$ satisfies
\begin{equation}\label{eq:simultaneous-strip-count}
 n_{U_\omega}(r,1/f_{\mathbf c})\sim c_\omega r^{\beta(\omega)}e^{ar},
 \qquad c_\omega>0,
\end{equation}
simultaneously for all $\omega\in\Omega_a$.
For the global count, put
\[
 \beta_*=\max_{\omega\in\Omega_a}\beta(\omega),\qquad
 \Sigma_\varepsilon=
 \{z\ne0:|\arg(\omega z)|<\varepsilon
                 \text{ for some }\omega\in\Omega_a\},
\]
where arguments are principal. However small $\varepsilon>0$ is,
the tail of each fixed strip lies in the corresponding sector.
For initial values outside a union of at most two proper linear
subspaces, the next theorem shows that these sectors contain
asymptotically all zeros.

\begin{theorem}\label{thm:global-zero-distribution}\label{thm:global-zero-upper-bound}
 There is $C>0$, depending only on~\eqref{eq:reduced}, such that
 every nonzero solution satisfies
 \begin{equation}\label{eq:global-zero-upper-bound}
  n(r,1/f)\le C r^{\beta_*+1/2}e^{ar},\qquad r\ge1.
 \end{equation}
 There is also $N>0$, depending only on $n$ and the coefficients,
 such that for every $0<\varepsilon<\pi$ there are constants
 $C_\varepsilon>0$ and $0\le a_\varepsilon<a$, independent of
 the solution, for which
 \begin{equation}\label{eq:global-directional-bound}
  n_{\C\setminus\Sigma_\varepsilon}(r,1/f)
       \le C_\varepsilon r^N e^{a_\varepsilon r},\qquad r\ge1.
 \end{equation}

 There are complex linear subspaces $V_-,V_+\subsetneq\C^n$
 with $\dim V_\pm\le n-1$ such that every solution $f_{\mathbf c}$ with
 $\mathbf c\notin V_-\cup V_+$ satisfies
 \begin{equation}\label{eq:generic-exponent}
  \lim_{r\to\infty}\frac{\log n(r,1/f_{\mathbf c})}r=a,
 \end{equation}
 and, for every $0<\varepsilon<\pi$,
 \[
  n_{\C\setminus\Sigma_\varepsilon}(r,1/f_{\mathbf c})
        =o\bigl(n(r,1/f_{\mathbf c})\bigr).
 \]
\end{theorem}

The proof is given in \hyperref[subsec:proof-global-zero-distribution]{Section~\ref*{subsec:proof-global-zero-distribution}}.

Together with the upper bounds above, the lower bound from
Theorem~\ref{thm:generic-density} for one $\omega\in\Omega_a$
gives~\eqref{eq:generic-exponent} and the concentration statement.
Thus only one pair of exceptional subspaces is needed.
The model~\eqref{eq:model-original} gives sharpness of the power
$1/2$, with $a=1$ and $\beta_*=0$.

\subsection*{Global zero distribution and weak convergence}

To state the main conclusion precisely, put
\begin{equation}\label{eq:dominant-weighted-frequencies}
 \Omega_* =\{\omega\in\Omega_a:\beta(\omega)=\beta_*\}
 \subseteq\Omega_a.
\end{equation}
The inclusion may be strict. Only these directions contribute to the
leading term of the global zero-counting function. Put
$D(r)=(ar)^{\beta_*+1/2}e^{ar}$.

The next theorem also describes how the zeros concentrate near the
disk boundary and near the rays $\arg z=-\arg\omega$,
$\omega\in\Omega_*$.

\begin{samepage}
Let $X=[0,\infty)\times\R$. For $x\in X$, the Dirac measure
$\delta_x$ satisfies $\delta_x(B)=1$ if $x\in B$, and $0$ otherwise,
for every Borel set $B\subset X$. Thus $\int_X g\,\dd\delta_x=g(x)$.
A point of multiplicity $m$ is recorded by $m\delta_x$; namely,
$\int_X g\,\dd(m\delta_x)=m\,g(x)$.
Let $\mu_r$ and $\mu$ be finite positive Borel measures on $X$.
Weak convergence $\mu_r\Longrightarrow\mu$ means that
\[
 \int_X g\,\dd\mu_r\longrightarrow\int_X g\,\dd\mu
\]
for every bounded continuous function
$g:[0,\infty)\times\R\to\R$.
\end{samepage}


\begin{theorem}\label{thm:exact-global}
There are positive constants $C_\omega$, $\omega\in\Omega_*$,
depending only on the equation, and a complete pluripolar cone
$\mathcal E\subset\C^n$ such that the following limits hold for
every $\mathbf c\notin\mathcal E$.
\begin{equation}\label{eq:exact-global-count}
 n(r,1/f_{\mathbf c})\sim
 \left(\sum_{\omega\in\Omega_*}C_\omega\right)D(r).
\end{equation}
Choose $\varepsilon>0$ sufficiently small that the sectors
$|\arg(\omega z)|<\varepsilon$, $\omega\in\Omega_*$, are disjoint.
For a zero $\zeta\ne0$ in $|\zeta|\le r$, set
$u=a(r-|\zeta|)$, $\tau=\sqrt{ar}\arg(\omega \zeta)$, and let $m_{\zeta}$
be its multiplicity.
\begin{equation}\label{eq:joint-global-law}
 \frac1{D(r)}
 \sum_{\substack{f_{\mathbf c}(\zeta)=0,\ 0<|\zeta|\le r\\
                  |\arg(\omega \zeta)|<\varepsilon}}
 m_{\zeta}\,\delta_{(u,\tau)}
 \ \Longrightarrow\
 C_\omega e^{-u}\frac{e^{-\tau^2/2}}{\sqrt{2\pi}}\dd u\dd\tau
\end{equation}
for every $\omega\in\Omega_*$. Zeros outside these sectors
contribute $o(D(r))$.
\end{theorem}

Each constant $C_\omega$ is explicitly computable from the perimeter
of the convex hull of the corresponding limiting characteristic roots.
The proof is given in \hyperref[subsec:proof-exact-global]{Section~\ref*{subsec:proof-exact-global}}.

\begin{figure}[!ht]
\centering
\begin{tikzpicture}[x=1cm,y=1cm,font=\small,
  axis/.style={thin,gray!65,-{Stealth[length=4pt]}},
  boundary/.style={thin,gray!70},
  selected/.style={thin,densely dashed,gray!75},
  profile/.style={blue!75!black,semithick},
  coordinate/.style={blue!75!black,semithick,-{Stealth[length=4pt]}},
  measure/.style={thin,{Stealth[length=3pt]}-{Stealth[length=3pt]}}]
  \begin{scope}[shift={(2.5,2.2)}]
    \def\RR{1.75}
    \def\RI{1.30}
    \def\TH{13}
    \node at (0,2.50) {$z$-plane};
    \fill[gray!12,even odd rule] (0,0) circle (\RR) (0,0) circle (\RI);
    \fill[blue!7] (0,0)--(-\TH:\RR) arc (-\TH:\TH:\RR)--cycle;
    \fill[blue!7] (0,0)--({-90-\TH}:\RR)
       arc ({-90-\TH}:{-90+\TH}:\RR)--cycle;
    \fill[blue!40] (-\TH:\RI)--(-\TH:\RR)
       arc (-\TH:\TH:\RR)--(\TH:\RI) arc (\TH:-\TH:\RI)--cycle;
    \fill[blue!40] ({-90-\TH}:\RI)--({-90-\TH}:\RR)
       arc ({-90-\TH}:{-90+\TH}:\RR)--({-90+\TH}:\RI)
       arc ({-90+\TH}:{-90-\TH}:\RI)--cycle;
    \draw[boundary] (0,0) circle (\RR);
    \draw[boundary,densely dotted] (0,0) circle (\RI);
    \draw[axis] (-2.05,0)--(2.02,0);
    \draw[axis] (0,-2.02)--(0,2.05);
    \draw[selected] (0,0)--(0:\RR);
    \draw[selected] (0,0)--(-90:\RR);
    \foreach \aa in {-13,13,-103,-77}
      {\draw[boundary] (0,0)--(\aa:\RR);}
    \fill (0,0) circle (0.018);
    \node[below left=1pt] at (0,0) {$0$};
    \draw[coordinate] (1.68,0)--(1.37,0);
    \node[blue!75!black,above=2pt] at (1.53,0) {$u$};
    \draw[coordinate] (-11:0.87) arc (-11:11:0.87);
    \node[blue!75!black,anchor=west] at (0.96,-0.30) {$\tau$};
    \draw[thin] (8:0.92)--(1.16,0.66)--(1.36,0.66);
    \node[anchor=west] at (1.39,0.66) {$T/\sqrt{ar}$};
    \draw[measure] (145:\RI)--(145:\RR);
    \draw[thin] (145:1.52)--(-1.64,1.12);
    \node[anchor=east] at (-1.68,1.12) {$h$};
    \draw[thin] (55:\RR)--(1.18,2.02);
    \node[anchor=west] at (1.21,2.02) {$|z|=r$};
    \draw[thin] (105:\RI)--(-0.68,1.83);
    \node[anchor=east] at (-0.71,1.83) {$|z|=r-h$};
    \draw[thin] (1.78,0)--(1.96,-0.23);
    \node[anchor=west] at (2.00,-0.23) {$\arg z=0$};
    \draw[thin] (0,-1.78)--(0.25,-1.88);
    \node[anchor=west] at (0.29,-1.88) {$\arg z=-\pi/2$};
    \node[align=center] at (0,-2.43)
      {$u=a(r-|z|)$,\quad $\tau=\sqrt{ar}\arg(\omega z)$};
  \end{scope}
  \begin{scope}[shift={(7.0,2.8)},x=1.2cm,y=1.5cm]
    \node at (2.15,1.39) {Radial density};
    \fill[blue!16] (0,0)
      -- plot[domain=0:1.1,samples=35] (\x,{exp(-\x)})
      --(1.1,0)--cycle;
    \draw[axis] (-0.10,0)--(4.8,0) node[right] {$u$};
    \draw[axis] (0,0)--(0,1.16);
    \draw[profile] plot[domain=0:4.5,samples=80] (\x,{exp(-\x)});
    \draw[thin,gray!70,densely dashed] (1.1,0)--(1.1,{exp(-1.1)});
    \draw[thin,gray!70] (-0.04,1)--(0.04,1);
    \node[anchor=east] at (-0.08,1) {$1$};
    \node[below=2pt] at (0,0) {$0$};
    \node[below=2pt] at (1.1,0) {$ah$};
    \node[blue!75!black,anchor=west] at (2.25,0.40) {$e^{-u}$};
  \end{scope}
  \begin{scope}[shift={(9.9,0.10)},x=0.85cm,y=3cm]
    \node at (0,0.65) {Angular density};
    \fill[blue!16] (-1.2,0)
      -- plot[domain=-1.2:1.2,samples=50]
         (\x,{exp(-\x*\x/2)/sqrt(2*pi)})
      --(1.2,0)--cycle;
    \draw[axis] (-3.55,0)--(3.72,0) node[right] {$\tau$};
    \draw[axis] (0,0)--(0,0.48);
    \draw[profile] plot[domain=-3.4:3.4,samples=100]
       (\x,{exp(-\x*\x/2)/sqrt(2*pi)});
    \foreach \tt in {-1.2,1.2}
      {\draw[thin,gray!70,densely dashed]
         (\tt,0)--(\tt,{exp(-\tt*\tt/2)/sqrt(2*pi)});}
    \node[below=2pt] at (-1.2,0) {$-T$};
    \node[below=2pt] at (0,0) {$0$};
    \node[below=2pt] at (1.2,0) {$T$};
    \node[blue!75!black,anchor=west] at (1.45,0.38)
      {$\dfrac{e^{-\tau^2/2}}{\sqrt{2\pi}}$};
  \end{scope}
\end{tikzpicture}
\caption{Geometry and limiting densities for $\Omega_*=\{1,i\}$.}
\label{fig:radial-angular-concentration}
\end{figure}
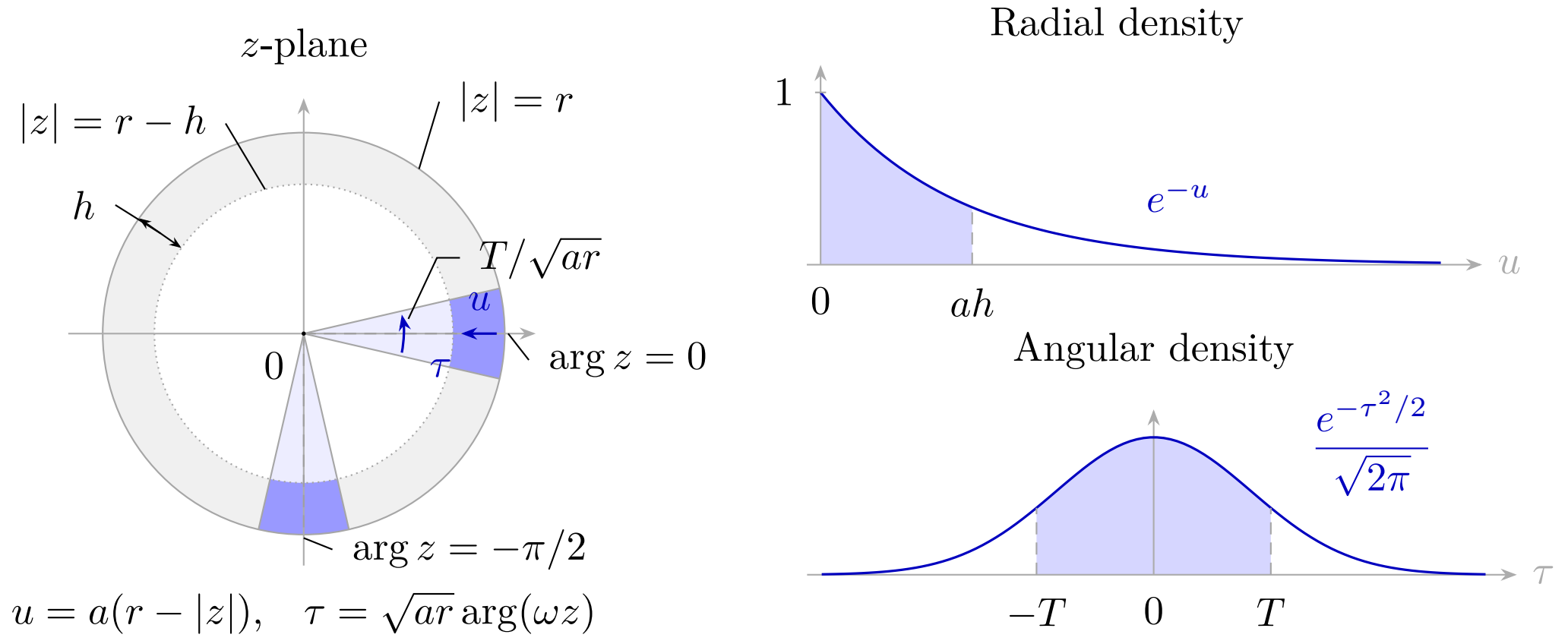

Figure~\ref{fig:radial-angular-concentration} shows the limiting radial
and angular densities after division by $C_\omega$.

Denote the measures on the left and right
of~\eqref{eq:joint-global-law} by $\mu_r$ and $\mu$.
For a Borel set $B\subset X$, $D(r)\mu_r(B)$ counts, with
multiplicities, the zeros in the chosen sector with $0<|\zeta|\le r$
whose coordinates $(u,\tau)$ lie in $B$. For fixed $\omega$, the limiting measure
$\mu$ is common to all $\mathbf c\notin\mathcal E$. The limit
$\mu_r(X)\to C_\omega$ gives only the normalized count in the
chosen sector. For every rectangle $B\subset X$, we also have
\[
 \mu_r(B)\longrightarrow\mu(B)
 =C_\omega\iint_B e^{-u}\frac{e^{-\tau^2/2}}{\sqrt{2\pi}}
 \dd u\dd\tau.
\]
These limits describe the radial and angular distribution.

In particular, the limiting fraction of zeros near the ray
$\arg z=-\arg\omega$ is $C_\omega/\sum_{\nu\in\Omega_*}C_\nu$.
For fixed $h>0$, formula~\eqref{eq:exact-global-count} gives
\[
 \frac{n(r-h,1/f_{\mathbf c})}{n(r,1/f_{\mathbf c})}\longrightarrow e^{-ah}.
\]
Reducing the disk radius by $h$ decreases the zero count
asymptotically by the exponential factor $e^{-ah}$.

Thus the annulus $r-h\le|z|\le r$ contains a limiting fraction
$1-e^{-ah}$ of the zeros. For fixed $T>0$, the sectors
$|\arg(\omega z)|\le T/\sqrt{ar}$, $\omega\in\Omega_*$,
correspond to $|\tau|\le T$ in~\eqref{eq:joint-global-law}.
Equal intervals of $\tau$ cut this annulus into curvilinear
quadrilaterals. For each fixed ray, their limiting fractions of
zeros are proportional to the Gaussian integrals over these
intervals and decrease symmetrically away from $\tau=0$.
For sufficiently large fixed $h,T>0$, the part of the annulus
inside these sectors contains a limiting fraction
arbitrarily close to $1$; see Figure~\ref{fig:radial-angular-concentration}.
These intersections have fixed radial thickness and transverse width of
order $\sqrt r$. This is the scale on which the model~\eqref{eq:model-original}
has rows contributing comparable numbers of zeros.

For example, consider
\[
 f^{(3)}+z^3(e^{3z}+8e^{3iz})f=0,\qquad z\in\C.
\]
For initial vectors outside the exceptional cone, the proportions of zeros
in $|z|\le r$ lying in fixed small sectors about $\arg z=0$ and
$\arg z=-\pi/2$ tend to $1/3$ and $2/3$, respectively
(see Example~\ref{ex:two-direction-weights}).

Put $\mathbf F=(f_1,\ldots,f_n)$. Unless otherwise specified, $\|\cdot\|$
denotes the Euclidean $2$-norm for vectors and its induced operator
norm for matrices. Then $f_{\mathbf c}=\mathbf F\mathbf c$ and, for each fixed $z$,
$\|\mathbf F(z)\|=\max_{\|\mathbf c\|=1}|f_{\mathbf c}(z)|$.
The proof of Theorem~\ref{thm:exact-global} first estimates
$\log\|\mathbf F\|$. Cancellation in a fixed linear combination $f_{\mathbf c}=\mathbf F\mathbf c$
can make $\log|f_{\mathbf c}|$ much smaller. Weighted averages in $z$ of
$\log|f_{\mathbf c}|-\log\|\mathbf F\|$ define a plurisubharmonic potential on the
initial-value space. At initial vectors where this potential is finite, the logarithmic estimates transfer to the fixed solution \(f_{\mathbf c}\).
The one-variable Lelong--Poincar\'e formula
\cite[Chapter~III, (2.15), p.~143]{Demailly} then converts the resulting
logarithmic limits into convergence of zero-counting measures.
The \(-\infty\)-set of the potential has the following geometry.

A set is \emph{locally complete pluripolar} if, in local complex
coordinates near each point, it is the $-\infty$-set of a
plurisubharmonic function that is not identically $-\infty$.
In one complex dimension, we say \emph{locally complete polar}.

\begin{remark}\phantomsection\label{rem:exceptional-cone-geometry}
The defining potential $\Psi\in\operatorname{PSH}(\C^n)$ for the
cone $\mathcal E$ in Theorem~\ref{thm:exact-global} can be normalized so that
\[
 \mathcal E=\{\mathbf c:\Psi(\mathbf c)=-\infty\},\qquad
 \Psi(\lambda \mathbf c)=\Psi(\mathbf c)+\log|\lambda|\quad(\lambda\ne0).
\]
Local integrability gives Lebesgue measure zero and a dense complement.
For every complex affine line $L\not\subset\mathcal E$, the restriction
$\Psi|_L$ is a nontrivial subharmonic function, so $\mathcal E\cap L$
is complete polar in $L$.

In a projective chart $c_j\ne0$, normalize $c_j=1$ and restrict $\Psi$
to this affine hyperplane. The restriction is plurisubharmonic and
not identically $-\infty$: otherwise homogeneity would put the open
set $\{c_j\ne0\}$ inside $\mathcal E$, contradicting measure zero.
Its $-\infty$-set consists exactly of the exceptional directions in
the chart. These directions therefore form a locally complete
pluripolar subset of $\mathbb P^{n-1}(\C)$; when $n=2$, they are
locally complete polar on the Riemann sphere.
\end{remark}

For the periodic subclass, Proposition~\ref{prop:periodic-finite-exception}
identifies the actual exceptional cone as a finite union of proper
complex linear subspaces.

Lebesgue measure zero alone does not control intersections with
complex lines.
For comparison, when $n\ge2$, the complex cone
\[
 K=\{\mathbf c\in\C^n:|c_1|=|c_2|\}
\]
has Lebesgue measure zero by Fubini's theorem. The complex affine
line $L=\{(t,1,0,\ldots,0):t\in\C\}$ is not contained in $K$,
but $K\cap L$ is the unit circle $|t|=1$ in the coordinate $t$,
which is not polar in $L$.

\subsection*{Stability under analytic perturbations}

The global zero-counting asymptotic formula and the joint
radial--angular limit in Theorem~\ref{thm:exact-global} remain
valid under entire coefficient perturbations satisfying the
growth condition below. These perturbations need not be
exponential sums. For an entire function $f$, define
\[
 M(r,f)=\max_{|z|=r}|f(z)|,\qquad
 \tau_1(f)=\limsup_{r\to\infty}
       \frac{\log\max\{1,M(r,f)\}}r.
\]
In particular, $\tau_1(0)=0$. The subscript $1$ refers to the power
$r^1$ in the denominator; functions of order less than one have
$\tau_1(f)=0$~\cite[Chapter~I, \S20, p.~84]{Levin}.

\begin{theorem}\label{thm:perturbation-stability}
Let $A_j$ be the coefficients in~\eqref{eq:reduced}, and let
$R_0,\ldots,R_{n-2}$ be entire functions satisfying
\begin{equation}\label{eq:perturbation-type-gap}
 \max_{0\le j\le n-2}\frac{\tau_1(R_j)}{n-j}<a.
\end{equation}
Fix $z_0\in\C$, and let $g_{\mathbf c}$ solve
\[
 g^{(n)}+\sum_{j=0}^{n-2}(A_j(z)+R_j(z))g^{(j)}=0,
 \qquad z\in\C,
\]
with $(g_{\mathbf c}(z_0),\ldots,g_{\mathbf c}^{(n-1)}(z_0))=\mathbf c\in\C^n$.
With the notation and constants of Theorem~\ref{thm:exact-global},
there is a complete pluripolar cone
$\widetilde{\mathcal E}\subset\C^n$ such that the following
limits hold for every $\mathbf c\notin\widetilde{\mathcal E}$.
\begin{equation}\label{eq:perturbed-global-count}
 n(r,1/g_{\mathbf c})\sim
 \left(\sum_{\omega\in\Omega_*}C_\omega\right)D(r).
\end{equation}
For every $\omega\in\Omega_*$,
\begin{equation}\label{eq:perturbed-joint-law}
 \frac1{D(r)}
 \sum_{\substack{g_{\mathbf c}(\zeta)=0,\ 0<|\zeta|\le r\\
                   |\arg(\omega\zeta)|<\varepsilon}}
 m_\zeta\,\delta_{(u,\tau)}
 \ \Longrightarrow\
 C_\omega e^{-u}\frac{e^{-\tau^2/2}}{\sqrt{2\pi}}\dd u\dd\tau.
\end{equation}
Zeros outside these sectors contribute $o(D(r))$.
\end{theorem}

The proof is given in \hyperref[subsec:proof-perturbation-stability]{Section~\ref*{subsec:proof-perturbation-stability}}.

\begin{corollary}\label{cor:subunit-order-perturbations}
The conclusions of Theorem~\ref{thm:perturbation-stability} hold for
arbitrary entire perturbations $R_j$ of exponential type strictly less
than $2a$, where $a$ is defined in~\eqref{eq:weighted}.
\end{corollary}
\begin{proof}
Since $n-j\ge2$, we have
$\tau_1(R_j)/(n-j)\le\tau_1(R_j)/2<a$, so~\eqref{eq:perturbation-type-gap} holds.
\end{proof}

Example~\ref{ex:subunit-order-perturbations} applies
Corollary~\ref{cor:subunit-order-perturbations} to a canonical product
and a hypergeometric function.
The exceptional cone $\widetilde{\mathcal E}$ may depend on the
perturbations. Proposition~\ref{prop:periodic-finite-exception} gives a
more precise description of the actual exceptional set for a periodic
subclass. Multiplying each $R_j$ by any fixed constant preserves
\eqref{eq:perturbation-type-gap}; exceptional solutions and individual
strip counts require separate analysis.

\subsection*{Background and related work}

For polynomial coefficients, the positive rational orders of
transcendental solutions can be found from Newton--Puiseux
diagrams~\cite{HelmrathNikolaus} or directly from the coefficient
degrees~\cite{GSW}.
When at least one coefficient is transcendental, Frei's
theorem~\cite{Frei} shows that finite-order solutions form a proper
subspace of the solution space.
This subspace can even be trivial: Langley~\cite[Theorem~2]{Langley1986}
settled a problem of M.~Ozawa by proving that every nonzero solution of
$f''+e^{-z}f'+P(z)f=0$ on $\C$ has infinite order whenever $P$ is a
nonconstant polynomial.

However, infinite-order solutions can still have a finite exponent of convergence of zeros or even be zero-free, as
in~\eqref{eq:model-original}. For $f''+Af=0$ on $\C$ with $A$ transcendental entire, the Bank--Laine conjecture
predicted zeros of infinite exponent of convergence for the product
of any two independent solutions when $A$ is 
of finite noninteger order. Bergweiler and Eremenko~\cite{BE2017,BE2019}
used quasiconformal methods to construct counterexamples.
For transcendental $A$ and a real solution basis whose product has
finite order and infinitely many zeros, all real,
Langley~\cite{LangleyB,LangleyReal} obtained sharp lower bounds for
the exponent of convergence of these zeros.
Bergweiler, Eremenko and Rempe~\cite{BER} refined these bounds to
discrete restrictions on the coefficient order and raised the
question of extensions to higher-order equations
\cite[Remark~1.7]{BER}.
For further background on complex LDEs,
see the growth estimates in~\cite{HKR,CGHRgrowth}, the survey
\cite{HITW} and the problems in~\cite{GundersenQuestions}.

A related question is where the zeros lie in the complex plane.
For second-order equations with coefficients
$\sum_j Q_j(z)e^{jP(z)}$, Bank~\cite[Theorems~1 and~3(A)]{Bank1988}
constructs two linearly independent solutions for each selected band
$|\Ima P(z)-\kappa|<\delta$ in a suitable sector.
Every linear combination with both coefficients nonzero has zeros
of infinite exponent of convergence in that band.
Langley~\cite[\S3]{LangleyB}
uses logarithmic tracts of a coefficient in the Eremenko--Lyubich
class $\mathcal B$ to obtain sectorial asymptotics after a conformal
change of variable. For higher-order equations with periodic
coefficients, Bank and Langley~\cite{BL} and Shimomura~\cite{Shimomura}
study solutions with sparse zeros; Shimomura also transfers
sectorial estimates to strips by $t=e^z$
\cite[Proposition~5.1 and Lemma~6.2]{Shimomura}.
A global asymptotic formula requires uniform estimates when
combining the contributions from the different strips.

The exceptional initial values are also related to classical value
distribution. With $\mathbf F$ as above, put, for $\mathbf c\ne0$ and $r\ge1$,
\[
 N_{\mathbf c}(r)=\int_1^r n(t,1/f_{\mathbf c})\,\frac{\mathrm dt}{t},\qquad
 T_{\mathbf F}(r)=\frac1{2\pi}\int_0^{2\pi}
       \log\|\mathbf F(re^{i\theta})\|\,\mathrm d\theta-\log\|\mathbf F(0)\|.
\]
The first quantity is the usual integrated zero count, up to an additive
constant; the second measures the growth of the whole solution basis.
Sadullaev~\cite[Proposition~3.4 and Theorem~3.5]{Sadullaev} implies
that $N_{\mathbf c}(r)\sim T_{\mathbf F}(r)$ outside a pluripolar set of coefficient
directions. His result is stated for hyperplanes and holomorphic maps.
Here the hyperplane $\sum c_jX_j=0$ contains
$[f_1(z):\cdots:f_n(z)]$ exactly when $f_{\mathbf c}(z)=0$, so its
intersection parameters are precisely the zeros of the corresponding
solution. The projective curve is well defined because the fundamental
solutions have no common zero. This explains the connection with
exceptional initial values. To obtain a concrete counting formula by
this route, one must still determine $T_{\mathbf F}(r)$ for the given equation.

There is a complementary formula for averages over the initial values.
Give the coefficients $c_j$ independent standard complex Gaussian
distributions, each with density $\pi^{-1}e^{-|w|^2}\,\mathrm dA(w)$,
and write
\[
 \nu_{\mathbf c}=\sum_{f_{\mathbf c}(\zeta)=0}m_\zeta\delta_\zeta
\]
for the zero-counting measure of $f_{\mathbf c}$. Thus
$\mathbb E\nu_{\mathbf c}(B)$ is the expected number of zeros in a region $B$,
counted with multiplicity. The Edelman--Kostlan formula, as proved in
Sodin~\cite[Theorem~1]{Sodin}, gives
\[
 \mathbb E\nu_{\mathbf c}=\frac1{2\pi}\Delta\log\|\mathbf F\|.
\]
Here $\Delta=\partial_x^2+\partial_y^2$, with $z=x+iy$, is understood
in the sense of distributions. The formula shows why it is natural
to estimate $\log\|\mathbf F\|$ before choosing a particular initial vector:
its Laplacian determines the average zero distribution of the solutions.
Obtaining the same asymptotic distribution for each fixed vector
outside a specified exceptional set requires a further argument.
Related uses of plurisubharmonic functions to control exceptional
parameters appear in value distribution and other counting problems;
see Russakovskii and Shiffman~\cite[Theorem~1.2 and \S6]{RS} and
Dinh and Vu~\cite[Theorem~1.2]{DinhVu}.

\bigskip
\noindent\textbf{Organization of the paper.}
Section~\ref{sec:coefficients} introduces the conformal coordinates and
states the uniform growth estimates. Section~\ref{sec:sectorial}
gives the sectorial zero counts and solution-space tools.
Section~\ref{sec:global} gives uniform zero bounds, states the
perimeter average in moving rectangles, and develops the pluripotential
argument for the common global asymptotic formula.
Section~\ref{sec:auxiliary-proofs} collects the longer, relatively
self-contained auxiliary proofs of
Lemmas~\ref{lem:moving-inverse}--\ref{lem:uniform-fundamental-row},
\ref{lem:slower-growth-subspace} and~\ref{lem:physical-window}.
Section~\ref{sec:main-proofs} then uses the preceding tools to prove
the main results of the introduction, namely
Theorems~\ref{thm:intro}--\ref{thm:perturbation-stability}.
Section~\ref{sec:questions} estimates zeros near critical rays,
classifies the exceptional initial values for periodic coefficients,
and gives examples of entire perturbations and an application to
osculating hyperplanes before concluding with further questions.

\section{Conformal coordinates and uniform growth}
\label{sec:coefficients}

We construct coordinates in which the equation tends to a
constant-coefficient equation, then state the estimates for the
inverse maps and fundamental solutions used below.

For positive quantities $u,v$, we write $u\asymp v$ if
$cv\le u\le Cv$ for constants $c,C>0$.
Subscripts on $O$ and $\asymp$ indicate parameters on which the
bounding constants may depend.
\phantomsection\label{par:positive-negative-parts}
For $x\in\R$, write $x^+=\max\{x,0\}$ and $x^-=\max\{-x,0\}$
for its positive and negative parts. We use $\log^+t=(\log t)^+$
for $t>0$, with $\log^+0=0$.

\subsection{Conformal coordinates and the limiting equation}\label{subsec:limiting-equation}

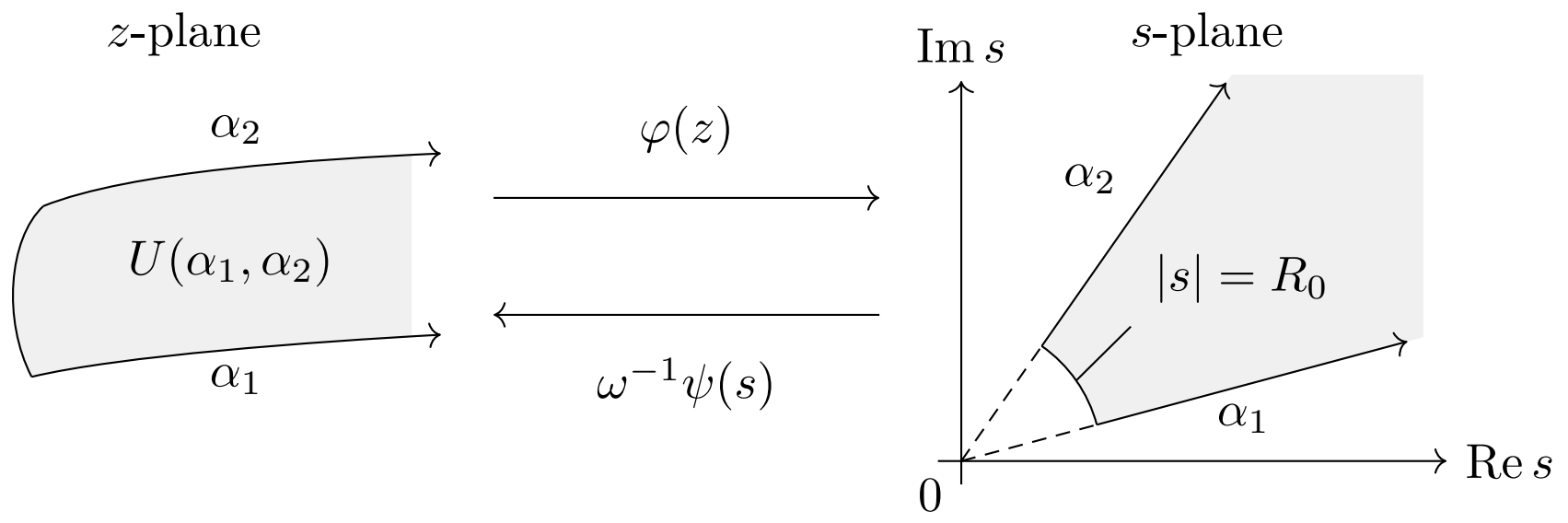
\begin{figure}[!htb]
\centering
\begin{tikzpicture}[x=0.84cm,y=0.84cm,font=\small,
  boundary/.style={thin},axis/.style={->,thin}]
  \node at (2.25,3.65) {$z$-plane};
  \fill[gray!12]
    (1.05,2.18) .. controls (1.8,2.48) and (3.35,2.58) .. (4.2,2.62)
    -- (4.2,1.07)
    .. controls (3.35,1.02) and (1.7,0.9) .. (0.95,0.72)
    .. controls (0.7,1.2) and (0.75,1.9) .. (1.05,2.18) -- cycle;
  \draw[boundary,->]
    (1.05,2.18) .. controls (1.8,2.48) and (3.35,2.58) .. (4.45,2.63);
  \draw[boundary,->]
    (0.95,0.72) .. controls (1.7,0.9) and (3.35,1.02) .. (4.45,1.08);
  \draw[boundary]
    (0.95,0.72) .. controls (0.7,1.2) and (0.75,1.9) .. (1.05,2.18);
  \node at (2.65,1.7) {$U(\alpha_1,\alpha_2)$};
  \node[above] at (2.7,2.55) {$\alpha_2$};
  \node[below] at (2.7,0.98) {$\alpha_1$};

  \draw[->,thin] (4.9,2.25) -- (8.2,2.25)
    node[midway,above=5pt] {$\varphi(z)$};
  \draw[<-,thin] (4.9,1.25) -- (8.2,1.25)
    node[midway,below=5pt] {$\omega^{-1}\psi(s)$};

  \begin{scope}[shift={(8.9,0)}]
    \node at (2.1,3.65) {$s$-plane};
    \begin{scope}
      \clip (-0.1,-0.1) rectangle (3.95,3.3);
      \fill[gray!12] (15:1.2) -- (15:10)
        arc[start angle=15,end angle=55,radius=10]
        -- (55:1.2) arc[start angle=55,end angle=15,radius=1.2] -- cycle;
    \end{scope}
    \draw[axis] (-0.2,0) -- (4.15,0) node[right] {$\Rea s$};
    \draw[axis] (0,-0.2) -- (0,3.25) node[above] {$\Ima s$};
    \node[below left] at (0,0) {$0$};
    \draw[densely dashed,thin] (0,0) -- (15:1.2);
    \draw[densely dashed,thin] (0,0) -- (55:1.2);
    \draw[boundary,->] (15:1.2) -- (15:3.95);
    \draw[boundary,->] (55:1.2) -- (55:3.95);
    \draw[boundary] (15:1.2)
      arc[start angle=15,end angle=55,radius=1.2];
    \node[below] at (15:2.5) {$\alpha_1$};
    \node[above left] at (55:2.6) {$\alpha_2$};
    \draw[thin] (35:1.2) -- (1.45,1.15)
      node[above right=1pt] {$|s|=R_0$};
  \end{scope}
\end{tikzpicture}
\caption{Pullback of a sector.}
\label{fig:domain-pullback}
\end{figure}

Fix $\omega\in\Omega_a$ and put $\beta=\beta(\omega)$, with
$\Omega_a$ and $\beta(\omega)$ as in~\eqref{eq:weighted}
and~\eqref{eq:beta}. For
$\Rea(\omega z)>0$, define
\begin{equation}\label{eq:phase}
 \varphi(z)=\Phi(\omega z)
 =C+\int_1^{\omega z}\xi^\beta e^\xi\dd\xi,
\end{equation}
where $C\in\C$ is fixed, the path lies in $\Rea\xi>0$, and
$\xi^\beta$ uses $|\arg\xi|<\pi/2$. For an angular interval
$J$ of length less than $\pi$ and an integer $\ell$, the map
$\varphi$ has an inverse branch $\omega^{-1}\psi$ on
$S_J=\{s:|s|>R_0,\ \arg s\in J\}$ for sufficiently large $R_0$, with
\[
 \psi(s)=\log s-\beta\log\log s+2\pi i\ell+o(1).
\]
Here $\log s=\log|s|+i\arg s$ with $\arg s\in J$, and
$\log\log s$ is principal. For each fixed $\ell$, the expansion
is uniform on closed subsectors. In particular, for $s=\rho e^{i\theta}$,
\[
 \omega z=\log\rho-\beta\log\log\rho
           +i(\theta+2\pi\ell)+o(1).
\]
Thus the strips extend in the direction $-\arg\omega$; the branch
index determines their transverse position, and $J$ determines
their width. Lemma~\ref{lem:phase-equation} proves these inverse
estimates. Figure~\ref{fig:domain-pullback} shows one sector and
its pullback.

For the set $\mathcal F$ of weighted frequencies used
in~\eqref{eq:weighted}, put
\[
 \mathcal H=\conv(\{0\}\cup\mathcal F).
\]

The hull $\mathcal H$ lies in the complex frequency plane, and
the fixed $\omega\in\Omega_a$ is a vertex of maximal modulus.
For a frequency $\lambda\in\Gamma_j$ of $A_j$
in~\eqref{eq:reduced},
\begin{equation}\label{eq:frequency-separation}
 \Rea(\lambda/\omega)\le n-j,
\end{equation}
with equality only when $\lambda=(n-j)\omega$.
Indeed, $|\lambda/\omega|\le n-j$ by the maximality of
$|\omega|$, and equality in the real-part inequality also requires
$\lambda/\omega$ to be positive real. In the coordinate $t=\omega z$,
the ratio of any other exponential to $e^{(n-j)t}$ therefore
decays exponentially as $\Rea t\to+\infty$ with $\Ima t$ bounded.
Among the terms with frequency $(n-j)\omega$, the polynomial
degrees determine the limiting coefficients.

\begin{figure}[htbp]
\centering
\begin{tikzpicture}[x=0.85cm,y=0.72cm,font=\small,
 dot/.style={circle,fill=black,inner sep=1.5pt}]
 \draw[->,thin] (-.2,0) -- (6.7,0) node[right] {$j$};
 \draw[->,thin] (0,-.15) -- (0,4.65) node[above] {$d$};
 \filldraw[fill=gray!10,thin] (0,4) -- (1,1) -- (6,0) -- cycle;
 \draw[very thick] (0,4) -- (6,0);
 \draw[densely dashed,thin] (-.3,4.2) -- (0,4);
 \foreach \p in {(0,4),(1,1),(2,1),(3,2),(4,1),(6,0)}
   \node[dot] at \p {};
 \node[above right] at (3,2) {$(j,d_j)$};
 \node[below] at (6,0) {$(n,0)$};
 \node at (3.7,3.7) {$d=\beta(n-j)$};
\end{tikzpicture}
\caption{The supporting edge of $\mathcal N_\omega$.}
\label{fig:newton-edge}
\end{figure}
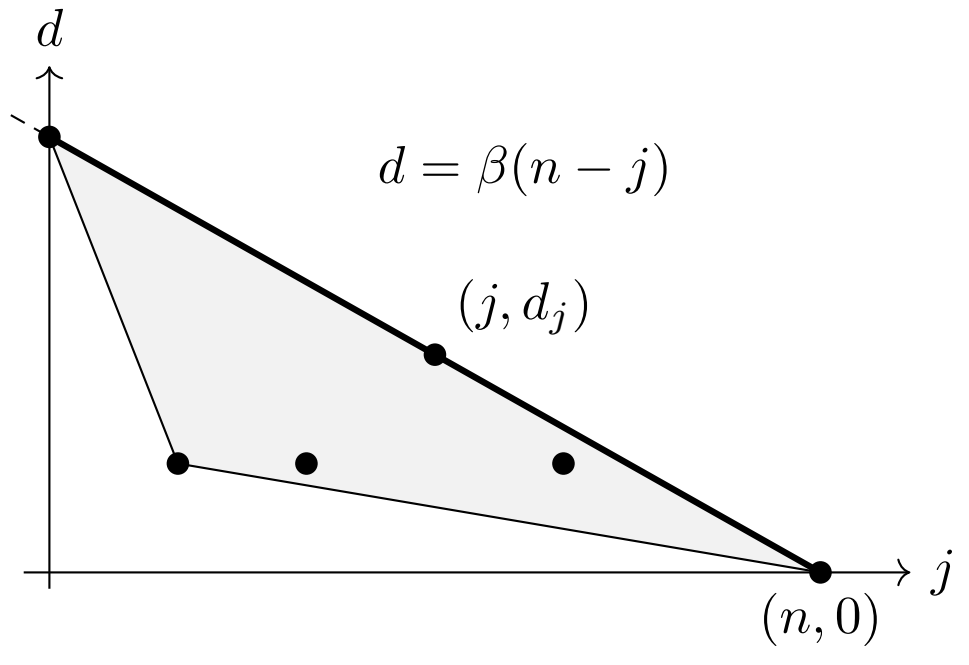

For a selected term write
\[
 d_j=\deg P_{j,(n-j)\omega},\qquad
 a_j=\operatorname{lc}(P_{j,(n-j)\omega}).
\]
Here $\operatorname{lc}$ denotes the leading coefficient.
Only indices for which this term exists are used. The corresponding
Newton polygon is
\[
 \mathcal N_\omega=\conv\bigl(\{(n,0)\}\cup
       \{(j,d_j):(n-j)\omega\in\Gamma_j\}\bigr).
\]
Its coordinates are the derivative order and the polynomial degree.
By~\eqref{eq:beta}, the upper supporting edge
incident with $(n,0)$ lies on $d=\beta(n-j)$ and has slope $-\beta$.
The \emph{edge indices} are the coefficient indices $j$ for which
$(j,d_j)$ lies on this supporting edge. They form the set
\[
 I_\omega=\{j:(n-j)\omega\in\Gamma_j,
                        \ d_j/(n-j)=\beta\}.
\]
This set is nonempty; see
Figure~\ref{fig:newton-edge}.

Put $t=\omega z$, $u(t)=f(t/\omega)$, and
$B_j(t)=\omega^{j-n}A_j(t/\omega)$. Since $\omega\ne0$, the
functions $u$ and $B_j$ are entire, and~\eqref{eq:reduced} becomes
\begin{equation}\label{eq:normalized}
 u^{(n)}+\sum_{j=0}^{n-2}B_j(t)u^{(j)}=0,\qquad t\in\C.
\end{equation}
Indeed,
\begin{equation}\label{eq:normalized-coefficients}
 B_j(t)=\omega^{j-n}\sum_{\lambda\in\Gamma_j}
 P_{j,\lambda}(t/\omega)e^{(\lambda/\omega)t},\qquad t\in\C,
\end{equation}
so each $B_j$ is again an exponential polynomial in $t$.
At the selected frequency $\lambda=(n-j)\omega$, the highest-degree part is
\[
 a_j\omega^{j-n-d_j}t^{d_j}e^{(n-j)t}.
\]
Fix an inverse branch $\psi:S_J\to\{t:\Rea t>0\}$ of $\Phi$;
thus $\Phi(\psi(s))=s$ for $s\in S_J$. Set
\[
 v(s)=u(\psi(s))\quad(s\in S_J),\qquad
 \Phi'(t)=t^\beta e^t\quad(\Rea t>0),
\]
with the power branch in~\eqref{eq:phase}.

For $t\in\psi(S_J)$, repeated differentiation of
$u(t)=v(\Phi(t))$ gives, for $m\ge0$,
\begin{equation}\label{eq:chain-rule-proof}
 u^{(m)}(t)=\sum_{k=0}^m C_{m,k}(t)v^{(k)}(\Phi(t)).
\end{equation}
Differentiation of~\eqref{eq:chain-rule-proof} gives the recurrence
\begin{equation}\label{eq:chain-coefficients}
 C_{0,0}=1,\qquad
 C_{m+1,k}=C'_{m,k}+\Phi' C_{m,k-1},\qquad C_{m,m}=(\Phi')^m.
\end{equation}
The functions $C_{m,k}$ are analytic on $\Rea t>0$.
Here $C_{m,k}=0$ for $k<0$ or $k>m$. In particular,
$C_{m,0}=0$ for $m\ge1$.

Substitute~\eqref{eq:chain-rule-proof} into~\eqref{eq:normalized}
and divide by $(\Phi'(t))^n$ on $\psi(S_J)$. Since $C_{n,n}=(\Phi')^n$, we obtain
\begin{equation}\label{eq:phase-equation}
 v^{(n)}(s)+\sum_{k=0}^{n-1}b_k(s)v^{(k)}(s)=0,\qquad s\in S_J,
\end{equation}
where $S_J=\{s:|s|>R_0,\ \arg s\in J\}$. The coefficients are
analytic on $S_J$. Since $C_{j,k}=0$ for $j<k$ and
$C_{k,k}=(\Phi')^k$, for $s\in S_J$ and $t=\psi(s)$ we have
$b_{n-1}(s)=C_{n,n-1}(t)/(\Phi'(t))^n$ and
\begin{equation}\label{eq:phase-coefficients}
 b_k(s)=\frac{B_k(t)}{(\Phi'(t))^{n-k}}
       +\frac{C_{n,k}(t)+\sum_{j=k+1}^{n-2}
                    B_j(t)C_{j,k}(t)}{(\Phi'(t))^n},
       \qquad 0\le k\le n-2.
\end{equation}
Define
\begin{equation}\label{eq:limiting-coefficients}
 c_k=\begin{cases}
 a_k\omega^{k-n-d_k},&k\in I_\omega,\\
 0,&k\notin I_\omega,
 \end{cases}
 \quad 0\le k\le n-2,\qquad c_{n-1}=0.
\end{equation}

\begin{lemma}\label{lem:phase-coefficient-limit}
For the chosen inverse branch $\psi$, the coefficients $b_k$ of
\eqref{eq:phase-equation} and the constants $c_k$ in
\eqref{eq:limiting-coefficients} satisfy
\[
 b_k(s)\longrightarrow c_k,\qquad 0\le k<n,
\]
as $|s|\to\infty$, uniformly on closed subsectors of $S_J$.
\end{lemma}

\begin{proof}
We first consider the ratios $B_k/(\Phi')^{n-k}$
in~\eqref{eq:phase-coefficients}.
Using~\eqref{eq:normalized-coefficients} and $\Phi'(t)=t^\beta e^t$,
we obtain, for $\Rea t>0$,
\[
 \frac{B_k(t)}{(\Phi'(t))^{n-k}}
 =\omega^{k-n}\sum_{\lambda\in\Gamma_k}
   t^{-\beta(n-k)}P_{k,\lambda}(t/\omega)
   e^{(\lambda/\omega-(n-k))t},\qquad 0\le k\le n-2.
\]
Take $\Rea t\to+\infty$, with $\Ima t$ in a fixed bounded interval.
For $\lambda\ne(n-k)\omega$,~\eqref{eq:frequency-separation} gives
$\Rea(\lambda/\omega)<n-k$, so the exponential factor decays;
the other factors grow at most as a power of $|t|$. Each such
summand therefore tends to zero, uniformly in this interval.
If $(n-k)\omega\notin\Gamma_k$, all summands are of this kind,
and the ratio tends to zero.

If $(n-k)\omega\in\Gamma_k$, its exponential factor is $1$.
The summand at this frequency is
\[
 \begin{aligned}
 \omega^{k-n}t^{-\beta(n-k)}P_{k,(n-k)\omega}(t/\omega)=a_k\omega^{k-n-d_k}t^{d_k-\beta(n-k)}
       +O\!\left(|t|^{d_k-1-\beta(n-k)}\right).
 \end{aligned}
\]
By the definition of $\beta$, $d_k-\beta(n-k)\le0$.
Consequently the error term tends to zero. The leading term
tends to $a_k\omega^{k-n-d_k}$ when $k\in I_\omega$, and to zero
otherwise. Together with the other frequencies, this proves
\[
 \frac{B_k(t)}{(\Phi'(t))^{n-k}}\longrightarrow c_k,
 \qquad 0\le k\le n-2,
\]
uniformly on each fixed horizontal strip.

It remains to show that the other terms in~\eqref{eq:phase-coefficients}
tend to zero. By Fa\`a di Bruno's formula~\cite[Eq.~(1.4.13)]{DLMF},
each $C_{m,k}$ is a finite linear combination of products of $k$ positive-order derivatives
of $\Phi$; see also~\cite[proof of Theorem~1]{CGHR}.
On each fixed horizontal strip, as $\Rea t\to+\infty$,
\[
 \Phi',\Phi'',\ldots,\Phi^{(n)}=O(\Phi').
\]
Hence $C_{m,k}=O((\Phi')^k)$.
Since $|\Phi'(t)|\to\infty$ and the ratios $B_j/(\Phi')^{n-j}$
are bounded, the remaining terms satisfy
$C_{n,k}/(\Phi')^n=O((\Phi')^{k-n})=o(1)$ for $0\le k<n$
and $B_jC_{j,k}/(\Phi')^n=O((\Phi')^{k-j})=o(1)$
for $0\le k<j\le n-2$.
In particular,
\[
 b_{n-1}(s)=\binom n2\left.\frac{\Phi''(t)}{(\Phi'(t))^2}
      \right|_{t=\psi(s)}=O(s^{-1}).
\]
Lemma~\ref{lem:phase-equation} verifies that the inverse branch
maps each closed subsector into such a horizontal strip, so the
above estimates are uniform on closed subsectors.
\end{proof}

Lemma~\ref{lem:phase-coefficient-limit} gives the constant-coefficient
limiting equation
\[
 y^{(n)}+\sum_{j=0}^{n-1}c_jy^{(j)}=0,\qquad s\in\C.
\]
Define its \textit{characteristic polynomial} by
\begin{equation}\label{eq:characteristic}
 p(\mu)=p_\omega(\mu)=\mu^n+\sum_{k=0}^{n-1}c_k\mu^k.
\end{equation}
The limiting equation has an exponential solution $y(s)=e^{\mu s}$
exactly when $p(\mu)=0$, since $y^{(j)}=\mu^je^{\mu s}$.
Since $I_\omega\ne\varnothing$, at least one $c_k$ with $k\le n-2$
is nonzero.
Write $\Lambda=\{\mu\in\C:p(\mu)=0\}$ for its distinct roots.

Thus the vertices $\omega$ of $\mathcal H$ with maximal modulus
determine the main rays $\arg z=-\arg\omega$. For each such $\omega$,
the Newton polygon $\mathcal N_\omega$ is formed using the degrees of
the polynomial factors $P_{k,(n-k)\omega}$ that occur in the coefficients;
its upper supporting edge through $(n,0)$ determines $\beta$ and
the set $I_\omega$ of indices with nonzero limiting coefficients $c_k$.
We will use the roots of $p$ to determine the possible
angular zero densities in the $s$-plane.

We first obtain a radial growth estimate for the derivative vector.
Fixed equivalent norms give the same exponential growth rates.
For a nonzero solution $v$ of~\eqref{eq:phase-equation}, set
\begin{equation}\label{eq:phase-derivative-vector}
 \mathbf Y(s)=(v(s),v'(s),\ldots,v^{(n-1)}(s))^{\mathsf T}.
\end{equation}
We write vectors and matrices in boldface. The scalar equation is equivalent
to the first-order system
\begin{equation}\label{eq:phase-companion-system}
 \frac{\mathrm d\mathbf Y}{\mathrm ds}=\mathbf A(s)\mathbf Y,\qquad s\in S_J,
 \qquad
 \mathbf A(s)=\begin{pmatrix}
 0&1&0&\cdots&0\\
 0&0&1&\cdots&0\\
 \vdots&\vdots&\vdots&\ddots&\vdots\\
 0&0&0&\cdots&1\\
 -b_0(s)&-b_1(s)&-b_2(s)&\cdots&-b_{n-1}(s)
 \end{pmatrix}.
\end{equation}
Replacing $b_j(s)$ by $c_j$ gives a constant matrix $\mathbf A_0$
with $\mathbf A(s)\to\mathbf A_0$. Expanding the determinant along
the last row gives
\begin{equation}\label{eq:companion-characteristic}
 \det(\mu\mathbf I-\mathbf A_0)
 =\mu^n+\sum_{j=0}^{n-1}c_j\mu^j=p(\mu).
\end{equation}
The following consequence of Coppel's theorem applies to this system.

\begin{lemma}\label{lem:radial-system-growth}
Let $S_J=\{s:|s|>R_0,\ \arg s\in J\}$ be a sector with a
continuous argument, where $R_0>0$ and $J$ is an open interval of
width less than $2\pi$. Suppose that $\mathbf A:S_J\to\C^{n\times n}$
is analytic and converges to a constant matrix $\mathbf A_0$ as $|s|\to\infty$,
uniformly on each closed subsector. Let $\Lambda$ be the set of
eigenvalues of $\mathbf A_0$. Every nonzero solution of
\[
 \frac{\mathrm d\mathbf Y}{\mathrm ds}=\mathbf A(s)\mathbf Y,\qquad s\in S_J,
\]
satisfies
\begin{equation}\label{eq:radial-system-growth}
 \lim_{r\to\infty}\frac{\log\|\mathbf Y(re^{i\theta})\|}{r}
 \in\{\Rea(\mu e^{i\theta}):\mu\in\Lambda\},\qquad\theta\in J.
\end{equation}
\end{lemma}

\begin{proof}
The vector $\mathbf Y$ never vanishes, by uniqueness for the linear system.
Fix $\theta\in J$ and $r_0>R_0$. Put
$\mathbf B_\theta(r)=e^{i\theta}\mathbf A(re^{i\theta})$. With $r$ as the real
independent variable, the system becomes
\[
 \frac{\mathrm d}{\mathrm dr}\mathbf Y(re^{i\theta})
 =\mathbf B_\theta(r)\mathbf Y(re^{i\theta}),\qquad r\ge r_0.
\]
Since $\mathbf B_\theta(r)\to e^{i\theta}\mathbf A_0$, the number
$M=\sup_{r\ge r_0}\|\mathbf B_\theta(r)\|$ is finite. Choose
$\kappa>M$ and set $\mathbf w(r)=e^{-\kappa r}\mathbf Y(re^{i\theta})$.
Gronwall's inequality gives, for $r\ge r_0$,
\[
 \|\mathbf Y(re^{i\theta})\|
 \le \|\mathbf Y(r_0e^{i\theta})\|e^{M(r-r_0)}.
\]
Multiplying by $e^{-\kappa r}$, we obtain
\begin{equation}\label{eq:perron-bounded-solution}
 \|\mathbf w(r)\|\le\|\mathbf w(r_0)\|
 e^{-(\kappa-M)(r-r_0)}\longrightarrow0.
\end{equation}
Thus $\mathbf w$ remains in a fixed ball in $\C^n$. Differentiating yields
\[
 \frac{\mathrm d\mathbf w}{\mathrm dr}
 =(e^{i\theta}\mathbf A_0-\kappa \mathbf I)\mathbf w+\mathbf E_\theta(r)\mathbf w,
 \qquad r\ge r_0,
\]
where $\mathbf E_\theta(r)=\mathbf B_\theta(r)-e^{i\theta}\mathbf A_0$ is continuous
and tends to zero. The first matrix is constant, and the perturbation satisfies
\begin{equation}\label{eq:perron-perturbation-bound}
 \|\mathbf E_\theta(r)\mathbf x\|\le\|\mathbf E_\theta(r)\|\,\|\mathbf x\|,
 \qquad r\ge r_0,\quad \mathbf x\in\C^n,
\end{equation}
as well as
\begin{equation}\label{eq:perron-mean-smallness}
 \int_r^{r+1}\|\mathbf E_\theta(t)\|\dd t
 \le\sup_{t\ge r}\|\mathbf E_\theta(t)\|\longrightarrow0.
\end{equation}
By~\eqref{eq:perron-bounded-solution}--\eqref{eq:perron-mean-smallness},
we may apply Coppel's theorem
\cite[Ch.~IV, Theorem~5, p.~97]{Coppel} to $\mathbf w$.
Since $\mathbf Y$ never vanishes, $\mathbf w$ cannot be eventually zero.
The eigenvalues of $e^{i\theta}\mathbf A_0-\kappa \mathbf I$ are
$e^{i\theta}\mu-\kappa$, with $\mu\in\Lambda$.
Applying the theorem and restoring the factor $e^{\kappa r}$ proves
\eqref{eq:radial-system-growth}.
\end{proof}

For the limiting companion matrix $\mathbf A_0$, $\Lambda$ is the root set of $p$.
Lemma~\ref{lem:indicator-limit} identifies the vector growth rate with
the scalar indicator of $v$ and gives an $L^1$ limit for its normalized
logarithm. Proposition~\ref{prop:sectorial-counting} then converts this
limit into the sectorial zero-counting formula.
For the classical exponential factors at an irregular singularity,
see~\cite[\S2.7(ii)]{DLMF}.

For the roots of $p_\omega$, put
\[
 K_\omega=\conv\{\overline\mu:p_\omega(\mu)=0\}.
\]
Let $P_\omega$ denote the perimeter of $K_\omega$, taken to be twice
its length when $K_\omega$ is a segment, and set
\begin{equation}\label{eq:global-perimeter-constant}
 C_\omega=\frac{P_\omega}{2\pi\sqrt{2\pi}}.
\end{equation}
Since $p_\omega$ has no $\mu^{n-1}$ term and has a nonzero lower coefficient, its
roots cannot all coincide. Thus $P_\omega>0$ and $C_\omega>0$.
 Perimeters also occur in the classical
zero-counting formulas for exponential polynomials
\cite[Theorem~2.1 and (2.6)]{HW}.

\subsection{The inverse map}

For fixed $\omega$ and $J$, the inverse branches of $\varphi$
map $S_J$ to asymptotically parallel strips at different transverse
positions: along each fixed branch,
$\Ima(\omega z)=2\pi\ell+\arg s+o(1)$ as $|s|\to\infty$.
Thus, as $r$ increases, the disk $|z|\le r$ meets strips with
increasingly large $|\ell|$.
The constants in estimates for one branch may depend on the strip's
transverse position. To combine the strip counts, the next lemma
constructs conformal inverse branches on entire exterior sectors
and estimates the inverse maps and transformed coefficients
uniformly in $\ell$.

\begin{lemma}\label{lem:phase-equation}\label{lem:moving-inverse}
Fix an argument interval $J$ of width less than $\pi$.
There are positive constants $L,C,N,\sigma,\delta_0$, independent of
$\ell\in\mathbb Z$, such that $\Phi$ has a conformal inverse branch
$\psi_\ell$ on
\[
 |s|>\exp(L(1+|\ell|)),\qquad \arg s\in J.
\]
Each $\psi_\ell$ takes values in $\{t:\Rea t>0\}$ and satisfies
$\Phi(\psi_\ell(s))=s$.
Put
\begin{equation}\label{eq:inverse-logarithmic-variables}
 X=\log|s|,\qquad W_\ell=\log_Js+2\pi i\ell.
\end{equation}
The coefficients $b_{k,\ell}$ of~\eqref{eq:phase-equation}
are analytic and satisfy
\begin{equation}\label{eq:moving-coefficient-limit}
 |b_{k,\ell}(s)-c_k|
 \le C\{X^{-\sigma}+X^Ne^{-\delta_0X}+|s|^{-1}\},\qquad 0\le k<n.
\end{equation}
For each fixed $\ell$, putting $L_s=\log_Js$ gives
\begin{equation}\label{eq:inverse-proof}
 \psi_\ell(s)=L_s-\beta\log L_s+2\pi i\ell
       +O\left(\frac{\log|L_s|+1}{|L_s|}\right),
\end{equation}
and $b_{k,\ell}(s)\to c_k$ uniformly on closed subsectors.
\end{lemma}

The proof is given in \hyperref[subsec:proof-moving-inverse]{Section~\ref*{subsec:proof-moving-inverse}}.

\begin{figure}[!b]
\centering
\resizebox{0.96\linewidth}{!}{%
\input{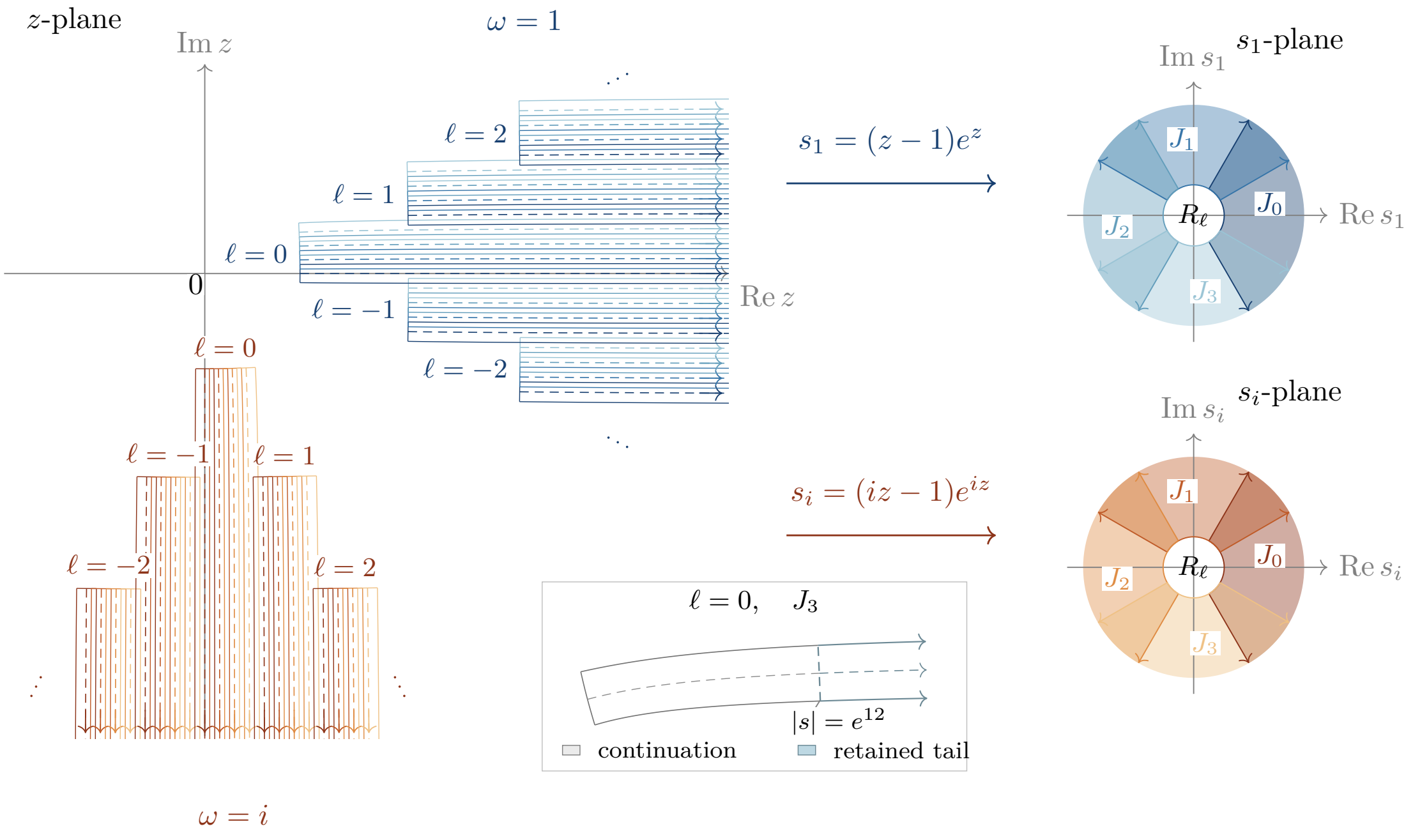}%
}
\caption{Curvilinear pullbacks for the four intervals $J_j$.}
\label{fig:multiple-directions}
\end{figure}

For each fixed branch, increasing the inner radius and slightly
enlarging the argument interval gives an analytic extension across
every finite boundary point of $S_J$. Its image is unbounded,
with bounded imaginary part in the $t$-plane; hence the corresponding
domain $U$ in the $z$-plane is a curvilinear half-strip.

To cover every argument of $s$, we may use the four intervals
\[
 J_j=\left(\frac{j\pi}{2}-\frac{\pi}{3},
           \frac{j\pi}{2}+\frac{\pi}{3}\right),\qquad j=0,1,2,3.
\]
Each has width $2\pi/3<\pi$, and together they cover all arguments
modulo $2\pi$. Lemma~\ref{lem:moving-inverse} applies separately
to each interval; for fixed $\omega$, all the resulting strips
still extend in the direction $-\arg\omega$.
The inner radii supplied by the lemma grow as $\exp(L(1+|\ell|))$,
so the starting points in the $t$-plane advance approximately
linearly with $|\ell|$.

For example, the equation $f''+z^2(e^{2z}+e^{2iz})f=0$ on $\C$ has
$\Omega_a=\{1,i\}$, $a=1$ and $\beta(1)=\beta(i)=1$.
Taking $C=0$ in~\eqref{eq:phase} gives $\Phi(t)=(t-1)e^t$, so
\[
 s_1=(z-1)e^z,\qquad s_i=(iz-1)e^{iz},\qquad z\in\C.
\]

For each $J_j$, we use the inverse branches on $|s|>R_\ell$,
where $R_\ell=\exp(L(1+|\ell|))$ with a common sufficiently large $L$.
For $\Rea t>1$, the side boundaries satisfy
$\Ima t+\arg(t-1)=2\pi\ell+j\pi/2\pm\pi/3$, $\ell\in\Z$.
The two groups of strip tails extend to the right and downwards,
respectively. Figure~\ref{fig:multiple-directions} distinguishes
$J_0,\ldots,J_3$ by four cool shades for $\omega=1$ and four warm
shades for $\omega=i$.

\FloatBarrier

\subsection{Uniform growth of the fundamental solutions}\label{subsec:uniform-row}

On each inverse branch supplied by Lemma~\ref{lem:moving-inverse},
the limiting polynomial $p_\omega$ in~\eqref{eq:characteristic}
is common to all solutions, but different linear combinations may
have different zero counts.
For the introductory model~\eqref{eq:model-original},
we have $\omega=1$, $a=1$, $\beta=0$, and
$p(\mu)=\mu^2-1/4$. With $s=e^z$, the exact transformed equation
on a fixed exterior sector $S_J$ is
\[
 v''+\frac1s v'-\left(\frac14+\frac1{4s^2}\right)v=0.
\]
Its solutions on $S_J$ are
$v(s)=s^{-1/2}(c_+e^{s/2}+c_-e^{-s/2})$.
When both constants are nonzero, the two exponential terms can cancel.
When one constant vanishes, the remaining term is zero-free.

We therefore first estimate the norm of a row of fundamental
solutions on the branches controlled by Lemma~\ref{lem:moving-inverse}.
Fix $z_0\in\C$ and let $\mathbf F=(f_1,\ldots,f_n)$ be the row vector of
fundamental solutions whose initial vectors at $z_0$ are the
standard basis of $\C^n$, so that
\begin{equation}\label{eq:fundamental-initial-matrix}
 \bigl(f_j^{(k-1)}(z_0)\bigr)_{k,j=1}^{n}=\mathbf I_n.
\end{equation}
For each fixed $z\in\C$,
\begin{equation}\label{eq:fundamental-row-norm}
 \|\mathbf F(z)\|=\left(\sum_{j=1}^n|f_j(z)|^2\right)^{1/2}
 =\max_{\|\mathbf c\|=1}\left|\sum_{j=1}^n c_jf_j(z)\right|.
\end{equation}
Indeed, the Cauchy--Schwarz inequality gives
\[
 |f_{\mathbf c}(z)|=|\mathbf F(z)\mathbf c|\le\|\mathbf F(z)\|,\qquad \|\mathbf c\|=1.
\]
Since the fundamental matrix is invertible, $\mathbf F(z)\ne0$, and
equality is attained by the unit vector with
$c_j=\overline{f_j(z)}/\|\mathbf F(z)\|$, $1\le j\le n$.
The maximizing initial vector may depend on $z$.
To obtain zero counts for a fixed solution $f_{\mathbf c}$, we must also
control cancellation in $\mathbf F(z)\mathbf c$.

We next estimate the growth of the fundamental solutions uniformly
across different inverse branches. We first compare their values
and derivatives on each branch with the fixed initial values at
$z_0$.
Fix a maximal weighted frequency $\omega\in\Omega_a$ and an
argument interval $J$ of width less than $\pi$ in the $s$-plane.
On the domain of $\psi_\ell$ in Lemma~\ref{lem:moving-inverse},
define the row vector
\begin{equation}\label{eq:branch-fundamental-row}
 \mathbf F_\ell(s)=\mathbf F(\psi_\ell(s)/\omega)
 =\bigl(f_1(\psi_\ell(s)/\omega),\ldots,
        f_n(\psi_\ell(s)/\omega)\bigr).
\end{equation}
Its values and derivatives of orders $1,\ldots,n-1$ form the matrix
\begin{equation}\label{eq:branch-fundamental-matrix}
 \mathbf X_\ell(s)=\begin{pmatrix}
 f_1(\psi_\ell(s)/\omega)&\cdots&f_n(\psi_\ell(s)/\omega)\\[2pt]
 \dfrac{\mathrm d}{\mathrm ds}f_1(\psi_\ell(s)/\omega)
 &\cdots&\dfrac{\mathrm d}{\mathrm ds}f_n(\psi_\ell(s)/\omega)\\
 \vdots&&\vdots\\
 \dfrac{\mathrm d^{n-1}}{\mathrm ds^{n-1}}f_1(\psi_\ell(s)/\omega)
 &\cdots&\dfrac{\mathrm d^{n-1}}{\mathrm ds^{n-1}}f_n(\psi_\ell(s)/\omega)
 \end{pmatrix}.
\end{equation}

On this branch, write $b_{k,\ell}(s)$ for the coefficients of
\eqref{eq:phase-equation}. The matrices $\mathbf A_\ell(s)$ and $\mathbf A_0$ are
obtained from~\eqref{eq:phase-companion-system} by replacing
$b_k(s)$ with $b_{k,\ell}(s)$ and $c_k$, respectively, where the
constants $c_k$ are defined in~\eqref{eq:limiting-coefficients}.
Thus $\mathbf X_\ell'=\mathbf A_\ell \mathbf X_\ell$.

\begin{lemma}\label{lem:branch-matrix-bounds}
Fix a closed interval $J'\subset J$. For every $\varepsilon>0$,
there are constants $L_\varepsilon,C_\varepsilon>0$, independent of
$\ell$ and $\theta$, such that, with
$\rho_{\ell,\varepsilon}=\exp(L_\varepsilon(1+|\ell|))$,
\begin{equation}\label{eq:uniform-companion-closeness}
 \|\mathbf A_\ell(s)-\mathbf A_0\|<\varepsilon,\qquad
 |s|\ge\rho_{\ell,\varepsilon},\quad \arg s\in J',
\end{equation}
and
\begin{equation}\label{eq:global-entrance}
 \log^+\bigl\|\mathbf X_\ell(\rho_{\ell,\varepsilon}e^{i\theta})^{\pm1}\bigr\|
 \le C_\varepsilon e^{C_\varepsilon(1+|\ell|)},
\end{equation}
for every $\ell\in\Z$ and $\theta\in J'$.
\end{lemma}

The proof is given in \hyperref[subsec:proof-branch-matrix-bounds]{Section~\ref*{subsec:proof-branch-matrix-bounds}}.

Lemma~\ref{lem:branch-matrix-bounds} provides the radii
$\rho_{\ell,\varepsilon}$ and error bounds for comparison with the
limiting constant system.
We use $\zeta=s/T=\Phi(t)/T$, with $t=\omega z$ and $T>0$;
the compact set $K$ below is fixed in the $\zeta$-plane.
Recall that $\mathbf F_\ell$ in~\eqref{eq:branch-fundamental-row} is the
row vector of fundamental solutions evaluated on branch $\ell$.
The next lemma determines the leading growth of its norm from
$p_\omega$ in~\eqref{eq:characteristic}, uniformly over a large number of inverse branches. 
\begin{lemma}\label{lem:uniform-fundamental-row}
For $\mathbf F_\ell$ defined in~\eqref{eq:branch-fundamental-row}, put
\begin{equation}\label{eq:fundamental-row-support-function}
 H_\omega(\zeta)=\max_{p_\omega(\mu)=0}\Rea(\mu\zeta),
 \qquad h_\omega(\theta)=H_\omega(e^{i\theta}).
\end{equation}
For every nonnegative function $m(T)=o(\log T)$ and compact set
$K\subset\{\zeta\ne0:\arg\zeta\in J\}$,
\begin{equation}\label{eq:uniform-fundamental-potential}
 \sup_{|\ell|\le m(T)}
 \left\|T^{-1}\log\|\mathbf F_\ell(T\zeta)\|-H_\omega(\zeta)
 \right\|_{L^1(K)}\longrightarrow0
 \qquad(T\to\infty).
\end{equation}
For these branches, $T^{-1}\log\|\mathbf F_\ell(T\zeta)\|$ is uniformly
bounded above and below on $K$ for all sufficiently large $T$.
The convergence is uniform on $K$ if $K$ avoids the rays
$\Rea((\mu-\nu)\zeta)=0$ for distinct roots $\mu,\nu$ of $p_\omega$.
\end{lemma}

The proof is given in \hyperref[subsec:proof-uniform-fundamental-row]{Section~\ref*{subsec:proof-uniform-fundamental-row}}.

\begin{remark}
Let $\mu_1,\ldots,\mu_d$ be the distinct roots of $p_\omega$.
If a compact set $K$ avoids the rays
$\Rea((\mu_j-\mu_k)\zeta)=0$, $j<k$, in the $\zeta$-plane, then
\[
 \delta_K:=\min_{\zeta\in K}\min_{j<k}
 \bigl|\Rea((\mu_j-\mu_k)\zeta)\bigr|>0.
\]
Thus the largest of the values $\Rea(\mu_j\zeta)$ exceeds all the
others by at least $\delta_K$. This separation is used for the
uniform lower growth estimate.

The $L^1$ conclusion of Lemma~\ref{lem:uniform-fundamental-row} permits
$K$ to meet these rays. Remove small angular neighbourhoods and
use uniform convergence on the remaining compact set.
The uniform boundedness stated in the lemma, together with the
continuity of $H_\omega$, makes the integral of the absolute error
over the removed neighbourhoods tend to zero with their area.
\end{remark}

Section~\ref{sec:sectorial} next uses the fixed-branch estimates to
count zeros in one strip. The uniform estimate
\eqref{eq:uniform-fundamental-potential} will be used in
Section~\ref{sec:global}, where the number of branches grows with
the radius.

\newpage
\section{Sectorial counting and solution spaces}\label{sec:sectorial}

We first obtain a zero-counting formula for $v(s)$ in a sector
of the $s$-plane, then transfer this formula to $f(z)$ in a
curvilinear strip in the $z$-plane, and identify the exceptional
initial values relevant to positive strip density.
Section~\ref{subsec:nonmaximal-frequencies} extends the strip counts
to certain weighted frequencies of smaller modulus.

\subsection{Sectorial counting, transfer and localization}

Let $R_*>0$ and let $S=\{s:|s|>R_*,\ \arg s\in J_0\}$ be an
exterior sector in the $s$-plane, with a continuous argument and
$|J_0|<2\pi$. Consider on $S$ the equation
\begingroup
\predisplaypenalty=-500\relax
\begin{equation}\label{eq:model-equation}
 v^{(n)}+\sum_{j=0}^{n-1}b_j(s)v^{(j)}=0,
 \qquad b_j(s)\longrightarrow c_j
\end{equation}
\endgroup
where the $b_j$ are analytic in $S$ and the convergence is as
$|s|\to\infty$, uniformly on each closed subsector.
As in~\eqref{eq:characteristic}, let
$p(\mu)=\mu^n+\sum_{j=0}^{n-1}c_j\mu^j$ be the characteristic
polynomial of the limiting equation, and put
$\Lambda=\{\mu:p(\mu)=0\}$. For a nonzero solution $v$, define its indicator
\cite[Chapter~I, \S15, p.~51]{Levin} by
\begin{equation}\label{eq:solution-indicator}
 h_v(\theta)=\limsup_{R\to\infty}
                    \frac{\log|v(Re^{i\theta})|}{R},\qquad \theta\in J_0.
\end{equation}
The estimates below show that $v$ grows at most exponentially
on closed subsectors, explaining the normalization by $R^1=R$.

For a finite nonempty set $E\subset\C$, let
$K=\conv\{\overline\mu:\mu\in E\}$.
Its supporting face with outward normal $e^{i\theta}$ consists
of the points $\nu\in K$ maximizing $\Rea(\nu e^{-i\theta})$.
Let $d_K(\theta)$ be the Euclidean length of this face, with
$d_K(\theta)=0$ if it is a point. For a segment, its full length
occurs at each of the two perpendicular normal directions.
These conventions give
$\log|e^{\mu s}|=|s|\Rea(\overline\mu e^{-i\arg s})$.

For example, the exponential sum $g(s)=e^s-e^{-s}$ has simple zeros
\[
 s_k=\pi ik,\quad k\in\Z,\qquad |s_{k+1}-s_k|=\pi.
\]
The segment joining the conjugate frequencies is $[-1,1]$,
of length $2$, and its normal rays are the positive and negative
imaginary axes. The classical connection between frequency hulls
and zero counts in their normal directions is reviewed by
Heittokangas et al.~\cite[\S2.2 and Theorem~2.5]{HITW}.

Let $R>0$ be a scaling parameter tending to infinity and set
$\zeta=s/R$; the disk $|s|\le R$ then becomes $|\zeta|\le1$.
In this example, the points $s_k/R$ are $\pi/R$ apart, so a
fixed segment of length $L>0$ on the imaginary axis contains
$RL/\pi+O(1)$ points. With weight $1/R$ at each point,
the total mass is $L/\pi+O(R^{-1})$, giving limiting density
$1/\pi$ per unit length.

\begin{figure}[htbp]
\centering
\begin{tikzpicture}[x=0.85cm,y=0.85cm,font=\small,
  vertex/.style={circle,fill=black,inner sep=1.4pt},
  axis/.style={->,thin}]
  \draw[axis] (-0.2,0) -- (5.5,0) node[right] {$\Rea\nu$};
  \draw[axis] (0,-0.2) -- (0,3.8) node[above] {$\Ima\nu$};
  \filldraw[fill=gray!12,thin]
    (0.7,0.6) -- (3.7,1) -- (2.7,2.732)
    -- (1,2.9) -- (0.4,1.8) -- cycle;
  \draw[densely dashed,thin] (4.25,0.047) -- (2.2,3.598);
  \draw[line width=1pt] (3.7,1) -- (2.7,2.732);
  \foreach \point in {(0.7,0.6),(3.7,1),(2.7,2.732),(1,2.9),(0.4,1.8)}
    \node[vertex] at \point {};
  \node at (1.2,1.9) {$K$};
  \draw[thin] (3.7,1) -- (3.31,0.775);
  \draw[thin] (2.7,2.732) -- (2.31,2.507);
  \draw[<->,thin] (3.31,0.775) -- (2.31,2.507)
    node[midway,left=5pt] {$d_K(\theta)$};
  \draw[->,thin] (3.2,1.866) -- (4.5,2.617)
    node[above right] {$e^{i\theta}$};
  \node at (2.5,4.1) {$\nu$-plane};
  \begin{scope}[shift={(7.6,0)}]
    \node at (2.5,4.1) {$\zeta$-plane};
    \draw[axis] (-0.2,0) -- (5.2,0) node[right] {$\Rea\zeta$};
    \draw[axis] (0,-0.2) -- (0,3.2) node[above] {$\Ima\zeta$};
    \node[below left] at (0,0) {$0$};
    \draw[->,line width=1pt] (0,0) -- (30:4.6);
    \draw[thin] (0:0.8) arc[start angle=0,end angle=30,radius=0.8];
    \node at (15:1.12) {$\theta$};
    \node[above left=3pt] at (30:3.0) {$\dfrac{d_K(\theta)}{2\pi}$};
    \node[below right=2pt] at (30:4.4) {$\arg\zeta=\theta$};
  \end{scope}
\end{tikzpicture}
\caption{A supporting face and its normal ray.}
\label{fig:supporting-face}
\end{figure}
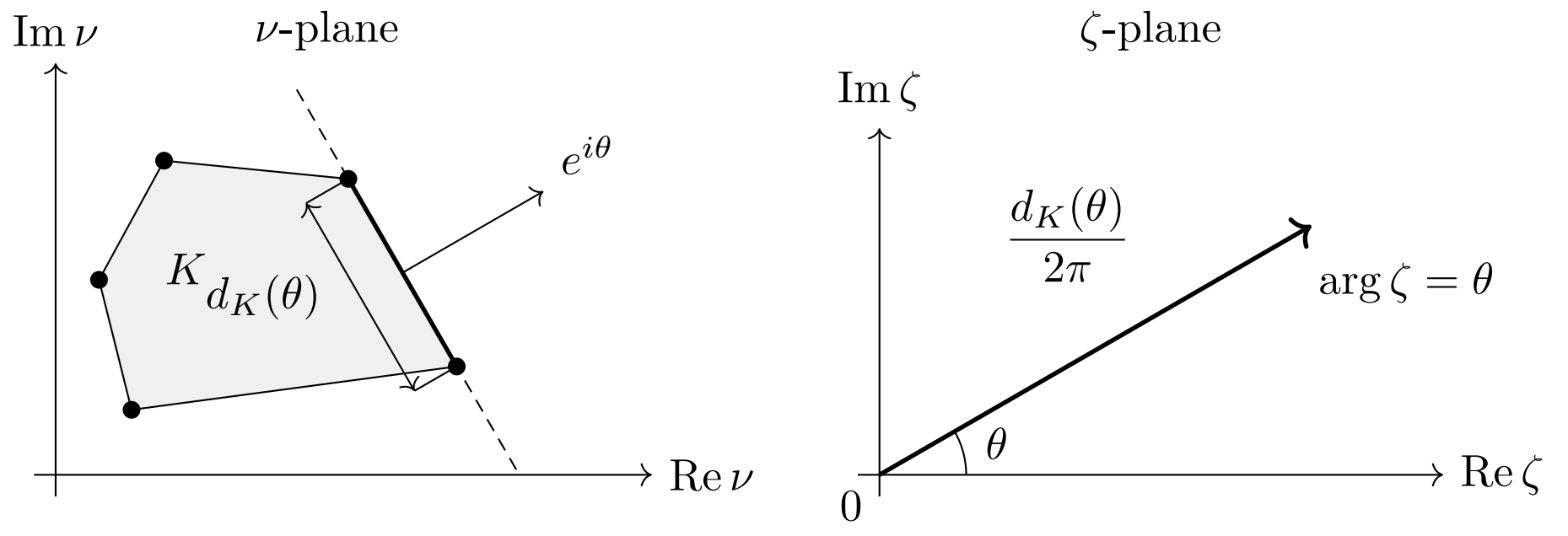

The proposition below
associates such a hull with each nonzero solution
of~\eqref{eq:model-equation}.
After rescaling the zeros by $\zeta=s/R$ and normalizing their
counting measure by $1/R$, an edge of length
$d_K(\theta)$ gives a limiting zero measure with density
$d_K(\theta)/(2\pi)$ along its normal ray; see
Proposition~\ref{prop:sectorial-counting} and
Figure~\ref{fig:supporting-face}.

Under the scaling $\zeta=s/R$, the domain of $v(R\zeta)$ is
$\{\zeta:|\zeta|>R_*/R,\ \arg\zeta\in J_0\}$.
For the limit, we use the fixed sectors in the $\zeta$-plane
\begin{equation}\label{eq:scaled-sectors}
 \mathcal S=\{\zeta\ne0:\arg\zeta\in J_0\},\qquad
 \mathcal S_J=\{\zeta\ne0:\arg\zeta\in J\},
\end{equation}
where $J\subset J_0$ is an open interval. Since $R_*/R\to0$,
$v(R\zeta)$ is defined on every compact subset of $\mathcal S$
for all sufficiently large $R$. Put
\begin{equation}\label{eq:homogeneous-indicator}
 H_v(\zeta)=|\zeta|h_v(\arg\zeta),\qquad \zeta\in\mathcal S.
\end{equation}

For the classical relation between indicators and angular zero
densities, see Levin~\cite[Chapter~III, \S3, Theorem~3]{Levin}.
We derive this zero-counting formula by combining
Lemma~\ref{lem:radial-system-growth} with convergence of
subharmonic functions.

\begin{samepage}
\begin{proposition}\label{prop:sectorial-counting}
Let $v\not\equiv0$ satisfy~\eqref{eq:model-equation}.
For every open interval $J\Subset J_0$ with $|J|<\pi$, there
is a nonempty set $\Lambda_v(J)\subseteq\Lambda$ such that
\begin{equation}\label{eq:support-limit}
 \frac1R\log|v(R\zeta)|\longrightarrow
 H_v(\zeta)=\max_{\mu\in\Lambda_v(J)}\Rea(\mu\zeta)
 \quad\text{in }L^1_{\mathrm{loc}}(\mathcal S_J).
\end{equation}
Put $K_v=\conv\{\overline\mu:\mu\in\Lambda_v(J)\}$, and let
$d_{K_v}(\theta)$ be the Euclidean length of the supporting face
of $K_v$ with outward normal $e^{i\theta}$.
For $R_0>R_*$ and $\alpha_1<\alpha_2$ in $J$ with
$d_{K_v}(\alpha_1)=d_{K_v}(\alpha_2)=0$,
\begin{equation}\label{eq:angular-count}
 n(\{R_0<|s|\le R,\ \alpha_1<\arg s<\alpha_2\},1/v)
 =\frac R{2\pi}\sum_{\alpha_1<\theta<\alpha_2}d_{K_v}(\theta)+o(R).
\end{equation}
The set $\Lambda_v(J)$ consists of those roots $\mu$ for which
$h_v(\theta)=\Rea(\mu e^{i\theta})$ on some nonempty open
subinterval of $J$. 
\end{proposition}
\end{samepage}

The proof is given \hyperref[proof:sectorial-counting]{after
Lemmas~\ref*{lem:indicator-limit} and~\ref*{lem:local-root-selection}}.

The equation determines the full root set $\Lambda$, whereas
$\Lambda_v(J)$ may depend on the solution $v$. Consequently,
different solutions may have different hulls $K_v$ and different
zero-counting asymptotics.
If the selected hull is a single point, the leading coefficient
in~\eqref{eq:angular-count} is zero; this does not imply that the
solution has only finitely many zeros in the sector.

Lemma~\ref{lem:radial-system-growth} gives the radial growth rate
of the vector $\mathbf Y=(v,v',\ldots,v^{(n-1)})^{\mathsf T}$
associated with the solution $v$. Lemma~\ref{lem:indicator-limit} shows
that the scalar indicator $h_v$ equals this rate and proves the
$L^1_{\mathrm{loc}}$ convergence of $R^{-1}\log|v(R\zeta)|$.
Lemma~\ref{lem:local-root-selection} then expresses the limit using
the selected roots $\Lambda_v(J)$. The limiting zero measure is \((2\pi)^{-1}\Delta_\zeta H_v\), with the Laplacian understood in the sense of distributions.

\begin{lemma}\label{lem:indicator-limit}
Let $v\not\equiv0$ satisfy~\eqref{eq:model-equation} on $S$.
The indicator $h_v$ in~\eqref{eq:solution-indicator} is finite and
continuous, and
\begin{equation}\label{eq:indicator-root-values}
 h_v(\theta)\in
 \{\Rea(\mu e^{i\theta}):\mu\in\Lambda\},\qquad \theta\in J_0.
\end{equation}
The function $H_v$ defined in~\eqref{eq:homogeneous-indicator} is
continuous and subharmonic on $\mathcal S$, and
\begin{equation}\label{eq:scalar-limit}
 \frac1R\log|v(R\zeta)|\longrightarrow H_v(\zeta)
 \quad\text{in }L^1_{\mathrm{loc}}(\mathcal S)
 \quad\text{as }R\to\infty.
\end{equation}
\end{lemma}

\begin{proof}
Set
\begin{equation}\label{eq:indicator-jet-system}
 \mathbf Y(s)=(v(s),v'(s),\ldots,v^{(n-1)}(s))^{\mathsf T},
 \qquad \mathbf Y'(s)=\mathbf A(s)\mathbf Y(s),\qquad s\in S.
\end{equation}
The matrix $\mathbf A(s)$ is the companion matrix
in~\eqref{eq:phase-companion-system}, with the coefficients of
\eqref{eq:model-equation}. Replacing $b_j(s)$ by $c_j$ gives its
limit $\mathbf A_0$, with $\det(\mu\mathbf I-\mathbf A_0)=p(\mu)$
by~\eqref{eq:companion-characteristic}.
The vector $\mathbf Y$ never vanishes: a zero value would give zero initial
values for the linear system, and uniqueness would imply
$v\equiv0$.

Fix a closed interval $J'\Subset J_0$. For sufficiently large
$r_0>R_*$, the matrices $\mathbf A(re^{i\theta})$ are uniformly
bounded for $r\ge r_0$, $\theta\in J'$. The values
$\mathbf Y(r_0e^{i\theta})$ are uniformly bounded as well, since
$\{r_0e^{i\theta}:\theta\in J'\}$ is a compact subset of $S$.
The radial Gronwall estimate used in the proof of
Lemma~\ref{lem:radial-system-growth} therefore gives
$\|\mathbf Y(re^{i\theta})\|\le Ce^{Cr}$ for $r\ge r_0$ and
$\theta\in J'$, with $C$ independent of $r$ and $\theta$.
Hence $\log^+|v(re^{i\theta})|=O(r)$ uniformly for $\theta\in J'$.
Lemma~\ref{lem:radial-system-growth} gives
\begin{equation}\label{eq:jet-ray-limit}
 \ell_v(\theta):=
 \lim_{r\to\infty}\frac{\log\|\mathbf Y(re^{i\theta})\|}{r}
 \in\{\Rea(\mu e^{i\theta}):\mu\in\Lambda\},
 \qquad \theta\in J_0.
\end{equation}

Fix $\theta\in J_0$. We next show that
$h_v(\theta)=\ell_v(\theta)$, with $h_v$ defined in~\eqref{eq:solution-indicator}.
Let $\mathbf X_r(\tau)$ be the solution of
\begin{equation}\label{eq:indicator-short-ray-system}
 \mathbf X_r'(\tau)=e^{i\theta}\mathbf A((r+\tau)e^{i\theta})\mathbf X_r(\tau),
 \qquad \mathbf X_r(0)=\mathbf I_n,\qquad 0\le\tau\le1.
\end{equation}
Uniqueness gives
$\mathbf Y((r+\tau)e^{i\theta})=\mathbf X_r(\tau)\mathbf Y(re^{i\theta})$.
Put $\mathbf X_0(\tau)=\exp(\tau e^{i\theta}\mathbf A_0)$, so that
$\mathbf X_0'=e^{i\theta}\mathbf A_0\mathbf X_0$ and
$\mathbf X_0(0)=\mathbf I_n$.
Integrating~\eqref{eq:indicator-short-ray-system} and the equation
for $\mathbf X_0$ gives
\begin{equation}\label{eq:indicator-integral-equations}
 \begin{aligned}
 \mathbf X_r(\tau)
 &=\mathbf I_n+\int_0^\tau e^{i\theta}
   \mathbf A((r+t)e^{i\theta})\mathbf X_r(t)\dd t,\\
 \mathbf X_0(\tau)
 &=\mathbf I_n+\int_0^\tau e^{i\theta}
   \mathbf A_0\mathbf X_0(t)\dd t.
 \end{aligned}
\end{equation}
Set $\mathbf D_r=\mathbf X_r-\mathbf X_0$. The difference of the
integrands can be decomposed using
\[
 \begin{aligned}
 \mathbf A((r+t)e^{i\theta})\mathbf X_r(t)
       -\mathbf A_0\mathbf X_0(t)=\mathbf A((r+t)e^{i\theta})\mathbf D_r(t)
       +\bigl(\mathbf A((r+t)e^{i\theta})-\mathbf A_0\bigr)
          \mathbf X_0(t).
 \end{aligned}
\]
Define $\varepsilon_r=\sup_{0\le t\le1}
\|\mathbf A((r+t)e^{i\theta})-\mathbf A_0\|$.
Coefficient convergence gives $\varepsilon_r\to0$ and a constant
$M>0$ such that $\|\mathbf A((r+t)e^{i\theta})\|\le M$ for all
sufficiently large $r$ and $0\le t\le1$.
Also put $K=\max_{0\le t\le1}\|\mathbf X_0(t)\|<\infty$;
both $M$ and $K$ are independent of $r$.
Subtracting~\eqref{eq:indicator-integral-equations}, taking norms,
and using $|e^{i\theta}|=1$ and $\tau\le1$, we obtain
\[
 \|\mathbf D_r(\tau)\|
 \le M\int_0^\tau\|\mathbf D_r(t)\|\dd t+K\varepsilon_r,
 \qquad 0\le\tau\le1.
\]
Gronwall's inequality yields
$\sup_{0\le\tau\le1}\|\mathbf D_r(\tau)\|
\le Ke^M\varepsilon_r\to0$, and hence
\begin{equation}\label{eq:indicator-short-ray-limit}
 \mathbf X_r(\tau)\longrightarrow\mathbf X_0(\tau)
 \quad\text{uniformly for }0\le\tau\le1.
\end{equation}

Let $\mathbf e_1,\ldots,\mathbf e_n$ be the standard coordinate vectors,
and write $\mathbf w=(w_1,\ldots,w_n)^{\mathsf T}$.
We claim that a constant $a_\theta>0$ satisfies
\begin{equation}\label{eq:companion-observation-bound}
m_\theta(\mathbf{w})= \max_{0\le\tau\le1}|\mathbf e_1^{\mathsf T}\mathbf X_0(\tau)\mathbf w|
 \ge a_\theta\|\mathbf w\|,\qquad \mathbf w\in\C^n.
\end{equation}
First restrict $\mathbf w$ to the unit sphere
$\{\mathbf w\in\C^n:\|\mathbf w\|=1\}$.
If the maximum were zero, all derivatives at $\tau=0$ of the scalar
function $\mathbf e_1^{\mathsf T}\mathbf X_0(\tau)\mathbf w$ would vanish.
The companion form gives
$\mathbf e_1^{\mathsf T}\mathbf A_0^j=\mathbf e_{j+1}^{\mathsf T}$, so
\[
 \left.\frac{\mathrm d^j}{\mathrm d\tau^j}
       \bigl(\mathbf e_1^{\mathsf T}\mathbf X_0(\tau)\mathbf w\bigr)\right|_{\tau=0}
 =e^{ij\theta}w_{j+1},\qquad 0\le j<n.
\]
All coordinates of $\mathbf w$ would then be zero, contradicting
$\|\mathbf w\|=1$. The triangle inequality and
$\|\mathbf X_0(\tau)\|\le K$ give, for $\mathbf w,\mathbf z\in\C^n$,
\[
 |m_\theta(\mathbf w)-m_\theta(\mathbf z)|
 \le\max_{0\le\tau\le1}
       |\mathbf e_1^{\mathsf T}\mathbf X_0(\tau)(\mathbf w-\mathbf z)|
 \le K\|\mathbf w-\mathbf z\|.
\]
Thus $m_\theta$ is a continuous, strictly positive function on the compact
unit sphere $\|\mathbf w\|=1$ and attains a positive minimum $a_\theta$ there. Namely,
\[
a_\theta
:=\min_{\|\mathbf w\|=1}m_\theta(\mathbf w)>0.
\]
For $\mathbf w\ne0$, applying this bound to
$\mathbf w/\|\mathbf w\|$ and using homogeneity gives
\eqref{eq:companion-observation-bound} for arbitrary $\mathbf w$;
the case $\mathbf w=0$ is immediate.

By~\eqref{eq:indicator-short-ray-limit}, for all sufficiently large
$r$ we have $\sup_{0\le\tau\le1}
\|\mathbf X_r(\tau)-\mathbf X_0(\tau)\|\le a_\theta/2$.
The triangle inequality therefore gives, for every $\mathbf w\in\C^n$,
\begin{equation}\label{eq:companion-maximum-comparison}
 \begin{aligned}
 \left|\max_{0\le\tau\le1}
       |\mathbf e_1^{\mathsf T}\mathbf X_r(\tau)\mathbf w|
       -m_\theta(\mathbf w)\right|
 &\le\max_{0\le\tau\le1}
       |\mathbf e_1^{\mathsf T}
        (\mathbf X_r(\tau)-\mathbf X_0(\tau))\mathbf w|\\
 &\le\left(\sup_{0\le\tau\le1}
        \|\mathbf X_r(\tau)-\mathbf X_0(\tau)\|\right)\|\mathbf w\|
       \le\frac{a_\theta}{2}\|\mathbf w\|.
 \end{aligned}
\end{equation}
Combining~\eqref{eq:companion-maximum-comparison} with
\eqref{eq:companion-observation-bound}, we obtain
\[
 \max_{0\le\tau\le1}
 |\mathbf e_1^{\mathsf T}\mathbf X_r(\tau)\mathbf w|
 \ge m_\theta(\mathbf w)-\frac{a_\theta}{2}\|\mathbf w\|
 \ge\frac{a_\theta}{2}\|\mathbf w\|.
\]
Recall from~\eqref{eq:indicator-jet-system} that
$v(s)=\mathbf e_1^{\mathsf T}\mathbf Y(s)$.
Uniqueness for~\eqref{eq:indicator-short-ray-system} gives
\[
 \mathbf e_1^{\mathsf T}\mathbf X_r(\tau)\mathbf Y(re^{i\theta})
 =\mathbf e_1^{\mathsf T}\mathbf Y((r+\tau)e^{i\theta})
 =v((r+\tau)e^{i\theta}).
\]
Thus, taking $\mathbf w=\mathbf Y(re^{i\theta})$ yields
\begin{equation}\label{eq:scalar-jet-segment-comparison}
 \max_{0\le\tau\le1}|v((r+\tau)e^{i\theta})|
 \ge \frac{a_\theta}{2}\|\mathbf Y(re^{i\theta})\|.
\end{equation}
Choose a maximizing $\tau=\tau_r$. Dividing the logarithm of
\eqref{eq:scalar-jet-segment-comparison} by $r+\tau_r$ and using
$(r+\tau_r)/r\to1$ shows that
$h_v(\theta)\ge\ell_v(\theta)$.
The inequality $|v|\le\|\mathbf Y\|$ gives the reverse inequality.
Thus $h_v=\ell_v$, and~\eqref{eq:jet-ray-limit} proves
\eqref{eq:indicator-root-values}.

\bigskip

We now compare the scalar and vector logarithms on
$\mathcal S$, the sector of the $\zeta$-plane defined
in~\eqref{eq:scaled-sectors}.
For $R$ large enough on each compact subset of $\mathcal S$, put
\begin{equation}\label{eq:scaled-scalar-and-jet-potentials}
 U_R(\zeta)=\frac1R\log|v(R\zeta)|,\qquad
 V_R(\zeta)=\frac1R\log\|\mathbf Y(R\zeta)\|,
\end{equation}
with $\log0=-\infty$. Recall from~\eqref{eq:homogeneous-indicator}
that $H_v(\rho e^{i\theta})=\rho h_v(\theta)$.
Fix $\zeta=\rho e^{i\theta}\in\mathcal S$, with $\rho>0$.
Since $R\rho\to\infty$,~\eqref{eq:jet-ray-limit} and
$\ell_v(\theta)=h_v(\theta)$ give the pointwise convergence
\[
 \begin{aligned}
 V_R(\zeta)
 &=\rho\,\frac{\log\|\mathbf Y(R\rho e^{i\theta})\|}{R\rho}\\
 &\longrightarrow\rho\ell_v(\theta)
   =\rho h_v(\theta)=H_v(\zeta).
 \end{aligned}
\]
For $|\xi|=1$, define the real directional derivative by
$D_\xi V_R(\zeta)=
\left.\frac{\mathrm d}{\mathrm dh}V_R(\zeta+h\xi)\right|_{h=0}$,
where $h\in\R$. Write $\mathbf Y=(Y_1,\ldots,Y_n)^{\mathsf T}$.
The chain rule gives $D_\xi Y_j(R\zeta)=R\xi Y_j'(R\zeta)$,
with $Y_j'=\mathrm dY_j/\mathrm ds$; the factor $R$ cancels
$1/R$ in~\eqref{eq:scaled-scalar-and-jet-potentials}.
Since $\mathbf Y\ne0$, we can differentiate
$\log\|\mathbf Y\|=\frac12\log\sum_{j=1}^n|Y_j|^2$.
Cauchy--Schwarz and~\eqref{eq:indicator-jet-system} then give
\[
 \begin{aligned}
 D_\xi V_R(\zeta)
 &=\frac{\Rea\!\left(\xi\sum_{j=1}^n
       Y_j'(R\zeta)\overline{Y_j(R\zeta)}\right)}
       {\|\mathbf Y(R\zeta)\|^2},\\
 |D_\xi V_R(\zeta)|
 &\le\frac{\|\mathbf Y'(R\zeta)\|}{\|\mathbf Y(R\zeta)\|}
 \le\|\mathbf A(R\zeta)\|.
 \end{aligned}
\]
The right-hand side is bounded independently of large $R$ on
each compact subset of $\mathcal S$. The functions $V_R$ are
therefore locally equicontinuous. Pointwise convergence, applied
on a finite covering of each compact set, gives
\begin{equation}\label{eq:jet-potential-local-limit}
 V_R\longrightarrow H_v
 \quad\text{locally uniformly on }\mathcal S.
\end{equation}
In particular, $H_v$ is continuous. The logarithm of the norm of
an analytic vector is subharmonic, so $V_R$ and the limit $H_v$
are subharmonic.

It remains to show that $U_R$ has the same limit. We have
$U_R\le V_R$, hence $U_R$ is locally uniformly bounded above.
Fix $\delta>0$ such that
$\{\eta:|\eta-\zeta|\le\delta\}\subset\mathcal S$.
For large $R$, the corresponding disk of radius $R\delta$ around
$R\zeta$ lies in $S$. Cauchy's inequalities give
\[
 |v^{(j)}(R\zeta)|
 \le\frac{j!}{(R\delta)^j}
       \sup_{|\eta-\zeta|\le\delta}|v(R\eta)|,
 \qquad 0\le j<n.
\]
Since $R\ge1$, summing the squares gives a constant $C_\delta$,
independent of $R$, such that
$\|\mathbf Y(R\zeta)\|\le C_\delta
\sup_{|\eta-\zeta|\le\delta}|v(R\eta)|$.
Taking logarithms gives
\begin{equation}\label{eq:cauchy-log-comparison}
 V_R(\zeta)\le
 \sup_{|\eta-\zeta|\le\delta}U_R(\eta)
       +\frac{\log C_\delta}{R}.
\end{equation}

\phantomsection\label{par:indicator-subharmonic-compactness}
By subharmonic compactness on an exhaustion of $\mathcal S$
\cite[Theorem~2.4.8, p.~28]{DinhSibonyCurrents}, every sequence
$(U_{R_k})$ with $R_k\to\infty$ has a subsequence converging in
$L^1_{\mathrm{loc}}(\mathcal S)$ to a subharmonic function
$\widetilde H\not\equiv-\infty$.
The alternative of locally uniform convergence to $-\infty$ is
excluded by~\eqref{eq:cauchy-log-comparison} and
\eqref{eq:jet-potential-local-limit}.
Relabel this subsequence as $(U_{R_k})$.
To prove $\widetilde H\le H_v$ at every point, fix
$\zeta\in\mathcal S$ and $\delta>0$ with
$\{\eta:|\eta-\zeta|\le\delta\}\subset\mathcal S$.
The submean inequality for $\widetilde H$
\cite[Theorem~2.4.5, p.~28]{DinhSibonyCurrents},
the $L^1_{\mathrm{loc}}$ convergence of $U_{R_k}$,
$U_{R_k}\le V_{R_k}$, and~\eqref{eq:jet-potential-local-limit} give
\begin{equation}\label{eq:scalar-limit-upper-comparison}
 \begin{aligned}
 \widetilde H(\zeta)
 &\le\frac1{\pi\delta^2}
      \int_{|\eta-\zeta|<\delta}\widetilde H(\eta)\,\mathrm dA(\eta)\\
 &=\lim_{k\to\infty}\frac1{\pi\delta^2}
      \int_{|\eta-\zeta|<\delta}U_{R_k}(\eta)\,\mathrm dA(\eta)\\
 &\le\frac1{\pi\delta^2}
      \int_{|\eta-\zeta|<\delta}H_v(\eta)\,\mathrm dA(\eta).
 \end{aligned}
\end{equation}
Letting $\delta\to0$ and using continuity of $H_v$ proves
$\widetilde H\le H_v$ on $\mathcal S$.
For every closed disk $\{\eta:|\eta-\zeta|\le\delta\}$ in
$\mathcal S$, Hartogs' lemma
\cite[Theorem~2.4.9, p.~28]{DinhSibonyCurrents}, together
with~\eqref{eq:cauchy-log-comparison}, gives
\[
 H_v(\zeta)\le
 \limsup_{k\to\infty}\sup_{|\eta-\zeta|\le\delta}U_{R_k}(\eta)
 \le\sup_{|\eta-\zeta|\le\delta}\widetilde H(\eta).
\]
Letting $\delta\to0$ and using upper semicontinuity yields
$H_v\le\widetilde H$. Together with the reverse inequality proved
using~\eqref{eq:scalar-limit-upper-comparison}, this gives
$\widetilde H=H_v$ on $\mathcal S$.

To obtain convergence of the full family, suppose
that~\eqref{eq:scalar-limit} fails. There would then be a compact
set $K\subset\mathcal S$, a number $\varepsilon>0$, and a sequence
$R_k\to\infty$ such that
\[
 \int_K|U_{R_k}-H_v|\,\mathrm dA\ge\varepsilon,
 \qquad k\ge1.
\]
Applying the same subharmonic compactness theorem
\cite[Theorem~2.4.8, p.~28]{DinhSibonyCurrents} again to this
sequence $(U_{R_k})$ gives a subsequence converging in
$L^1_{\mathrm{loc}}(\mathcal S)$. Its limit must equal $H_v$
by the preceding identification. Along this subsequence, the
integrals over $K$ tend to zero, a contradiction.
This proves~\eqref{eq:scalar-limit} for the full family as $R\to\infty$.
\end{proof}

\begin{remark}
The convergence in~\eqref{eq:scalar-limit} means that, for every
compact set $K\subset\mathcal S$,
\[
 \int_K\left|\frac1R\log|v(R\zeta)|-H_v(\zeta)\right|
 \,\mathrm dA(\zeta)\longrightarrow0,
\]
where $\mathrm dA$ denotes planar area measure.
Locally uniform convergence is stronger than
$L^1_{\mathrm{loc}}$ convergence. For the vector logarithms
$V_R(\zeta)=R^{-1}\log\|\mathbf Y(R\zeta)\|$, where
$\mathbf Y=(v,v',\ldots,v^{(n-1)})^{\mathsf T}$, we have
\[
 \int_K|V_R-H_v|\,\mathrm dA
 \le\operatorname{area}(K)\sup_{\zeta\in K}|V_R(\zeta)-H_v(\zeta)|.
\]
Thus locally uniform convergence implies $L^1_{\mathrm{loc}}$
convergence; the converse need not hold.
The scalar logarithm takes the value $-\infty$ at zeros of
$v(R\zeta)$, but its logarithmic singularities are locally integrable.

The growth theorem used in Lemma~\ref{lem:radial-system-growth}
allows nontrivial Jordan blocks. The comparison with the scalar
solution uses only the companion form of $\mathbf A_0$.
Repeated characteristic roots are therefore allowed in
Lemma~\ref{lem:indicator-limit}.
\end{remark}

The limit in Lemma~\ref{lem:indicator-limit} describes the growth
of $v$ through its logarithm.\phantomsection\label{par:distributional-laplacian} We recall the planar potential-theoretic
facts that convert this limit into a zero-counting measure.

For $s=x+iy$, write $\mathrm dA=\mathrm dx\,\mathrm dy$.
The \textit{Laplacian} is $\Delta_s=\partial_x^2+\partial_y^2$.
A $C^2$ real-valued function is \textit{harmonic} if its Laplacian is zero;
in particular, the real part of a holomorphic function is harmonic.
For functions with logarithmic singularities, the Laplacian is
understood in the sense of distributions. Let $D\subset\C$ be a
domain. A function $u$ is \textit{locally integrable}, written
$u\in L^1_{\mathrm{loc}}(D)$, if $\int_K|u|\,\mathrm dA<\infty$
for every compact $K\subset D$.
The space $C_c^\infty(D)$ consists of smooth functions whose
supports are compact subsets of $D$. Its elements are called
\textit{test functions}. For $u\in L^1_{\mathrm{loc}}(D)$, its
\textit{distributional Laplacian} is defined by
\[
 \langle\Delta_s u,\phi\rangle
 :=\int_Du\,\Delta_s\phi\,\mathrm dA,\qquad
 \phi\in C_c^\infty(D).
\]
The brackets denote the action of a distribution on a test function.
The integral is finite because $\Delta_s\phi$ is bounded and
compactly supported. When $u$ is $C^2$, two integrations by parts,
with no boundary terms, give $\int_D(\Delta_su)\phi\,\mathrm dA$;
thus this definition extends the usual Laplacian
\cite[Chapter~I, (2.5)--(2.6), p.~15]{Demailly}.

An upper semicontinuous function $u:D\to[-\infty,\infty)$,
not identically $-\infty$, is \textit{subharmonic} if it satisfies the
disk mean inequality
$u(b)\le(\pi r^2)^{-1}\int_{|s-b|<r}u(s)\,\mathrm dA(s)$
whenever the closed disk lies in $D$.
Such a function is locally integrable, and its distributional
Laplacian is a positive measure \cite[Chapter~I, Corollary~4.19, p.~35]{Demailly}.
In particular, $\langle\Delta_su,\phi\rangle\ge0$ for
nonnegative test functions. This applies to the subharmonic
limit $H_v$ in
Lemma~\ref{lem:indicator-limit}.

The \textit{unit point mass} $\delta_b$ assigns mass $1$ to a set containing
$b$ and mass $0$ otherwise; equivalently,
\begin{equation}\label{eq:point-mass-pairing}
 \langle\delta_b,\phi\rangle=\int_D\phi\,\mathrm d\delta_b=\phi(b).
\end{equation}
Thus $\sum_bm_b\delta_b$ counts points with weights $m_b$.
The logarithmic fundamental solution gives
\[
 \Delta_s\log|s-b|=2\pi\delta_b
 \qquad\text{on }\C
\]
in the distributional sense
\cite[Chapter~I, (3.22) and the following calculation,
pp.~25--26]{Demailly}.
The function $\log|s-b|$ is locally integrable even at $b$, since
$\int_0^r t|\log t|\,\dd t<\infty$, and is harmonic away from $b$.
Consequently, the identity means
$\int_{\C}\log|s-b|\Delta_s\phi\,\mathrm dA=2\pi\phi(b)$,
rather than a pointwise equality at $b$.

For a holomorphic function $v$ on $D$ that is not identically
zero, let $m_b$ be the multiplicity of its zero $b$. The one-variable
form of the Lelong--Poincar\'e formula
\cite[Chapter~III, (2.15), p.~143]{Demailly} is
\begin{equation}\label{eq:planar-zero-measure}
 \frac1{2\pi}\Delta_s\log|v(s)|
 =\sum_{\substack{b\in D\\v(b)=0}}m_b\delta_b
 \qquad\text{on }D.
\end{equation}
For this Laplacian formulation and its proof, see
Hough--Krishnapur--Peres--Vir\'ag
\cite[\S2.4.1, (2.4.3)--(2.4.4), p.~24]{HKPV}.
The zero set is locally finite, so only finitely many terms
contribute for each test function. Equivalently, for every
$\phi\in C_c^\infty(D)$,
\[
 \frac1{2\pi}\int_D\log|v|\,\Delta_s\phi\,\mathrm dA
 =\sum_{v(b)=0}m_b\phi(b).
\]
To see why the coefficients are the zero multiplicities, write
$v(s)=(s-b)^{m_b}g(s)$ near $b$, with $g$ holomorphic and
nonvanishing. Then
\begin{equation}\label{eq:local-log-factorization}
 \log|v(s)|=m_b\log|s-b|+\log|g(s)|.
\end{equation}
The last term is harmonic, being locally the real part of a
holomorphic logarithm of $g$. Thus $\log|v|$ is locally integrable.
Applying the fundamental identity
therefore gives $2\pi m_b\delta_b$ near $b$.
Away from the zeros the Laplacian is zero. Splitting a test
function into a finite sum of smooth functions supported in
these neighborhoods and adding the local identities gives
\eqref{eq:planar-zero-measure} on $D$.

Finally, the distributional Laplacian respects
$L^1_{\mathrm{loc}}$ limits: if $u_k\to u$ in
$L^1_{\mathrm{loc}}(D)$, then, for each $\phi\in C_c^\infty(D)$,
\[
 |\langle\Delta_s(u_k-u),\phi\rangle|
 \le\sup_D|\Delta_s\phi|
       \int_{\operatorname{supp}\phi}|u_k-u|\,\mathrm dA
 \longrightarrow0.
\]
Recall $U_R(\zeta)=R^{-1}\log|v(R\zeta)|$ from
\eqref{eq:scaled-scalar-and-jet-potentials} in the proof of
Lemma~\ref{lem:indicator-limit}.
By~\eqref{eq:scalar-limit}, $U_R\to H_v$ in
$L^1_{\mathrm{loc}}(\mathcal S)$, where $\mathcal S$ is the
scaled sector in~\eqref{eq:scaled-sectors}.
On each fixed relatively compact subdomain of $\mathcal S$, for all
sufficiently large $R$, formula~\eqref{eq:planar-zero-measure}
identifies $(2\pi)^{-1}\Delta_\zeta U_R$ with the measure placing
mass $m_b/R$ at each rescaled zero $b/R$.
Computing $\Delta_\zeta H_v$ will therefore determine the
limiting zero measure.

Formula~\eqref{eq:indicator-root-values} shows that, on each ray,
$h_v(\theta)$ equals $\Rea(\mu e^{i\theta})$ for some characteristic
root $\mu$, which may depend on $\theta$.
To construct the convex hull used in Proposition~\ref{prop:sectorial-counting},
we need a fixed subset of roots whose maximum real projection gives
$h_v$ throughout an argument interval $J$. The next lemma proves
this representation when $|J|<\pi$ and identifies the vertices of
the resulting hull.

\begin{samepage}
\begin{lemma}\label{lem:local-root-selection}
Under the assumptions of Lemma~\ref{lem:indicator-limit}, let
$J\subset J_0$ be an open interval with $|J|<\pi$. Define
$\Lambda_v(J)$ to consist of the roots $\mu\in\Lambda$ for which
$h_v(\theta)=\Rea(\mu e^{i\theta})$ on a nonempty open subinterval
of $J$. Then $\Lambda_v(J)$ is nonempty and
\begin{equation}\label{eq:local-support-indicator}
 h_v(\theta)=
 \max_{\mu\in\Lambda_v(J)}\Rea(\mu e^{i\theta}),
 \qquad \theta\in J.
\end{equation}
The points $\overline\mu$, where $\mu\in\Lambda_v(J)$, are exactly the
vertices of their convex hull.
\end{lemma}
\end{samepage}

\begin{proof}
We first identify the roots $\mu$ for which
$h_v(\theta)=\Rea(\mu e^{i\theta})$ on an open subinterval of $J$.
We then use the subharmonicity of $H_v$ to prove that
$\Rea(\mu\zeta)\le H_v(\zeta)$ for every such root and every
$\zeta\in\mathcal S_J$. These two steps give the maximum
representation in~\eqref{eq:local-support-indicator}.

For distinct roots $\mu,\nu$, the equation
$\Rea((\mu-\nu)e^{i\theta})=0$ has two solutions modulo $2\pi$.
Consequently, only finitely many angles $\theta\in J$ satisfy
$\Rea(\mu e^{i\theta})=\Rea(\nu e^{i\theta})$ for some distinct roots
$\mu,\nu$. Removing these angles from $J$ leaves one or more open
intervals. On each interval, the values
in~\eqref{eq:indicator-root-values} are pairwise distinct. Continuity of $h_v$, supplied by
Lemma~\ref{lem:indicator-limit}, forces $h_v$ to equal one fixed
function $\Rea(\mu e^{i\theta})$ throughout that interval.
Indeed, a switch from $\mu$ to $\nu\ne\mu$ would, by continuity,
require $\Rea(\mu e^{i\theta})=\Rea(\nu e^{i\theta})$ at the
switching angle, but all such angles have been removed.
The set $\Lambda_v(J)$ consists of the roots selected on these
intervals, with each root included only once.
There is at least one such interval, so the set is nonempty.

Recall that $H_v(\rho e^{i\theta})=\rho h_v(\theta)$.
The preceding root selection makes $H_v$ continuous and piecewise
linear on $\mathcal S_J$ in~\eqref{eq:scaled-sectors}.
The sector $\mathcal S_J$ is convex because $|J|<\pi$.

Identify $\C$ with $\mathbb R^2$ and write $\zeta=x+iy$.
Where $H_v$ is differentiable, its gradient is
$\nabla H_v=(\partial_xH_v,\partial_yH_v)$.
On a subsector where $H_v(\zeta)=\Rea(\mu\zeta)
=(\Rea\mu)x-(\Ima\mu)y$, this gives
$\nabla H_v=(\Rea\mu,-\Ima\mu)$, the vector corresponding to
$\overline\mu$.

Fix a \textit{separating ray} $\arg\zeta=\theta_0$, across which
the linear expression for $H_v$ changes. For this local
calculation, rotate it to the positive real axis by
$\zeta=e^{i\theta_0}w$ and set
$\widetilde H(w)=H_v(e^{i\theta_0}w)$ on
$e^{-i\theta_0}\mathcal S_J$. Without loss of generality, write $w=x+iy$ for the rotated
coordinate. Near this ray, write
\begin{equation}\label{eq:root-indicator-linear-pieces}
 \widetilde H(x+iy)=
 \begin{cases}
 A_+x+B_+y,&y>0,\\
 A_-x+B_-y,&y<0,
 \end{cases}
 \qquad A_\pm,B_\pm\in\R.
\end{equation}
Continuity on $y=0$, $x>0$, gives $A_+x=A_-x$, hence
$A_+=A_-$. Thus the two $x$-derivatives agree, and their upward
$y$-derivatives have the constant jump $\kappa=B_+-B_-$.

In the $w$-plane, choose a disk $D$ centred on the positive real
axis, with $\overline D\subset e^{-i\theta_0}\mathcal S_J$ and
small enough to meet no other separating ray.
The chosen ray divides $D$ into upper and lower half-disks,
on which the $y>0$ and $y<0$ expressions in
\eqref{eq:root-indicator-linear-pieces} hold, respectively.
For $\phi\in C_c^\infty(D)$, apply Green's formula
\cite[Chapter~I, (4.4), p.~30]{Demailly} to the two harmonic
linear expressions on the respective half-disks.
The circular-arc terms vanish because $\phi$ and its first
derivatives vanish near $\partial D$.
On the common diameter, the outward normals of the upper and
lower half-disks point downwards and upwards, respectively.
The boundary integrands are therefore
$-A_+x\phi_y+B_+\phi$ and $A_-x\phi_y-B_-\phi$,
with $\phi$ and $\phi_y$ evaluated at $(x,0)$.
Using $A_+=A_-$ and adding the two identities gives
\begin{equation}\label{eq:root-indicator-interface}
 \int_D\widetilde H\,\Delta\phi\,\mathrm dA
 =\int_{D\cap\R}(B_+-B_-)\phi(x,0)\,\dd x
 =\kappa\int_{D\cap\R}\phi(x,0)\,\dd x.
\end{equation}
Rotation preserves subharmonicity. The
\hyperref[par:distributional-laplacian]{positivity of the distributional Laplacian}
therefore gives $\langle\Delta\widetilde H,\phi\rangle\ge0$
for $\phi\ge0$ \cite[Chapter~I, Corollary~4.19, p.~35]{Demailly}.
Choose a nonnegative $\phi$ that is positive near the centre
of $D$. Then $\int_{D\cap\R}\phi(x,0)\,\dd x>0$, so
\eqref{eq:root-indicator-interface} and the preceding
nonnegativity imply $\kappa\ge0$.

Since $A_+=A_-$ and $B_+-B_-=\kappa\ge0$, the difference
between the two linear expressions in
\eqref{eq:root-indicator-linear-pieces} is $\kappa y$.
Thus the upper expression is the larger one for $y>0$ and
the lower expression is the larger one for $y<0$, giving
\begin{equation}\label{eq:root-local-maximum}
 \widetilde H(x+iy)
 =\max\{A_+x+B_+y,\ A_-x+B_-y\},
 \qquad x+iy\in D.
\end{equation}
This local maximum representation will control the left and
right derivatives along a line segment.

For any $\zeta_1,\zeta_2\in\mathcal S_J$, the segment joining
them lies in $\mathcal S_J$ because this sector is convex.
Put
\[
 g(t)=H_v((1-t)\zeta_1+t\zeta_2),\qquad 0\le t\le1.
\]
On each open subsector considered above,
$H_v(\zeta)=\Rea(\nu\zeta)$ for a fixed $\nu\in\Lambda_v(J)$.
Thus the finitely many separating rays give a partition
$0=t_0<t_1<\cdots<t_m=1$ and roots
$\nu_j\in\Lambda_v(J)$, $1\le j\le m$, such that
\begin{equation}\label{eq:root-segment-linear-pieces}
 g(t)=\Rea(\nu_j\zeta_1)
       +t\,\Rea\bigl(\nu_j(\zeta_2-\zeta_1)\bigr),
 \qquad t_{j-1}<t<t_j.
\end{equation}
Continuity of $H_v$ extends~\eqref{eq:root-segment-linear-pieces}
to $t=t_{j-1}$ and $t=t_j$.
For each fixed $j$, the root $\nu_j$ is constant, so $g$ is a
polynomial of degree at most one in $t$ on $[t_{j-1},t_j]$.
Thus $g$ is continuous and piecewise linear, with constant
derivative on each $(t_{j-1},t_j)$.
If the segment is contained in a separating ray, take $m=1$.

\begin{figure}[!htb]
\centering
\resizebox{0.98\linewidth}{!}{\input{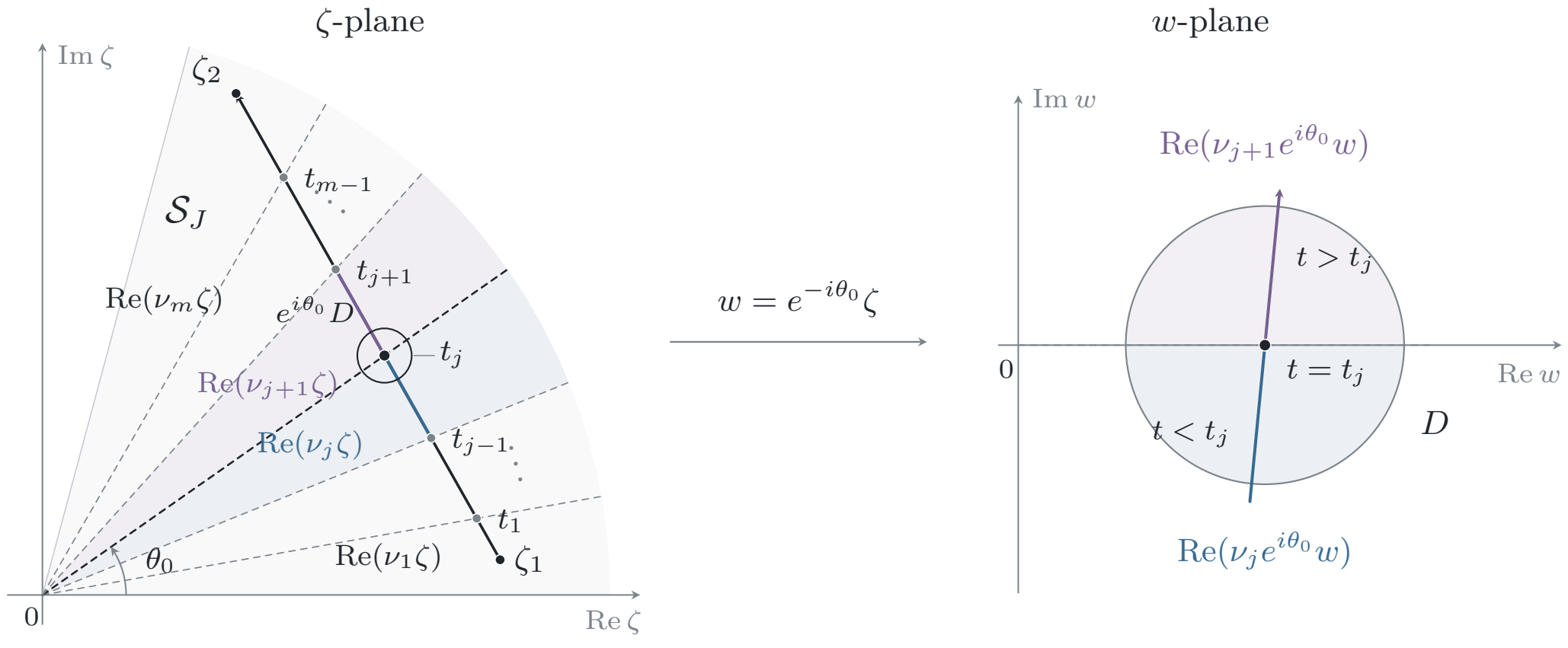}}
\caption{Segment decomposition and local rotation.}
\label{fig:segment-partition-rotation}
\end{figure}

For $1\le j<m$, apply~\eqref{eq:root-local-maximum} at the
separating ray of argument
$\theta_0=\arg((1-t_j)\zeta_1+t_j\zeta_2)$, as illustrated in
Figure~\ref{fig:segment-partition-rotation}.
Choose the disk $D$ in this calculation centred at
$e^{-i\theta_0}((1-t_j)\zeta_1+t_j\zeta_2)>0$.
Recall that $\widetilde H(w)=H_v(e^{i\theta_0}w)$, so
\[
 \widetilde H(e^{-i\theta_0}\zeta)
 =H_v(e^{i\theta_0}e^{-i\theta_0}\zeta)
 =H_v(\zeta),\qquad \zeta\in e^{i\theta_0}D.
\]
The adjacent roots are $\nu_j$ and $\nu_{j+1}$.
For either $\nu\in\{\nu_j,\nu_{j+1}\}$, its linear expression
in the rotated coordinates is $\Rea(\nu e^{i\theta_0}w)$.
For $t$ sufficiently close to $t_j$,
$w=e^{-i\theta_0}((1-t)\zeta_1+t\zeta_2)$ lies in $D$, and
\[
 \Rea(\nu e^{i\theta_0}w)
 =\Rea\bigl(\nu((1-t)\zeta_1+t\zeta_2)\bigr)
 =\Rea(\nu\zeta_1)+t\,\Rea\bigl(\nu(\zeta_2-\zeta_1)\bigr).
\]
Thus~\eqref{eq:root-local-maximum} gives, for $t$ sufficiently
close to $t_j$,
\[
 g(t)=\max_{\nu\in\{\nu_j,\nu_{j+1}\}}
 \bigl\{\Rea(\nu\zeta_1)
       +t\,\Rea\bigl(\nu(\zeta_2-\zeta_1)\bigr)\bigr\}.
\]
These are exactly the two polynomials in $t$ in
\eqref{eq:root-segment-linear-pieces} on the intervals adjacent
to $t_j$. Their maximum is convex, so
$g'_+(t_j)\ge g'_-(t_j)$, where $g'_-$ and $g'_+$ are the
left and right derivatives.
\(g\) is convex on \([0,1]\). In particular,
$g(t)\le(1-t)g(0)+tg(1)$ for $0\le t\le1$.
Equivalently,
\begin{equation}\label{eq:root-indicator-convexity}
 H_v((1-t)\zeta_1+t\zeta_2)
 \le(1-t)H_v(\zeta_1)+tH_v(\zeta_2),
 \qquad 0\le t\le1.
\end{equation}
Since $\zeta_1,\zeta_2$ are arbitrary, this proves that $H_v$
is convex on $\mathcal S_J$.

Fix $\mu\in\Lambda_v(J)$. By definition, choose a nonempty open
interval $I_\mu\subset J$ on which
$h_v(\theta)=\Rea(\mu e^{i\theta})$.
Then $H_v(\zeta)=\Rea(\mu\zeta)$ on the open subsector
$\{\zeta\in\mathcal S_J:\arg\zeta\in I_\mu\}$.
Choose $\zeta_0$ in this subsector. Since $H_v$ agrees with the
linear function $\Rea(\mu\zeta)$ near $\zeta_0$, it is
differentiable there with
$\nabla H_v(\zeta_0)=(\Rea\mu,-\Ima\mu)$.
Write $\mu=a+ib$ and $\zeta-\zeta_0=\Delta x+i\Delta y$.
Direct multiplication gives
\begin{equation}\label{eq:gradient-real-part}
 \nabla H_v(\zeta_0)\cdot(\Delta x,\Delta y)
 =a\Delta x-b\Delta y
 =\Rea\bigl((a+ib)(\Delta x+i\Delta y)\bigr)
 =\Rea\bigl(\mu(\zeta-\zeta_0)\bigr).
\end{equation}
The supporting-plane inequality for the convex function $H_v$,
using~\eqref{eq:root-indicator-convexity}, therefore gives
\[
 H_v(\zeta)\ge H_v(\zeta_0)
       +\Rea(\mu(\zeta-\zeta_0))
       =\Rea(\mu\zeta),\qquad \zeta\in\mathcal S_J.
\]
Since $\mu$ is arbitrary,
$H_v(\zeta)\ge\max_{\mu\in\Lambda_v(J)}\Rea(\mu\zeta)$.

For the reverse inequality, fix $\zeta\in\mathcal S_J$.
This point lies in the closure of an open subsector on which
$H_v(\eta)=\Rea(\nu\eta)$ for a fixed $\nu\in\Lambda_v(J)$.
Letting $\eta\to\zeta$ within that subsector and using
continuity gives $H_v(\zeta)=\Rea(\nu\zeta)$, and hence
$H_v(\zeta)\le\max_{\mu\in\Lambda_v(J)}\Rea(\mu\zeta)$.
Together these inequalities give equality with the maximum.
Setting $\zeta=e^{i\theta}$ proves~\eqref{eq:local-support-indicator}.

Set $K_v=\conv\{\overline\nu:\nu\in\Lambda_v(J)\}$.
By~\eqref{eq:gradient-real-part}, each linear piece
$\Rea(\nu\zeta)$ has gradient $(\Rea\nu,-\Ima\nu)$,
represented by $\overline\nu$ in $\C$.
Thus $K_v$ is the convex hull of the gradients of these pieces.
Suppose, for a contradiction, that $\overline\mu$ is not a
vertex of $K_v$ for some $\mu\in\Lambda_v(J)$.
Taking conjugates in its representation as a convex combination
of the other points gives
\[
 \mu=\sum_{\nu\in\Lambda_v(J)\setminus\{\mu\}}a_\nu\nu,
 \qquad a_\nu\ge0,\qquad \sum_{\nu\ne\mu}a_\nu=1.
\]
For $\theta\in I_\mu$, we have $h_v(\theta)=\Rea(\mu e^{i\theta})$,
and hence
\[
 0=h_v(\theta)-\Rea(\mu e^{i\theta})
  =\sum_{\nu\ne\mu}a_\nu
       \bigl(h_v(\theta)-\Rea(\nu e^{i\theta})\bigr).
\]
Every summand is nonnegative by~\eqref{eq:local-support-indicator}.
Choose $\nu$ with $a_\nu>0$. Then
$\Rea((\mu-\nu)e^{i\theta})=0$ throughout $I_\mu$.
Differentiating this identity on the open interval $I_\mu$
also gives $\Ima((\mu-\nu)e^{i\theta})=0$ there.
Thus $\mu=\nu$, contradicting $\nu\ne\mu$.
It follows that every $\overline\mu$, $\mu\in\Lambda_v(J)$,
is a vertex of $K_v$.
Conversely, every vertex of the finite convex hull $K_v$
belongs to $\{\overline\nu:\nu\in\Lambda_v(J)\}$.
\end{proof}

\phantomsection\label{proof:sectorial-counting}
\begin{proof}[Proof of {\hyperref[prop:sectorial-counting]{Proposition~\ref*{prop:sectorial-counting}}}]
Lemmas~\ref{lem:indicator-limit} and~\ref{lem:local-root-selection}
establish~\eqref{eq:support-limit}. It remains to prove
\eqref{eq:angular-count}.
Set $H=H_v$, $h=h_v$ and $K=K_v$.
With $R_0>R_*$ fixed, write $n_{\alpha_1,\alpha_2}(R)$ for
the $s$-plane zero count in~\eqref{eq:angular-count}.
By~\eqref{eq:homogeneous-indicator},
$H(re^{i\theta})=rh(\theta)$ for $r>0$, $\theta\in J$.

On each open angular interval containing no edge normal of $K$,
a fixed vertex $\nu$ realizes the maximum, so
\[
 h(\theta)=\Rea(\nu e^{-i\theta})
 =(\Rea\nu)\cos\theta+(\Ima\nu)\sin\theta,
 \qquad h''(\theta)=-h(\theta).
\]
At an edge normal $\theta_0\in J$, let $\nu_-$ and $\nu_+$ be
the maximizing endpoints for $\theta<\theta_0$ and
$\theta>\theta_0$, respectively, with $\theta$ sufficiently
close to $\theta_0$.
Since $\nu_+-\nu_-=d_K(\theta_0)i e^{i\theta_0}$,
\begin{equation}\label{eq:angular-support-jump}
 h'(\theta_0+)-h'(\theta_0-)
 =\Rea\!\left(-i(\nu_+-\nu_-)e^{-i\theta_0}\right)
 =d_K(\theta_0).
\end{equation}
See Levin~\cite[Chapter~I, \S19, pp.~79--80]{Levin}.
For a segment $K$, each perpendicular normal contributes its
full length.
All angular sums below are over the finitely many edge normal
angles $\theta_0\in J$ with $d_K(\theta_0)>0$.

Take $\phi\in C_c^\infty(\mathcal S_J)$, using the same symbol
for its polar-coordinate expression $\phi(r,\theta)=\phi(re^{i\theta})$.
Its polar expression has compact support in $(0,\infty)\times J$.
For fixed $r>0$, split $J$ at the edge normal angles.
On each resulting interval $(a,b)$, two integrations by parts
and $h''=-h$ give
\[
 \int_a^b h(\phi+\phi_{\theta\theta})\dd\theta
 =[h\phi_\theta-h'\phi]_a^b.
\]
Continuity of $h$ cancels the terms involving
$\phi_\theta$, and compact support removes the outer
boundary terms. Summing and using~\eqref{eq:angular-support-jump}
gives
\begin{equation}\label{eq:angular-support-pairing}
 \int_J h(\theta)
       \bigl(\phi(r,\theta)+\phi_{\theta\theta}(r,\theta)\bigr)
       \dd\theta
 =\sum_{\theta_0}d_K(\theta_0)\phi(r,\theta_0).
\end{equation}
For fixed $\theta\in J$, radial integration by parts gives
\begin{equation}\label{eq:radial-test-integration}
 \int_0^\infty
 (r^2\phi_{rr}+r\phi_r-\phi)\dd r
 =[r^2\phi_r-r\phi]_0^\infty=0.
\end{equation}
Using $\Delta=\partial_r^2+r^{-1}\partial_r+r^{-2}\partial_\theta^2$
and $\mathrm dA=r\dd r\dd\theta$, then
\eqref{eq:radial-test-integration} and~\eqref{eq:angular-support-pairing},
we obtain
\[
 \begin{aligned}
 \langle\Delta H,\phi\rangle
 &=\int_0^\infty\!\int_J h(\theta)
 (r^2\phi_{rr}+r\phi_r+\phi_{\theta\theta})
 \dd\theta\dd r\\
 &=\int_0^\infty\!\int_J
 h(\theta)(\phi+\phi_{\theta\theta})
 \dd\theta\dd r\\
 &=\sum_{\theta_0}d_K(\theta_0)
       \int_0^\infty\phi(re^{i\theta_0})\dd r.
 \end{aligned}
\]
Let $\sigma_{\theta_0}$ denote Euclidean length measure on the ray
$\{re^{i\theta_0}:r>0\}$. Then
\begin{equation}\label{eq:ray-length-measure}
 \int\phi\,\mathrm d\sigma_{\theta_0}
 =\int_0^\infty\phi(re^{i\theta_0})\dd r.
\end{equation}
Since $\phi$ is arbitrary, we obtain the distributional identity
\begin{equation}\label{eq:laplacian-support}
 \Delta H=\sum_{\theta_0}d_K(\theta_0)\,\sigma_{\theta_0}
 \qquad\text{on }\mathcal S_J.
\end{equation}

Recall $U_R(\zeta)=R^{-1}\log|v(R\zeta)|$ from
\eqref{eq:scaled-scalar-and-jet-potentials}.
If $b$ is a zero of $v$ of multiplicity $m_b$, then $v(R\zeta)$
has a zero at $\zeta=b/R$ of the same multiplicity. Define
\[
 \nu_R=\frac1R
   \sum_{\substack{v(b)=0,\ |b|>R_0\\ \arg b\in J}}m_b\delta_{b/R}.
\]
On every fixed relatively compact subdomain of $\mathcal S_J$,
\eqref{eq:planar-zero-measure} gives
$\nu_R=(2\pi)^{-1}\Delta_\zeta U_R$ for all sufficiently large $R$,
since the inner radius $R_0/R$ tends to zero.
Thus~\eqref{eq:scalar-limit} gives
\begin{equation}\label{eq:scaled-zero-test-convergence}
 \int\phi\,\mathrm d\nu_R
 =\frac1{2\pi}\int_{\operatorname{supp}\phi}U_R\,\Delta\phi\,\mathrm dA
 \longrightarrow
 \frac1{2\pi}\int_{\operatorname{supp}\phi}H\,\Delta\phi\,\mathrm dA.
\end{equation}
For a compact set $B\subset\mathcal S_J$, the definition of $\nu_R$ gives
\[
 \nu_R(B)=\frac1R
 \sum_{\substack{v(b)=0,\ |b|>R_0\\ \arg b\in J,\ b/R\in B}}m_b.
\]
Choose a smooth cutoff function satisfying
\[
 \phi\in C_c^\infty(\mathcal S_J),\qquad
 0\le\phi\le1,\qquad \phi=1\ \text{on }B.
\]
Since $\nu_R\ge0$, \eqref{eq:scaled-zero-test-convergence} gives
$0\le\nu_R(B)=\int_B1\,\mathrm d\nu_R\le\int\phi\,\mathrm d\nu_R=O(1)$.
\phantomsection\label{par:sectorial-smooth-approximation}
Put $\nu=(2\pi)^{-1}\Delta H$ and fix $f\in C_c(\mathcal S_J)$.
By smooth approximation
\cite[Proposition~2.1.1, p.~21]{DinhSibonyCurrents}, for each
$\varepsilon>0$ choose $\phi_\varepsilon\in C_c^\infty(\mathcal S_J)$
with $\|f-\phi_\varepsilon\|_\infty<\varepsilon$ and
$\operatorname{supp}f\cup\operatorname{supp}\phi_\varepsilon\subset B$,
where the compact set $B\subset\mathcal S_J$ is independent of
$\varepsilon$. For sufficiently large $R$, the mass bound
$\nu_R(B)\le C_B$ gives
\[
 \begin{aligned}
 \left|\int f\,\mathrm d\nu_R-\int f\,\mathrm d\nu\right|
 \le{}\varepsilon\bigl(C_B+\nu(B)\bigr)+\left|\int\phi_\varepsilon\,\mathrm d\nu_R
             -\int\phi_\varepsilon\,\mathrm d\nu\right|.
 \end{aligned}
\]
For fixed $\varepsilon$, the last term tends to zero
by~\eqref{eq:scaled-zero-test-convergence}.
Since $\nu(B)<\infty$, letting first $R\to\infty$ and then
$\varepsilon\to0$ proves that $\nu_R$ converges weakly to $\nu$
against continuous compactly supported test functions
\cite[Definition~1.3.5, p.~16]{DinhSibonyCurrents}.

Put $c=(2\pi)^{-1}
\sum_{\alpha_1<\theta_0<\alpha_2}d_K(\theta_0)$ and let
\[
 A=\{\zeta:\tfrac12<|\zeta|\le1,\
                   \alpha_1<\arg\zeta<\alpha_2\}.
\]
Its closure is compact in $\mathcal S_J$.
Recall from~\eqref{eq:ray-length-measure} that $\sigma_{\theta_0}$
is Euclidean length measure on $\{re^{i\theta_0}:r>0\}$.
Thus, for $\alpha_1<\theta_0<\alpha_2$,
\[
 \sigma_{\theta_0}(A)=\int_{1/2}^1\dd r=\frac12.
\]
Hence~\eqref{eq:laplacian-support} gives $\nu(A)=c/2$.
Each of the two circular boundary arcs meets each ray in at most
one point and hence has zero $\nu$-measure. The angular sides also
have zero $\nu$-measure because \mbox{$d_K(\alpha_1)=d_K(\alpha_2)=0$}.
Thus $\nu(\partial A)=0$.
Since $\nu_R(A)=\int\mathbf1_A\,\mathrm d\nu_R$, we use
continuous compactly supported functions bounding the indicator
$\mathbf1_A$ from below and above and agreeing with it outside
arbitrarily small neighbourhoods of $\partial A$.
By weak convergence
\cite[Exercise~1.3.12, p.~17]{DinhSibonyCurrents},
\begin{equation}\label{eq:weak-measure-set-bounds}
 \nu(\operatorname{int}(A))
 \le\liminf_{R\to\infty}\nu_R(A)
 \le\limsup_{R\to\infty}\nu_R(A)
 \le\nu(\overline A).
\end{equation}
Both outer terms equal $\nu(A)$, so
$\nu_R(A)\to\nu(A)=c/2$.
Since $R\nu_R(A)=n_{\alpha_1,\alpha_2}(R)
-n_{\alpha_1,\alpha_2}(R/2)$ for $R>2R_0$,
\begin{equation}\label{eq:annular-count-proof}
 n_{\alpha_1,\alpha_2}(R)
 -n_{\alpha_1,\alpha_2}(R/2)=cR/2+o(R).
\end{equation}

Finally, write $n(t)=n_{\alpha_1,\alpha_2}(t)$ for brevity.
Fix $\varepsilon>0$ and choose $T=T(\varepsilon)>2R_0$ so that
the error in~\eqref{eq:annular-count-proof} is at most
$\varepsilon t$ for every $t\ge T$.
Keeping $T$ fixed, for each $R\ge T$ choose an integer $m=m(R)\ge0$
such that $R/2^m\ge T>R/2^{m+1}$.
We sum~\eqref{eq:annular-count-proof} at the radii $R/2^j\ge T$,
$0\le j\le m$, using
\[
 \sum_{j=0}^{m}\bigl(n(R/2^j)-n(R/2^{j+1})\bigr)
 =n(R)-n(R/2^{m+1}).
\]
The main terms sum to $cR(1-2^{-(m+1)})$.
The remaining count satisfies $n(R/2^{m+1})\le n(T)<\infty$,
since the fixed truncated sector of outer radius $T$ has compact
closure in the domain of $v$. Also, $cR/2^{m+1}\le cT$. Therefore
\[
 |n(R)-cR|
 \le\varepsilon\sum_{j=0}^{m}\frac R{2^j}+n(T)+cT
 \le2\varepsilon R+O_T(1).
\]
Let $R\to\infty$ and then $\varepsilon\to0$
to obtain~\eqref{eq:angular-count}.
\end{proof}

We now transfer a sectorial count in the $s$-plane to the original
disk $|z|\le r$. Use the sector $S_J$ from~\eqref{eq:phase-equation}
and let $\chi:S_J\to T$ be conformal, extending analytically across
every finite boundary point of $S_J$. Suppose that
$f(\chi(s))=G(\chi(s))v(s)$ on $S_J$, where $f$ is nonzero and
entire, $G$ is analytic and zero-free on $T$, and $v$ is analytic
on $S_J$.

Fix $[\alpha_1,\alpha_2]\Subset J$ and put
\[
 U=\chi\{s\in S_J:\alpha_1<\arg s<\alpha_2\}.
\]
Write $N_v(R)$ for the sectorial zero count on the left-hand side
of~\eqref{eq:angular-count}.
\phantomsection\label{par:radial-transfer-monotonicity}
Assume that, for some constant $R_1>R_0$,
$\rho\mapsto|\chi(\rho e^{i\theta})|$ is strictly increasing on
$[R_1,\infty)$ for every $\theta\in[\alpha_1,\alpha_2]$.
For our inverse branches $\chi=\omega^{-1}\psi$, this condition
is verified in the \hyperref[par:inverse-radial-monotonicity]{proof
of Theorem~\ref*{thm:intro}}.

\begin{lemma}\label{lem:radial-transfer}
In this setting, let $a>0$ and $Q(r)\to\infty$.
Suppose that
\begin{equation}\label{eq:transfer-sector-count}
 N_v(R)=c_fR+o(R),\qquad R\to\infty,
\end{equation}
and, for every fixed $c>0$,
\begin{equation}\label{eq:transfer-radius}
 |\chi(cQ(r)e^{i\theta})|
 =r+\frac{\log c}{a}+o(1),\qquad r\to\infty,
\end{equation}
uniformly for $\alpha_1\le\theta\le\alpha_2$. Then
\[
 n_U(r,1/f)=c_fQ(r)+o(Q(r)),\qquad r\to\infty.
\]
\end{lemma}

\begin{proof}
Conformality of $\chi$ and the absence of zeros of $G$ preserve
zero multiplicities. We compare the conditions $|s|\le R$
and $|\chi(s)|\le r$.

Fix $0<\varepsilon<1$. Apply~\eqref{eq:transfer-radius} with
$c=1-\varepsilon$ and $c=1+\varepsilon$. Since the remainder is
uniform in $\theta$, $a>0$, and
$\log(1-\varepsilon)<0<\log(1+\varepsilon)$, there exists
$r_\varepsilon>0$ such that
\begin{equation}\label{eq:transfer-cutoff-bracket}
 |\chi((1-\varepsilon)Q(r)e^{i\theta})|
 <r<
 |\chi((1+\varepsilon)Q(r)e^{i\theta})|
\end{equation}
for all $r\ge r_\varepsilon$ and $\theta\in[\alpha_1,\alpha_2]$.
Recall $R_1>R_0$ from the
\hyperref[par:radial-transfer-monotonicity]{radial monotonicity assumption}.
Since $Q(r)\to\infty$, we may increase $r_\varepsilon$ so that
$(1-\varepsilon)Q(r)>R_1$ for all $r\ge r_\varepsilon$.

For $r\ge r_\varepsilon$ and $s\in S_J$ with
$\alpha_1<\arg s<\alpha_2$, the monotonicity assumption and
\eqref{eq:transfer-cutoff-bracket} imply
\[
 \begin{aligned}
 R_1\le |s|\le(1-\varepsilon)Q(r)
 &\ \Longrightarrow\ |\chi(s)|<r,\\
 |s|\ge R_1,\quad |\chi(s)|\le r
 &\ \Longrightarrow\ |s|<(1+\varepsilon)Q(r).
 \end{aligned}
\]

The part with $R_0<|s|\le R_1$ has image with compact closure:
$\chi$ extends across the inner boundary, and
$[\alpha_1,\alpha_2]\Subset J$.
That compact set contains only finitely many zeros of the
nonzero entire function $f$. Their total multiplicity is bounded
by a constant $C_0$, independent of $r$.
The correspondence of zeros therefore gives
\begin{equation}\label{eq:transfer-count-bracket}
 N_v((1-\varepsilon)Q(r))-C_0
 \le n_U(r,1/f)
 \le N_v((1+\varepsilon)Q(r))+C_0.
\end{equation}

Insert~\eqref{eq:transfer-sector-count} into
\eqref{eq:transfer-count-bracket} and divide by $Q(r)$.
Since $Q(r)\to\infty$, the constant $C_0$ disappears, and
\[
 c_f(1-\varepsilon)
 \le\liminf_{r\to\infty}\frac{n_U(r,1/f)}{Q(r)}
 \le\limsup_{r\to\infty}\frac{n_U(r,1/f)}{Q(r)}
 \le c_f(1+\varepsilon).
\]
Letting $\varepsilon\to0$ proves the conclusion.
\end{proof}

The sectorial formula places the limiting zero measure on finitely
many rays in the $s$-plane. We now estimate the number of zeros
outside shrinking neighbourhoods of their inverse images in a
fixed strip.
Fix $\omega\in\Omega_a$, put $\beta=\beta(\omega)$ as
in~\eqref{eq:beta}, and choose one inverse branch $\psi$ supplied by
Lemma~\ref{lem:phase-equation}. Write $\varphi:T\to S_J$ and
$\varphi^{-1}(s)=\omega^{-1}\psi(s)$.
For a nonzero solution $f$ of~\eqref{eq:reduced}, put
$v(s)=f(\omega^{-1}\psi(s))$ on $S_J$ and let
$K_f=\conv\{\overline\mu:\mu\in\Lambda_v(J)\}$, with the root
set supplied by Proposition~\ref{prop:sectorial-counting}.
For a supporting normal $\theta$ of $K_f$, put
\[
 \Gamma_\theta=
       \{\omega^{-1}\psi(\rho e^{i\theta}):\rho>R_0\}.
\]
These disjoint analytic curves in the $z$-plane tend to infinity
within the fixed strip. The angle $\theta$ is measured in the
$s$-plane; all the curves extend in direction $-\arg\omega$
in the $z$-plane.

\begin{proposition}\label{prop:shrinking-tubes}
Let $f$, the inverse branch and $K_f$ be as above. Fix $\alpha_1<\alpha_2$ in $J$ with
$d_{K_f}(\alpha_1)=d_{K_f}(\alpha_2)=0$, and put
\[
 U=U(\alpha_1,\alpha_2),\qquad
 \Theta=\{\theta\in(\alpha_1,\alpha_2):d_{K_f}(\theta)>0\},
 \qquad \Gamma=\bigcup_{\theta\in\Theta}\Gamma_\theta.
\]
If $\Theta\ne\varnothing$, there is a nonincreasing function
$\varepsilon(r)>0$, tending to zero, such that
\[
 n\bigl(\{z\in U:|z|\le r,\
             \operatorname{dist}(z,\Gamma)>\varepsilon(r)\},1/f\bigr)
       =o\bigl(r^\beta e^{|\omega|r}\bigr).
\]
\end{proposition}

\begin{proof}
Put $Q(r)=(|\omega|r)^\beta e^{|\omega|r}$ and
$Z(R,\theta)=\omega^{-1}\psi(Re^{i\theta})$.
Recall $\Phi$ from~\eqref{eq:phase}.
By~\eqref{eq:phase-expansion-proof},
$\Phi(t)/\Phi'(t)=1+O(1/t)$ along the fixed inverse branch, so
\[
 \frac{\partial Z}{\partial\theta}(R,\theta)
 =\frac{i}{\omega}\frac{\Phi(\psi(Re^{i\theta}))}
                            {\Phi'(\psi(Re^{i\theta}))}
 =\frac{i}{\omega}+O\!\left(\frac1{\log R}\right)
\]
uniformly for $\alpha_1\le\theta\le\alpha_2$.
Thus there are constants $C_1>0$ and $R_1>R_0$ such that
$|\partial_\theta Z(R,\theta)|\le C_1$ for $R\ge R_1$ and
$\alpha_1\le\theta\le\alpha_2$.
For $\theta_0\in\Theta$ with $|\theta-\theta_0|\le\delta$,
integration in $\theta$ gives
\[
 \operatorname{dist}(Z(R,\theta),\Gamma)
 \le |Z(R,\theta)-Z(R,\theta_0)|
 \le C_1|\theta-\theta_0|\le C_1\delta,
\]
since $Z(R,\theta_0)\in\Gamma_{\theta_0}\subset\Gamma$.

Choose a decreasing sequence $\delta_m>0$ with $\delta_m\to0$ so that the angular neighbourhoods
of radius $\delta_m$ about the points of $\Theta$ are disjoint
and contained in $(\alpha_1,\alpha_2)$. Their complement in
$[\alpha_1,\alpha_2]$ is a finite union of closed angular intervals
on which $d_{K_f}(\theta)=0$, by the definition of $\Theta$ and the
endpoint conditions $d_{K_f}(\alpha_1)=d_{K_f}(\alpha_2)=0$.
Applying~\eqref{eq:angular-count} on slightly larger intervals
where $d_{K_f}(\theta)=0$, followed by Lemma~\ref{lem:radial-transfer},
gives an $o(Q(r))$ bound for zeros in the pullbacks of the original
closed intervals.
The part with $R_0<|s|\le R_1$ contributes only finitely many
zeros, as in the proof of Lemma~\ref{lem:radial-transfer}.
Thus, for every fixed $m$, the number of zeros at
distance greater than $C_1\delta_m$ from $\Gamma$ is $o(Q(r))$.
Choose a strictly increasing sequence $r_m\to\infty$ so that this count is at most
$Q(r)/m$ for all $r\ge r_m$. Taking
$\varepsilon(r)=C_1\delta_m$ for $r_m\le r<r_{m+1}$ proves
the assertion.
\end{proof}

\subsection{Exceptional initial values}
\label{subsec:exceptional-subspaces}

Theorem~\ref{thm:generic-density} guarantees a positive leading
coefficient in a suitably chosen strip for initial vectors outside
the union of two proper complex linear subspaces.
For this purpose, we will choose a ray on which one characteristic
root has a strictly larger real projection than the others. We therefore need to identify the initial values that
remove this maximal growth. The next lemma constructs those
initial values as a linear subspace by describing the solutions
that remain bounded after an exponential rescaling. This is the
constant-limit case of the stability of exponential dichotomies;
see \cite[Lecture~4]{CoppelDichotomies}.

\begin{lemma}\label{lem:slower-growth-subspace}
Let $\mathbf B:[r_*,\infty)\to\C^{n\times n}$ be continuous and satisfy
$\mathbf B(r)\to \mathbf B_0$. Suppose that $\mathbf B_0$ has exactly one eigenvalue with
maximal real part $\sigma$, and that its algebraic multiplicity is
$m<n$.
Every nonzero solution $\mathbf x:[r_*,\infty)\to\C^n$ of
$\mathrm d\mathbf x/\mathrm dr=\mathbf B(r)\mathbf x$ in the real variable $r$
has a growth rate
\[
 \lambda(\mathbf x)=\lim_{r\to\infty}\frac{\log\|\mathbf x(r)\|}{r}.
\]
The initial values $\mathbf x(r_*)$ for which $\lambda(\mathbf x)<\sigma$,
together with $0$, form a complex linear subspace $V$ with
$\dim V=n-m$. Every solution with $\mathbf x(r_*)\notin V$ has
$\lambda(\mathbf x)=\sigma$.
\end{lemma}

The proof is given in \hyperref[subsec:proof-slower-growth-subspace]{Section~\ref*{subsec:proof-slower-growth-subspace}}.

\begin{remark}
The conclusion of Lemma~\ref{lem:slower-growth-subspace} also holds
when several distinct eigenvalues have maximal real part, provided
$m<n$ denotes their total algebraic multiplicity.
\end{remark}

\begin{example}\label{ex:fourth-order-subspaces}
In the second-order model~\eqref{eq:model-original}, the two
exceptional subspaces are the hyperplanes $A=0$ and $B=0$ in the
coordinates $f=Af_++Bf_-$. For an example of codimension two,
consider on $\C$ the equation
\[
 f^{(4)}-\left(2e^{2z}+\frac52\right)f''-4e^{2z}f'
 +\left(e^{4z}-\frac32e^{2z}+\frac9{16}\right)f=0.
\]
Fix an analytic logarithm on a simply connected domain
$D\subset\C\setminus\{0\}$. The substitutions $z=\log s$ and
$v(s)=e^{3\log s/2}f(\log s)$ give $v^{(4)}-2v''+v=0$ on $D$.
This constant-coefficient equation and its solutions extend to
all $s\in\C$. Hence the general solution of the original equation on $\C$ is
\[
 f=e^{-3z/2}\bigl((A+Be^z)e^{e^z}+(C+De^z)e^{-e^z}\bigr).
\]
Fix $0<\varepsilon<\pi/2$. Along
$\arg s=\pi/2-\varepsilon$, the factor $e^s$ grows and
$e^{-s}$ decays, so removing the growing part requires
$A=B=0$. Along $\arg s=\pi/2+\varepsilon$, it instead
requires $C=D=0$. Thus the two subspaces in coefficient
coordinates are
\[
 \{(A,B,C,D):A=B=0\},\qquad
 \{(A,B,C,D):C=D=0\}.
\]
Both have complex dimension $2$, whereas a hyperplane in
$\C^4$ has dimension $3$. The invertible linear map from
coefficient vectors to initial vectors preserves these dimensions.

Multiplication by a zero-free analytic factor preserves zeros
and their multiplicities. On a conformal branch
$\omega^{-1}\psi:S_J\to T$, let $H$ be analytic on $T$ and set
\[
 v(s)=e^{-H(\omega^{-1}\psi(s))}f(\omega^{-1}\psi(s)),\qquad s\in S_J.
\]
If $v$ satisfies the hypotheses of Proposition~\ref{prop:sectorial-counting},
Lemma~\ref{lem:radial-transfer}, with $G=e^H$, transfers its
sectorial zero count to the original disk $|z|\le r$.
In the fourth-order example above, $\omega=1$, $s=e^z$ and
$H(z)=-3z/2$.
\end{example}

\subsection{Zero counts for nonmaximal frequencies}\label{subsec:nonmaximal-frequencies}

The strip-counting conclusions of Theorems~\ref{thm:intro}
and~\ref{thm:generic-density} extend to certain nonmaximal weighted
frequencies. The proof of Lemma~\ref{lem:phase-coefficient-limit}
uses $\Rea(\lambda/\omega)<n-j$ for $\lambda\ne(n-j)\omega$.
Writing $\nu=\lambda/(n-j)$ as in~\eqref{eq:weighted}, these
inequalities for $\omega\in\mathcal F\setminus\{0\}$ become
\begin{equation}\label{eq:own-direction-dominance}
 \Rea(\nu/\omega)<1\qquad
 (\nu\in\{0\}\cup\mathcal F,\ \nu\ne\omega).
\end{equation}
In the weighted-frequency plane, the line
$L_\omega=\{\nu:\Rea(\nu\overline\omega)=|\omega|^2\}$ passes through
$\omega$ and is perpendicular to the segment from $0$ to $\omega$.
Condition~\eqref{eq:own-direction-dominance} says that all other
frequencies lie strictly on the side of $L_\omega$ containing $0$.
Thus $L_\omega$ meets the frequency hull $\mathcal H$ only at $\omega$.

\begin{proposition}\label{prop:nonmaximal-strip-count}
Let $\omega\in\mathcal F\setminus\{0\}$ satisfy
\eqref{eq:own-direction-dominance}.
Define $\beta(\omega)$ by~\eqref{eq:beta} and $\varphi$
by~\eqref{eq:phase}, with its fixed inverse branches as in
Lemma~\ref{lem:phase-equation}. Then that lemma holds for
this choice of $\omega$. The conclusions of Theorems~\ref{thm:intro}
and~\ref{thm:generic-density} hold with $a$ replaced by $|\omega|$.
In particular, outside the union of two proper complex linear
subspaces of initial values, there is a fixed strip with
\[
 n_U(r,1/f)\sim c_\omega r^{\beta(\omega)}e^{|\omega|r},
 \qquad c_\omega>0.
\]
\end{proposition}

The proof is given at the end of
\hyperref[par:proof-nonmaximal-strip-count]{\S\ref*{subsec:proof-generic-density}}.

\begin{example}\label{ex:nonmaximal-weighted-frequency}
Put $\eta=\tfrac35+\tfrac32 i$ and consider
\[
 f''+\bigl(ze^{2z}+e^{2\eta z}+e^{4iz}\bigr)f=0,\qquad z\in\C.
\]
Dividing the three frequencies by $2$ gives
$\mathcal F=\{1,\eta,2i\}$. The frequency $2i$ has maximal modulus.
The nonmaximal frequency $1$ also satisfies~\eqref{eq:own-direction-dominance},
since $\Rea\eta=3/5$ and $\Rea(2i)=0$. Its supporting line
$L_1$ is the vertical line $x=1$ in Figure~\ref{fig:weighted-frequency-support}.
The dashed line $x+y=21/10$ supports
$\mathcal H$ at $\eta$ alone, so $\eta$ is a vertex. However,
the point $2i$ lies beyond $L_\eta$, since $\Rea(2i/\eta)=100/87>1$.

\begin{figure}[htbp]
\centering
\begin{tikzpicture}[
  x=1.55cm,y=1.55cm,font=\small,
  wf axis/.style={-{Stealth[length=3.8pt]},draw=black!55,line width=.45pt},
  wf edge/.style={draw=black!60,line width=.6pt},
  wf point/.style={circle,fill=black!80,inner sep=1.6pt},
  wf selected/.style={circle,inner sep=2pt},
  wf normal/.style={-{Stealth[length=4.5pt]},line width=.8pt},
  wf line/.style={line width=.95pt},
  pics/wf hull/.style={code={
    \filldraw[wf edge,fill=black!6] (0,0)--(1,0)--(.6,1.5)--(0,2)--cycle;
  }},
  pics/wf axes/.style={code={
    \draw[wf axis] (-.20,0)--(1.72,0) node[right,black] {$\Re\nu$};
    \draw[wf axis] (0,-.13)--(0,2.27) node[left,black] {$\Im\nu$};
  }}]
  \definecolor{wfTeal}{RGB}{24,109,116}
  \definecolor{wfRust}{RGB}{157,62,40}
  \definecolor{wfSlate}{RGB}{60,80,126}

  \begin{scope}
    \pic {wf hull};
    \pic {wf axes};
    \draw[wf line,wfTeal] (1,-.12)--(1,2.11)
      node[below right,wfTeal] {$L_1$};
    \draw[wf normal,wfTeal] (0,0)--(.93,0);
    \draw[thin,wfTeal] (.88,0)--(.88,.12)--(1,.12);
    \node[wf point,label=below left:$0$] at (0,0) {};
    \node[wf selected,fill=wfTeal,label=below:$1$] at (1,0) {};
    \node[wf point,label=above right:$\eta$] at (.6,1.5) {};
    \node[wf point,label=left:$2i$] at (0,2) {};
    \node at (.82,-.45) {(a) $\omega=1$};
  \end{scope}

  \begin{scope}[shift={(3.78,0)}]
    \pic {wf hull};
    \fill[wfRust!23] (0,1.74)--(.6,1.5)--(0,2)--cycle;
    \draw[wf edge] (0,2)--(.6,1.5);
    \pic {wf axes};
    \draw[wfSlate,dashed,line width=.85pt] (-.06,2.16)--(1.68,.42);
    \draw[wf line,wfRust] (-.17,1.808)--(1.68,1.068)
      node[above,wfRust] {$L_\eta$};
    \draw[wf normal,wfRust] (0,0)--(.57,1.425);
    \node[wf point,label=below left:$0$] at (0,0) {};
    \node[wf point,label=below:$1$] at (1,0) {};
    \node[wf selected,fill=wfRust,label=above right:$\eta$] at (.6,1.5) {};
    \node[wf point,label=left:$2i$] at (0,2) {};
    \node at (.82,-.45) {(b) $\omega=\eta$};
  \end{scope}
\end{tikzpicture}
\caption{Supporting lines of the weighted frequency hull.}
\label{fig:weighted-frequency-support}
\end{figure}
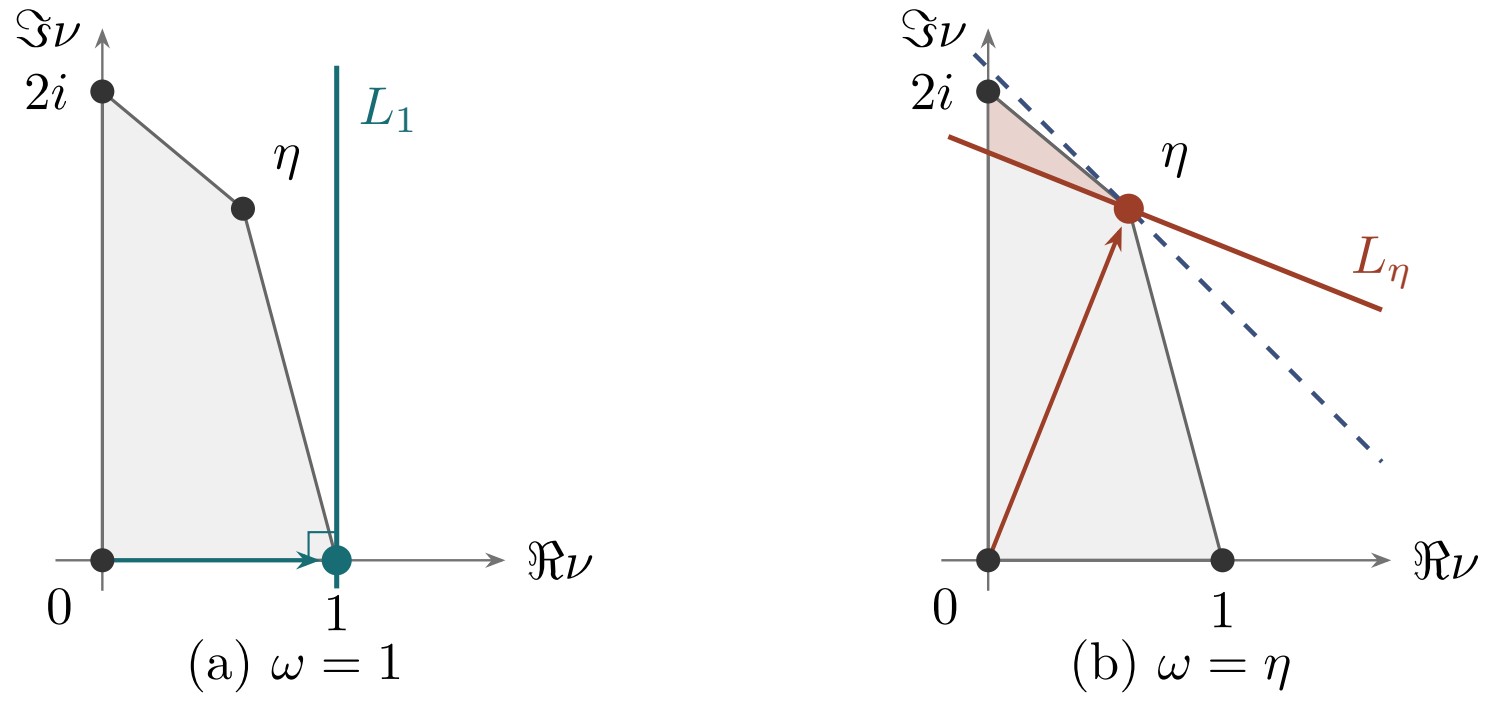

Since $n=2$ and $j=0$, definition~\eqref{eq:beta} gives
$\beta(\omega)=\tfrac12\deg P_{0,2\omega}$. The multipliers of
$e^{2z}$ and $e^{4iz}$ are $z$ and $1$, respectively, so
$\beta(1)=1/2$ and $\beta(2i)=0$.
By~\eqref{eq:limiting-coefficients} and~\eqref{eq:characteristic},
the limiting polynomials are $p_1(\mu)=\mu^2+1$ and
$p_{2i}(\mu)=\mu^2-1/4$, with root segments of lengths $2$ and $1$,
respectively. Proposition~\ref{prop:nonmaximal-strip-count} thus gives
fixed strips $U_1$ extending to the right and $U_{2i}$ extending
downwards for which
\[
 n_{U_1}(r,1/f_{\mathbf c})\sim\frac{\sqrt r\,e^r}{\pi},
 \qquad
 n_{U_{2i}}(r,1/f_{\mathbf c})\sim\frac{e^{2r}}{2\pi}.
\]
Both formulas hold whenever the initial vector
$\mathbf c=(f_{\mathbf c}(0),f_{\mathbf c}'(0))$ lies outside
the union of at most four complex lines in $\C^2$, two for each
frequency.
\end{example}

Recall from Section~\ref{subsec:limiting-equation} that
$\mathcal H=\conv(\{0\}\cup\mathcal F)$. In convex geometry, a point of
$\mathcal H$ is an \textit{exposed point} if some supporting line meets
$\mathcal H$ there alone~\cite[\S3.2.4, p.~29]{Tropp}. In
Example~\ref{ex:nonmaximal-weighted-frequency}, the question is whether
the vector $\omega$ itself exposes the vertex $\omega$.

\begin{corollary}\label{cor:exposed-frequency-vertices}
The frequencies satisfying the hypothesis of
Proposition~\ref{prop:nonmaximal-strip-count} are precisely the nonzero
vertices $\omega$ of $\mathcal H$ for which the line through $\omega$
perpendicular to $[0,\omega]$ supports $\mathcal H$ at $\omega$ alone.
\end{corollary}

\begin{proof}
Since $\Rea(\nu\overline\omega)=|\omega|^2\Rea(\nu/\omega)$,
condition~\eqref{eq:own-direction-dominance} means that $\omega$ uniquely
maximizes $\nu\mapsto\Rea(\nu\overline\omega)$ on $\{0\}\cup\mathcal F$,
or equivalently on its convex hull $\mathcal H$. Thus the level line
$L_\omega$ exposes $\omega$ precisely when this condition holds.
\end{proof}

Even for a vertex outside the scope of Proposition~\ref{prop:nonmaximal-strip-count},
an inverse branch of the phase map may take an exterior sector conformally
onto a curvilinear strip. The following numerical examples suggest that a row of zeros
initially lying in such a strip for $\omega=1$ may leave it and continue
to infinity in a neighbouring strip associated with $q=2+i$.

\begin{figure}[htbp]
\centering
\includegraphics[width=.98\linewidth]{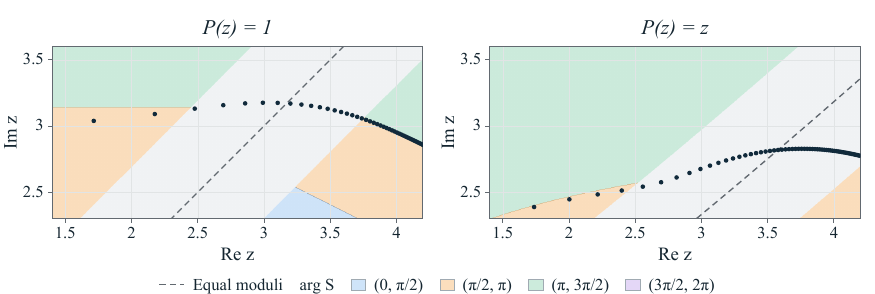}
\caption{Computed zeros near the change of dominant frequency.}
\label{fig:numerical-strip-transition}
\end{figure}

\begin{example}\label{ex:numerical-strip-transition}
Consider the two equations
\[
 f''+\bigl(P(z)e^{2z}+e^{2qz}\bigr)f=0,\qquad z\in\C,\qquad
 P(z)=1\ \text{or}\ P(z)=z,
\]
with $f(0)=0$ and $f'(0)=1$. In both cases $1$ is a vertex of the
frequency hull and fails~\eqref{eq:own-direction-dominance} since
$\Rea(q)=2>1$. Figure~\ref{fig:numerical-strip-transition} shows the
computed zeros on either side of the curve $|P(z)e^{2z}|=|e^{2qz}|$.

For frequency $1$ we use $S=e^z$ or $S=\int_0^z t^{1/2}e^t\dd t$,
respectively, on $\Rea z>0$. For $q$ we use $S=e^{qz}/q$.
The shading indicates the quadrant of $S$ where the corresponding
coefficient term is at least four times the other in modulus and $|S|\ge3$.
\end{example}

\FloatBarrier

\section{Global estimates and transfer to fixed solutions}\label{sec:global}

To count zeros in $|z|\le r$, we use different estimates near and
away from the rays $\arg z=-\arg\omega$, where $\omega$ ranges over
the set $\Omega_*$ defined in~\eqref{eq:dominant-weighted-frequencies}.
Lemma~\ref{lem:coefficient-integral-zero-bound} bounds the zero count
in any bounded region by an integral determined by the coefficients.
Applied to the whole disk, it gives the upper bound in
Theorem~\ref{thm:global-zero-distribution}. Applied outside rectangles
centred at $ar/\omega$, it shows that, after division by $D(r)$,
the contribution there can be made arbitrarily small by enlarging the
rectangles in the rescaled coordinates.
To determine the limiting zero measure inside the rectangles, we use
the weak limit of the rescaled $\log\|\mathbf F\|$ from
Lemma~\ref{lem:physical-window}, where
$\mathbf F=(f_1,\ldots,f_n)$ is the row of fundamental solutions
normalized at $z_0$ by~\eqref{eq:fundamental-initial-matrix}. Along the radii
$r=\sqrt m$ used in the proof of Theorem~\ref{thm:exact-global},
Lemma~\ref{lem:pluripolar-projection} gives the same limit for the
rescaled $\log|f_{\mathbf c}|$, with $f_{\mathbf c}=\mathbf F\mathbf c$,
whenever the initial vector $\mathbf c$
lies outside one complete pluripolar cone.

\subsection{Zero bounds and perimeter averages in moving rectangles}

For the coefficients of~\eqref{eq:reduced}, define
\begin{equation}\label{eq:weight-function}
 W(z)=1+\sum_{j=0}^{n-2}\sum_{\lambda\in\Gamma_j}
 (1+|z|)^{\deg P_{j,\lambda}/(n-j)}
       \exp\!\left(\frac{\Rea(\lambda z)}{n-j}\right).
\end{equation}
The next lemma bounds the number of zeros in $\overline B(z,1)$ by
a constant multiple of $W(z)$. Integrating this bound over the centres
gives a zero-counting bound for any bounded Borel set.
The disk bound used below follows from Yakovenko's estimate for
the variation of argument.

\begin{lemma}\label{lem:coefficient-integral-zero-bound}
Let $W$ be given by~\eqref{eq:weight-function}. There is a
constant $C>0$, depending only on $n$ and the coefficients
$A_0,\ldots,A_{n-2}$, such that every nonzero solution $f$
of~\eqref{eq:reduced} satisfies
\begin{equation}\label{eq:unit-disk-zero-bound}
 n(\overline B(z,1),1/f)\le C W(z),\qquad z\in\C.
\end{equation}
For every bounded Borel set $E\subset\C$, it also satisfies
\begin{equation}\label{eq:coefficient-integral-zero-bound}
 n(E,1/f)\le C\int_{\overline B(E,1)}W(z)\,\mathrm dA(z),
 \qquad \overline B(E,1)=\{z:\operatorname{dist}(z,E)\le1\}.
\end{equation}
\end{lemma}

\begin{proof}
We first work with a monic equation whose coefficients are bounded
on a disk. Let $R\ge1$, and suppose that $g\not\equiv0$ and the
coefficients $b_0,\ldots,b_{n-1}$ are analytic on a neighbourhood
of $|\zeta|\le5R$, with $|b_j|\le B$ there, and
\[
 g^{(n)}+\sum_{j=0}^{n-1}b_j(\zeta)g^{(j)}=0,
 \qquad |\zeta|<5R.
\]
Choose $r\in(R,2R)$ such that $g$ has no zeros on the boundary
of the square
$D_r=\{\zeta:|\Rea\zeta|<r,\ |\Ima\zeta|<r\}$.
It contains $|\zeta|\le R$ and lies in $|\zeta|<3R$.
Yakovenko's rectangle estimate
\cite[Corollary~2.7, p.~507; rectangle example, p.~508]{Yakovenko1999},
applied with coefficient bound $\max\{1,B\}$ and perimeter
$8r<16R$, gives
\begin{equation}\label{eq:bounded-coefficient-disk-count}
 n(\{|\zeta|\le R\},1/g)\le C(n,B)R.
\end{equation}
To see why $C(n,B)$ must depend on the coefficient bound $B$,
take $R=1$ and $g(\zeta)=\sin(K\zeta)$ with $K>0$.
This function solves $g''+K^2g=0$ on $\C$ and has zeros at
$\zeta=m\pi/K$, $m\in\mathbb Z$. Their number in $|\zeta|\le1$ is
$2\lfloor K/\pi\rfloor+1$, which tends to infinity as $K\to\infty$.

We now choose the scale for the original equation.
Fix $z\in\C$ and let $|w-z|\le5$. For each term of $A_j$
in~\eqref{eq:reduced}, write $P_{j,\lambda}(w)=\sum_{q=0}^d a_qw^q$.
Since $1+|w|\le6(1+|z|)$, we have
\[
 |P_{j,\lambda}(w)|
 \le \Bigl(\sum_{q=0}^d|a_q|\Bigr)(1+|w|)^d
 \le 6^d\Bigl(\sum_{q=0}^d|a_q|\Bigr)(1+|z|)^d.
\]
Also, $\Rea(\lambda w)\le\Rea(\lambda z)+5|\lambda|$. Hence, with
$C_{j,\lambda}=6^d e^{5|\lambda|}\sum_{q=0}^d|a_q|$,
\[
 |P_{j,\lambda}(w)e^{\lambda w}|
 \le C_{j,\lambda}
 \left((1+|z|)^{d/(n-j)}e^{\Rea(\lambda z)/(n-j)}\right)^{n-j}
 \le C_{j,\lambda}W(z)^{n-j}.
\]
Summing over $\lambda\in\Gamma_j$ and choosing $C_0>0$ at least as
large as each of the finitely many sums $\sum_{\lambda\in\Gamma_j}C_{j,\lambda}$
gives
\begin{equation}\label{eq:weight-local-coefficient-bound}
 |A_j(w)|\le C_0 W(z)^{n-j},
 \qquad |w-z|\le5,\quad 0\le j\le n-2.
\end{equation}
With $z$ fixed, use the local coordinate $\zeta=W(z)(w-z)$ and write
$w=z+\zeta/W(z)$. Then $g(\zeta)=f(w)$ is entire and satisfies
\[
 g^{(n)}+\sum_{j=0}^{n-2}
 \frac{A_j(z+\zeta/W(z))}{W(z)^{n-j}}g^{(j)}=0,
 \qquad \zeta\in\C.
\]
If $|\zeta|\le5W(z)$, then $|w-z|=|\zeta|/W(z)\le5$.
By~\eqref{eq:weight-local-coefficient-bound}, each rescaled coefficient
therefore satisfies $|A_j(w)|/W(z)^{n-j}\le C_0$ on this disk,
with the same $C_0$ for every centre $z$.
Apply~\eqref{eq:bounded-coefficient-disk-count} with $R=W(z)\ge1$.
The disk $|\zeta|\le W(z)$ corresponds to $|w-z|\le1$, and the
affine change preserves zero multiplicities. This proves
\eqref{eq:unit-disk-zero-bound}.

Let $D=\overline B(E,1)$ and denote the multiplicity of each zero $w$
of $f$ by $m_w$. For each centre $z$,
\[
 n(\overline B(z,1),1/f)
 =\sum_{w:\,f(w)=0}m_w\,\mathbf{1}_{\{|z-w|\le1\}}.
\]
The summands are nonnegative, so Tonelli's theorem
\cite[Corollary~1.4.46, p.~109]{TaoMeasure} gives
\[
 \begin{aligned}
 \int_D n(\overline B(z,1),1/f)\,\mathrm dA(z)
 &=\sum_{w:\,f(w)=0}m_w\int_D\mathbf{1}_{\{|z-w|\le1\}}\,\mathrm dA(z)\\
 &=\sum_{w:\,f(w)=0}m_w\,\operatorname{area}\bigl(D\cap\overline B(w,1)\bigr).
 \end{aligned}
\]
For $w\in E$ we have $\overline B(w,1)\subset D$, so each such zero
contributes $\pi m_w$. Keeping these terms and applying
\eqref{eq:unit-disk-zero-bound} gives
\[
 \pi n(E,1/f)
 \le\int_D n(\overline B(z,1),1/f)\,\mathrm dA(z)
 \le C\int_D W(z)\,\mathrm dA(z).
\]
Absorbing the factor $\pi$ into the constant proves
\eqref{eq:coefficient-integral-zero-bound}.
\end{proof}

We next determine the averaged logarithmic growth of the fundamental
row in moving rectangles near each main direction.
The row $\mathbf F=(f_1,\ldots,f_n)$ has the
initial matrix $\mathbf I_n$ in~\eqref{eq:fundamental-initial-matrix}.
For $\omega\in\Omega_*=\{\omega\in\Omega_a:\beta(\omega)=\beta_*\}$
as defined in~\eqref{eq:dominant-weighted-frequencies}, we use the phase map $\Phi$
from~\eqref{eq:phase}, with $\beta=\beta_*$, and the support
functions $H_\omega,h_\omega$
from~\eqref{eq:fundamental-row-support-function}.
The constant $P_\omega$ is the perimeter of
$K_\omega=\conv\{\overline\mu:p_\omega(\mu)=0\}$, with twice
the length used when $K_\omega$ is a segment.
The coordinates below describe rectangles near $z=ar/\omega$,
with bounded width along the main ray $\arg z=-\arg\omega$ and
width of order $\sqrt r$ perpendicular to it.
As this latter width increases, the rectangles meet more inverse
branches, so estimates uniform in the branch index are needed.

\begin{lemma}\label{lem:physical-window}
Fix $\omega\in\Omega_*$, and put $R=ar$, $Q=R^{\beta_*}e^R$ for $r>0$.
As $r\to\infty$, the functions
\[
 U_r(\xi,\eta)=Q^{-1}
 \log\left\|\mathbf F\left((R+\xi+i\sqrt R\eta)/\omega\right)\right\|,
 \qquad (\xi,\eta)\in\R^2,
\]
are uniformly bounded in absolute value on each compact subset of
$\R^2$. For every continuous compactly supported test function
$\phi$ on $\R^2$,
\begin{equation}\label{eq:physical-potential-average}
 \int_{\R^2}U_r(\xi,\eta)\phi(\xi,\eta)\,\dd\xi\dd\eta
 \longrightarrow
 \frac{P_\omega}{2\pi}
 \int_{\R^2}e^\xi\phi(\xi,\eta)\,\dd\xi\dd\eta,
 \qquad r\to\infty.
\end{equation}
\end{lemma}

The proof is given in
\hyperref[subsec:proof-physical-window]{Section~\ref*{subsec:proof-physical-window}}.

\begin{figure}[htbp]
\centering
\resizebox{0.80\linewidth}{!}{%
\input{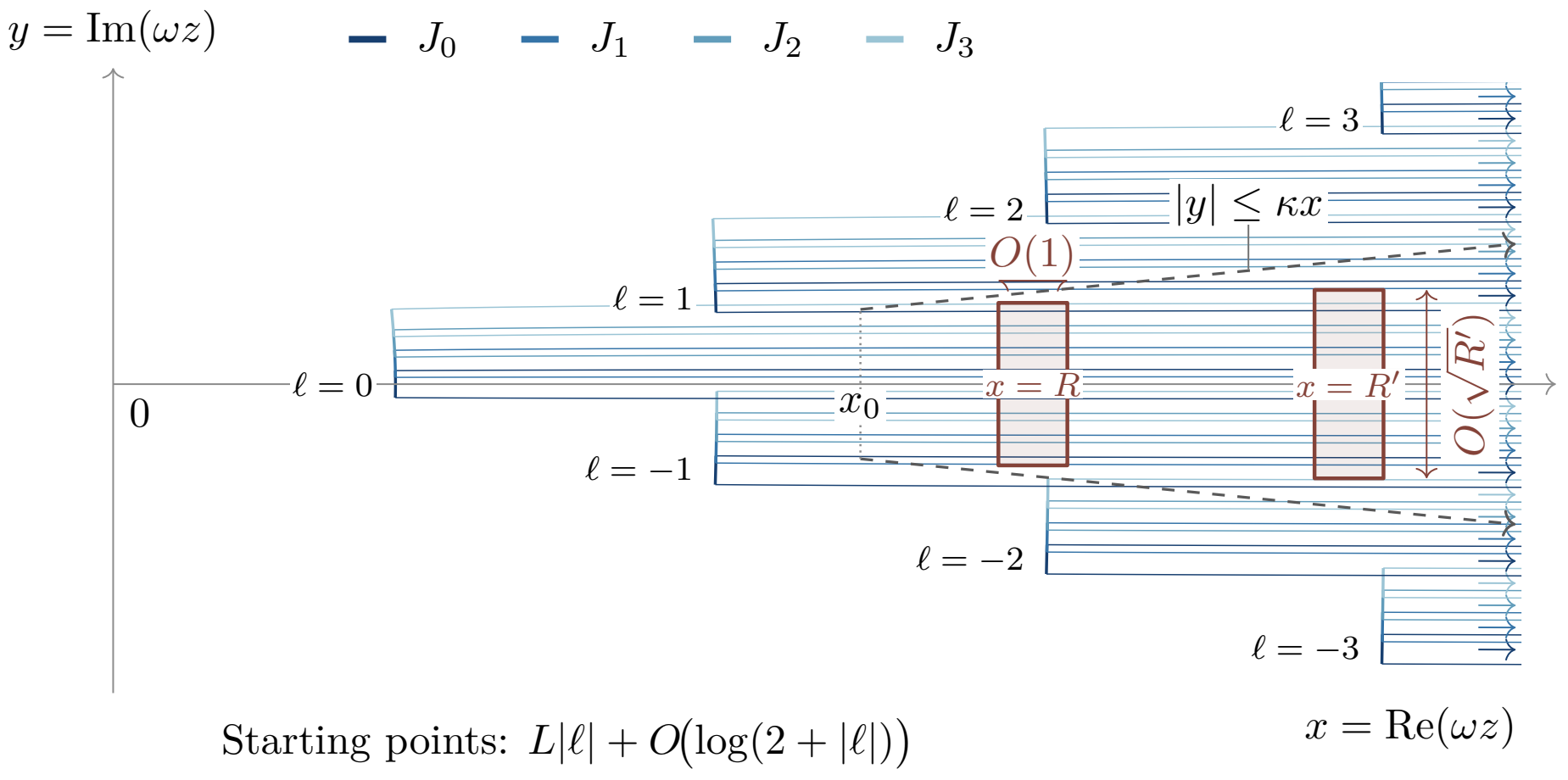}%
}
\caption{Strip tails used for uniform estimates.}
\label{fig:controlled-strip-windows}
\end{figure}

\enlargethispage{2\baselineskip}
In the coordinates $t=\omega z=x+iy=R+\xi+i\sqrt R\eta$, a fixed
square $|\xi|,|\eta|\le L$ corresponds to the rectangle
\begin{equation}\label{eq:physical-window-rectangle}
 |x-R|\le L,\qquad |y|\le L\sqrt R.
\end{equation}
Its points satisfy $|y|/x=O_L(R^{-1/2})\to0$, so for large $R$
it lies in the fixed region $x>0$, $|y|\le\kappa x$, with
$\kappa>0$ sufficiently small for the uniform estimates in
Lemma~\ref{lem:moving-inverse}.
Figure~\ref{fig:controlled-strip-windows} shows the rectangles
and the strip tails covered by these estimates.

For fixed $R$, the relation $y=\sqrt R\eta$ means that an increment
$2\pi/\sqrt R$ in $\eta$ increases $y$ by $2\pi$. The proof approximates $U_r(\xi,\eta)$ in
$L^1_{\mathrm{loc}}$ by $e^\xi h_\omega(\sqrt R\eta)$.
Since $h_\omega$ is $2\pi$-periodic, this approximating expression
has exact period $2\pi/\sqrt R$ in $\eta$, tending to zero as
$R\to\infty$. Its mean over each period is
\[
 \frac{e^\xi}{2\pi}\int_0^{2\pi}h_\omega(\theta)\dd\theta
 =\frac{P_\omega}{2\pi}e^\xi.
\]
This explains the perimeter factor in~\eqref{eq:physical-potential-average}.

\begin{remark}
Convergence in $L^1_{\mathrm{loc}}(\R^2)$ is stronger than the weak
convergence in~\eqref{eq:physical-potential-average}. Indeed, if
$U_r\to U$ in $L^1_{\mathrm{loc}}(\R^2)$, then every continuous
compactly supported test function $\phi$, with $K=\operatorname{supp}\phi$,
satisfies
\[
 \left|\int_{\R^2}(U_r-U)\phi\,\dd\xi\dd\eta\right|
 \le \|\phi\|_\infty\int_K|U_r-U|\,\dd\xi\dd\eta
 \longrightarrow0.
\]
Formula~\eqref{eq:physical-potential-average} allows continuous test
functions, whereas distributional convergence uses smooth ones.
Under the local uniform bounds of Lemma~\ref{lem:physical-window},
uniform approximation of continuous test functions by smooth ones
makes these two formulations equivalent.
\end{remark}

\begin{example}\label{ex:two-direction-weights}
For the equation $f^{(3)}+z^3(e^{3z}+8e^{3iz})f=0$ in the
introduction, \eqref{eq:characteristic} gives
$p_1(\mu)=\mu^3+1$ and $p_i(\mu)=\mu^3-8$.
The two root hulls are equilateral triangles of circumradii $1$ and $2$,
so their perimeters in~\eqref{eq:physical-potential-average} are
$P_1=3\sqrt3$ and $P_i=6\sqrt3$.
Hence $C_1=3\sqrt3/(2\pi\sqrt{2\pi})$ and $C_i=2C_1$.
Since $a=\beta_*=1$ and $\Omega_*=\{1,i\}$,
Theorem~\ref{thm:exact-global} gives, outside the exceptional cone,
\[
 n(r,1/f_{\mathbf c})\sim (C_1+C_i)r^{3/2}e^r
 =\frac{9\sqrt3}{2\pi\sqrt{2\pi}}r^{3/2}e^r.
\]
The two directions thus account for $1/3$ and $2/3$ of the zeros,
respectively.
\end{example}

\FloatBarrier

\subsection{Growth of solutions outside a complete pluripolar cone}

Lemma~\ref{lem:physical-window} gives the averaged growth of
$\log\|\mathbf F\|$. For a fixed initial vector $\mathbf c$,
the terms in $f_{\mathbf c}(z)=\sum_{j=1}^n c_jf_j(z)$ may
cancel one another, so $|f_{\mathbf c}(z)|$ can be small even
when $\|\mathbf F(z)\|$ is large.
The next lemma controls the difference between
$\log|f_{\mathbf c}|$ and $\log\|\mathbf F\|$ in normalized integrals.
The same complete pluripolar cone of initial vectors applies to
all prescribed sequences $(z_{m,\omega})_{m\ge1}$ and every compact
rectangle $K\subset\R^2$ in the $(\xi,\eta)$-plane.

Identifying $\C$ with $\R^2$, an \textit{invertible real affine map}
is an invertible real linear transformation followed by a translation.
For any positive sequence $r_m\to\infty$, the rectangle coordinates
of Lemma~\ref{lem:physical-window} give the real affine maps
\[
 z_{m,\omega}(\xi,\eta)
 =\omega^{-1}(R_m+\xi+i\sqrt{R_m}\eta),
 \qquad R_m=ar_m,\quad \omega\in\Omega_*.
\]
In the coordinates $t=\omega z$, the image of $[-L,L]^2$ is
the rectangle in~\eqref{eq:physical-window-rectangle} with $R=R_m$.

For the notation $\operatorname{PSH}(\C^n)$ and complete pluripolar sets,
see the \hyperref[par:pluripotential-definitions]{definitions in the
introduction on p.~\pageref*{par:pluripotential-definitions}}.

In the next lemma, let $\mathbf F=(f_1,\ldots,f_n)$ be a holomorphic
row on $\C$ whose components are linearly independent and have no
common zero. The fundamental row normalized
by~\eqref{eq:fundamental-initial-matrix} satisfies these assumptions.
Its fundamental matrix is invertible at every $z\in\C$, so its
first row $\mathbf F(z)$ cannot equal $\mathbf0$. Individual
components and linear combinations may still have zeros.

\begin{lemma}\label{lem:pluripolar-projection}
Let $\mathbf F$ be as above, and let $\Omega$ be a finite nonempty
index set. For $m\ge1$ and
$\omega\in\Omega$, let $z_{m,\omega}:\R^2\to\C$ be an invertible
real affine map. Suppose that $Q_m>0$ and
$\sum_{m=1}^\infty Q_m^{-1}<\infty$.
There are a function $\Psi\in\operatorname{PSH}(\C^n)$,
$\Psi\not\equiv-\infty$, and a constant $\kappa>0$ such that
\begin{equation}\label{eq:projection-homogeneity}
 \Psi(\lambda \mathbf c)=\Psi(\mathbf c)+\kappa\log|\lambda|,
 \qquad \mathbf c\in\C^n,\quad \lambda\in\C\setminus\{0\}.
\end{equation}
The set $E=\{\mathbf c:\Psi(\mathbf c)=-\infty\}$ is a complete pluripolar cone,
and, for every $\mathbf c\notin E$, every $\omega\in\Omega$ and every
compact rectangle $K\subset\R^2$,
\begin{equation}\label{eq:pluripolar-projection-loss}
 \frac1{Q_m}\int_K
 \left|\log\frac{|\mathbf F(z_{m,\omega}(\xi,\eta))\mathbf c|}
 {\|\mathbf F(z_{m,\omega}(\xi,\eta))\|}\right|
 \dd\xi\dd\eta\longrightarrow0
 \qquad(m\to\infty).
\end{equation}
\end{lemma}

\begin{proof}
We construct a plurisubharmonic function $\Psi$ on the initial-value
space $\C^n$, using $\{\Psi=-\infty\}$ as a common exceptional set.
Choose nonnegative functions $\phi_L\in C_c^\infty(\R^2)$ with
$\phi_L\ge1$ on $[-L,L]^2$ for $L\ge1$, and put
$M_L=\int_{\R^2}\phi_L\dd\xi\dd\eta$.
For $m\ge1$, $\omega\in\Omega$ and $L\ge1$, define
\begin{equation}\label{eq:projection-log-average}
 v_{m,\omega,L}(\mathbf c)=\frac1{Q_m}\int_{\R^2}\phi_L(\xi,\eta)
 \log\frac{|\mathbf F(z_{m,\omega}(\xi,\eta))\mathbf c|}
 {\|\mathbf F(z_{m,\omega}(\xi,\eta))\|}\dd\xi\dd\eta,
 \qquad \mathbf c\in\C^n.
\end{equation}
The denominator is positive because $\mathbf F$ has no common zero.
For $\mathbf c\ne\mathbf0$, linear independence makes $\mathbf F\mathbf c$ a nonzero
holomorphic function, so $\log|\mathbf F\mathbf c|$ is locally integrable
by the local decomposition~\eqref{eq:local-log-factorization}.
The invertibility of $z_{m,\omega}$ therefore makes
\eqref{eq:projection-log-average} finite for $\mathbf c\ne\mathbf0$.
At $\mathbf c=\mathbf0$ its value is $-\infty$.

Each $v_{m,\omega,L}$ is plurisubharmonic in $\mathbf c$.
For fixed $z$ and $\delta>0$, the map
$\mathbf c\mapsto(\mathbf F(z)\mathbf c,\sqrt\delta)\in\C^2$
is holomorphic: its first component is complex linear and its
second is constant. Its Euclidean norm is positive, and we put
\[
 u_{\delta,z}(\mathbf c)
 =\log\|(\mathbf F(z)\mathbf c,\sqrt\delta)\|
  -\log\|\mathbf F(z)\|
 =\frac12\log\frac{|\mathbf F(z)\mathbf c|^2+\delta}
                         {\|\mathbf F(z)\|^2}.
\]
To check plurisubharmonicity, we restrict to an arbitrary complex
line $\mathbf c=\mathbf c_0+\lambda\mathbf a$ in the initial-value
space, where $\lambda\in\C$ and $\mathbf a\ne\mathbf0$.
Put $g(\lambda)=\mathbf F(z)(\mathbf c_0+\lambda\mathbf a)$.
Since $g'=\mathbf F(z)\mathbf a$ and
$\partial_{\bar\lambda}g=0$, the identity
$\Delta_\lambda=4\partial_{\bar\lambda}\partial_\lambda$ gives
\[
 \Delta_\lambda u_{\delta,z}(\mathbf c_0+\lambda\mathbf a)
 =2\partial_{\bar\lambda}
   \left(\frac{g'(\lambda)\overline{g(\lambda)}}{|g(\lambda)|^2+\delta}\right)
 =\frac{2\delta|\mathbf F(z)\mathbf a|^2}
          {\bigl(|g(\lambda)|^2+\delta\bigr)^2}\ge0.
\]
Thus the restriction to every complex line is subharmonic, so
$u_{\delta,z}$ is plurisubharmonic.
Since $\phi_L\ge0$ and $Q_m>0$, its averages in
\eqref{eq:projection-log-average} are also plurisubharmonic by
integrating the submean inequality on complex lines
\cite[Chapter~I, (5.2), p.~39]{Demailly}.
For fixed $\delta$, compact support ensures that these averages
are continuous in $\mathbf c$.
As $\delta\to0$, these averages decrease to
\eqref{eq:projection-log-average}. The limit is plurisubharmonic
or identically $-\infty$
\cite[Chapter~I, Theorem~5.4, p.~39]{Demailly}; finiteness for
$\mathbf c\ne\mathbf0$ excludes the latter alternative.
Using the Cauchy--Schwarz estimate in~\eqref{eq:fundamental-row-norm}
and the nonnegative weight $\phi_L$ in~\eqref{eq:projection-log-average}
gives, for every $\mathbf c\in\C^n$,
\begin{equation}\label{eq:projection-upper-bound}
 v_{m,\omega,L}(\mathbf c)
 \le\frac1{Q_m}\int_{\R^2}\phi_L(\xi,\eta)\log^+\|\mathbf c\|
       \dd\xi\dd\eta
 =\frac{M_L}{Q_m}\log^+\|\mathbf c\|.
\end{equation}

Define
\begin{equation}\label{eq:projection-potential}
 \Psi(\mathbf c)=\sum_{\omega\in\Omega}\sum_{L,m\ge1}
                b_Lv_{m,\omega,L}(\mathbf c),
\end{equation}
where $b_L=2^{-L}/(1+M_L)$ for $L\ge1$.
By~\eqref{eq:projection-upper-bound}, $b_LM_L\le2^{-L}$ and
$\sum_mQ_m^{-1}<\infty$, the positive parts of this series are
summable for every $\mathbf c$. Thus $\Psi$ is well defined with
values in $[-\infty,\infty)$.
To prove $\Psi\not\equiv-\infty$, we integrate in $\mathbf c$ against the Gaussian
probability measure
$\dd\gamma_n(\mathbf c)=\pi^{-n}e^{-\|\mathbf c\|^2}\,\mathrm dV(\mathbf c)$ on $\C^n$.
Here $\mathrm dV$ is Lebesgue measure. For fixed $z$, choose a unitary
matrix $\mathbf U_z$ whose first row is the unit vector
$\mathbf F(z)/\|\mathbf F(z)\|$. With $\mathbf w=\mathbf U_z\mathbf c$
and $\mathbf c=\mathbf U_z^*\mathbf w$, we have
\[
 \frac{\mathbf F(z)\mathbf c}{\|\mathbf F(z)\|}
 =\frac{\mathbf F(z)\mathbf U_z^*\mathbf w}{\|\mathbf F(z)\|}
 =w_1.
\]
The standard complex Gaussian measure is invariant under unitary
transformations~\cite[\S2.1, Exercise~2.1.1(i), p.~13]{HKPV},
so this substitution preserves $\gamma_n$. Its independent coordinates
have product density
$\dd\gamma_n(\mathbf w)=\prod_{j=1}^n
(\pi^{-1}e^{-|w_j|^2}\,\mathrm dA(w_j))$
\cite[\S2.1, p.~13]{HKPV}.
Since $|\log|w_1||$ is nonnegative, Tonelli's theorem
\cite[Theorem~1.7.18, p.~202]{TaoMeasure} allows us to integrate
first in $w_2,\ldots,w_n$, each contributing a factor $1$.
Thus, for every $z\in\C$,
\begin{equation}\label{eq:projection-gaussian-identity}
 \begin{aligned}
 \int_{\C^n}\left|\log\frac{|\mathbf F(z)\mathbf c|}{\|\mathbf F(z)\|}\right|
       \dd\gamma_n(\mathbf c)
       &=\int_{\C^n}|\log|w_1||\dd\gamma_n(\mathbf w)\\
       &=\frac1\pi\int_{\C}|\log|w_1||e^{-|w_1|^2}\,\mathrm dA(w_1)=C_0,
 \end{aligned}
\end{equation}
where polar coordinates give
$C_0=2\int_0^\infty|\log r|r e^{-r^2}\dd r<\infty$,
independently of $z$.

Next apply the triangle inequality to~\eqref{eq:projection-log-average}.
Since $\phi_L\ge0$, the resulting integrand is nonnegative, so the
same Tonelli theorem permits exchanging the integrals in $\mathbf c$
and $(\xi,\eta)$. Writing $z=z_{m,\omega}(\xi,\eta)$ inside the integral,
\eqref{eq:projection-gaussian-identity} gives
\begin{equation}\label{eq:projection-gaussian-bound}
 \begin{aligned}
 &\int_{\C^n}|v_{m,\omega,L}(\mathbf c)|\dd\gamma_n(\mathbf c)\\
 &\quad\le\frac1{Q_m}\int_{\R^2}\phi_L(\xi,\eta)
 \left[\int_{\C^n}\left|\log\frac{|\mathbf F(z)\mathbf c|}{\|\mathbf F(z)\|}\right|
 \dd\gamma_n(\mathbf c)\right]\dd\xi\dd\eta
 =\frac{C_0M_L}{Q_m}.
 \end{aligned}
\end{equation}

We now use Tonelli's theorem for sums and integrals
\cite[Corollary~1.4.46, p.~109]{TaoMeasure} together with
\eqref{eq:projection-gaussian-bound} to obtain
\[
 \begin{aligned}
 \int_{\C^n}\sum_{\omega\in\Omega}\sum_{L,m\ge1}
 b_L|v_{m,\omega,L}(\mathbf c)|\dd\gamma_n(\mathbf c)
 &=
 \sum_{\omega\in\Omega}\sum_{L,m\ge1}b_L
 \int_{\C^n}|v_{m,\omega,L}(\mathbf c)|\dd\gamma_n(\mathbf c)\\
 &\le C_0|\Omega|\left(\sum_{L\ge1}b_LM_L\right)
                  \left(\sum_{m\ge1}Q_m^{-1}\right)<\infty.
 \end{aligned}
\]
Hence the nonnegative series inside the integral is finite for
$\gamma_n$-a.e.\ $\mathbf c$, so
\eqref{eq:projection-potential} converges absolutely there.
Since the Gaussian density is strictly positive, every open ball
has positive $\gamma_n$-measure and therefore contains a point where
this series converges absolutely. In particular, $\Psi\not\equiv-\infty$.

\bigskip

We next prove that $\Psi$ is plurisubharmonic on $\C^n$.
On each ball $\|\mathbf c\|<B$, with $B\ge1$,
\eqref{eq:projection-upper-bound} bounds every summand above by
$b_LM_LQ_m^{-1}\log B$, and these upper bounds have finite sum.
Define the partial sums
\[
 S_N(\mathbf c)=\sum_{\omega\in\Omega}\sum_{L,m=1}^{N}
 \left(b_Lv_{m,\omega,L}(\mathbf c)
       -\frac{b_LM_L}{Q_m}\log B\right).
\]
Each $S_N$ is plurisubharmonic, since finite sums of plurisubharmonic
functions remain plurisubharmonic
\cite[Chapter~I, Theorem~5.6, p.~39]{Demailly}.
Every summand is nonpositive on this ball, so
$S_{N+1}(\mathbf c)\le S_N(\mathbf c)$ there for every $N\ge1$.
Their limit is therefore plurisubharmonic or identically $-\infty$
\cite[Chapter~I, Theorem~5.4, p.~39]{Demailly}.
At an absolute convergence point in the ball, $S_N$ has a finite
limit because the subtracted constants have finite sum, excluding
the latter alternative. Adding back the finite sum of upper bounds
shows that~\eqref{eq:projection-potential} defines a nontrivial
plurisubharmonic function on $\C^n$.

\bigskip

We now verify~\eqref{eq:projection-homogeneity} and specify its constant
$\kappa$. For $\lambda\ne0$, linearity of $\mathbf F(z)\mathbf c$
in $\mathbf c$ and~\eqref{eq:projection-log-average} give
\[
 v_{m,\omega,L}(\lambda\mathbf c)
 =v_{m,\omega,L}(\mathbf c)
   +\frac{\log|\lambda|}{Q_m}\int_{\R^2}\phi_L\,\dd\xi\dd\eta
 =v_{m,\omega,L}(\mathbf c)+\frac{M_L}{Q_m}\log|\lambda|.
\]
Define
\begin{equation}\label{eq:projection-homogeneity-constant}
 \kappa:=\sum_{\omega\in\Omega}\sum_{L,m\ge1}\frac{b_LM_L}{Q_m}
 =|\Omega|\left(\sum_{L\ge1}\frac{2^{-L}M_L}{1+M_L}\right)
          \left(\sum_{m\ge1}Q_m^{-1}\right).
\end{equation}
Since $M_L>0$, $0<b_LM_L\le2^{-L}$ and $\Omega\ne\varnothing$,
we have $0<\kappa\le|\Omega|\sum_{m\ge1}Q_m^{-1}<\infty$.
Summing the scaling identity over finite index sets in
\eqref{eq:projection-potential} and passing to the limit proves
\eqref{eq:projection-homogeneity}. By~\eqref{eq:projection-homogeneity-constant},
the coefficients of $\log|\lambda|$ in these finite sums tend to
$\kappa<\infty$, so~\eqref{eq:projection-homogeneity} also holds when
$\Psi(\mathbf c)=-\infty$.
Set $E:=\{\Psi=-\infty\}$. Then
\[
 \mathbf c\notin E
 \quad\Longleftrightarrow\quad\Psi(\mathbf c)>-\infty,
 \qquad \mathbf c\in\C^n.
\]
Since $\Psi\in\operatorname{PSH}(\C^n)$ and $\Psi\not\equiv-\infty$,
$E$ is complete pluripolar by the
\hyperref[par:pluripotential-definitions]{definition in the introduction}.
Moreover, $\kappa<\infty$ and~\eqref{eq:projection-homogeneity} give
$\mathbf c\in E\iff\lambda\mathbf c\in E$ for $\lambda\ne0$.
Together with $\Psi(\mathbf0)=-\infty$, this shows that $E$ is a complex cone.

Fix $\mathbf c\notin E$. Recall that $x=x^+-x^-$ and
$|x|=x^++x^-$.
By~\eqref{eq:projection-upper-bound} and
\eqref{eq:projection-homogeneity-constant},
\[
 \sum_{\omega\in\Omega}\sum_{L,m\ge1}
 b_Lv_{m,\omega,L}(\mathbf c)^+
 \le\kappa\log^+\|\mathbf c\|<\infty.
\]
The negative parts are also summable: otherwise,
\eqref{eq:projection-potential} would give $\Psi(\mathbf c)=-\infty$,
contrary to $\mathbf c\notin E$. Thus $|x|=x^++x^-$ gives
\[
 \sum_{\omega\in\Omega}\sum_{L,m\ge1}
 b_L|v_{m,\omega,L}(\mathbf c)|<\infty.
\]
For each fixed $\omega$ and $L$, positivity of $b_L$ implies
$\sum_{m\ge1}|v_{m,\omega,L}(\mathbf c)|<\infty$, and hence
$v_{m,\omega,L}(\mathbf c)\to0$ as $m\to\infty$.
The positive part of the logarithm in~\eqref{eq:projection-log-average}
is at most $\log^+\|\mathbf c\|$.
Using $|x|=2x^+-x$ and $M_L=\int_{\R^2}\phi_L\dd\xi\dd\eta$, we obtain
\begin{equation}\label{eq:projection-weighted-loss}
 \begin{aligned}
 &\frac1{Q_m}\int_{\R^2}\phi_L(\xi,\eta)
 \left|\log\frac{|\mathbf F(z_{m,\omega}(\xi,\eta))\mathbf c|}
 {\|\mathbf F(z_{m,\omega}(\xi,\eta))\|}\right|\dd\xi\dd\eta\\
 &\qquad=-v_{m,\omega,L}(\mathbf c)
 +\frac2{Q_m}\int_{\R^2}\phi_L(\xi,\eta)
 \log^+\frac{|\mathbf F(z_{m,\omega}(\xi,\eta))\mathbf c|}
 {\|\mathbf F(z_{m,\omega}(\xi,\eta))\|}\dd\xi\dd\eta\\
 &\qquad\le-v_{m,\omega,L}(\mathbf c)
          +\frac{2M_L}{Q_m}\log^+\|\mathbf c\|
       \longrightarrow0.
 \end{aligned}
\end{equation}
Here $Q_m^{-1}\to0$ follows from its summability.
Every compact rectangle $K$ lies in $[-L,L]^2$ for some $L$,
where $\phi_L\ge1$. Thus~\eqref{eq:projection-weighted-loss}
proves~\eqref{eq:pluripolar-projection-loss} for every $K$,
with the same exceptional cone $E$ for all $\omega\in\Omega$.
\end{proof}

\begin{remark}
The Gaussian estimate~\eqref{eq:projection-gaussian-bound} and
Tonelli's theorem give
\[
 \int_{\C^n}\sum_{\omega\in\Omega}\sum_{L,m\ge1}
 b_L|v_{m,\omega,L}(\mathbf c)|\dd\gamma_n(\mathbf c)
 \le C_0\kappa<\infty.
\]
The nonnegative series is therefore finite $\gamma_n$-a.e.
Since the Gaussian density is strictly positive, this also holds
Lebesgue-a.e. By~\eqref{eq:projection-upper-bound},
the positive parts are summable everywhere, so $E=\{\Psi=-\infty\}$
is exactly the set where~\eqref{eq:projection-potential} fails to
converge absolutely. It has Lebesgue measure zero and contains
$\mathbf0$; it may also contain nonzero vectors, although each
individual $v_{m,\omega,L}$ is finite there.
\end{remark}

\clearpage
\section{Proofs of the auxiliary lemmas}\label{sec:auxiliary-proofs}

\subsection{Proof of \hyperref[lem:moving-inverse]{Lemma~\ref*{lem:moving-inverse}}}
\label{subsec:proof-moving-inverse}

\begin{proof}
Fix $\omega\in\Omega_a$ and put $\beta=\beta(\omega)$.
Recall from~\eqref{eq:phase} that $\Phi'(t)=t^\beta e^t$ for
$\Rea t>0$. For each frequency $\lambda\in\Gamma_j$ with
$\lambda\ne(n-j)\omega$, put
$\delta_{j,\lambda}=n-j-\Rea(\lambda/\omega)>0$
by~\eqref{eq:frequency-separation}.
Since there are only finitely many frequencies, we may choose
$\kappa>0$ such that, for every $0\le j\le n-2$ and
$\lambda\in\Gamma_j$ with $\lambda\ne(n-j)\omega$,
\begin{equation}\label{eq:inverse-kappa-choice}
 \kappa\left|\Ima(\lambda/\omega)\right|
 \le\frac{\delta_{j,\lambda}}2
 =\frac{n-j-\Rea(\lambda/\omega)}2.
\end{equation}
For $t=x+iy$ with $x>0$ and $|y|\le\kappa x$, this gives
$|e^{(\lambda/\omega-(n-j))t}|\le e^{-\delta_{j,\lambda}x/2}$
since $\Rea((\lambda/\omega-(n-j))t)
=-\delta_{j,\lambda}x-\Ima(\lambda/\omega)y$.

Integrating by parts twice in~\eqref{eq:phase}, along the path from $1$ to $1+iy$ and then to
$t=x+iy$, gives
\[
 \Phi(t)=e^t(t^\beta-\beta t^{\beta-1})
 +\beta(\beta-1)\int_1^t\xi^{\beta-2}e^\xi\dd\xi+O(1).
\]
On the vertical segment $\xi=1+iv$, with $v$ running from $0$ to
$y$, we have $|e^\xi|=e$. Thus the vertical remainder integral is
bounded by $e\int_0^{|y|}(1+v^2)^{(\beta-2)/2}\,\mathrm dv$,
which grows at most as a power of $1+|y|$.

On the horizontal segment $\xi=u+iy$, $1\le u\le x$, split the
integral at $u=x/2$, taking $x>2$. The first part is exponentially
smaller, since $e^u\le e^{x/2}$ and the other factor grows at most
as a power of $x$. On the second part, $|\xi|\asymp|t|$, so its
contribution is
$O(|t|^{\beta-2}\int_{x/2}^x e^u\,\mathrm du)
 =O(e^x|t|^{\beta-2})$.
Since $|t|\asymp x$ when $|y|\le\kappa x$, the vertical part and
the constant term are also absorbed by this bound. Hence, uniformly
as $x\to\infty$ with $|y|\le\kappa x$,
\begin{equation}\label{eq:phase-expansion-proof}
 \Phi(t)=\Phi'(t)\{1-\beta/t+O(t^{-2})\}.
\end{equation}

Put $h=\Phi/\Phi'$. Then $h=1-\beta/t+O(t^{-2})$.
Differentiating the quotient gives
$h'=1-(1+\beta/t)h=O(t^{-2})$.
For sufficiently large $x$ with $|y|\le\kappa x$, we have
$|h(t)-1|<1/2$. Thus the principal logarithm defines an analytic
function $\log h(t)$ on this region, with
$\log h(t)=O(t^{-1})\to0$ as $x\to\infty$.

To solve $\Phi(t)=s$, use $\Phi(t)=t^\beta e^t h(t)$ and take
logarithms. For the branch indexed by $\ell$, this leads, with
$W_\ell$ defined in~\eqref{eq:inverse-logarithmic-variables}, to
\[
 t+\beta\log t+\log h(t)=W_\ell.
\]
Here $\Rea W_\ell=X$.
Choose $0<\kappa_0<\kappa$ sufficiently small and then $L$
sufficiently large that $X>L(1+|\ell|)$ ensures
$|\Ima W_\ell|=|\arg_Js+2\pi\ell|\le\kappa_0X$.

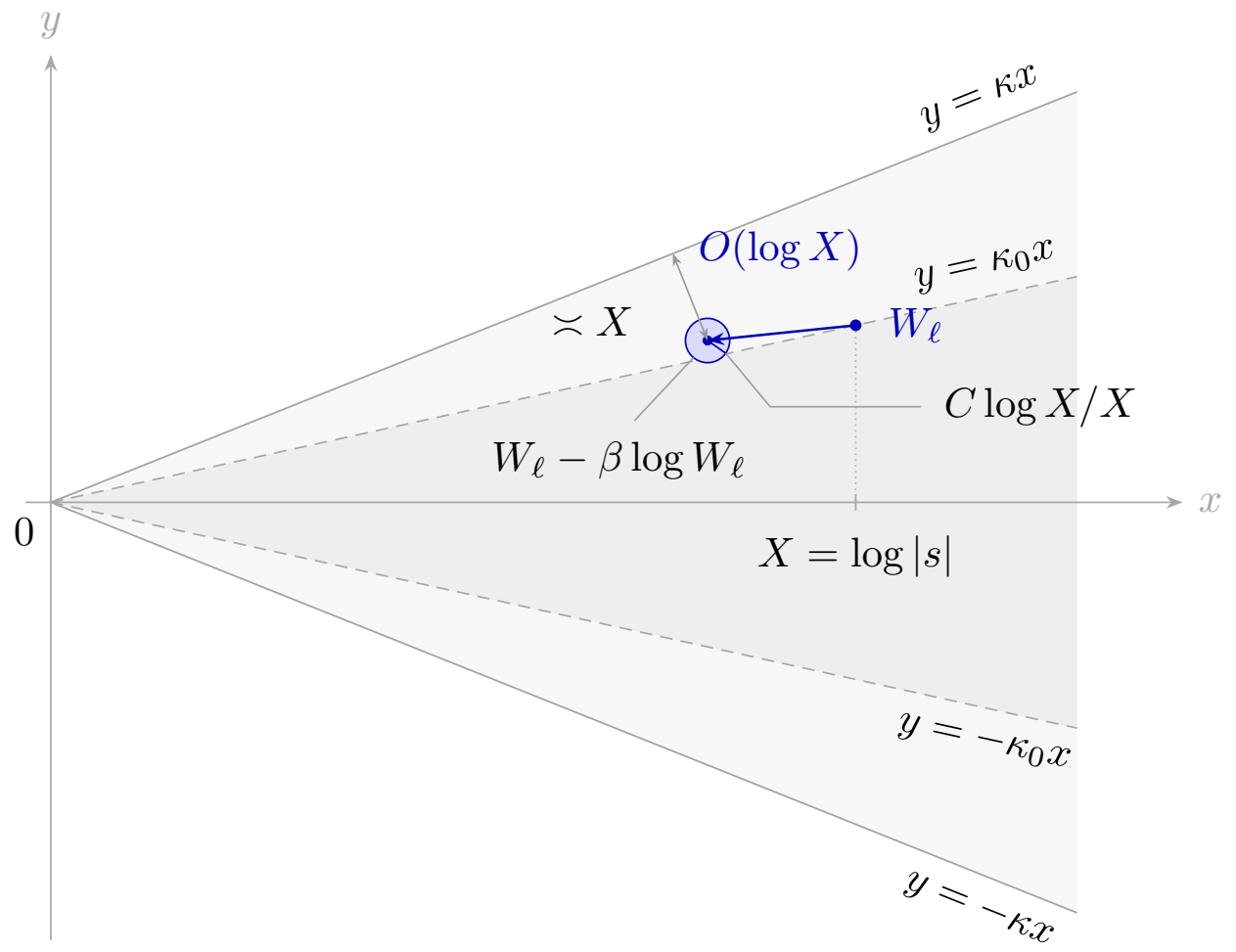
\begin{figure}[!htb]
\centering
\scalebox{0.8}{%
\begin{tikzpicture}[x=0.85cm,y=0.85cm,font=\small,
 axis/.style={thin,gray!70,-{Stealth[length=4pt]}},
 boundary/.style={thin,gray!75},
 inner/.style={thin,densely dashed,gray!70},
 measure/.style={thin,gray!80,{Stealth[length=3pt]}-{Stealth[length=3pt]}}]
 \fill[gray!6] (0,0)--(10.2,4.08)--(10.2,-4.08)--cycle;
 \fill[gray!13] (0,0)--(10.2,2.244)--(10.2,-2.244)--cycle;
 \draw[boundary] (0,0)--(10.2,4.08);
 \draw[boundary] (0,0)--(10.2,-4.08);
 \draw[inner] (0,0)--(10.2,2.244);
 \draw[inner] (0,0)--(10.2,-2.244);
 \draw[axis] (-0.25,0)--(11.25,0) node[right] {$x$};
 \draw[axis] (0,-4.35)--(0,4.45) node[above] {$y$};
 \node[below left] at (0,0) {$0$};
 \node[above,rotate=21.8] at (9.35,3.74) {$y=\kappa x$};
 \node[below,rotate=-21.8] at (9.35,-3.74) {$y=-\kappa x$};
 \node[above,rotate=12.4] at (9.35,2.057) {$y=\kappa_0x$};
 \node[below,rotate=-12.4] at (9.35,-2.057) {$y=-\kappa_0x$};

 \coordinate (W) at (8,1.76);
 \coordinate (centre) at (6.5278481,1.6084148);
 \coordinate (foot) at (6.1820811,2.4728324);
 \draw[thin,densely dotted,gray!65] (W)--(8,0);
 \draw[thin,gray!70] (8,0.08)--(8,-0.08);
 \node[below=4pt] at (8,0) {$X=\log|s|$};

 \filldraw[fill=blue!14,draw=blue!70!black,thin]
   (centre) circle[radius=0.22];
 \draw[blue!70!black,semithick,-{Stealth[length=4pt]}]
   (W)--(centre);
 \node[blue!70!black,above=3pt] at (7.25,2.02) {$O(\log X)$};
 \fill[blue!70!black] (W) circle[radius=1.4pt];
 \fill[blue!70!black] (centre) circle[radius=1.2pt];
 \node[blue!70!black,right=4pt] at (W) {$W_\ell$};

 \draw[measure] (centre)--(foot);
 \node[anchor=east] at (5.95,1.80) {$\asymp X$};
 \draw[thin,gray!75] (6.3878481,1.4384148)--(5.80,0.81);
 \node[anchor=north] at (5.65,0.76)
   {$W_\ell-\beta\log W_\ell$};
 \draw[thin,blue!70!black] (centre)--(6.7080481,1.4822148);
 \draw[thin,gray!75] (6.7080481,1.4822148)--(7.15,0.95)--(8.65,0.95);
 \node[anchor=west] at (8.72,0.95) {$C\log X/X$};
\end{tikzpicture}
}
\caption{The disk used in the inverse construction.}
\label{fig:inverse-disk}
\end{figure}

For fixed $s$ and $\ell$, consider the map
\begin{equation}\label{eq:inverse-contraction-map}
 G_{s,\ell}(t)=W_\ell-\beta\log t-\log h(t).
\end{equation}
We apply the contraction mapping theorem to construct $\psi_\ell(s)$
satisfying
\[
 \left|\psi_\ell(s)-(W_\ell-\beta\log W_\ell)\right|
 \le C\frac{\log X}{X},
\]
with $C$ independent of $\ell$. Write
$a_s=W_\ell-\beta\log W_\ell$ and $r_s=C\log X/X$;
we work on the closed disk $|t-a_s|\le r_s$.
On this disk, $x=X+O(\log X)$ and $y=\Ima W_\ell+O(1)$.
Since $\kappa_0<\kappa$, for large $X$ the disk lies in
$x>0$, $|y|\le\kappa x$, and $|t|\asymp X$ there; see
Figure~\ref{fig:inverse-disk}. The derivative of the map
in~\eqref{eq:inverse-contraction-map} is
$G'_{s,\ell}(t)=-\beta/t-h'(t)/h(t)=O(X^{-1})$.

At the centre $a_s$, the map in~\eqref{eq:inverse-contraction-map}
satisfies $|G_{s,\ell}(a_s)-a_s|=O(\log X/X)$, since
$\log a_s-\log W_\ell=O(\log X/X)$ and
$\log h(a_s)=O(X^{-1})$.
Let $K_0$ be a constant in this bound, so that
$|G_{s,\ell}(a_s)-a_s|\le K_0\log X/X$.
Since $a_s$ does not depend on $C$, we may choose $C\ge4K_0$;
then $|G_{s,\ell}(a_s)-a_s|\le r_s/4$ by the definition of $r_s$.
With this $C$ fixed, for large $X$ we have $|t|\ge X/2$,
$|h(t)|\ge1/2$ and $|h'(t)|\le M/|t|^2$, with $M$ independent
of $s$ and $\ell$. Hence
\[
 |G'_{s,\ell}(t)|
 \le\frac{|\beta|}{|t|}+\frac{|h'(t)|}{|h(t)|}
 \le\frac{2|\beta|}{X}+\frac{8M}{X^2}
 \le\frac{K_1}{X},
\]
where $K_1=2|\beta|+8M$.
Increase $L$ so that these estimates hold on all the prescribed
domains and $L>\max\{1,2K_1\}$.
Then $X>L(1+|\ell|)\ge L$ implies $K_1/X<1/2$.

For $|t-a_s|\le r_s$, the line segment
$a_s+u(t-a_s)$, $0\le u\le1$, stays in the disk.
Integrating the derivative along this segment gives
\[
 \begin{aligned}
 G_{s,\ell}(t)-G_{s,\ell}(a_s)
 &=(t-a_s)\int_0^1G'_{s,\ell}(a_s+u(t-a_s))\,\mathrm du,\\
 |G_{s,\ell}(t)-G_{s,\ell}(a_s)|
 &\le \frac{K_1}{X}|t-a_s|\le\tfrac12|t-a_s|.
 \end{aligned}
\]
The triangle inequality therefore yields
\[
 \begin{aligned}
 |G_{s,\ell}(t)-a_s|
 &\le |G_{s,\ell}(t)-G_{s,\ell}(a_s)|
      +|G_{s,\ell}(a_s)-a_s|\\
 &\le \tfrac12|t-a_s|+\tfrac14r_s\le \tfrac12r_s+\tfrac14r_s=\tfrac34r_s<r_s.
 \end{aligned}
\]
Thus the map in~\eqref{eq:inverse-contraction-map} sends the closed
disk into its interior. The same line-segment estimate for any
two points in the disk gives contraction factor $K_1/X<1/2$.
The contraction mapping theorem gives a unique fixed point
$t=\psi_\ell(s)$ satisfying
\begin{equation}\label{eq:inverse-fixed-point}
 G_{s,\ell}(\psi_\ell(s))=\psi_\ell(s),\qquad
 |\psi_\ell(s)-a_s|\le3r_s/4<r_s.
\end{equation}

The same constants work for every $s$ in the prescribed sector.
Fix $\ell$ and a point $s_0$ in this sector.
At the corresponding fixed point, the derivative with respect to
$t$ of the logarithmic equation is $1-G'_{s_0,\ell}(t)\ne0$.
The analytic implicit function theorem gives an analytic solution
for $s$ near $s_0$. Since this solution lies strictly inside the
disk at $s_0$, continuity of the solution, centre $a_s$ and radius
$r_s$ keeps it inside the corresponding disk for nearby $s$.
Uniqueness of the fixed point then identifies it with $\psi_\ell(s)$.
As $s_0$ was arbitrary, $\psi_\ell$ is analytic throughout the sector.
Exponentiating the logarithmic equation gives
$\Phi(\psi_\ell(s))=s$. This identity makes $\psi_\ell$ injective,
so it is a conformal inverse branch.
The radius of the disk gives
\begin{equation}\label{eq:moving-inverse}
 \psi_\ell(s)=W_\ell-\beta\log W_\ell
                   +O(\log X/X),
\end{equation}
uniformly in $\ell$ on the stated domains, where $\log W_\ell$
is the principal logarithm. For fixed $\ell$, recall that
$L_s=\log_Js$ and $W_\ell=L_s+2\pi i\ell$. Thus
\[
 \log W_\ell-\log L_s
 =\log\left(1+\frac{2\pi i\ell}{L_s}\right)
 =O(|L_s|^{-1}).
\]
Substitution into~\eqref{eq:moving-inverse} gives
\eqref{eq:inverse-proof}.

\smallskip
The argument in Lemma~\ref{lem:phase-coefficient-limit} treats a
fixed branch, and its error bounds may depend on the strip height.
We now estimate the transformed coefficients with constants
independent of $\ell$, so that zero counts over an increasing
number of strips can be combined. Put $t=\psi_\ell(s)$.
Since $\Rea W_\ell=X$ and $|\Ima W_\ell|\le\kappa_0X$, we have
$|W_\ell|\asymp X$ and $|\log W_\ell|=O(\log X)$.
Equation~\eqref{eq:moving-inverse} therefore gives
$t=W_\ell+O(\log X)$, so $|t|\asymp X$. Taking real parts also gives
\[
 x=X-\beta\log|W_\ell|+O(\log X/X)=X+O(\log X).
\]
The constants in these bounds do not depend on $\ell$, provided
$X>L(1+|\ell|)$ and $\arg s\in J$.
If $(n-j)\omega\in\Gamma_j$, write
$P_{j,(n-j)\omega}(z)=\sum_{d=0}^{d_j}p_{j,d}z^d$ for $z\in\C$,
where $p_{j,d_j}=a_j$.
Then~\eqref{eq:normalized-coefficients} and $\Phi'(t)=t^\beta e^t$ give
\[
 \begin{aligned}
 \frac{\omega^{j-n}P_{j,(n-j)\omega}(t/\omega)e^{(n-j)t}}
      {(\Phi'(t))^{n-j}}
 =\omega^{j-n}t^{-\beta(n-j)}P_{j,(n-j)\omega}(t/\omega)
 =\sum_{d=0}^{d_j}p_{j,d}\omega^{j-n-d}t^{d-\beta(n-j)}.
 \end{aligned}
\]
By~\eqref{eq:beta}, $d-\beta(n-j)\le d_j-\beta(n-j)\le0$
for $0\le d\le d_j$. If $j\in I_\omega$, then $d_j=\beta(n-j)$, and the term with
$d=d_j$ is exactly
\[
 p_{j,d_j}\omega^{j-n-d_j}t^0
 =a_j\omega^{j-n-d_j}=c_j
\]
by~\eqref{eq:limiting-coefficients}; the terms with $d<d_j$
have exponents $d-d_j<0$.
If $(n-j)\omega\in\Gamma_j$ but $j\notin I_\omega$, then
$d_j<\beta(n-j)$, so every term has negative exponent,
in agreement with $c_j=0$.
Only finitely many negative exponents occur, so they are all
at most $-\sigma$ for some $\sigma>0$.
Since $|t|\asymp X$, the sum of these terms is $O(X^{-\sigma})$.

For $\lambda\ne(n-j)\omega$, the exponential estimate above,
together with $x=X+O(\log X)$, bounds each term by
$O(X^Ne^{-\delta_0X})$. The constants $N,\delta_0>0$ can be
chosen independently of $\ell$, since the frequencies and
polynomial degrees form finite sets. Combining the two cases gives
\begin{equation}\label{eq:normalized-ratio-error}
 \frac{B_j(t)}{(\Phi'(t))^{n-j}}-c_j
 =O(X^{-\sigma})+O(X^Ne^{-\delta_0X}),\qquad 0\le j\le n-2.
\end{equation}
In particular, for $X$ sufficiently large,
\begin{equation}\label{eq:normalized-ratio-uniform-bound}
 \left|\frac{B_j(t)}{(\Phi'(t))^{n-j}}\right|\le B,
 \qquad 0\le j\le n-2,
\end{equation}
with $B>0$ independent of $s$ and $\ell$.

It remains to estimate the chain-rule terms
in~\eqref{eq:phase-coefficients}.
Since $\Phi''/\Phi'=1+\beta/t$, the recurrence
\eqref{eq:chain-coefficients} shows that
$C_{m,k}/(\Phi')^k$ is a polynomial in $1/t$ for
$1\le k\le m\le n$. Hence
\begin{equation}\label{eq:chain-coefficient-uniform-bound}
 |C_{m,k}(t)|\le A|\Phi'(t)|^k,\qquad 1\le k\le m\le n,
\end{equation}
for $t=x+iy$ with $x$ sufficiently large and $|y|\le\kappa x$,
where $A>0$ is independent of $s$ and $\ell$.
To express these bounds in terms of $s$, substitute $t=\psi_\ell(s)$
in~\eqref{eq:phase-expansion-proof}. Using $\Phi(\psi_\ell(s))=s$
and $|\psi_\ell(s)|\asymp X$, we obtain, for $X$ sufficiently large,
\begin{align}
 \frac{\Phi'(\psi_\ell(s))}{s}
 &=\left\{1-\frac{\beta}{\psi_\ell(s)}+O(X^{-2})\right\}^{-1}
 =1+O(X^{-1}),\notag\\
 |\Phi'(t)|&\ge\frac{|s|}{2},\qquad t=\psi_\ell(s).
 \label{eq:phase-derivative-lower-bound}
\end{align}
Combining~\eqref{eq:normalized-ratio-uniform-bound},
\eqref{eq:chain-coefficient-uniform-bound} and
\eqref{eq:phase-derivative-lower-bound}, we obtain
\begin{equation}\label{eq:chain-rule-error-bounds}
 \begin{aligned}
 \frac{C_{n,k}(t)}{(\Phi'(t))^n}
 &=O(|s|^{k-n}),\qquad 1\le k<n,\\
 \frac{B_j(t)C_{j,k}(t)}{(\Phi'(t))^n}
 &=O(|s|^{k-j}),\qquad 1\le k<j\le n-2.
 \end{aligned}
\end{equation}
Both bounds in~\eqref{eq:chain-rule-error-bounds} are $O(|s|^{-1})$.
For $k=0$, all correction terms vanish because $C_{m,0}=0$ for
$m\ge1$. For $k=n-1$,
$b_{n-1}(s)=C_{n,n-1}(t)/(\Phi'(t))^n=O(|s|^{-1})$,
while $c_{n-1}=0$.
Combining~\eqref{eq:normalized-ratio-error} and
\eqref{eq:chain-rule-error-bounds} with~\eqref{eq:phase-coefficients}
gives~\eqref{eq:moving-coefficient-limit}.
\end{proof}

\begin{remark}
For $\Phi(t)=e^t$, we have $\beta=0$ and $h=1$, so the map
in~\eqref{eq:inverse-contraction-map} reduces to $G_{s,\ell}(t)=W_\ell$.
For fixed $s$ and $\ell$, it sends the whole disk to its unique
fixed point $\psi_\ell(s)=W_\ell=\log_Js+2\pi i\ell$.
As $s$ varies over $S_J$, these points form a horizontal half-strip
in the $t$-plane. For general $\Phi$, the preceding contraction
argument constructs the inverse branches with estimates uniform
in $\ell$.
\end{remark}

\subsection{Proof of \hyperref[lem:branch-matrix-bounds]{Lemma~\ref*{lem:branch-matrix-bounds}}}
\label{subsec:proof-branch-matrix-bounds}

\begin{proof}
The first $n-1$ rows of $\mathbf A_\ell(s)$ and $\mathbf A_0$ coincide; their last
rows are $(-b_{0,\ell}(s),\ldots,-b_{n-1,\ell}(s))$ and
$(-c_0,\ldots,-c_{n-1})$, respectively.
Thus, for $\mathbf x=(x_1,\ldots,x_n)^{\mathsf T}\in\C^n$,
\begin{equation}\label{eq:companion-difference-action}
 (\mathbf A_\ell(s)-\mathbf A_0)\mathbf x
 =\left(0,\ldots,0,\sum_{k=0}^{n-1}
 (c_k-b_{k,\ell}(s))x_{k+1}\right)^{\mathsf T}.
\end{equation}
Taking the maximum over $\|\mathbf x\|=1$ in
\eqref{eq:companion-difference-action} and applying the same
Cauchy--Schwarz argument as in~\eqref{eq:fundamental-row-norm} gives
\[
 \|\mathbf A_\ell(s)-\mathbf A_0\|
 =\left(\sum_{k=0}^{n-1}|b_{k,\ell}(s)-c_k|^2\right)^{1/2}.
\]
By the uniform coefficient estimate~\eqref{eq:moving-coefficient-limit},
we may choose $L_\varepsilon>\max\{L,1\}$, where $L$ is the constant
in Lemma~\ref{lem:moving-inverse}, so that
\eqref{eq:uniform-companion-closeness} holds.

To prove~\eqref{eq:global-entrance}, put
$s=\rho_{\ell,\varepsilon}e^{i\theta}$ on the circle
$|s|=\rho_{\ell,\varepsilon}$, where $\theta\in J'$.
Then $X=\log|s|=L_\varepsilon(1+|\ell|)$, and
\eqref{eq:inverse-logarithmic-variables} gives
$W_\ell=X+i(\theta+2\pi\ell)$.
Since $\theta$ is bounded and $\Rea W_\ell=X$, we have
$|W_\ell|\asymp_\varepsilon(1+|\ell|)$ and
$|\log W_\ell|=O_\varepsilon(\log(2+|\ell|))$.
Since $X\ge L_\varepsilon>1$, we have $0\le\log X/X\le1$.
The remainder in~\eqref{eq:moving-inverse} is therefore bounded
independently of $\ell$ and $\theta$, and we obtain
\begin{equation}\label{eq:inverse-endpoint-bound}
 |\psi_\ell(s)|\le C_\varepsilon(1+|\ell|).
\end{equation}
Let $\mathbf C(z)$ be the companion matrix of~\eqref{eq:reduced}
in the variable $z$, and define the original fundamental matrix by
\begin{equation}\label{eq:original-fundamental-matrix}
 \mathbf H(z)=\bigl(f_j^{(k-1)}(z)\bigr)_{k,j=1}^n,\qquad z\in\C.
\end{equation}
Then $\mathbf H'=\mathbf C\mathbf H$ and
$\mathbf H(z_0)=\mathbf I_n$ by~\eqref{eq:fundamental-initial-matrix}.
Differentiating $\mathbf H^{-1}\mathbf H=\mathbf I_n$, and substituting
$\mathbf H'=\mathbf C\mathbf H$, we obtain
\[
 \begin{aligned}
 (\mathbf H^{-1})'
 =-\mathbf H^{-1}\mathbf H'\mathbf H^{-1}
 =-\mathbf H^{-1}(\mathbf C\mathbf H)\mathbf H^{-1}
 =-\mathbf H^{-1}\mathbf C.
 \end{aligned}
\]
The superdiagonal part of $\mathbf C(z)$ has norm $1$; the
remaining part has only its last row nonzero.
The triangle inequality and the row-norm calculation in
\eqref{eq:fundamental-row-norm}, together with the finite exponential
type of the coefficients $A_j$, give a constant $C>0$, depending
only on $n$ and the coefficients of~\eqref{eq:reduced}, such that
\begin{equation}\label{eq:original-companion-growth}
 \|\mathbf C(z)\|
 \le1+\left(\sum_{j=0}^{n-2}|A_j(z)|^2\right)^{1/2}
 \le Ce^{C|z|},\qquad z\in\C.
\end{equation}

Put $t=\psi_\ell(s)$ and $z_*=t/\omega$.
To estimate $\mathbf H(z_*)$ from its value $\mathbf I_n$ at $z_0$,
parametrize the straight segment by
$\gamma(\tau)=z_0+\tau(z_*-z_0)$, $0\le\tau\le1$.
Since $z_0$ and $\omega$ are fixed,~\eqref{eq:inverse-endpoint-bound}
shows that $\gamma$ has length $O_\varepsilon(1+|\ell|)$ and that
$|z|\le C_\varepsilon(1+|\ell|)$ for every $z\in\gamma$.
Thus~\eqref{eq:original-companion-growth} gives
$\|\mathbf C(z)\|\le C_\varepsilon e^{C_\varepsilon(1+|\ell|)}$
on $\gamma$.
Set $\mathbf U(\tau)=\mathbf H(\gamma(\tau))$, $0\le\tau\le1$.
Since $\gamma'(\tau)=z_*-z_0$, the chain rule and the equations
for $\mathbf H$ and $\mathbf H^{-1}$ give
\[
 \begin{aligned}
 \frac{\mathrm d\mathbf U}{\mathrm d\tau}
 &=(z_*-z_0)\mathbf C(\gamma(\tau))\mathbf U(\tau),
 &&\mathbf U(0)=\mathbf I_n,\\
 \frac{\mathrm d(\mathbf U^{-1})}{\mathrm d\tau}
 &=-(z_*-z_0)\mathbf U(\tau)^{-1}\mathbf C(\gamma(\tau)),
 &&\mathbf U(0)^{-1}=\mathbf I_n.
 \end{aligned}
\]
Integrating these equations and using $\|\mathbf I_n\|=1$ and
submultiplicativity gives
\[
 \|\mathbf U(\tau)^{\pm1}\|
 \le1+\int_0^\tau |z_*-z_0|\,\|\mathbf C(\gamma(u))\|
                  \|\mathbf U(u)^{\pm1}\|\dd u.
\]
Since $\mathbf U(1)=\mathbf H(z_*)$, the integral form of Gronwall's
inequality gives
\begin{equation}\label{eq:original-matrix-gronwall}
 \begin{aligned}
 \log^+\|\mathbf H(z_*)^{\pm1}\|
 &\le\int_\gamma\|\mathbf C(z)\|\,|\mathrm dz|\\
 &\le\operatorname{length}(\gamma)
       \sup_{z\in\gamma}\|\mathbf C(z)\|
 \le C_\varepsilon(1+|\ell|)e^{C_\varepsilon(1+|\ell|)}.
 \end{aligned}
\end{equation}
The factor $1+|\ell|$ can be absorbed into the exponential by
increasing $C_\varepsilon$.

We now relate the derivatives in $z$ to those in $s$.
For $1\le j\le n$, define $u_j(t)=f_j(t/\omega)$, $t\in\C$,
and $v_j(s)=u_j(\psi_\ell(s))$ on the domain of $\psi_\ell$.
Then $u_j(t)=v_j(\Phi(t))$ on the image of $\psi_\ell$.
Differentiating $u_j(t)=f_j(t/\omega)$ and applying
\eqref{eq:chain-rule-proof} to each pair $u_j,v_j$, then evaluating
at $t=\psi_\ell(s)=\omega z_*$, gives
\[
 f_j^{(m)}(z_*)=\omega^m u_j^{(m)}(t)
 =\sum_{k=0}^m\omega^m C_{m,k}(t)v_j^{(k)}(s),\qquad 0\le m<n.
\]
Define the $n\times n$ matrix
$\mathbf T(t)=\bigl(\omega^m C_{m,k}(t)\bigr)_{m,k=0}^{n-1}$.
By~\eqref{eq:chain-coefficients}, it has the lower triangular form
\[
 \mathbf T(t)=\begin{pmatrix}
 1&0&0&\cdots&0\\
 0&\omega\Phi'(t)&0&\cdots&0\\
 0&\omega^2\Phi''(t)&\omega^2(\Phi'(t))^2&\cdots&0\\
 \vdots&\vdots&\vdots&\ddots&\vdots\\
 0&\omega^{n-1}\Phi^{(n-1)}(t)&\omega^{n-1}C_{n-1,2}(t)
 &\cdots&(\omega\Phi'(t))^{n-1}
 \end{pmatrix}.
\]
Collecting these identities for all $j$, using the matrix definitions
in~\eqref{eq:original-fundamental-matrix} and
\eqref{eq:branch-fundamental-matrix}, gives
\begin{equation}\label{eq:derivative-change-matrix}
 \mathbf H(z_*)=\mathbf T(t)\mathbf X_\ell(s).
\end{equation}
 $\mathbf T(t)$ is invertible.
Also, recall that $s=\rho_{\ell,\varepsilon}e^{i\theta}$, \(W_\ell=X+i(\theta+2\pi\ell)\), and
$X=L_\varepsilon(1+|\ell|)$. Here $|W_\ell|\asymp X$ gives
$|\log W_\ell|=O(\log X)$, so~\eqref{eq:moving-inverse} yields
$$t-W_\ell=-\beta\log W_\ell+O(\log X/X)=O(\log X).$$
Hence $|t|\ge X/2\ge1$ after increasing $L_\varepsilon$ if necessary.
Since $\beta\ge0$, $\Rea t>0$ and~\eqref{eq:inverse-endpoint-bound}
gives $|t|\le C_\varepsilon(1+|\ell|)$, we have
\begin{equation}\label{eq:phase-derivative-endpoint-bounds}
 1\le|t^\beta e^t|=|\Phi'(t)|
 \le C_\varepsilon e^{C_\varepsilon(1+|\ell|)}.
\end{equation}

The bounds~\eqref{eq:chain-coefficient-uniform-bound}, together
with $C_{0,0}=1$ and $C_{m,0}=0$ for $m\ge1$, show that
$C_{m,k}(t)/(\Phi'(t))^k$ is uniformly bounded on these branches.
The corresponding matrix factorization is
\begin{equation}\label{eq:derivative-change-factorization}
 \begin{aligned}
 \mathbf T(t)
 &=\operatorname{diag}(1,\omega,\ldots,\omega^{n-1})
 \begin{pmatrix}
 1&0&0&\cdots&0\\
 0&1&0&\cdots&0\\
 0&\dfrac{C_{2,1}(t)}{\Phi'(t)}&1&\cdots&0\\
 \vdots&\vdots&\vdots&\ddots&\vdots\\
 0&\dfrac{C_{n-1,1}(t)}{\Phi'(t)}
 &\dfrac{C_{n-1,2}(t)}{(\Phi'(t))^2}&\cdots&1
 \end{pmatrix}\\
 &\quad{}\times\operatorname{diag}
 \bigl(1,\Phi'(t),\ldots,(\Phi'(t))^{n-1}\bigr).
 \end{aligned}
\end{equation}
The lower triangular matrix in~\eqref{eq:derivative-change-factorization}
has diagonal entries $1$ and uniformly bounded entries.
Since its determinant is $1$, the entries of its inverse are
cofactors, hence polynomials in the bounded matrix entries;
they are therefore uniformly bounded.
Consequently,~\eqref{eq:phase-derivative-endpoint-bounds} and the
factorization~\eqref{eq:derivative-change-factorization} give
\begin{equation}\label{eq:derivative-change-norm-bounds}
 \|\mathbf T(t)\|\le C|\Phi'(t)|^{n-1}
 \le C_\varepsilon e^{C_\varepsilon(1+|\ell|)},\qquad
 \|\mathbf T(t)^{-1}\|\le C.
\end{equation}
Thus the logarithms of both norms are bounded above by
$C_\varepsilon(1+|\ell|)$.
From~\eqref{eq:derivative-change-matrix},
$\mathbf X_\ell(s)=\mathbf T(t)^{-1}\mathbf H(z_*)$ and
$\mathbf X_\ell(s)^{-1}=\mathbf H(z_*)^{-1}\mathbf T(t)$.
The submultiplicative property of the matrix norm,
\eqref{eq:original-matrix-gronwall} and
\eqref{eq:derivative-change-norm-bounds} therefore give
\[
 \log^+\|\mathbf X_\ell(s)^{\pm1}\|
 \le\log^+\|\mathbf H(z_*)^{\pm1}\|
       +C_\varepsilon(1+|\ell|)
 \le C_\varepsilon e^{C_\varepsilon(1+|\ell|)}.
\]
This proves~\eqref{eq:global-entrance}.
\end{proof}

\subsection{Proof of \hyperref[lem:uniform-fundamental-row]{Lemma~\ref*{lem:uniform-fundamental-row}}}
\label{subsec:proof-uniform-fundamental-row}

\begin{proof}
By~\eqref{eq:branch-fundamental-matrix}, $\mathbf F_\ell(s)$ is the first
row of $\mathbf X_\ell(s)$.
The columns of $\mathbf X_\ell$ solve the system with coefficient matrix
$\mathbf A_\ell(s)$, whose limit is the companion matrix $\mathbf A_0$ in
\eqref{eq:phase-companion-system}.
We first consider a compact set $K$ avoiding the rays
$\Rea((\mu-\nu)\zeta)=0$ for distinct roots $\mu,\nu$ of $p_\omega$.
We will prove that, for every sufficiently small $\varepsilon>0$,
\begin{equation}\label{eq:quantified-row-growth}
 \sup_{\zeta\in K}
 \left|\frac{\log\|\mathbf F_\ell(T\zeta)\|}{T}-H_\omega(\zeta)\right|
 \le C_K\varepsilon+
       \frac{C_{\varepsilon,K}e^{C_{\varepsilon,K}(1+|\ell|)}}{T},
\end{equation}
provided $T\min_{\zeta\in K}|\zeta|\ge\rho_{\ell,\varepsilon}$.
Here the radii have the form specified in
Lemma~\ref{lem:branch-matrix-bounds}, on a closed interval
$J'\subset J$ containing the arguments of points in $K$.
The constants in~\eqref{eq:quantified-row-growth} are independent
of $T$ and $\ell$.

To handle repeated roots, we choose coordinates in which $\mathbf A_0$
differs from a diagonal matrix by an arbitrarily small fixed
matrix. Write the distinct roots of $p_\omega$ as
$\mu_1,\ldots,\mu_d$, with multiplicities $m_1,\ldots,m_d$,
so that $\sum_{j=1}^d m_j=n$.
For the column vector
$\mathbf v(\mu)=(1,\mu,\ldots,\mu^{n-1})^{\mathsf T}$,
the companion form gives
\[
 \mathbf A_0\mathbf v(\mu)=\mu\mathbf v(\mu)-p_\omega(\mu)\mathbf e_n,
 \qquad
 \mathbf w_{j,k}=\frac1{k!}
 \left.\frac{\mathrm d^k\mathbf v}{\mathrm d\mu^k}\right|_{\mu=\mu_j},
 \quad 0\le k<m_j,
\]
where $\mathbf e_n$ is the last standard basis vector.
Differentiating the first identity $k$ times at $\mu_j$, and
using the vanishing of $p_\omega$ and its first $m_j-1$
derivatives there, gives
$\mathbf A_0\mathbf w_{j,k}=\mu_j\mathbf w_{j,k}+\mathbf w_{j,k-1}$, with $\mathbf w_{j,-1}=0$.
The matrix with columns $\mathbf w_{j,k}$ is invertible.

Indeed, suppose a row vector $\mathbf a=(a_0,\ldots,a_{n-1})$ satisfies
$\mathbf a\mathbf w_{j,k}=0$ for all $j$ and $0\le k<m_j$.
The polynomial $P(\mu)=\mathbf a\mathbf v(\mu)=\sum_{r=0}^{n-1}a_r\mu^r$
then satisfies
\[
 \frac{P^{(k)}(\mu_j)}{k!}=\mathbf a\mathbf w_{j,k}=0,
 \qquad 0\le k<m_j.
\]
Thus each $\mu_j$ is a zero of $P$ of multiplicity at least $m_j$.
Since $\sum_jm_j=n$ and $P$ has degree less than $n$, we have
$P\equiv0$, and hence $\mathbf a=0$.

For $0<\tau<1$, rescale each column $\mathbf w_{j,k}$ by
$\tau^k$ to obtain the $n\times n$ matrix
\begin{equation}\label{eq:uniform-row-change-of-basis}
 \mathbf V_\tau=
 \begin{pmatrix}
 \mathbf w_{1,0}&\cdots&
 \tau^{m_1-1}\mathbf w_{1,m_1-1}&\cdots&
 \mathbf w_{d,0}&\cdots&
 \tau^{m_d-1}\mathbf w_{d,m_d-1}
 \end{pmatrix}.
\end{equation}
Each scaling factor is nonzero, so $\mathbf V_\tau$ is invertible. Then
\begin{equation}\label{eq:uniform-row-jordan-coordinates}
 \mathbf V_\tau^{-1}\mathbf A_0\mathbf V_\tau
 =\begin{pmatrix}
 \mu_1\mathbf I_{m_1}&0&\cdots&0\\
 0&\mu_2\mathbf I_{m_2}&\cdots&0\\
 \vdots&\vdots&\ddots&\vdots\\
 0&0&\cdots&\mu_d\mathbf I_{m_d}
 \end{pmatrix}+\tau \mathbf N,
\end{equation}
where
\[
 \mathbf N=\begin{pmatrix}
 \mathbf N_{m_1}&0&\cdots&0\\
 0&\mathbf N_{m_2}&\cdots&0\\
 \vdots&\vdots&\ddots&\vdots\\
 0&0&\cdots&\mathbf N_{m_d}
 \end{pmatrix},
 \qquad
 \mathbf N_m=\begin{pmatrix}
 0&1&&\\
 &0&\ddots&\\
 &&\ddots&1\\
 &&&0
 \end{pmatrix}\quad(m\ge2),
\]
with $\mathbf N_1=(0)$.
Each column of $\mathbf N$ contains at most one $1$ and is otherwise
zero, so its induced matrix norm satisfies $\|\mathbf N\|_1\le1$.
The first component of $\mathbf w_{j,0}=\mathbf v(\mu_j)$ is $1$.
For $1\le k<m_j$, the first component of $\mathbf w_{j,k}$ is
$\frac1{k!}\frac{\mathrm d^k}{\mathrm d\mu^k}1=0$.
Since $\tau^0=1$, the full first row of $\mathbf V_\tau$ in
\eqref{eq:uniform-row-change-of-basis} is therefore
\begin{equation}\label{eq:uniform-row-first-row}
 \left(
 \underbrace{1,0,\ldots,0}_{m_1\text{ columns}},
 \quad\ldots,\quad
 \underbrace{1,0,\ldots,0}_{m_d\text{ columns}}
 \right).
\end{equation}

Fix $0<\varepsilon<1$ and choose $\tau<\varepsilon/2$.
Since $\mathbf V_\tau$ is fixed, submultiplicativity and equivalence of
matrix norms allow us to apply Lemma~\ref{lem:branch-matrix-bounds}
with the matrix error in~\eqref{eq:uniform-companion-closeness}
small enough that
$\|\mathbf V_\tau^{-1}(\mathbf A_\ell(s)-\mathbf A_0)\mathbf V_\tau\|_1<\varepsilon/2$.
Writing the resulting constants as $L_\varepsilon,C_\varepsilon$,
set
\begin{equation}\label{eq:uniform-row-radius}
 \rho_0=\rho_{\ell,\varepsilon}
 =\exp(L_\varepsilon(1+|\ell|)).
\end{equation}
The inequality holds for $|s|\ge\rho_0$ and $\arg s\in J'$.

For a solution $f$ of~\eqref{eq:reduced}, let
$v(s)=f(\psi_\ell(s)/\omega)$; it satisfies
\eqref{eq:phase-equation} with coefficients $b_{k,\ell}$.
Set $\mathbf Y(s)=(v(s),v'(s),\ldots,v^{(n-1)}(s))^{\mathsf T}$.
Fix $\theta\in J'$. On the ray
$\{s=\rho e^{i\theta}:\rho\ge\rho_0\}$, define
\begin{equation}\label{eq:uniform-row-radial-coordinates}
 \mathbf u(\rho):=\mathbf V_\tau^{-1}\mathbf Y(\rho e^{i\theta}),\qquad \rho\ge\rho_0.
\end{equation}
The companion system~\eqref{eq:phase-companion-system}, the chain
rule and~\eqref{eq:uniform-row-jordan-coordinates} give
\begin{equation}\label{eq:uniform-row-radial-system}
 \begin{aligned}
 \frac{\mathrm d\mathbf u}{\mathrm d\rho}
 &=(\mathbf D_\theta+\mathbf B_{\ell,\theta}(\rho))\mathbf u,\qquad
 \mathbf D_\theta=\operatorname{diag}(e^{i\theta}\mu_j\mathbf I_{m_j}),\\
 \mathbf B_{\ell,\theta}(\rho)
 &=e^{i\theta}\left\{\tau \mathbf N+
 \mathbf V_\tau^{-1}(\mathbf A_\ell(\rho e^{i\theta})-\mathbf A_0)\mathbf V_\tau\right\}.
 \end{aligned}
\end{equation}
By the choice of $\tau$ and~\eqref{eq:uniform-row-radius}, we have
$\|\mathbf B_{\ell,\theta}(\rho)\|_1\le\varepsilon$ for all $\rho\ge\rho_0$.
For the initial vector $\mathbf c$ of $f$ at $z_0$,
\eqref{eq:branch-fundamental-matrix} gives
$\mathbf Y(s)=\mathbf X_\ell(s)\mathbf c$.
Hence~\eqref{eq:uniform-row-radial-coordinates} yields
\begin{equation}\label{eq:uniform-row-initial-coordinates}
 \mathbf u(\rho_0)
 =\mathbf V_\tau^{-1}\mathbf X_\ell(\rho_0e^{i\theta})\mathbf c.
\end{equation}
The initial matrix is $\mathbf V_\tau^{-1}\mathbf X_\ell(s)$,
with inverse $\mathbf X_\ell(s)^{-1}\mathbf V_\tau$, evaluated
at $s=\rho_0e^{i\theta}$. Submultiplicativity gives
\[
 \log^+\|\mathbf V_\tau^{-1}\mathbf X_\ell(s)\|
 \le\log^+\|\mathbf X_\ell(s)\|+\log^+\|\mathbf V_\tau^{-1}\|.
\]
The inverse matrix is estimated in the same way.
Since $\mathbf V_\tau^{\pm1}$ are fixed,~\eqref{eq:global-entrance}
and~\eqref{eq:uniform-row-radius} bound the logarithmic norms of both matrices by
$C_\varepsilon e^{C_\varepsilon(1+|\ell|)}$, after increasing
$C_\varepsilon$.

\bigskip
We now establish upper and lower bounds showing that the norm of
the row vector $\mathbf F_\ell$ of fundamental solutions defined in
\eqref{eq:branch-fundamental-row} has leading exponential growth
rate $h_\omega(\theta)$ along $s=\rho e^{i\theta}$ as
$\rho\to\infty$, with $\theta\in J'$ fixed.

We begin with the upper bound.
Recall from~\eqref{eq:fundamental-row-support-function} that
$h=h_\omega(\theta)=\max_j\Rea(e^{i\theta}\mu_j)$.
Since $\mathbf D_\theta$ is diagonal with diagonal entries of
real part at most $h$, and
$\|\mathbf B_{\ell,\theta}(\rho)\|_1\le\varepsilon$,
the standard norm estimate for linear systems
\cite[Ch.~II, \S1, p.~41; Ch.~III, Theorem~3,
pp.~58--59]{Coppel} yields
\begin{equation}\label{eq:uniform-row-coordinate-upper-bound}
 \|\mathbf u(\rho)\|_1
 \le e^{(h+\varepsilon)(\rho-\rho_0)}\|\mathbf u(\rho_0)\|_1,
 \qquad \rho\ge\rho_0.
\end{equation}
We now transfer this bound to the row $\mathbf F_\ell$.
Take an arbitrary global initial vector $\mathbf c$ with
$\|\mathbf c\|=1$, so that $v(s)=\mathbf F_\ell(s)\mathbf c$.
Using $\|\mathbf x\|\le\|\mathbf x\|_1\le\sqrt n\,\|\mathbf x\|$
and~\eqref{eq:uniform-row-initial-coordinates}, we have
\[
 \|\mathbf u(\rho_0)\|_1
 \le\sqrt n\,\|\mathbf V_\tau^{-1}\|\,
       \|\mathbf X_\ell(\rho_0e^{i\theta})\|.
\]
By~\eqref{eq:uniform-row-radial-coordinates},
$\mathbf Y(\rho e^{i\theta})=\mathbf V_\tau\mathbf u(\rho)$.
Since $v(\rho e^{i\theta})$ is the first component of this vector,
\eqref{eq:uniform-row-coordinate-upper-bound} gives
\[
 \begin{aligned}
 |\mathbf F_\ell(\rho e^{i\theta})\mathbf c|
 &\le\|\mathbf Y(\rho e^{i\theta})\|
 \le\|\mathbf V_\tau\|\,\|\mathbf u(\rho)\|_1\\
 &\le\sqrt n\,\|\mathbf V_\tau\|\|\mathbf V_\tau^{-1}\|\,
 e^{(h+\varepsilon)(\rho-\rho_0)}
 \|\mathbf X_\ell(\rho_0e^{i\theta})\|.
 \end{aligned}
\]
Taking the maximum over $\|\mathbf c\|=1$ by
\eqref{eq:fundamental-row-norm}, and then taking logarithms, yields
\[
 \log\|\mathbf F_\ell(\rho e^{i\theta})\|
 \le(h+\varepsilon)(\rho-\rho_0)
   +\log^+\|\mathbf X_\ell(\rho_0e^{i\theta})\|+C_\varepsilon.
\]
Here $C_\varepsilon$ 
is independent of $\ell,\theta,\rho$.
By~\eqref{eq:global-entrance} and~\eqref{eq:uniform-row-radius}, this implies
\begin{equation}\label{eq:uniform-row-upper-growth}
 \log\|\mathbf F_\ell(\rho e^{i\theta})\|
 \le\rho h+\varepsilon\rho
       +C_\varepsilon e^{C_\varepsilon(1+|\ell|)}.
\end{equation}

For the lower bound, fix $\rho\ge\rho_0$ and construct a solution with
$|v(\rho e^{i\theta})|\ge1/2$.
Since $K$ is compact and avoids the rays
$\Rea((\mu_j-\mu_k)\zeta)=0$ for $j\ne k$, we have
\[
 \delta:=\min_{\substack{\zeta\in K\\1\le j<k\le d}}
 \left|\Rea\bigl((\mu_j-\mu_k)e^{i\arg\zeta}\bigr)\right|>0.
\]
For each ray $\arg s=\theta=\arg\zeta$, $\zeta\in K$, there is a unique
index $j_\theta$ such that
\begin{equation}\label{eq:uniform-row-root-gap}
 \Rea(e^{i\theta}\mu_{j_\theta})=h_\omega(\theta),\qquad
 \Rea(e^{i\theta}\mu_k)\le h_\omega(\theta)-\delta
 \quad(k\ne j_\theta).
\end{equation}
After permuting the coordinates in~\eqref{eq:uniform-row-radial-system},
write $\mathbf u=\binom{\mathbf a}{\mathbf b}$, where
$\mathbf a\in\C^{m_{j_\theta}}$ consists of the coordinates
corresponding to $\mu_{j_\theta}$ and
$\mathbf b\in\C^{n-m_{j_\theta}}$ contains the remaining coordinates.
Choose
$\mathbf u(\rho_0)=\binom{\mathbf a_0}{\mathbf 0}$ with
$\mathbf a_0\ne\mathbf 0$.
By~\eqref{eq:uniform-row-initial-coordinates}, this choice corresponds
to a unique global initial vector $\mathbf c$, since both matrices
are invertible.
Put $a(\rho)=\|\mathbf a(\rho)\|_1$ and $b(\rho)=\|\mathbf b(\rho)\|_1$.
To derive their differential inequalities, let $m=m_{j_\theta}$,
write $\mathbf D_\theta=\operatorname{diag}(\lambda_1,\ldots,\lambda_n)$,
and put
$\mathbf r=(r_1,\ldots,r_n)^{\mathsf T}
=\mathbf B_{\ell,\theta}(\rho)\mathbf u$.
The component equations in~\eqref{eq:uniform-row-radial-system}
are $u_i'=\lambda_i u_i+r_i$, $1\le i\le n$.
By~\eqref{eq:uniform-row-root-gap},
$\Rea\lambda_i=h$ for $i\le m$ and
$\Rea\lambda_i\le h-\delta$ for $i>m$. Moreover,
\begin{equation}\label{eq:uniform-row-forcing-sum}
 \sum_{i=1}^n|r_i|
 =\|\mathbf B_{\ell,\theta}(\rho)\mathbf u\|_1
 \le\varepsilon\|\mathbf u\|_1=\varepsilon(a+b).
\end{equation}
At points where $u_i\ne0$, differentiating
$|u_i|^2=u_i\overline{u_i}$ gives
\[
 (|u_i|)'=\frac{\Rea(u_i'\overline{u_i})}{|u_i|}
 =\Rea\lambda_i\,|u_i|
   +\Rea\!\left(\frac{r_i\overline{u_i}}{|u_i|}\right).
\]
Consequently,
\[
 \Rea\lambda_i\,|u_i|-|r_i|
 \le (|u_i|)'
 \le\Rea\lambda_i\,|u_i|+|r_i|.
\]
The components $u_i$ are real analytic, so their moduli are locally
Lipschitz and hence locally absolutely continuous.
If $u_i=0$ at a point where $|u_i|$ is differentiable, then
$(|u_i|)'=0$, and the inequality above still holds.
Since $a=\sum_{i=1}^m|u_i|$ and $b=\sum_{i=m+1}^n|u_i|$,
summing the lower bounds for $i\le m$ and the upper bounds for
$i>m$, and using~\eqref{eq:uniform-row-forcing-sum}, gives a.e.
\begin{equation}\label{eq:uniform-row-block-norms}
 \begin{aligned}
 a'&\ge ha-\sum_{i=1}^m|r_i|\ge ha-\varepsilon(a+b),\\
 b'&\le(h-\delta)b+\sum_{i=m+1}^n|r_i|
      \le(h-\delta)b+\varepsilon(a+b).
 \end{aligned}
\end{equation}
For $\varepsilon<\delta/8$, subtracting half the first inequality
in~\eqref{eq:uniform-row-block-norms} from the second gives,
a.e.,
\begin{equation}\label{eq:uniform-row-difference-inequality}
 \begin{aligned}
 \left(b-\frac a2\right)'
 &\le(h-\delta)b-\frac h2a+\frac{3\varepsilon}{2}(a+b)\\
 &=\left(h-\delta+\frac{3\varepsilon}{2}\right)
       \left(b-\frac a2\right)
   +\left(-\frac\delta2+\frac{9\varepsilon}{4}\right)a\\
 &\le\left(h-\delta+\frac{3\varepsilon}{2}\right)
       \left(b-\frac a2\right),
 \end{aligned}
\end{equation}
where $a\ge0$.
Since $a,b$ are locally absolutely continuous, we may multiply
\eqref{eq:uniform-row-difference-inequality} by
$e^{-(h-\delta+3\varepsilon/2)(\rho-\rho_0)}$ and integrate from
$\rho_0$ to $\rho$. Using $b(\rho_0)=0$ and
$a(\rho_0)=\|\mathbf a_0\|_1>0$, we obtain
\[
 b(\rho)-\frac{a(\rho)}2
 \le-\frac{a(\rho_0)}2
 e^{(h-\delta+3\varepsilon/2)(\rho-\rho_0)}<0,
 \qquad \rho\ge\rho_0.
\]
Together with $b\ge0$, this proves $a>0$ and $b<a/2$ on the
whole half-line $\rho\ge\rho_0$.
The first inequality in~\eqref{eq:uniform-row-block-norms} now
gives $a'\ge(h-3\varepsilon/2)a$ a.e. Namely,
\[
\frac{d}{d\rho}
\left[e^{-(h-3\varepsilon/2)(\rho-\rho_0)}a(\rho)\right]
\ge0
\qquad\text{a.e.},
\]
and hence, by absolute continuity and integration over $[\rho_0,\rho]$,
\begin{equation}\label{eq:uniform-row-leading-norm-growth}
 a(\rho)\ge e^{(h-3\varepsilon/2)(\rho-\rho_0)}\|\mathbf a_0\|_1,
 \qquad \rho\ge\rho_0.
\end{equation}
Recall that $m=m_{j_\theta}$.
Let $\mathbf a(\rho;\mathbf a_0)$ denote the first $m$ coordinates of the solution
of~\eqref{eq:uniform-row-radial-system} with initial values
$\mathbf a(\rho_0)=\mathbf a_0$ and $\mathbf b(\rho_0)=\mathbf0$.
For each fixed $\rho\ge\rho_0$, linearity of the system defines
an $m\times m$ matrix $\mathbf Z(\rho)$ by
\[
 \mathbf Z(\rho)\mathbf a_0=\mathbf a(\rho;\mathbf a_0),
 \qquad \mathbf a_0\in\C^m.
\]
By~\eqref{eq:uniform-row-leading-norm-growth},
$\mathbf a_0\ne\mathbf0$ implies $\mathbf Z(\rho)\mathbf a_0\ne\mathbf0$,
so $\mathbf Z(\rho)$ is injective.
Its domain and codomain both have dimension $m$; hence it is
surjective and therefore invertible.
The lower bound~\eqref{eq:uniform-row-leading-norm-growth}
also gives
\[
 \|\mathbf Z(\rho)^{-1}\|_1
 \le e^{-(h-3\varepsilon/2)(\rho-\rho_0)}.
\]

Let
$\mathbf e_0=(1,0,\ldots,0)^{\mathsf T}\in\C^m$ and take
$\mathbf a_0=\mathbf Z(\rho)^{-1}\mathbf e_0$.
Then $\mathbf a(\rho)=\mathbf e_0$ and $$\|\mathbf b(\rho)\|_1=b(\rho)<\frac{a(\rho)}{2}=\frac{\|\mathbf a(\rho)\|_1}{2}=\frac{1}{2}.$$
For $j\ne j_\theta$, let $b_{j,1}(\rho)$ denote the first coordinate
of the group in $\mathbf b(\rho)$ corresponding to $\mu_j$.
By the first row of $\mathbf V_\tau$ in~\eqref{eq:uniform-row-first-row},
the first component of $\mathbf Y=\mathbf V_\tau\mathbf u$ is
\[
 Y_1(\rho e^{i\theta})
 =a_1(\rho)+\sum_{j\ne j_\theta}b_{j,1}(\rho)
 =1+\sum_{j\ne j_\theta}b_{j,1}(\rho).
\]
The reverse triangle inequality gives
\[
 \begin{aligned}
 |Y_1(\rho e^{i\theta})|
 \ge1-\left|\sum_{j\ne j_\theta}b_{j,1}(\rho)\right|
 \ge1-\sum_{j\ne j_\theta}|b_{j,1}(\rho)|
 \ge1-\|\mathbf b(\rho)\|_1>\frac12.
 \end{aligned}
\]
Recall that we chose
$\mathbf u(\rho_0)=\binom{\mathbf a_0}{\mathbf0}$.
By~\eqref{eq:uniform-row-initial-coordinates}, the unique global initial
vector corresponding to this choice is
$\mathbf c=\mathbf X_\ell(\rho_0e^{i\theta})^{-1}\mathbf V_\tau\binom{\mathbf a_0}{\mathbf 0}$,
so
\[
 \|\mathbf c\|\le
 \|\mathbf X_\ell(\rho_0e^{i\theta})^{-1}\|\,\|\mathbf V_\tau\|
 e^{-(h-3\varepsilon/2)(\rho-\rho_0)}.
\]
Recall from~\eqref{eq:branch-fundamental-row} that
$v(s)=f_{\mathbf c}(\psi_\ell(s)/\omega)=\mathbf F_\ell(s)\mathbf c$,
and $Y_1(s)=v(s)$ by~\eqref{eq:phase-derivative-vector}.
For each $\rho\ge\rho_0$, write $\mathbf c=\mathbf c_\rho\ne\mathbf0$.
Applying the row-norm identity~\eqref{eq:fundamental-row-norm} to the unit vector
$\mathbf c_\rho/\|\mathbf c_\rho\|$ gives
\begin{equation}\label{eq:uniform-row-lower-growth}
 \|\mathbf F_\ell(\rho e^{i\theta})\|
 \ge\frac{|\mathbf F_\ell(\rho e^{i\theta})\mathbf c_\rho|}{\|\mathbf c_\rho\|}
 \ge\frac{e^{(h-3\varepsilon/2)(\rho-\rho_0)}}
 {2\|\mathbf X_\ell(\rho_0e^{i\theta})^{-1}\|\,\|\mathbf V_\tau\|}.
\end{equation}
If there is only one distinct root, take $\mathbf a=\mathbf u$ with no $\mathbf b$ block;
the same lower estimate follows directly from
$a'\ge(h-\varepsilon)a$.
By~\eqref{eq:global-entrance} and~\eqref{eq:uniform-row-radius},
\eqref{eq:uniform-row-lower-growth} gives
\[
 \log\|\mathbf F_\ell(\rho e^{i\theta})\|
 \ge(h-3\varepsilon/2)\rho
      -C_\varepsilon e^{C_\varepsilon(1+|\ell|)}.
\]
Combining this lower bound with~\eqref{eq:uniform-row-upper-growth}
and enlarging the constants, we obtain
\[
 \left|\log\|\mathbf F_\ell(\rho e^{i\theta})\|
       -\rho h_\omega(\theta)\right|
 \le\frac32\varepsilon\rho
      +C_{\varepsilon,K}e^{C_{\varepsilon,K}(1+|\ell|)}.
\]
For $\zeta\in K$, put $\theta=\arg\zeta$ and $\rho=T|\zeta|$.
Then $\rho e^{i\theta}=T\zeta$ and
$\rho h_\omega(\theta)/T=H_\omega(\zeta)$
by~\eqref{eq:fundamental-row-support-function}.
Dividing by $T$ and taking the supremum over $K$ proves
\eqref{eq:quantified-row-growth}, with
$C_K=\frac32\max_{\zeta\in K}|\zeta|$.
Since $m(T)=o(\log T)$ and $\min_{\zeta\in K}|\zeta|>0$,
\eqref{eq:uniform-row-radius} implies that, for fixed $\varepsilon>0$
and all sufficiently large $T$, we have
$$T\min_{\zeta\in K}|\zeta|\ge\rho_{\ell,\varepsilon}=\rho_0$$
for every $|\ell|\le m(T)$.
Thus $\rho=T|\zeta|\ge\rho_0$ for all $\zeta\in K$ and $|\ell|\le m(T)$.

\bigskip
We now obtain bounds at every angle and pass to $L^1$ convergence.
Let $K$ be an arbitrary compact set as in the statement and choose
a closed interval $J'\subset J$ containing $\arg K$.
Apply Lemma~\ref{lem:branch-matrix-bounds} with $\varepsilon=1$.
By~\eqref{eq:uniform-companion-closeness},
$\|\mathbf A_\ell(s)\|\le\|\mathbf A_0\|+1$ whenever
$|s|\ge\rho_{\ell,1}$ and $\arg s\in J'$.
For all sufficiently large $T$,
$T\min_{\zeta\in K}|\zeta|\ge\rho_{\ell,1}$ for every $|\ell|\le m(T)$.
Applying Gronwall along the ray $\arg s=\arg\zeta$, from
$s=\rho_{\ell,1}e^{i\arg\zeta}$ to $s=T\zeta$, and
using~\eqref{eq:global-entrance}, gives
\begin{equation}\label{eq:uniform-row-coarse-matrix-growth}
 \log^+\|\mathbf X_\ell(T\zeta)^{\pm1}\|
 \le C_KT+C_Ke^{C_K(1+|\ell|)}.
\end{equation}
Recall that $\mathbf F_\ell$ is the first row of $\mathbf X_\ell$, so
$\mathbf F_\ell \mathbf X_\ell^{-1}=\mathbf e_1^{\mathsf T}$.
Consequently $\|\mathbf F_\ell\|\le\|\mathbf X_\ell\|$ and
$1\le\|\mathbf F_\ell\|\,\|\mathbf X_\ell^{-1}\|$.
Together with~\eqref{eq:uniform-row-coarse-matrix-growth}, these
inequalities give
\begin{equation}\label{eq:uniform-row-coarse-potential-bound}
 \left|\frac1T\log\|\mathbf F_\ell(T\zeta)\|\right|
 \le C_K+\frac{C_Ke^{C_K(1+|\ell|)}}T.
\end{equation}
For every fixed constant $C$, the condition $m(T)=o(\log T)$ gives
\begin{equation}\label{eq:uniform-row-branch-remainder}
 \sup_{|\ell|\le m(T)}\frac{e^{C(1+|\ell|)}}T
 \le\exp\{C(1+m(T))-\log T\}\longrightarrow0.
\end{equation}
Thus the normalized logarithms are uniformly bounded above and
below on $K$.

For $\alpha>0$, remove from $K$ the angular $\alpha$-neighbourhoods
of the finitely many rays in the $\zeta$-plane satisfying
$\Rea((\mu_j-\mu_k)\zeta)=0$, $1\le j<k\le d$,
where $\mu_1,\ldots,\mu_d$ are the distinct roots of $p_\omega$.
Write $K_\alpha$ for the remaining compact set.
The removed area is $O_K(\alpha)$.
On $K_\alpha$, the largest of the values
$\Rea(\mu_j e^{i\arg\zeta})$ is separated from all the others
by a constant $\delta_\alpha>0$, as in~\eqref{eq:uniform-row-root-gap}.
Fix $0<\varepsilon<\min\{1,\delta_\alpha/8\}$ and apply
\eqref{eq:quantified-row-growth} with $K_\alpha$ in place of $K$.
For these fixed $\alpha,\varepsilon$,~\eqref{eq:uniform-row-branch-remainder}
gives
\[
 \sup_{|\ell|\le m(T)}
 \frac{C_{\varepsilon,K_\alpha}
 e^{C_{\varepsilon,K_\alpha}(1+|\ell|)}}{T}
 \longrightarrow0\qquad(T\to\infty).
\]
Since the $L^1$ norm is the integral of the absolute error,
\eqref{eq:quantified-row-growth} gives
\[
 \begin{aligned}
 \left\|T^{-1}\log\|\mathbf F_\ell(T\zeta)\|-H_\omega(\zeta)
 \right\|_{L^1(K_\alpha)}
\le\operatorname{area}(K_\alpha)
 \left(C_{K_\alpha}\varepsilon+
 \frac{C_{\varepsilon,K_\alpha}
 e^{C_{\varepsilon,K_\alpha}(1+|\ell|)}}{T}\right).
 \end{aligned}
\]
Taking the supremum over $|\ell|\le m(T)$ and then the limsup as
$T\to\infty$, and absorbing $\operatorname{area}(K_\alpha)$ into the
constant, we obtain
\[
 \limsup_{T\to\infty}\sup_{|\ell|\le m(T)}
 \left\|T^{-1}\log\|\mathbf F_\ell(T\zeta)\|-H_\omega(\zeta)
 \right\|_{L^1(K_\alpha)}
 \le C_{K_\alpha}\varepsilon.
\]
By~\eqref{eq:uniform-row-coarse-potential-bound} and
\eqref{eq:uniform-row-branch-remainder}, the normalized logarithm
is uniformly bounded on all of $K$ for sufficiently large $T$
and $|\ell|\le m(T)$.
By~\eqref{eq:fundamental-row-support-function},
$|H_\omega(\zeta)|\le|\zeta|\max_{p_\omega(\mu)=0}|\mu|$,
so $H_\omega$ is bounded on the compact set $K$.
Their absolute difference is therefore at most $C_K$,
independently of $\alpha,\varepsilon$.
Consequently,
\[
 \left\|T^{-1}\log\|\mathbf F_\ell(T\zeta)\|-H_\omega(\zeta)
 \right\|_{L^1(K\setminus K_\alpha)}
 \le C_K\operatorname{area}(K\setminus K_\alpha)
 \le C_K\alpha.
\]
Hence the limsup in~\eqref{eq:uniform-fundamental-potential} is at most
$C_K\alpha+C_{K_\alpha}\varepsilon$.
Letting first $\varepsilon\to0$ and then
$\alpha\to0$ proves~\eqref{eq:uniform-fundamental-potential}.
If $K$ avoids these rays, apply~\eqref{eq:quantified-row-growth}
directly on $K$.
For each fixed sufficiently small $\varepsilon>0$, its second term
tends to zero uniformly for $|\ell|\le m(T)$ as $T\to\infty$,
by~\eqref{eq:uniform-row-branch-remainder}.
Letting $\varepsilon\to0$ then proves uniform convergence on $K$.
\end{proof}

\subsection{Proof of \hyperref[lem:slower-growth-subspace]{Lemma~\ref*{lem:slower-growth-subspace}}}
\label{subsec:proof-slower-growth-subspace}

\begin{proof}
Set $\mathbf E(r)=\mathbf B(r)-\mathbf B_0$.
Put $M=\sup_{r\ge r_*}\|\mathbf B(r)\|<\infty$ and choose $\kappa>M$.
Gronwall's inequality gives
\[
 \|\mathbf x(r)\|\le\|\mathbf x(r_*)\|e^{M(r-r_*)},\qquad r\ge r_*.
\]
Thus $e^{-\kappa(r-r_*)}\mathbf x(r)$ is bounded.
Applying Coppel's theorem
\cite[Ch.~IV, Theorem~5, p.~97]{Coppel} as in the proof of
Lemma~\ref{lem:radial-system-growth}, and restoring the exponential
factor, gives
\begin{equation}\label{eq:continuous-system-growth-values}
 \lambda(\mathbf x)\in
 \{\Rea\mu:\det(\mu \mathbf I-\mathbf B_0)=0\}
\end{equation}
for every nonzero solution; uniqueness excludes the eventually
zero alternative.

If all eigenvalues of $\mathbf B_0$ have the same real part $\sigma$,
\eqref{eq:continuous-system-growth-values} already gives
$\lambda(\mathbf x)=\sigma$ for every nonzero solution. In that case
the subspace of initial values for slower growth is $\{\mathbf 0\}$.

Under the hypotheses of the lemma, the unique eigenvalue with maximal
real part has multiplicity $m<n$, so the remaining eigenvalues have
a largest real part $\sigma_2<\sigma$. In Theorem~\ref{thm:generic-density}
this gap holds on each of the two chosen rays because one root has
strictly larger real projection than every other distinct root. Choose
$\sigma_2<\gamma<\sigma$ and put $\mathbf D=\mathbf B_0-\gamma \mathbf I$.
The generalized eigenspaces of $\mathbf D$ split according to the signs
of the real parts. Let $\mathbf P_s,\mathbf P_u$ be the corresponding spectral
projections, of ranks $n-m$ and $m$, respectively. Choose
$0<\delta<\min\{\gamma-\sigma_2,\sigma-\gamma\}$.
For a Jordan block $\lambda\mathbf I+\mathbf N$ of $\mathbf D$
of size $k$, one has $\mathbf N^k=0$ and
\[
 e^{t(\lambda\mathbf I+\mathbf N)}
 =e^{\lambda t}\sum_{j=0}^{k-1}\frac{t^j\mathbf N^j}{j!}.
\]
For a block of size \(k\), the finite sum is a matrix polynomial in \(t\) of degree at most \(k-1\le n-1\).
We call the blocks with $\Rea\lambda<0$ \textit{stable} and those
with $\Rea\lambda>0$ \textit{unstable}. On the stable blocks,
$\Rea\lambda\le\sigma_2-\gamma$, whereas on the unstable blocks,
$\Rea\lambda=\sigma-\gamma$.
Combining the blocks and returning to the original coordinates gives,
for some constant $C>0$ independent of $t$,
\[
 \begin{aligned}
 \|e^{t\mathbf D}\mathbf P_s\|
 &\le C(1+t^{n-1})e^{-(\gamma-\sigma_2)t},&&t\ge0,\\
 \|e^{t\mathbf D}\mathbf P_u\|
 &\le C(1+|t|^{n-1})e^{(\sigma-\gamma)t},&&t\le0.
 \end{aligned}
\]
The functions $(1+t^{n-1})e^{-(\gamma-\sigma_2-\delta)t}$ for $t\ge0$
and $(1+|t|^{n-1})e^{-(\sigma-\gamma-\delta)|t|}$ for $t\le0$
are bounded
.
Hence there is a constant $K>0$ such that
\begin{equation}\label{eq:stable-unstable-semigroup-bounds}
 \|e^{t\mathbf D}\mathbf P_s\|\le Ke^{-\delta t}\quad(t\ge0),\qquad
 \|e^{t\mathbf D}\mathbf P_u\|\le Ke^{\delta t}\quad(t\le0).
\end{equation}

Choose $r_0\ge r_*$ so large that
\begin{equation}\label{eq:stable-contraction-constant}
 q:=\frac{2K}{\delta}\sup_{r\ge r_0}\|\mathbf E(r)\|<1.
\end{equation}
For a solution $\mathbf x$, the substitution
$\mathbf w(r)=e^{-\gamma(r-r_0)}\mathbf x(r)$ gives
$\mathbf w'=(\mathbf D+\mathbf E(r))\mathbf w$ on $[r_0,\infty)$.
We seek the initial vectors for which this rescaled solution is bounded.
Write $\operatorname{im}\mathbf P_s=\{\mathbf P_s\mathbf y:\mathbf y\in\C^n\}$
for the image of $\mathbf P_s$. Fix $\boldsymbol{\xi}\in\operatorname{im}\mathbf P_s$.
We work in the space of bounded continuous functions
$\mathbf w:[r_0,\infty)\to\C^n$, with norm
$\|\mathbf w\|_\infty:=\sup_{r\ge r_0}\|\mathbf w(r)\|$.
Following Coppel's integral construction
\cite[Lecture~4, (4), pp.~30--31]{CoppelDichotomies}, with
unperturbed fundamental matrix $e^{\mathbf D(r-r_0)}$, we construct
a bounded solution with $\mathbf P_s\mathbf w(r_0)=\boldsymbol{\xi}$
by solving
\begin{equation}\label{eq:stable-integral-equation}
 \begin{aligned}
 \mathbf w(r)={}&e^{\mathbf D(r-r_0)}\boldsymbol{\xi}
 +\int_{r_0}^r e^{\mathbf D(r-t)}\mathbf P_s\mathbf E(t)\mathbf w(t)\,\dd t\\
 &-\int_r^\infty e^{\mathbf D(r-t)}\mathbf P_u\mathbf E(t)\mathbf w(t)\,\dd t,
 \qquad r\ge r_0.
 \end{aligned}
\end{equation}
Since $\mathbf P_s\boldsymbol{\xi}=\boldsymbol{\xi}$,
\eqref{eq:stable-unstable-semigroup-bounds} gives
$\|e^{\mathbf D(r-r_0)}\boldsymbol{\xi}\|\le K\|\boldsymbol{\xi}\|$.
The same estimate yields, for $r\ge r_0$,
\[
 \begin{aligned}
 \int_{r_0}^r\|e^{\mathbf D(r-t)}\mathbf P_s\|\,\dd t
 &\le K\int_{r_0}^r e^{-\delta(r-t)}\,\dd t\le\frac K\delta,\\
 \int_r^\infty\|e^{\mathbf D(r-t)}\mathbf P_u\|\,\dd t
 &\le K\int_r^\infty e^{-\delta(t-r)}\,\dd t=\frac K\delta.
 \end{aligned}
\]
Both integrals in~\eqref{eq:stable-integral-equation} are therefore
well defined for every function in this space.
Let $\mathcal F$ denote the right-hand side
of~\eqref{eq:stable-integral-equation}.
The bound $\|\mathcal F(\mathbf w)\|_\infty
\le K\|\boldsymbol{\xi}\|+q\|\mathbf w\|_\infty$ shows that
$\mathcal F$ maps this space into itself. For any two such functions,
\[
 \begin{aligned}
 \|\mathcal F(\mathbf w_1)-\mathcal F(\mathbf w_2)\|_\infty
 &\le\left(\frac K\delta+\frac K\delta\right)
      \sup_{t\ge r_0}\|\mathbf E(t)\|\,
      \|\mathbf w_1-\mathbf w_2\|_\infty\\
 &=q\|\mathbf w_1-\mathbf w_2\|_\infty.
 \end{aligned}
\]
Since $q<1$ by~\eqref{eq:stable-contraction-constant},
$\mathcal F$ is a contraction and has a unique fixed point
$\mathbf w_{\boldsymbol{\xi}}$.
Differentiating~\eqref{eq:stable-integral-equation} verifies
$\mathbf w_{\boldsymbol{\xi}}'=(\mathbf D+\mathbf E)\mathbf w_{\boldsymbol{\xi}}$.

Conversely, let $\mathbf y$ be any bounded solution of
$\mathbf y'=(\mathbf D+\mathbf E)\mathbf y$ and put
$\boldsymbol{\xi}=\mathbf P_s\mathbf y(r_0)$.
The spectral projections commute with $\mathbf D$, so variation of
constants gives, for $T>r\ge r_0$,
\[
 \mathbf P_u\mathbf y(r)
 =e^{\mathbf D(r-T)}\mathbf P_u\mathbf y(T)
 -\int_r^T e^{\mathbf D(r-t)}
   \mathbf P_u\mathbf E(t)\mathbf y(t)\,\dd t.
\]
For fixed $r$,~\eqref{eq:stable-unstable-semigroup-bounds} gives
\[
 \|e^{\mathbf D(r-T)}\mathbf P_u\mathbf y(T)\|
 \le Ke^{-\delta(T-r)}\sup_{t\ge r_0}\|\mathbf y(t)\|
 \longrightarrow0,\qquad T\to\infty.
\]
Letting $T\to\infty$ and adding the stable component obtained by
variation of constants from $r_0$ to $r$ yields
\eqref{eq:stable-integral-equation} with $\mathbf w=\mathbf y$.
Uniqueness in~\eqref{eq:stable-integral-equation} then gives
$\mathbf y=\mathbf w_{\boldsymbol{\xi}}$.

Let $\boldsymbol{\xi}_1,\boldsymbol{\xi}_2\in\operatorname{im}\mathbf P_s$
and $a,b\in\C$. By linearity of~\eqref{eq:stable-integral-equation},
the bounded function $a\mathbf w_{\boldsymbol{\xi}_1}
+b\mathbf w_{\boldsymbol{\xi}_2}$ satisfies the same integral equation
with parameter $a\boldsymbol{\xi}_1+b\boldsymbol{\xi}_2$.
Uniqueness therefore gives
\[
 \mathbf w_{a\boldsymbol{\xi}_1+b\boldsymbol{\xi}_2}(r)
 =a\mathbf w_{\boldsymbol{\xi}_1}(r)
  +b\mathbf w_{\boldsymbol{\xi}_2}(r),\qquad r\ge r_0.
\]
Evaluating this identity at \(r=r_0\) shows that the map from the parameter \(\boldsymbol\xi\) to the initial vector \(\mathbf w_{\boldsymbol\xi}(r_0)\) is complex linear. Applying $\mathbf P_s$ to~\eqref{eq:stable-integral-equation}
at $r=r_0$, using $\mathbf P_s\mathbf P_u=0$ and commutation with
$e^{t\mathbf D}$, gives
$\mathbf P_s\mathbf w_{\boldsymbol{\xi}}(r_0)=\boldsymbol{\xi}$.
Thus $\mathbf w_{\boldsymbol{\xi}_1}(r_0)
=\mathbf w_{\boldsymbol{\xi}_2}(r_0)$ implies
$\boldsymbol{\xi}_1=\boldsymbol{\xi}_2$, so this linear map is injective.
Since $\dim\operatorname{im}\mathbf P_s
=\operatorname{rank}\mathbf P_s=n-m$, its image is a complex linear
subspace of dimension $n-m$. The converse proved above shows that
\[
 \{\mathbf w_{\boldsymbol{\xi}}(r_0):
   \boldsymbol{\xi}\in\operatorname{im}\mathbf P_s\}
\]
consists exactly of the initial vectors at $r_0$ of all bounded solutions.

Finally, the list of rates in
\eqref{eq:continuous-system-growth-values} shows that every
nonzero solution $\mathbf x$ has either $\lambda(\mathbf x)=\sigma$ or
$\lambda(\mathbf x)\le\sigma_2$.
Since $\sigma_2<\gamma<\sigma$, the rescaled solution
$e^{-\gamma(r-r_0)}\mathbf x(r)$ is bounded exactly when
$\lambda(\mathbf x)<\sigma$.
At $r=r_0$, the factor in
$\mathbf w(r)=e^{-\gamma(r-r_0)}\mathbf x(r)$ equals $1$, so
$\mathbf w(r_0)=\mathbf x(r_0)$.
The solution map from $r_*$ to $r_0$ is an invertible linear map,
so transporting this subspace back to $r_*$ proves the statement.
\end{proof}

\subsection{Proof of \hyperref[lem:physical-window]{Lemma~\ref*{lem:physical-window}}}
\label{subsec:proof-physical-window}

\begin{proof}
Fix $L>0$ and consider the rectangle
\eqref{eq:physical-window-rectangle} in the variable $t=\omega z=x+iy$,
with $R=ar$ and $Q=R^{\beta_*}e^R$ as in the statement. 
	We want to verify some of the curvilinear strips in Lemma~\ref{lem:moving-inverse} cover this rectangle and the pull-back functions in Lemma~\ref{lem:moving-inverse} are well-defined there.
	
Fix $\kappa>0$ satisfying~\eqref{eq:inverse-kappa-choice}, which ensures
exponential decay of the terms with $\lambda\ne(n-j)\omega$ in
$x>0$, $|y|\le\kappa x$.
Since $x\asymp R$ and $|y|/x=O_L(R^{-1/2})\to0$, the rectangle
lies in $x>0$, $|y|\le\kappa x$ for all sufficiently large $R$.
Figure~\ref{fig:controlled-strip-windows} shows these rectangles
and the domains of the controlled inverse branches.

On this rectangle, $|x-R|\le L$ and $|y|\le L\sqrt R$, so
$|t/R-1|\le L/R+L/\sqrt R=O_L(R^{-1/2})$ and $|t|\asymp R$.
For the power branch specified after~\eqref{eq:phase}, the Taylor
expansion $(1+w)^{\beta_*}=1+\beta_*w+O(|w|^2)$ at $w=0$ gives
\[
 t^{\beta_*}
 =R^{\beta_*}\left(1+\frac{t-R}{R}\right)^{\beta_*}
 =R^{\beta_*}(1+O_L(R^{-1/2})).
\]
Using~\eqref{eq:phase-expansion-proof} and
$R^{\beta_*}e^t=Qe^{x-R+iy}$, we obtain, uniformly on the rectangle,
\begin{equation}\label{eq:window-phase}
 \begin{aligned}
 \Phi(t)=t^{\beta_*}e^t(1+O_L(R^{-1}))=Qe^{x-R+iy}(1+O_L(R^{-1/2})).
 \end{aligned}
\end{equation}
Thus $\Phi$ has no zeros on the rectangle for large $R$, and
$|\Phi(t)|/Q$ stays between two positive constants depending only
on $L$.

For each fixed $x$, choose a continuous branch of
$\arg\Phi(x+iy)$ along the vertical segment.
Taking the reciprocal in~\eqref{eq:phase-expansion-proof} gives
$\Phi'(t)/\Phi(t)=1+\beta_*/t+O(t^{-2})=1+O_L(R^{-1})$.
Since the argument is locally $\Ima\log\Phi(x+iy)$, the chain rule
therefore gives, uniformly on the rectangle,
\begin{equation}\label{eq:window-argument-derivative}
 \partial_y\arg\Phi(x+iy)
 =\Ima\!\left(i\frac{\Phi'(t)}{\Phi(t)}\right)
 =\Rea\frac{\Phi'(t)}{\Phi(t)}
 =1+O_L(R^{-1}).
\end{equation}

We first determine which inverse branches \(\psi_\ell\) of \(\Phi\), supplied by Lemma 2.3, have images covering this rectangle.
Cover the unit circle in the $s$-plane by finitely many closed
arcs, each lying inside an open argument interval $J$ of width
less than $\pi$. On each arc use the corresponding logarithm
$\log_J$, and consider $s=\Phi(t)$ with $t$ in the rectangle
\eqref{eq:physical-window-rectangle} and $s/|s|$ on that arc.
Taking logarithms of
$\Phi(t)=t^{\beta_*}e^t(\Phi(t)/\Phi'(t))$ gives an integer $\ell$
such that
\begin{equation}\label{eq:window-logarithmic-identity}
 t+\beta_*\log t+\log(\Phi(t)/\Phi'(t))
 =\log_Js+2\pi i\ell.
\end{equation}
The logarithm of $\Phi/\Phi'$ is the one tending to zero as
$x\to\infty$ with $|y|\le\kappa x$, supplied
by~\eqref{eq:phase-expansion-proof}.
Taking imaginary parts in~\eqref{eq:window-logarithmic-identity},
with $\arg t$ and $\arg_Js$ bounded, gives
\[
 2\pi\ell=y+O_L(1).
\]
Since $|y|\le L\sqrt R$, it follows that $|\ell|\le C_L\sqrt R$.
For $s=\Phi(t)$, \eqref{eq:window-phase} gives
$\log|s|=R+\beta_*\log R+O_L(1)$ uniformly for $t$ in the rectangle
\eqref{eq:physical-window-rectangle}.
Choose a fixed $L_0>0$ large enough for Lemma~\ref{lem:moving-inverse}
on the finitely many chosen intervals $J$.
The bounds on $\log|s|$ and $|\ell|$ then give
$\log|s|>L_0(1+|\ell|)$ for all sufficiently large $R$.
Together with $\arg_Js\in J$, this ensures that the inverse branch
$\psi_\ell$ supplied by that lemma is defined at $s=\Phi(t)$.

By~\eqref{eq:moving-inverse} and the bounds for $X=\log|s|$ and
$\ell$ above, the segment joining $t$ and $\psi_\ell(s)$ lies in
$x=R+O_L(\log R)$, $|y|=O_L(\sqrt R)$.
For large $R$ it is therefore contained in $x>0$, $|y|\le\kappa x$.
The derivative estimate in Section~\ref{subsec:proof-moving-inverse}
for $G_{s,\ell}$ in~\eqref{eq:inverse-contraction-map} gives
$|G'_{s,\ell}|\le C_L/R$ along this segment.
By~\eqref{eq:inverse-fixed-point},
$G_{s,\ell}(\psi_\ell(s))=\psi_\ell(s)$.
Equation~\eqref{eq:window-logarithmic-identity} gives
$G_{s,\ell}(t)=t$, so $t$ is also a fixed point of $G_{s,\ell}$.
Integrating along the segment gives
\[
 |t-\psi_\ell(s)|
 =|G_{s,\ell}(t)-G_{s,\ell}(\psi_\ell(s))|
 \le\frac{C_L}{R}|t-\psi_\ell(s)|.
\]
For $R>C_L$, this implies $t=\psi_\ell(s)$. Hence the branches with
$|\ell|\le C_L\sqrt R$ cover the rectangle.

\bigskip
We can now apply the uniform growth estimate.
By~\eqref{eq:window-phase} and $|x-R|\le L$, the normalized points
$\zeta=s/Q$ in the $\zeta$-plane satisfy
\[
 \tfrac12e^{-L}\le|\zeta|\le2e^L
\]
for all sufficiently large $R$. 
Restricting the arguments to the chosen closed arcs gives fixed
compact subsets of the corresponding sectors in the $\zeta$-plane.
Put $T=Q$ and choose $m(T)=C_L\sqrt{\log T}$ for $T>1$,
with $C_L$ enlarged if necessary.
Since $\log Q\sim R$, the bound $|\ell|\le C_L\sqrt R$ gives
$|\ell|\le m(T)$, while $m(T)=o(\log T)$, as required by
Lemma~\ref{lem:uniform-fundamental-row}.
For every fixed coefficient tolerance $\varepsilon>0$, the radii
$\rho_{\ell,\varepsilon}$ in
Lemma~\ref{lem:branch-matrix-bounds} are at most
$\exp(C(1+\sqrt R))=o(Q)$.
For each such compact set $K$ and all sufficiently large $R$,
\[
 Q\min_{\zeta\in K}|\zeta|
 \ge\tfrac12e^{-L}Q\ge\rho_{\ell,\varepsilon}.
\]
This verifies the radius condition stated
after~\eqref{eq:quantified-row-growth}.
Since $t=\psi_\ell(s)$ and $s=Q\zeta$,
\eqref{eq:branch-fundamental-row} gives
$U_r(\xi,\eta)=Q^{-1}\log\|\mathbf F_\ell(Q\zeta)\|$.
The bound~\eqref{eq:uniform-row-coarse-potential-bound}, valid at
every angle in the chosen sectors, therefore gives
\[
 \begin{aligned}
 |U_r(\xi,\eta)|
\le C_K+
 \frac{C_Ke^{C_K(1+C_L\sqrt R)}}{R^{\beta_*}e^R}
 =C_K+o(1).
 \end{aligned}
\]
The finitely many chosen intervals $J$ cover all arguments in the
$s$-plane, so the maximum of their constants gives a uniform bound
for $|U_r|$ on $|\xi|,|\eta|\le L$.

To compare the leading terms, remove the intervals
$(\theta_0-\alpha,\theta_0+\alpha)$, modulo $2\pi$, around the
finitely many directions $\theta_0$ satisfying
$\Rea((\mu-\nu)e^{i\theta_0})=0$ for distinct roots $\mu,\nu$
of $p_\omega$.
For fixed $\alpha>0$, the largest of the values
$\Rea(\mu e^{i\theta})$ exceeds all the others by at least a
constant $\delta_\alpha>0$ on the remaining angles, as
in~\eqref{eq:uniform-row-root-gap}.
Choose $0<\varepsilon<\min\{1,\delta_\alpha/8\}$ and
apply~\eqref{eq:quantified-row-growth}.
The part of that estimate containing the branch index is bounded by
\[
 \frac{C_{\varepsilon,\alpha,L}
       e^{C_{\varepsilon,\alpha,L}(1+C_L\sqrt R)}}{R^{\beta_*}e^R}
 \longrightarrow0.
\]
Thus, away from the removed angles,
$Q^{-1}|\log\|\mathbf F(t/\omega)\|-H_\omega(\Phi(t))|$
is at most $C_{\alpha,L}\varepsilon+o_{\varepsilon,\alpha,L}(1)$.
Let \(B_{R,\alpha}\) be the part of \([-L,L]^2\) mapped to the removed angular regions in the \(s\)-plane. The uniform bound for $U_r$ gives
\[
 \int_{B_{R,\alpha}}|U_r|\,\mathrm d\xi\,\mathrm d\eta
 \le C_L\,\operatorname{area}(B_{R,\alpha}).
\]
The bound for $|\Phi(t)|/Q$ gives the same area estimate for the
error above, with $C_L$ enlarged independently of $\alpha$ and
$\varepsilon$.

For each fixed $x$, let $\Theta_x(y)$ be the continuous argument
of $\Phi(x+iy)$. By~\eqref{eq:window-argument-derivative},
$1/2\le\Theta_x'(y)\le2$ for large $R$.
As $y$ ranges over $[-L\sqrt R,L\sqrt R]$, the total increase is
$O_L(\sqrt R)$. Divide this argument range at the values $2\pi k$,
$k\in\mathbb Z$, into complete intervals $[2\pi k,2\pi(k+1)]$
and at most two shorter intervals at the endpoints.
The derivative lower bound gives
\[
 y_2-y_1\le2\bigl(\Theta_x(y_2)-\Theta_x(y_1)\bigr),
 \qquad -L\sqrt R\le y_1\le y_2\le L\sqrt R.
\]
Thus the preimage of any angular interval has length in $y$
at most twice its angular length.
Each complete interval contains removed angles of total length
$O(\alpha)$, and there are $O_L(\sqrt R)$ such intervals.
Their removed preimages therefore have total length
$O_L(\alpha\sqrt R)$.
Each of the two possible endpoint intervals has angular length
less than $2\pi$, so their preimages together have length at most
$8\pi=O(1)$.
Integrating over $|x-R|\le L$ gives area $O_L(\alpha\sqrt R+1)$
for the part of the rectangle mapping to the removed angles.
Integrating the bound from~\eqref{eq:quantified-row-growth} outside
the removed angles over an area at most $4L^2\sqrt R$, and using
the uniform bound on the remaining part, gives
\begin{equation}\label{eq:window-row-comparison}
 \begin{aligned}
 &\limsup_{r\to\infty}\frac1{Q\sqrt R}
 \int_{\substack{|x-R|\le L\\|y|\le L\sqrt R}}
 \left|\log\|\mathbf F(t/\omega)\|-H_\omega(\Phi(t))\right|
 \dd x\dd y\\
 &\quad\le\limsup_{r\to\infty}
 \left[4L^2\bigl(C_{\alpha,L}\varepsilon
 +o_{\varepsilon,\alpha,L}(1)\bigr)
 +C_L\bigl(\alpha+R^{-1/2}\bigr)\right]\\
 &\quad\le C_L\alpha+C_{\alpha,L}\varepsilon,
 \end{aligned}
\end{equation}
after enlarging the constants. First let $\varepsilon\to0$ for fixed $\alpha$, and then
let $\alpha\to0$. The normalized integral of~\eqref{eq:window-row-comparison} tends to zero.
Recall that $t=R+\xi+i\sqrt R\,\eta$ and
$U_r(\xi,\eta)=Q^{-1}\log\|\mathbf F(t/\omega)\|$.
Since $\dd x\dd y=\sqrt R\,\dd\xi\dd\eta$,
\eqref{eq:window-row-comparison} gives
\[
 \int_{[-L,L]^2}
 \left|U_r(\xi,\eta)-\frac{H_\omega(\Phi(t))}{Q}\right|
 \dd\xi\dd\eta\longrightarrow0.
\]
Moreover, \eqref{eq:window-phase} gives
$\Phi(t)/Q=e^{\xi+i\sqrt R\eta}+O_L(R^{-1/2})$ uniformly
on this square.
By~\eqref{eq:fundamental-row-support-function}, $H_\omega$ is
Lipschitz with constant $\max_{p_\omega(\mu)=0}|\mu|$ and
satisfies $H_\omega(q\zeta)=qH_\omega(\zeta)$ for $q>0$.
In particular, $H_\omega(e^{\xi+i\sqrt R\eta})
=e^\xi h_\omega(\sqrt R\eta)$, so
\[
 \begin{aligned}
 \left|\frac{H_\omega(\Phi(t))}{Q}
       -e^\xi h_\omega(\sqrt R\eta)\right|
 &=\left|H_\omega\!\left(\frac{\Phi(t)}Q\right)
        -H_\omega(e^{\xi+i\sqrt R\eta})\right|\\
 &\le\left(\max_{p_\omega(\mu)=0}|\mu|\right)
       \left|\frac{\Phi(t)}Q-e^{\xi+i\sqrt R\eta}\right|
 \le C_LR^{-1/2}.
 \end{aligned}
\]
The triangle inequality therefore gives
\[
 \begin{aligned}
 &\int_{[-L,L]^2}
 \left|U_r(\xi,\eta)-e^\xi h_\omega(\sqrt R\eta)\right|
 \dd\xi\dd\eta\\
 &\qquad\le
 \int_{[-L,L]^2}
 \left|U_r(\xi,\eta)-\frac{H_\omega(\Phi(t))}{Q}\right|
 \dd\xi\dd\eta+4L^2C_LR^{-1/2}\longrightarrow0.
 \end{aligned}
\]
Since $L>0$ is arbitrary, this proves the approximation
$U_r(\xi,\eta)-e^\xi h_\omega(\sqrt R\eta)\to0$
in $L^1_{\rm loc}(\R^2)$.

\bigskip
It remains to average the periodic function $h_\omega$.
Cauchy's perimeter formula~\cite[Ch.~I, \S19, (1.115), p.~80]{Levin} gives
\[
 \frac1{2\pi}\int_0^{2\pi}h_\omega(\theta)\dd\theta
 =\frac{P_\omega}{2\pi}.
\]
For the convex hull $K_\omega=\conv\{\overline\mu:p_\omega(\mu)=0\}$,
this formula also follows from~\eqref{eq:angular-support-jump}.
Between consecutive edge normal angles,
$h_\omega(\theta)=\Rea(\mu e^{i\theta})$ for a fixed root $\mu$,
so $h_\omega''=-h_\omega$.
At each edge normal angle $\theta_0$, the jump of $h_\omega'$
equals $d_{K_\omega}(\theta_0)$, the length of the corresponding
supporting face. Thus, as distributions on the angular circle
$\R/(2\pi\mathbb Z)$,
\[
 h_\omega+h_\omega''
 =\sum_{\theta_0}d_{K_\omega}(\theta_0)\delta_{\theta_0},
\]
where the sum is over the normal angles with positive face length.
Pairing with the constant function $1$ gives
$\langle h_\omega'',1\rangle=\langle h_\omega,1''\rangle=0$,
while~\eqref{eq:point-mass-pairing} gives
$\langle\delta_{\theta_0},1\rangle=1$. Hence
\[
 \int_0^{2\pi}h_\omega(\theta)\dd\theta
 =\sum_{\theta_0}d_{K_\omega}(\theta_0)=P_\omega.
\]
If $K_\omega$ is a segment, its two perpendicular normal
directions each contribute its full length, giving twice its length, as in the definition of \(P_\omega\).
Fix $\phi\in C_c(\R^2)$ and choose $L>0$ such that
$\operatorname{supp}\phi\subset[-L,L]^2$. Partition the $\eta$-axis into intervals
$I_k=[2\pi k/\sqrt R,\,2\pi(k+1)/\sqrt R]$, $k\in\mathbb Z$.
The substitution $\theta=\sqrt R\eta$ and the $2\pi$-periodicity
of $h_\omega$ give
\[
 \int_{I_k}h_\omega(\sqrt R\eta)\dd\eta
 =\frac1{\sqrt R}\int_{2\pi k}^{2\pi(k+1)}
 h_\omega(\theta)\dd\theta
 =\frac{P_\omega}{\sqrt R}.
\]
Since $|I_k|=2\pi/\sqrt R$, the mean on each full interval is
$P_\omega/(2\pi)$.
Put $\bar h=P_\omega/(2\pi)$, $\Delta_R=2\pi/\sqrt R$ and
$a_k=2\pi k/\sqrt R$.
Uniform continuity of $\phi$ gives $\varepsilon_R\to0$ such that
$|\phi(\xi,\eta)-\phi(\xi,a_k)|\le\varepsilon_R$
for $|\xi|\le L$ and $\eta\in I_k\subset[-L,L]$.
Since $\int_{I_k}(h_\omega(\sqrt R\eta)-\bar h)\dd\eta=0$,
\[
 \int_{I_k}\phi(\xi,\eta)
       \bigl(h_\omega(\sqrt R\eta)-\bar h\bigr)\dd\eta
 =
 \int_{I_k}\bigl(\phi(\xi,\eta)-\phi(\xi,a_k)\bigr)
       \bigl(h_\omega(\sqrt R\eta)-\bar h\bigr)\dd\eta.
\]
Let $M=\max_{0\le\theta\le2\pi}|h_\omega(\theta)-\bar h|$.
The last integral has modulus at most $M\varepsilon_R|I_k|$.
The full intervals contained in $[-L,L]$ have total length at
most $2L$, while the two possible remaining end intervals have
total length at most $2\Delta_R$. Consequently,
\[
 \left|\int_{-L}^{L}\phi(\xi,\eta)
       \bigl(h_\omega(\sqrt R\eta)-\bar h\bigr)\dd\eta\right|
 \le2LM\varepsilon_R
       +2\Delta_RM\|\phi\|_\infty\longrightarrow0,
\]
uniformly for $|\xi|\le L$.
Multiplying by $e^\xi$ and integrating in $\xi$ therefore gives
\[
 \int_{\R^2}e^\xi h_\omega(\sqrt R\eta)\phi(\xi,\eta)
 \dd\xi\dd\eta
 \longrightarrow
 \frac{P_\omega}{2\pi}\int_{\R^2}e^\xi\phi(\xi,\eta)\dd\xi\dd\eta.
\]
For the same test function $\phi$, the $L^1_{\rm loc}$ estimate
for the difference $U_r-e^\xi h_\omega(\sqrt R\eta)$ gives
\[
 \int_{\R^2}\bigl(U_r(\xi,\eta)-e^\xi h_\omega(\sqrt R\eta)\bigr)
 \phi(\xi,\eta)\dd\xi\dd\eta\longrightarrow0.
\]
Combining these two integral limits proves the weak convergence
in~\eqref{eq:physical-potential-average}.
\end{proof}

\section{Proofs of the main results}\label{sec:main-proofs}

\subsection{Proof of \hyperref[thm:intro]{Theorem~\ref*{thm:intro}}}
\label{subsec:proof-intro}

\begin{proof}
On the inverse branch supplied by Lemma~\ref{lem:phase-equation},
$v(s)=f(\omega^{-1}\psi(s))$ satisfies~\eqref{eq:phase-equation},
with limiting polynomial $p$ in~\eqref{eq:characteristic}.
The inverse construction extends this equation to a slightly
larger sector, so that Proposition~\ref{prop:sectorial-counting}
applies to $J$. Put
$\Lambda_f(J)=\Lambda_v(J)$ and
$K_f=\conv\{\overline\mu:\mu\in\Lambda_f(J)\}$.
To transfer the sectorial count to the disk $|z|\le r$, we
compare radii in the two coordinates. Put $Q(r)=(ar)^\beta e^{ar}$.
Since $\log\log(cQ(r))=\log(ar)+o(1)$ for fixed $c>0$,
\eqref{eq:inverse-proof} gives, uniformly on closed subintervals of $J$,
\[
 \left|\omega^{-1}\psi(cQ(r)e^{i\theta})\right|
       =r+\frac{\log c}{a}+o(1).
\]
On the fixed branch, the imaginary part of $\psi$ is bounded,
so replacing its real
part by its modulus adds $O(r^{-1})$.
\phantomsection\label{par:inverse-radial-monotonicity}
Moreover, differentiation of $\Phi(\psi(s))=s$ gives
\[
 \frac{\partial\psi(\rho e^{i\theta})}{\partial\log\rho}
       =\frac{\Phi(\psi(\rho e^{i\theta}))}
                    {\Phi'(\psi(\rho e^{i\theta}))}=1+o(1).
\]
Since $\psi/|\psi|\to1$ on these closed angular intervals,
differentiating the modulus gives
\[
 \frac{\partial}{\partial\log\rho}
       |\omega^{-1}\psi(\rho e^{i\theta})|
       =\frac1a+o(1)>0.
\]
Thus the modulus increases along each radial line in the
$s$-plane, uniformly on the chosen angular interval. Apply
Lemma~\ref{lem:radial-transfer} to the full map
$s\mapsto\omega^{-1}\psi(s)$, with $G=1$,
to obtain
\begin{equation}\label{eq:maincount}
 n_{U(\alpha_1,\alpha_2)}(r,1/f)
 =\frac{(ar)^\beta e^{ar}}{2\pi}
       \sum_{\alpha_1<\theta<\alpha_2}d_{K_f}(\theta)
       +o(r^\beta e^{ar}).
\end{equation}
Define
\[
 M_f(\alpha)=\frac{a^\beta}{2\pi}
       \sum_{\substack{\theta\in J\\\theta\le\alpha}}d_{K_f}(\theta).
\]
The sum is over outward normals to the faces of $K_f$ with
positive length. Formula~\eqref{eq:maincount} is then
\eqref{eq:intro-count}. Possible normal directions come from
pairs of characteristic roots, hence are fixed independently of $f$.

If $M_f$ is constant on a neighbourhood of $[\alpha_1,\alpha_2]$,
choose $\alpha_1'<\alpha_1<\alpha_2<\alpha_2'$ in that neighbourhood
with $[\alpha_1',\alpha_2']\subset J$. The strip
$U(\alpha_1',\alpha_2')$ contains $U(\alpha_1,\alpha_2)$ together
with both angular boundary curves. Its count is
$o(r^\beta e^{ar})$ by~\eqref{eq:maincount}, proving the assertion
for either convention on those curves.
\end{proof}

\begin{example}\label{ex:strip-polynomial-factor}
For a nonconstant polynomial multiplier, put $H(z)=(z-1)e^z$,
$B(z)=(z+3)e^z+z^2e^{2z}$, and
$f_1(z)=e^{-2z/3}(e^{H(z)}-1)$. Then
\[
f_1'''-(B+1/3)f_1'+(2/27-2B/3)f_1=0,\qquad z\in\C.
\]
Indeed, $H'=ze^z$, $H'''-2H''+H'=0$, and
$B=(H')^2+3H''-2H'$. Substituting $e^H-1$ in
$y'''-2y''+(1-B)y'=0$ on $\C$ and then $y=e^{2z/3}f_1$ gives the equation.
Its selected terms have $\omega=1$, $\beta=1$, and
$p(\mu)=\mu^3-\mu$. We may take $\varphi=H$.
The zeros satisfy $H(z)=2\pi i k$. On the row with
$0<\Ima z<\pi/2$ and $\Rea z>1$, the equation for their arguments is
\[
 y+\arctan\frac{y}{x-1}=\frac\pi2.
\]
For $x>1$, it has one solution $y(x)$, and
$y'=y/((x-1)^2+y^2+x-1)>0$, with $y(x)\to\pi/2$.
At the disk boundary $x=r+O(r^{-1})$ and $|z-1|=r+O(1)$.
Counting the points $2\pi i k$ therefore gives, in a fixed
pullback strip containing this row,
\[
 n_U(r,1/f_1)\sim\frac{r e^r}{2\pi}.
\]
Thus the factor $r^\beta$ in Theorem~\ref{thm:intro} cannot generally
be suppressed.
\end{example}

\subsection{Proof of \hyperref[thm:generic-density]{Theorem~\ref*{thm:generic-density}}}
\label{subsec:proof-generic-density}

\begin{proof}
Fix a maximal weighted frequency $\omega\in\Omega_a$, put
$\beta=\beta(\omega)$ as in~\eqref{eq:beta}, and let $p$ be
its limiting polynomial from~\eqref{eq:characteristic}. Its coefficient of $\mu^{n-1}$
vanishes, whereas at least one $c_j$, $j\le n-2$, is nonzero.
Thus $p$ cannot be a power $(\mu-c)^n$: comparison of the
$\mu^{n-1}$ coefficient gives $c=0$, contradicting
$p\ne\mu^n$. In particular, the convex hull
\[
 K=\conv\{\overline\mu:p(\mu)=0\}
\]
has a supporting face of positive length. Choose such a face with outward
unit normal $e^{i\theta_0}$ and label its endpoints $\overline{\mu_-}$,
$\overline{\mu_+}$ in the direction $i e^{i\theta_0}$ along the face.
Fix an open interval $J$ of width less than $\pi$ containing $\theta_0$
and an inverse branch $\psi_\ell$ from Lemma~\ref{lem:phase-equation} on $S_J$.
Choose $\delta_0>0$ so that $[\theta_0-\delta_0,\theta_0+\delta_0]\subset J$
contains no other angle at which two distinct roots have equal projections.
This is possible since, for each pair $\mu\ne\nu$, the equation
$\Rea((\mu-\nu)e^{i\theta})=0$ has only finitely many solutions in $J$.
At $\theta_0$ every conjugate root off the face has smaller projection
than either endpoint, so by continuity this ordering holds throughout
the chosen interval.
Along the oriented face, the projection $\nu\mapsto\Rea(\nu e^{-i\theta})$
has slope $\Rea(i e^{i\theta_0}e^{-i\theta})=\sin(\theta-\theta_0)$
with respect to arclength. Thus for every $0<\delta<\delta_0$ the
projection decreases along the face at $\theta_-=\theta_0-\delta$ and
increases at $\theta_+=\theta_0+\delta$. Hence $\mu_-$ and $\mu_+$
uniquely maximize $\Rea(\mu e^{i\theta})$ at $\theta_-$ and $\theta_+$,
respectively.

Let $\mathbf A_0$ be the companion matrix of $p$.
Choose $r_*$ large enough that $s_\pm=r_*e^{i\theta_\pm}$ lie in the
domain of $\psi_\ell$. On each ray the limiting matrix is
$e^{i\theta_\pm}\mathbf A_0$. By Lemma~\ref{lem:slower-growth-subspace},
the initial vectors at $s_\pm$ whose solutions have radial growth rate
strictly below $\Rea(e^{i\theta_\pm}\mu_\pm)$, together with the zero
vector, form a proper linear subspace $W_\pm\subset\C^n$.

For an initial vector $\mathbf c$ at $z_0$, put
$v(s)=f_{\mathbf c}(\psi_\ell(s)/\omega)$. By~\eqref{eq:branch-fundamental-matrix},
its vector of derivatives of orders $0,\ldots,n-1$ at $s_\pm$ is
$\mathbf X_\ell(s_\pm)\mathbf c$. If this vector is zero, uniqueness
for the transformed equation gives $v\equiv0$. Hence $f_{\mathbf c}$
vanishes on the corresponding strip and, by the identity theorem, on
$\C$, so $\mathbf c=0$. Thus $\mathbf X_\ell(s_\pm)$ is invertible.
The exceptional subspaces at $z_0$ are therefore
\[
 V_\pm=\mathbf X_\ell(s_\pm)^{-1}W_\pm,\qquad
 \mathbf c\in V_\pm\ \Longleftrightarrow\
 \mathbf X_\ell(s_\pm)\mathbf c\in W_\pm.
\]
Let $m_-$ and $m_+$ be the algebraic multiplicities of $\mu_-$ and
$\mu_+$ in $p$, respectively. By Lemma~\ref{lem:slower-growth-subspace}
and the invertibility of $\mathbf X_\ell(s_\pm)$,
$\dim V_\pm=\dim W_\pm=n-m_\pm$.

To compare the choices of ray, write $V_\pm(\delta)$ for the subspaces
just obtained. For $\mathbf c\ne\mathbf 0$, Lemma~\ref{lem:indicator-limit}
identifies the radial growth rate with $h_v(\theta_\pm)$, so
\[
 \mathbf c\notin V_\pm(\delta)
 \quad\Longleftrightarrow\quad
 h_v(\theta_\pm)=\Rea(\mu_\pm e^{i\theta_\pm})
 \quad\Longleftrightarrow\quad
 \mu_\pm\in\Lambda_v(J).
\]
Indeed, \eqref{eq:local-support-indicator} gives the maximum over
$\Lambda_v(J)$, and $\mu_\pm$ is the unique maximizing root at $\theta_\pm$.
For fixed $J$ and $\psi_\ell$, the set $\Lambda_v(J)$ depends only on
$\mathbf c$. Hence every $0<\delta<\delta_0$ gives the same $V_-$ on
the left and the same $V_+$ on the right. For $\mathbf c\notin V_-\cup V_+$,
both endpoints therefore belong to $K_f$, and
\[
 h_v'(\theta_0+)-h_v'(\theta_0-)
       =|\mu_+-\mu_-|>0.
\]
Take continuity directions $\alpha_1<\theta_0<\alpha_2$
close enough to include no other jump. Theorem~\ref{thm:intro}
then gives~\eqref{eq:generic-density} with
\[
 c_\omega=\frac{a^\beta}{2\pi}|\mu_+-\mu_-|>0.
\]
Apply this argument separately to each $\omega\in\Omega_a$,
always using the same point $z_0$ for the initial values.
The resulting exceptional subspaces form a finite family,
and their union has Lebesgue measure zero in $\C^n$.
Outside this union all the estimates hold simultaneously,
which proves~\eqref{eq:simultaneous-strip-count}.
\end{proof}

\begin{example}\label{ex:interior-selected-roots}
Apply the operator $(\mathrm d/\mathrm dz)^3-8e^{3z}$ to the left-hand
side of~\eqref{eq:model-original}. For $s=e^z$, the resulting fifth-order
equation has limiting polynomial $p(\mu)=(\mu^3-8)(\mu^2-\tfrac14)$,
so $K$ is the triangle with vertices $2,-1\pm i\sqrt3$.
For the infinite-order solution $f=f_+-f_-$ of~\eqref{eq:model-original},
we have $v(s)=s^{-1/2}(e^{s/2}-e^{-s/2})$. Hence on $J=(\pi/3,2\pi/3)$,
\[
 \Lambda_f(J)=\{-\tfrac12,\tfrac12\},\qquad
 K_f=[-\tfrac12,\tfrac12]\subset\operatorname{int}K.
\]
\end{example}

\begin{proof}[Proof of {\hyperref[prop:nonmaximal-strip-count]{Proposition~\ref*{prop:nonmaximal-strip-count}}}]
\leavevmode\phantomsection\label{par:proof-nonmaximal-strip-count}%
For $\nu=\lambda/(n-j)\ne\omega$, \eqref{eq:own-direction-dominance}
gives $\Rea(\lambda/\omega-(n-j))=(n-j)(\Rea(\nu/\omega)-1)<0$.
Since there are finitely many frequencies, all such factors
$e^{(\lambda/\omega-(n-j))t}$ decay exponentially in a common narrow
sector about the positive real axis as $\Rea t\to+\infty$.
This is where maximal modulus was used in the coefficient estimates
of Lemma~\ref{lem:phase-equation}. The inverse construction depends
only on $\beta$, so the same estimates hold uniformly in the branch
index, and the radius comparison in
\hyperref[subsec:proof-intro]{\S\ref*{subsec:proof-intro}} gives the
strip count with $a$ replaced by $|\omega|$. Since $c_{n-1}=0$ and
at least one lower coefficient in~\eqref{eq:limiting-coefficients}
is nonzero, the limiting polynomial has at least two distinct roots.
The supporting-face argument in the proof of
Theorem~\ref{thm:generic-density} therefore gives positive density
outside the union of two proper linear subspaces.
\end{proof}

\subsection{Proof of \hyperref[thm:global-zero-distribution]{Theorem~\ref*{thm:global-zero-distribution}}}
\label{subsec:proof-global-zero-distribution}

\begin{proof}
The global upper bound follows from the integral zero bound.
For $q\ne0$, $b\ge0$ and $R\ge1$, use the rotated polar coordinates
$z=e^{-i\arg q}\rho e^{i\theta}$, with $\rho=|z|$ and
$-\pi\le\theta\le\pi$. Here $\theta$ is the argument of $qz$, so
$\Rea(qz)=|q|\rho\cos\theta$. After factoring out $e^{|q|\rho}$,
the angular integral is estimated using
$1-\cos\theta\ge c\theta^2$ on $[-\pi,\pi]$ and the substitution
$u=\sqrt\rho\,\theta$, which gives, for $\rho\ge1$,
\[
 \int_{-\pi}^{\pi}e^{-|q|\rho(1-\cos\theta)}\dd\theta
 \le\rho^{-1/2}\int_{\R}e^{-c|q|u^2}\dd u
 =C_q\rho^{-1/2}.
\]
For $0\le\rho\le1$, the integral is at most $2\pi$; enlarging $C_q$
therefore gives the bound $C_q(1+\rho)^{-1/2}$ for all $\rho\ge0$.
Polar coordinates therefore give
\begin{equation}\label{eq:exponential-disk-integral}
 \begin{aligned}
 \int_{|z|\le R}(1+|z|)^b e^{\Rea(qz)}\,\mathrm dA(z)
 &\le C_q\int_0^R\rho(1+\rho)^{b-1/2}e^{|q|\rho}\,\mathrm d\rho\\
 &\le C_{q,b}R^{b+1/2}e^{|q|R}.
 \end{aligned}
\end{equation}
For $q=0$, direct integration gives $O(R^{b+2})$.
Lemma~\ref{lem:coefficient-integral-zero-bound}, applied with
$E=\{|z|\le r\}$, gives
\[
 n(r,1/f)\le C\int_{|z|\le r+1}W(z)\,\mathrm dA(z).
\]
Apart from the constant term, each summand of $W$
in~\eqref{eq:weight-function} has the form
$(1+|z|)^b e^{\Rea(qz)}$, where
$q=\lambda/(n-j)$ and $b=\deg P_{j,\lambda}/(n-j)$.
For each term with $q\ne0$, apply~\eqref{eq:exponential-disk-integral}
with $R=r+1$.
If $|q|=a$, then $b\le\beta_*$; if $|q|<a$, the exponential
gap absorbs every fixed power of $r$. The constant term and terms
with $q=0$ contribute only polynomial growth after integration.
Summing these finitely many bounds gives~\eqref{eq:global-zero-upper-bound}.

Outside the sectors $\Sigma_\varepsilon$, every exponential
term has a strictly smaller rate. More precisely, put
\[
 d_* =\max_{j,\lambda}\frac{\deg P_{j,\lambda}}{n-j},\qquad
 a_\varepsilon=\max\bigl(\{0,a\cos\varepsilon\}\cup
       \{|q|:q\in\mathcal F,\ |q|<a\}\bigr)<a.
\]
Then
\[
 \Rea(qz)\le a_\varepsilon r,\qquad
 q\in\mathcal F,\quad |z|\le r,\quad z\notin\Sigma_\varepsilon.
\]
For $E_r=\{|z|\le r\}\setminus\Sigma_\varepsilon$ and
$w\in \overline B(E_r,1)$, choose $z\in E_r$ with $|w-z|\le1$.
Then $|w|\le r+1$ and
$\Rea(qw)\le a_\varepsilon r+a$ for every
$q\in\mathcal F$. Thus
\[
 W(w)\le C(1+r)^{d_*}e^{a_\varepsilon r},\qquad
 w\in \overline B(E_r,1).
\]
Since $\overline B(E_r,1)\subset\{|w|\le r+1\}$ has area at most
$\pi(r+1)^2$,
Lemma~\ref{lem:coefficient-integral-zero-bound} yields
\[
 n(E_r,1/f)\le C_\varepsilon r^{d_*+2}e^{a_\varepsilon r},
\]
which proves~\eqref{eq:global-directional-bound} with $N=d_*+2$,
independent of $\varepsilon$.

For the lower bound, choose any one $\omega\in\Omega_a$. By
Theorem~\ref{thm:generic-density}, outside two proper linear
subspaces of initial values, a fixed strip gives
\[
 n(r,1/f)\ge c_f r^{\beta(\omega)}e^{ar}
\]
for all sufficiently large $r$. Combining this with the two
upper estimates proves~\eqref{eq:generic-exponent} and the
stated concentration in $\Sigma_\varepsilon$.
\end{proof}

\subsection{Proof of \hyperref[thm:exact-global]{Theorem~\ref*{thm:exact-global}}}
\label{subsec:proof-exact-global}

\begin{proof}
To obtain the limit in the radial and angular coordinates
$(u,\tau)$ of the theorem, we first determine the limiting zero
measure in the coordinates $(\xi,\eta)$ of
Lemma~\ref{lem:physical-window} and then change variables.
Let $\mathbf F=(f_1,\ldots,f_n)$ be the row vector of
fundamental solutions normalized by
the initial values at $z_0$. It has no common zero, since its
jet matrix is invertible. Put $r_m=\sqrt m$, $R_m=ar_m$, and
$Q_m=R_m^{\beta_*}e^{R_m}$. This choice makes
$\sum_m Q_m^{-1}<\infty$, as required by the parameter estimate,
while the gaps between successive radii tend to zero.
Apply Lemma~\ref{lem:pluripolar-projection} to the finite set
$\Omega_*$ and the affine maps
\[
 z_{m,\omega}(\xi,\eta)
       =\omega^{-1}(R_m+\xi+i\sqrt{R_m}\eta).
\]
Outside one complete pluripolar cone $\mathcal E\subset\C^n$,
the lemma applies simultaneously to all main directions
$\arg z=-\arg\omega$, $\omega\in\Omega_*$, and every fixed
compact rectangle in the $(\xi,\eta)$-plane.

Fix $\mathbf c\notin\mathcal E$ and $\omega\in\Omega_*$.
For a real test function $\phi\in C_c^\infty(\R^2)$, write
$R=ar$, $Q=R^{\beta_*}e^R$, so that $D(r)=Q\sqrt R$.
Set
\[
 \phi_r(z)=\phi\left(\Rea(\omega z)-R,
                         \frac{\Ima(\omega z)}{\sqrt R}\right).
\]
In these coordinates
$\Delta_z\phi_r=a^2(\phi_{\xi\xi}+R^{-1}\phi_{\eta\eta})$
and $\mathrm dA_z=a^{-2}\sqrt R\dd\xi\dd\eta$.
The zero-measure identity~\eqref{eq:planar-zero-measure},
applied to $f_{\mathbf c}$ in the $z$-plane, gives
\[
 \begin{aligned}
 \sum_{f_{\mathbf c}(b)=0}m_b\phi_r(b)
 =\frac1{2\pi}
   \bigl\langle\Delta_z\log|f_{\mathbf c}|,\phi_r\bigr\rangle
 =\frac1{2\pi}\int_{\C}
   \log|f_{\mathbf c}(z)|\,\Delta_z\phi_r(z)\,\mathrm dA_z.
 \end{aligned}
\]
Changing variables and dividing by $D(r)=Q\sqrt R$, we obtain
\[
 \frac1{D(r)}\sum_{f_{\mathbf c}(z)=0}m_z\phi_r(z)
 =\frac1{2\pi}\int_{\R^2}
   \frac{\log|f_{\mathbf c}(\omega^{-1}(R+\xi+i\sqrt R\eta))|}{Q}
   (\phi_{\xi\xi}+R^{-1}\phi_{\eta\eta})\dd\xi\dd\eta.
\]
For $r=r_m$, set
$V_m(\xi,\eta)=Q_m^{-1}\log|f_{\mathbf c}(z_{m,\omega}(\xi,\eta))|$.
The preceding identity becomes
\[
 I_m:=\frac1{D(r_m)}\sum_{f_{\mathbf c}(b)=0}m_b\phi_{r_m}(b)
 =\frac1{2\pi}\int_{\R^2}
 V_m(\phi_{\xi\xi}+R_m^{-1}\phi_{\eta\eta})\dd\xi\dd\eta.
\]
Let $K$ be a compact rectangle containing $\operatorname{supp}\phi$,
and let $U_{r_m}$ be the function defined in
Lemma~\ref{lem:physical-window}. Since $f_{\mathbf c}=\mathbf F\mathbf c$,
\eqref{eq:pluripolar-projection-loss} in
Lemma~\ref{lem:pluripolar-projection} gives
$\|V_m-U_{r_m}\|_{L^1(K)}\to0$.
Subtracting the integrals and using $R_m\ge1$ for large $m$, we obtain
\[
 \begin{aligned}
 &\left|I_m-\frac1{2\pi}\int_{\R^2}
 U_{r_m}(\phi_{\xi\xi}+R_m^{-1}\phi_{\eta\eta})\dd\xi\dd\eta\right|\\
 &\qquad\le
 \frac{\|\phi_{\xi\xi}\|_\infty+\|\phi_{\eta\eta}\|_\infty}{2\pi}
 \int_K|V_m-U_{r_m}|\dd\xi\dd\eta\longrightarrow0.
 \end{aligned}
\]
By Lemma~\ref{lem:physical-window}, $|U_{r_m}|\le C_K$ on $K$
for large $m$, so
$R_m^{-1}\bigl|\int_K U_{r_m}\phi_{\eta\eta}\dd\xi\dd\eta\bigr|
\le C_KR_m^{-1}\|\phi_{\eta\eta}\|_{L^1(K)}\to0$.
Thus
\[
 \begin{aligned}
 I_m=\frac1{2\pi}\int_{\R^2}
 U_{r_m}(\phi_{\xi\xi}+R_m^{-1}\phi_{\eta\eta})\dd\xi\dd\eta+o(1)=\frac1{2\pi}\int_{\R^2}U_{r_m}\phi_{\xi\xi}\dd\xi\dd\eta+o(1).
 \end{aligned}
\]
Lemma~\ref{lem:physical-window} and two integrations by parts
in $\xi$ then give for every real $\phi\in C_c^\infty(\R^2)$,
\begin{equation}\label{eq:cartesian-global-zeros}
 \frac1{D(r_m)}\sum_{f_{\mathbf c}(b)=0}m_b\phi_{r_m}(b)
 \longrightarrow
 \frac{P_\omega}{4\pi^2}
 \int_{\R^2}e^\xi\phi(\xi,\eta)\dd\xi\dd\eta.
\end{equation}

\bigskip
We next pass from this sequence to arbitrary radii. If
$r_m\le r<r_{m+1}$, then $r-r_m=O(r^{-1})$ and
$D(r)/D(r_m)=1+o(1)$. The union
$B_r:=\operatorname{supp}\phi_r\cup\operatorname{supp}\phi_{r_m}$
in the $z$-plane lies in $|z|\le r+C_\phi$, and
$\|\phi_r-\phi_{r_m}\|_\infty=O_\phi(r-r_m)$.
Theorem~\ref{thm:global-zero-distribution} gives
\begin{equation}\label{eq:global-support-zero-bound}
 n(B_r,1/f_{\mathbf c})\le n(r+C_\phi,1/f_{\mathbf c})
 =O_\phi(D(r)),
\end{equation}
uniformly in $\mathbf c$. Since both test
functions vanish outside $B_r$, we obtain
\[
 \begin{aligned}
 \frac1{D(r)}\left|\sum_{f_{\mathbf c}(b)=0}
 m_b\bigl(\phi_r(b)-\phi_{r_m}(b)\bigr)\right|
 &\le\frac{\|\phi_r-\phi_{r_m}\|_\infty}{D(r)}
 n(B_r,1/f_{\mathbf c})\\
 &=O_\phi(r-r_m)\longrightarrow0.
 \end{aligned}
\]
Together with $D(r)/D(r_m)\to1$, this proves that
\eqref{eq:cartesian-global-zeros} holds for all
real $r\to\infty$ for smooth compactly supported test functions.

For a Borel set $B\subset\R^2$, define the normalized zero measure
in the $(\xi,\eta)$-coordinates by
\begin{equation}\label{eq:global-coordinate-zero-measure}
 \nu_r(B):=\frac1{D(r)}
 \sum_{\substack{f_{\mathbf c}(b)=0\\
 (\Rea(\omega b)-R,\,\Ima(\omega b)/\sqrt R)\in B}}m_b.
\end{equation}
Thus $\nu_r(B)$ is the total multiplicity of zeros whose rescaled
coordinates lie in $B$, divided by $D(r)$.
Define the measure $\nu$ by its density
\begin{equation}\label{eq:global-limit-zero-density}
 \mathrm d\nu(\xi,\eta)
 =\frac{P_\omega}{4\pi^2}e^\xi\dd\xi\dd\eta.
\end{equation}
Equation~\eqref{eq:cartesian-global-zeros} now reads
$\int\phi\,\mathrm d\nu_r\to\int\phi\,\mathrm d\nu$
for every $\phi\in C_c^\infty(\R^2)$.

To extend this convergence, fix $\phi\in C_c(\R^2)$.
We use the same smooth approximation as in the
\hyperref[par:sectorial-smooth-approximation]{proof of Proposition~\ref*{prop:sectorial-counting}};
see~\cite[Proposition~2.1.1, p.~21]{DinhSibonyCurrents}.
Explicitly, choose a nonnegative $\chi\in C_c^\infty(\R^2)$
supported in the unit disk, with $\int_{\R^2}\chi\,\mathrm dA=1$.
Write $\mathbf x=(\xi,\eta)\in\R^2$ for a point in these coordinates
and let $\mathbf y\in\R^2$ denote the integration variable.
For $0<h<1$, convolution with the rescaled function $\chi$
gives a smooth approximation $\phi_h$ satisfying
\begin{equation}\label{eq:global-smooth-test-construction}
 \begin{aligned}
 \phi_h(\mathbf x)
 &=\int_{\R^2}\chi(\mathbf y)\phi(\mathbf x-h\mathbf y)
       \,\mathrm dA(\mathbf y),\\
 \|\phi_h-\phi\|_\infty
 &\le\sup_{\substack{\mathbf x\in\R^2\\\|\mathbf y\|\le1}}
       |\phi(\mathbf x-h\mathbf y)-\phi(\mathbf x)|
       \longrightarrow0\qquad(h\to0),
 \end{aligned}
\end{equation}
where the convergence follows from the uniform continuity of $\phi$.
Since $\chi$ is supported in the unit disk,
$\chi(\mathbf y)\ne0$ implies $\|\mathbf y\|\le1$, and hence
$\|\mathbf x-(\mathbf x-h\mathbf y)\|=h\|\mathbf y\|\le h$.
If $\operatorname{dist}(\mathbf x,\operatorname{supp}\phi)>h$,
these points $\mathbf x-h\mathbf y$ therefore lie outside
$\operatorname{supp}\phi$, so $\phi(\mathbf x-h\mathbf y)=0$
whenever $\chi(\mathbf y)\ne0$.
Thus the integrand in~\eqref{eq:global-smooth-test-construction}
vanishes for every $\mathbf y$, giving $\phi_h(\mathbf x)=0$.
Since $0<h<1$, this gives
\[
 \operatorname{supp}\phi_h
 \subset\{\mathbf x:\operatorname{dist}(\mathbf x,\operatorname{supp}\phi)\le h\}
 \subset\{\mathbf x:\operatorname{dist}(\mathbf x,\operatorname{supp}\phi)\le1\}.
\]
Choose a fixed compact rectangle $K$ in the $(\xi,\eta)$-plane
containing this closed unit neighbourhood of $\operatorname{supp}\phi$.
Under $z=\omega^{-1}(R+\xi+i\sqrt R\eta)$, it corresponds to a
moving rectangle in the $z$-plane. The unit-disk condition applies
only to $\operatorname{supp}\chi$; neither $\operatorname{supp}\phi$
nor $K$ is required to lie in that disk.
For each $\varepsilon>0$, choose $h=h(\varepsilon)$ sufficiently
small and set $\phi_\varepsilon:=\phi_{h(\varepsilon)}$.
Then $\|\phi-\phi_\varepsilon\|_\infty<\varepsilon$ and both
supports lie in $K$, independently of $\varepsilon$.

Choose a fixed $\vartheta\in C_c^\infty(\R^2)$ with
$0\le\vartheta\le1$ and $\vartheta=1$ on $K$, and put
$\vartheta_r(z)=\vartheta(\Rea(\omega z)-R,\Ima(\omega z)/\sqrt R)$.
By~\eqref{eq:global-coordinate-zero-measure} and
\eqref{eq:global-support-zero-bound}, applied with $\phi=\vartheta$,
we have, for large $r$,
\begin{equation}\label{eq:global-compact-zero-mass}
 \begin{aligned}
 \nu_r(K)\le\int\vartheta\,\mathrm d\nu_r
 =\frac1{D(r)}\sum_{f_{\mathbf c}(b)=0}m_b\vartheta_r(b)\le\frac{n(\operatorname{supp}\vartheta_r,1/f_{\mathbf c})}{D(r)}
 \le C_K.
 \end{aligned}
\end{equation}
Also,~\eqref{eq:global-limit-zero-density} gives
$\nu(K)=\frac{P_\omega}{4\pi^2}\int_K e^\xi\dd\xi\dd\eta<\infty$,
since $K$ has finite area and $e^\xi$ is bounded on $K$.
Hence
\begin{equation}\label{eq:global-continuous-test-approximation}
 \begin{aligned}
 \left|\int\phi\,\mathrm d\nu_r-\int\phi\,\mathrm d\nu\right|
 &\le\|\phi-\phi_\varepsilon\|_\infty
       \bigl(\nu_r(K)+\nu(K)\bigr)
       +\left|\int\phi_\varepsilon\,\mathrm d\nu_r
                   -\int\phi_\varepsilon\,\mathrm d\nu\right|\\
 &\le\varepsilon\bigl(C_K+\nu(K)\bigr)
       +\left|\int\phi_\varepsilon\,\mathrm d\nu_r
                   -\int\phi_\varepsilon\,\mathrm d\nu\right|.
 \end{aligned}
\end{equation}
For fixed $\varepsilon$, the last term tends to zero
by~\eqref{eq:cartesian-global-zeros}. Letting first $r\to\infty$
and then $\varepsilon\to0$ proves that $\nu_r$ converges weakly
to $\nu$ against continuous compactly supported test functions.

\bigskip
It remains to restrict the count to the disk and change
coordinates. For $z=\omega^{-1}(R+\xi+i\sqrt R\eta)$, write,
uniformly for $(\xi,\eta)$ in each fixed compact set,
\begin{equation}\label{eq:global-radial-angular-coordinates}
 u_r:=a(r-|z|)=-\xi-\eta^2/2+o(1),\qquad
 \tau_r:=\sqrt R\arg(\omega z)=\eta+o(1).
\end{equation}
These are the coordinates denoted by $(u,\tau)$
in~\eqref{eq:joint-global-law}; the subscript records their dependence on $r$.
In these coordinates, $|z|\le r$ is equivalent to
$|R+\xi+i\sqrt R\eta|\le R$. Thus the disk and its limiting region are
\[
 A_r=\{(\xi,\eta):\xi+\eta^2/2+\xi^2/(2R)\le0\},\qquad
 A=\{(\xi,\eta):\xi+\eta^2/2\le0\}.
\]
We have $A_r\subset A$, and the disk boundary tends locally to
the parabola $\xi=-\eta^2/2$.
To restrict~\eqref{eq:cartesian-global-zeros}
to the disk, we must control the normalized zero count near this
boundary, where the characteristic function of the disk is discontinuous.

For a fixed compact rectangle $K$, put
$B_{\delta,K}=K\cap\{|\xi+\eta^2/2|\le\delta\}$, where $0<\delta<1$.
Let $N_r(B_{\delta,K})$ count, with multiplicity, the zeros of
$f_{\mathbf c}$ whose $(\xi,\eta)$-coordinates lie in $B_{\delta,K}$.
By~\eqref{eq:global-coordinate-zero-measure},
$\nu_r(B_{\delta,K})=N_r(B_{\delta,K})/D(r)$.
The set $B_{\delta,K}$ is compact and has area $O_K(\delta)$.
The density of $\nu$ in~\eqref{eq:global-limit-zero-density}
is bounded on $K$, since $K$ is compact and $e^\xi$ is continuous.
Since $\nu_r$ converges weakly to $\nu$ and $B_{\delta,K}$ is compact,
the upper bound in~\eqref{eq:weak-measure-set-bounds} gives
$\limsup_{r\to\infty}\nu_r(B_{\delta,K})\le\nu(B_{\delta,K})$.
Hence
\begin{equation}\label{eq:global-boundary-strip-count}
 \limsup_{r\to\infty}\frac{N_r(B_{\delta,K})}{D(r)}
 \le\frac{P_\omega}{4\pi^2}
 \int_{B_{\delta,K}}e^\xi\dd\xi\dd\eta
 \le C_K\delta.
\end{equation}
Since $\xi^2/(2R)\to0$ uniformly on $K$, for each fixed $\delta>0$
there is $r_{\delta,K}$ such that
\begin{equation}\label{eq:global-parabolic-region-difference}
 K\cap(A\setminus A_r)\subset B_{\delta,K},
 \qquad r\ge r_{\delta,K}.
\end{equation}
As $A_r\subset A$, this is precisely where the two regions disagree
inside $K$; see Figure~\ref{fig:disk-parabolic-limit}.
By~\eqref{eq:global-boundary-strip-count}, the normalized
zero count in this difference has limit superior at most $C_K\delta$,
so letting $\delta\to0$ shows that this normalized count tends to
zero as $r\to\infty$.
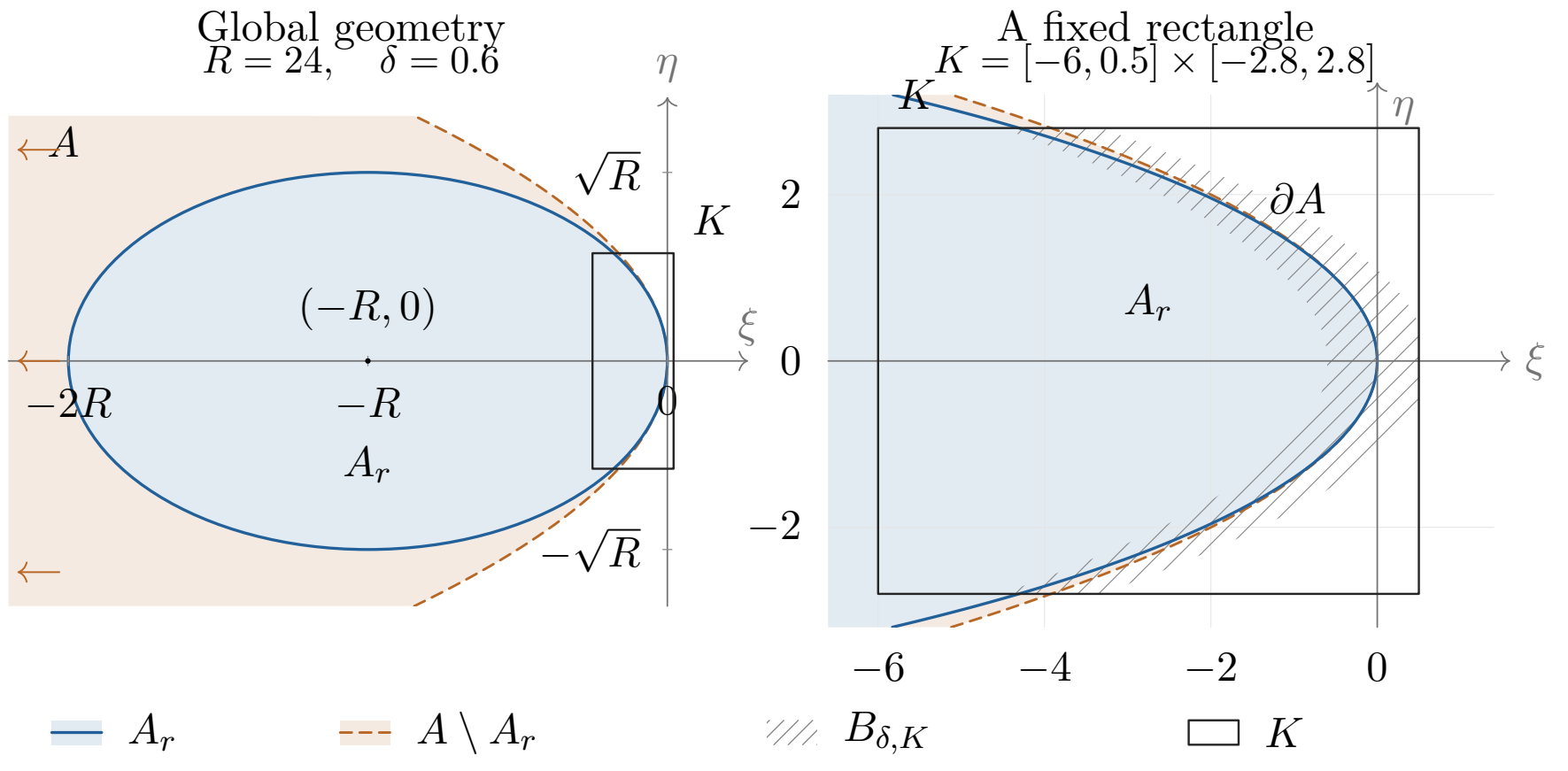
\begin{figure}[htbp]
\centering
\resizebox{0.76\linewidth}{!}{%
\begin{tikzpicture}[font=\normalsize,line join=round,line cap=round]
\definecolor{diskblue}{RGB}{34,96,154}
\definecolor{limitorange}{RGB}{184,104,38}
\tikzset{diskedge/.style={diskblue,line width=0.75pt},
  limitedge/.style={limitorange,line width=0.65pt,densely dashed},
  regionaxis/.style={black!55,line width=0.4pt,->}}

\node at (3.25,2.98) {Global geometry};
\node at (10.5,2.98) {A fixed rectangle};
\node[font=\small] at (3.25,2.69) {$R=24,\quad\delta=0.6$};
\node[font=\small] at (10.5,2.69)
  {$K=[-6,0.5]\times[-2.8,2.8]$};

\begin{scope}[shift={(6.1,0)},x=2.7cm,y=1.7cm]
  \begin{scope}
    \clip (-2.2,-1.3) rectangle (0.22,1.3);
    \fill[limitorange!14] (-2.2,-1.3)
      -- plot[domain=-1.3:1.3,samples=100,variable=\yy]
        ({-\yy*\yy/2},\yy) -- (-2.2,1.3) -- cycle;
    \fill[white] (-1,0) ellipse (1 and 1);
    \fill[diskblue!13] (-1,0) ellipse (1 and 1);
    \draw[limitedge] plot[domain=-1.3:1.3,samples=100,variable=\yy]
      ({-\yy*\yy/2},\yy);
    \draw[diskedge] (-1,0) ellipse (1 and 1);
  \end{scope}
  \draw[regionaxis] (-2.2,0)--(0.27,0) node[above] {$\xi$};
  \draw[regionaxis] (0,-1.3)--(0,1.4) node[above] {$\eta$};
  \foreach \xx/\lab in {-2/{-2R},-1/{-R},0/{0}}
    {\draw[black!45,line width=0.3pt] (\xx,-0.025)--(\xx,0.025);
     \node[below=3pt] at (\xx,0) {$\lab$};}
  \foreach \yy/\lab in {-1/{-\sqrt R},1/{\sqrt R}}
    {\draw[black!45,line width=0.3pt] (-0.015,\yy)--(0.015,\yy);
     \node[left=3pt] at (0,\yy) {$\lab$};}
  \fill (-1,0) circle[radius=0.7pt];
  \node[above=4pt] at (-1,0) {$(-R,0)$};
  \node at (-1,-0.55) {$A_r$};
  \node[anchor=west] at (-2.12,1.16) {$A$};
  \foreach \yy in {-1.12,0,1.12}
    \draw[limitorange,line width=0.45pt,->] (-2.03,\yy)--(-2.17,\yy);
  \draw[black!85,line width=0.55pt]
    (-0.25,{-2.8/sqrt(24)}) rectangle ({0.5/24},{2.8/sqrt(24)});
  \node[above right=1pt] at ({0.5/24},{2.8/sqrt(24)}) {$K$};
\end{scope}

\begin{scope}[shift={(12.5,0)},x=0.75cm,y=0.75cm]
  \begin{scope}
    \clip (-6.6,-3.2) rectangle (1.4,3.2);
    \fill[limitorange!14] (-6.6,-3.2)
      -- plot[domain=-3.2:3.2,samples=150,variable=\yy]
        ({-\yy*\yy/2},\yy) -- (-6.6,3.2) -- cycle;
    \fill[white] (-6.6,-3.2)
      -- plot[domain=-3.2:3.2,samples=150,variable=\yy]
        ({-24+sqrt(576-24*\yy*\yy)},\yy) -- (-6.6,3.2) -- cycle;
    \fill[diskblue!13] (-6.6,-3.2)
      -- plot[domain=-3.2:3.2,samples=150,variable=\yy]
        ({-24+sqrt(576-24*\yy*\yy)},\yy) -- (-6.6,3.2) -- cycle;
    \foreach \xx in {-6,-4,-2,0}
      \draw[black!10,line width=0.25pt] (\xx,-3.2)--(\xx,3.2);
    \foreach \yy in {-2,0,2}
      \draw[black!10,line width=0.25pt] (-6.6,\yy)--(1.4,\yy);
    \begin{scope}
      \clip (-6,-2.8) rectangle (0.5,2.8);
      \clip plot[domain=-2.8:2.8,samples=120,variable=\yy]
          ({-\yy*\yy/2-0.6},\yy)
        -- plot[domain=2.8:-2.8,samples=120,variable=\yy]
          ({-\yy*\yy/2+0.6},\yy) -- cycle;
      \foreach \xx in {-12,-11.75,...,3}
        \draw[black!48,line width=0.25pt] (\xx,-3.2)--({\xx+6.4},3.2);
    \end{scope}
    \draw[limitedge] plot[domain=-3.2:3.2,samples=150,variable=\yy]
      ({-\yy*\yy/2},\yy);
    \draw[diskedge] plot[domain=-3.2:3.2,samples=150,variable=\yy]
      ({-24+sqrt(576-24*\yy*\yy)},\yy);
  \end{scope}
  \draw[regionaxis] (-6.6,0)--(1.6,0) node[right] {$\xi$};
  \draw[regionaxis] (0,-3.2)--(0,3.35) node[below right] {$\eta$};
  \draw[black!85,line width=0.55pt] (-6,-2.8) rectangle (0.5,2.8);
  \foreach \xx in {-6,-4,-2,0}
    \node[below=3pt] at (\xx,-3.2) {$\xx$};
  \foreach \yy in {-2,0,2}
    \node[left=3pt] at (-6.6,\yy) {$\yy$};
  \node at (-2.75,0.7) {$A_r$};
  \node[anchor=west] at (-1.48,1.95) {$\partial A$};
  \node[above right=1pt] at (-6,2.8) {$K$};
\end{scope}

\fill[diskblue!13] (0.55,-3.46) rectangle (1.0,-3.24);
\draw[diskedge] (0.55,-3.35)--(1.0,-3.35);
\node[anchor=west] at (1.10,-3.35) {$A_r$};
\fill[limitorange!14] (3.15,-3.46) rectangle (3.6,-3.24);
\draw[limitedge] (3.15,-3.35)--(3.6,-3.35);
\node[anchor=west] at (3.70,-3.35) {$A\setminus A_r$};
\begin{scope}
  \clip (7.0,-3.46) rectangle (7.45,-3.24);
  \foreach \xx in {6.8,6.93,...,7.45}
    \draw[black!48,line width=0.3pt] (\xx,-3.46)--({\xx+0.22},-3.24);
\end{scope}
\node[anchor=west] at (7.55,-3.35) {$B_{\delta,K}$};
\draw[black!85,line width=0.55pt] (10.8,-3.46) rectangle (11.25,-3.24);
\node[anchor=west] at (11.35,-3.35) {$K$};
\end{tikzpicture}%
}
\caption{The rescaled disk, its parabolic limit, and the boundary strip.}
\label{fig:disk-parabolic-limit}
\end{figure}
We now justify convergence after restriction to $A$.
Fix $\phi\in C_c(\R^2)$ with $\operatorname{supp}\phi\subset K$,
put $q(\xi,\eta)=\xi+\eta^2/2$, and define the continuous function
\[
 \theta_\delta(\xi,\eta)=
 \begin{cases}
  1,&q(\xi,\eta)\le-\delta,\\
  \dfrac{\delta-q(\xi,\eta)}{2\delta},&|q(\xi,\eta)|<\delta,\\
  0,&q(\xi,\eta)\ge\delta.
 \end{cases}
\]
It takes values in $[0,1]$ and agrees with $\mathbf1_A$ outside
$\{|q|\le\delta\}$. Since $\operatorname{supp}\phi\subset K$,
we have the pointwise bound
\[
 |\phi(\mathbf1_A-\theta_\delta)|
 \le\|\phi\|_\infty\,\mathbf1_{B_{\delta,K}}.
\]
Integrating this bound against either positive measure
$\mu=\nu_r$ or $\mu=\nu$ gives
\begin{equation}\label{eq:global-parabolic-cutoff-error}
 \left|\int\phi\mathbf1_A\,\mathrm d\mu
       -\int\phi\theta_\delta\,\mathrm d\mu\right|
 \le\|\phi\|_\infty\,\mu(B_{\delta,K}).
\end{equation}
For fixed $\delta$, the product $\phi\theta_\delta$ is continuous
and compactly supported, so the weak convergence established
using~\eqref{eq:global-continuous-test-approximation} gives
$\int\phi\theta_\delta\,\mathrm d\nu_r
\to\int\phi\theta_\delta\,\mathrm d\nu$.
By~\eqref{eq:global-parabolic-cutoff-error} and
\eqref{eq:global-boundary-strip-count},
\[
 \limsup_{r\to\infty}
 \left|\int_A\phi\,\mathrm d\nu_r-\int_A\phi\,\mathrm d\nu\right|
 \le 2C_K\|\phi\|_\infty\delta.
\]
Letting $\delta\to0$ proves convergence on $A$.
Finally,~\eqref{eq:global-parabolic-region-difference} and
\eqref{eq:global-boundary-strip-count} give
\[
 \left|\int_{A_r}\phi\,\mathrm d\nu_r-\int_A\phi\,\mathrm d\nu_r\right|
 \le\|\phi\|_\infty\,
 \nu_r\bigl(K\cap(A\setminus A_r)\bigr)\longrightarrow0.
\]
Together with the density definition~\eqref{eq:global-limit-zero-density},
this yields
\begin{equation}\label{eq:global-disk-restricted-limit}
 \int_{A_r}\phi\,\mathrm d\nu_r
 \longrightarrow\int_A\phi\,\mathrm d\nu
 =\frac{P_\omega}{4\pi^2}
   \int_A e^\xi\phi(\xi,\eta)\dd\xi\dd\eta.
\end{equation}
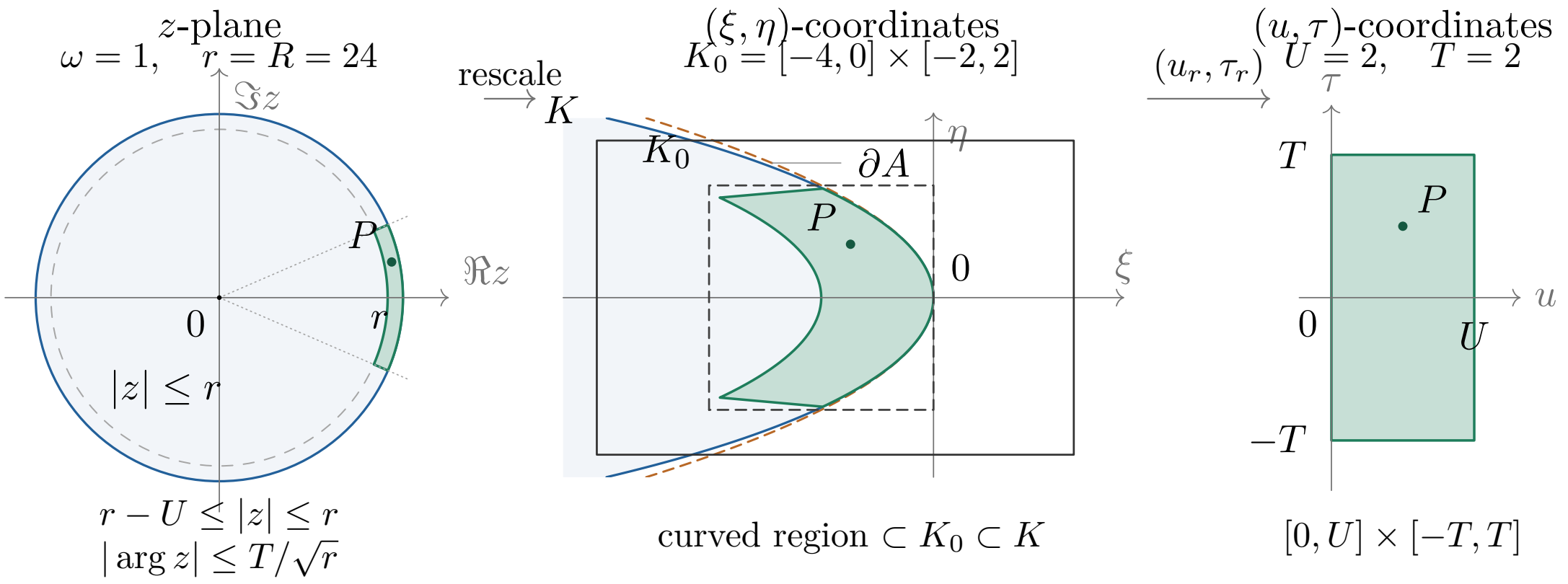
\begin{figure}[!b]
\centering
\resizebox{0.94\linewidth}{!}{%
\begin{tikzpicture}[font=\normalsize,line join=round,line cap=round]
\definecolor{diskblue}{RGB}{34,96,154}
\definecolor{limitorange}{RGB}{184,104,38}
\definecolor{windowgreen}{RGB}{30,125,91}
\tikzset{diskedge/.style={diskblue,line width=.65pt},
  limitedge/.style={limitorange,line width=.6pt,densely dashed},
  windowedge/.style={windowgreen,line width=.7pt},
  coordaxis/.style={black!55,line width=.35pt,->}}
\node at (2.1,2.65) {$z$-plane};
\node[font=\small] at (2.1,2.35) {$\omega=1,\quad r=R=24$};
\node at (8.3,2.65) {$(\xi,\eta)$-coordinates};
\node[font=\small] at (8.3,2.35) {$K_0=[-4,0]\times[-2,2]$};
\node at (13.7,2.65) {$(u,\tau)$-coordinates};
\node[font=\small] at (13.7,2.35) {$U=2,\quad T=2$};

\begin{scope}[shift={(2.1,0)},x=.075cm,y=.075cm]
  \fill[diskblue!6] (0,0) circle[radius=24];
  \draw[diskedge] (0,0) circle[radius=24];
  \draw[black!35,line width=.4pt,dashed] (0,0) circle[radius=22];
  \fill[windowgreen!25]
    plot[domain=-2:2,samples=100,variable=\tt]
      ({24*cos(deg(\tt/sqrt(24)))},{24*sin(deg(\tt/sqrt(24)))})
    --plot[domain=2:-2,samples=100,variable=\tt]
      ({22*cos(deg(\tt/sqrt(24)))},{22*sin(deg(\tt/sqrt(24)))})--cycle;
  \draw[windowedge]
    plot[domain=-2:2,samples=100,variable=\tt]
      ({24*cos(deg(\tt/sqrt(24)))},{24*sin(deg(\tt/sqrt(24)))})
    --plot[domain=2:-2,samples=100,variable=\tt]
      ({22*cos(deg(\tt/sqrt(24)))},{22*sin(deg(\tt/sqrt(24)))})--cycle;
  \foreach \sgn in {-1,1}
    \draw[black!35,densely dotted,line width=.35pt] (0,0)--
      ({27*cos(deg(\sgn*2/sqrt(24)))},{27*sin(deg(\sgn*2/sqrt(24)))});
  \draw[coordaxis] (-28,0)--(30,0) node[above right] {$\Re z$};
  \draw[coordaxis] (0,-28)--(0,30) node[below right] {$\Im z$};
  \fill (0,0) circle[radius=.65pt];
  \node[below left] at (0,0) {$0$};
  \node[below left] at (24,0) {$r$};
  \node at (-7,-12) {$|z|\le r$};
  \fill[windowgreen!70!black]
    ({23*cos(deg(1/sqrt(24)))},{23*sin(deg(1/sqrt(24)))}) circle[radius=1.3pt];
  \node[above left] at
    ({23*cos(deg(1/sqrt(24)))},{23*sin(deg(1/sqrt(24)))}) {$P$};
\end{scope}
\node[font=\small,align=center] at (2.1,-2.35)
  {$r-U\le|z|\le r$\\$|\arg z|\le T/\sqrt r$};

\begin{scope}[shift={(9.1,0)},x=.55cm,y=.55cm]
  \begin{scope}
    \clip (-6.6,-3.2) rectangle (3.2,3.2);
    \fill[diskblue!6] (-6.6,-3.2)
      --plot[domain=-3.2:3.2,samples=120,variable=\yy]
        ({-24+sqrt(576-24*\yy*\yy)},\yy)--(-6.6,3.2)--cycle;
    \draw[limitedge] plot[domain=-3.2:3.2,samples=120,variable=\yy]
      ({-\yy*\yy/2},\yy);
    \draw[diskedge] plot[domain=-3.2:3.2,samples=120,variable=\yy]
      ({-24+sqrt(576-24*\yy*\yy)},\yy);
    \fill[windowgreen!25]
      plot[domain=-2:2,samples=100,variable=\tt]
        ({24*cos(deg(\tt/sqrt(24)))-24},{sqrt(24)*sin(deg(\tt/sqrt(24)))})
      --plot[domain=2:-2,samples=100,variable=\tt]
        ({22*cos(deg(\tt/sqrt(24)))-24},{22/sqrt(24)*sin(deg(\tt/sqrt(24)))})--cycle;
    \draw[windowedge]
      plot[domain=-2:2,samples=100,variable=\tt]
        ({24*cos(deg(\tt/sqrt(24)))-24},{sqrt(24)*sin(deg(\tt/sqrt(24)))})
      --plot[domain=2:-2,samples=100,variable=\tt]
        ({22*cos(deg(\tt/sqrt(24)))-24},{22/sqrt(24)*sin(deg(\tt/sqrt(24)))})--cycle;
  \end{scope}
  \draw[coordaxis] (-6.6,0)--(3.4,0) node[above] {$\xi$};
  \draw[coordaxis] (0,-3.2)--(0,3.3) node[below right] {$\eta$};
  \draw[black!80,line width=.55pt] (-6,-2.8) rectangle (2.5,2.8);
  \draw[black!75,line width=.55pt,densely dashed] (-4,-2) rectangle (0,2);
  \node[above left=1pt] at (-6,2.8) {$K$};
  \node[above left=1pt] at (-4,2) {$K_0$};
  \node[above right=1pt] at (0,0) {$0$};
  \fill[windowgreen!70!black]
    ({23*cos(deg(1/sqrt(24)))-24},{23/sqrt(24)*sin(deg(1/sqrt(24)))}) circle[radius=1.3pt];
  \node[above left] at
    ({23*cos(deg(1/sqrt(24)))-24},{23/sqrt(24)*sin(deg(1/sqrt(24)))}) {$P$};
  \draw[black!45,line width=.3pt] (-2.88,2.4)--(-1.65,2.4);
  \node[anchor=west] at (-1.6,2.4) {$\partial A$};
\end{scope}
\node[font=\small] at (8.3,-2.35)
  {curved region $\subset K_0\subset K$};

\begin{scope}[shift={(13.0,0)},x=.70cm,y=.70cm]
  \fill[windowgreen!25] (0,-2) rectangle (2,2);
  \draw[windowedge] (0,-2) rectangle (2,2);
  \draw[coordaxis] (-.45,0)--(2.65,0) node[right] {$u$};
  \draw[coordaxis] (0,-2.7)--(0,2.7) node[above] {$\tau$};
  \node[below left] at (0,0) {$0$};
  \node[below=3pt] at (2,0) {$U$};
  \node[left=3pt] at (0,2) {$T$};
  \node[left=3pt] at (0,-2) {$-T$};
  \fill[windowgreen!70!black] (1,1) circle[radius=1.3pt];
  \node[above right] at (1,1) {$P$};
\end{scope}
\node[font=\small] at (13.7,-2.35)
  {$[0,U]\times[-T,T]$};

\draw[black!55,line width=.45pt,->] (4.70,1.95)--(5.20,1.95);
\node[above,font=\small] at (4.95,1.95) {rescale};
\draw[black!55,line width=.45pt,->] (11.20,1.95)--(12.40,1.95);
\node[above,font=\small] at (11.80,1.95) {$(u_r,\tau_r)$};
\end{tikzpicture}%
}
\caption{An annular sector in three coordinate systems.}
\label{fig:annular-sector-coordinates}
\end{figure}
The limiting coordinates $u=-q(\xi,\eta)$, $\tau=\eta$ map $A$
bijectively onto $X=[0,\infty)\times\R$, with inverse
$\xi=-u-\tau^2/2$, $\eta=\tau$.
By~\eqref{eq:global-radial-angular-coordinates},
$u_r=u+o(1)$ and $\tau_r=\tau+o(1)$ uniformly on each fixed compact set.
In particular,
\[
 \left|\det\frac{\partial(\xi,\eta)}{\partial(u,\tau)}\right|
 =\left|\det\begin{pmatrix}-1&-\tau\\0&1\end{pmatrix}\right|=1,
 \qquad e^\xi=e^{-u}e^{-\tau^2/2}.
\]
Fix $g\in C_c(X)$. The function
$\phi(\xi,\eta)=g(|q(\xi,\eta)|,\eta)$ belongs to $C_c(\R^2)$
and equals $g(-q(\xi,\eta),\eta)$ on $A$.
Taking this test function in~\eqref{eq:global-disk-restricted-limit}
and using $P_\omega/(4\pi^2)=C_\omega/\sqrt{2\pi}$
from~\eqref{eq:global-perimeter-constant} gives
\begin{equation}\label{eq:global-limiting-coordinate-test}
 \int_{A_r}g(-q,\eta)\,\mathrm d\nu_r
 \longrightarrow C_\omega\int_0^\infty\int_{\R}
 g(u,\tau)e^{-u}\frac{e^{-\tau^2/2}}{\sqrt{2\pi}}\dd\tau\dd u.
\end{equation}

To express the zero count as in~\eqref{eq:joint-global-law}, use
$u_r,\tau_r$ from~\eqref{eq:global-radial-angular-coordinates}.
If $\operatorname{supp}g\subset[0,U]\times[-T,T]$ and $R>U$,
the identity $\omega z=(R-u_r)e^{i\tau_r/\sqrt R}$ implies
$-U-T^2/2\le\xi\le0$ and $|\eta|\le T$ whenever
$g(u_r,\tau_r)\ne0$. Thus these points lie in the fixed rectangle
$K_0=[-U-T^2/2,0]\times[-T,T]$.
Choose a fixed compact rectangle $K$ containing
$K_0\cup\operatorname{supp}\phi$.
Figure~\ref{fig:annular-sector-coordinates} shows the corresponding
regions in the three coordinate systems.
By~\eqref{eq:global-radial-angular-coordinates}, uniform continuity
of $g$, and the bound~\eqref{eq:global-compact-zero-mass},
\[
 \left|\int_{A_r}\bigl(g(u_r,\tau_r)-g(-q,\eta)\bigr)\,\mathrm d\nu_r\right|
 \le\nu_r(K)\sup_{A_r\cap K}|g(u_r,\tau_r)-g(-q,\eta)|
 \longrightarrow0.
\]
Moreover, on this support,
$|\arg(\omega z)|\le T/\sqrt R<\varepsilon$ for large $r$.
Thus~\eqref{eq:global-coordinate-zero-measure} and
\eqref{eq:global-limiting-coordinate-test} give
\eqref{eq:joint-global-law} against continuous compactly supported
test functions on $X$.

\bigskip
To extend this convergence to bounded continuous test functions,
we control the zeros outside compact subsets of $X$. The
normalizing function is $D(r)=(ar)^{\beta_*}e^{ar}\sqrt{ar}$, and
only the frequencies in $\Omega_*$ contribute at this order.
We claim that there are constants $C,c>0$, depending only on the
equation, such that, for $U,T\ge1$, uniformly over nonzero solutions,
\begin{align}
 \limsup_{r\to\infty}\frac{n(r-U/a,1/f)}{D(r)}
       &\le Ce^{-U},\label{eq:global-radial-tail}\\
 \limsup_{r\to\infty}\frac{n(E_{r,T},1/f)}{D(r)}
       &\le Ce^{-cT^2},\label{eq:global-angular-tail}
\end{align}
where
\[
 E_{r,T}=\{0<|z|\le r:
  |\arg(\omega z)|>T/\sqrt{ar}\text{ for every }\omega\in\Omega_*\}.
\]
In particular, zeros outside any fixed neighbourhoods of the
directions $-\arg\omega$, $\omega\in\Omega_*$, contribute $o(D(r))$.

Apply~\eqref{eq:coefficient-integral-zero-bound} and integrate
each term of $W$. A term with weighted frequency $\nu=\lambda/(n-j)$ and
degree $b=\deg P_{j,\lambda}/(n-j)$ has full disk integral
$O(r^{b+1/2}e^{|\nu|r})$ if $\nu\ne0$, by
estimate~\eqref{eq:exponential-disk-integral}; zero frequencies
give only polynomial growth. Thus $|\nu|<a$ is exponentially
negligible after division by $D(r)$, and $|\nu|=a$, $b<\beta_*$
is $O(r^{b-\beta_*})=o(1)$. The remaining terms have
$\nu\in\Omega_*$ and $b=\beta_*$. Fix $U\ge1$.
The unit neighbourhood in~\eqref{eq:coefficient-integral-zero-bound}
enlarges the disk of radius $r-U/a$ to radius $r-U/a+1$, so
\[
 \begin{aligned}
 n(r-U/a,1/f)
 &\le C\int_{|w|\le r-U/a+1}W(w)\,\mathrm dA(w),\\
 \frac{n(r-U/a,1/f)}{D(r)}
 &\le C\left(\frac{r-U/a+1}{r}\right)^{\beta_*+1/2}
       e^{-U}+o_U(1).
 \end{aligned}
\]
Here we used~\eqref{eq:exponential-disk-integral} and
$e^{a(r-U/a+1)}/e^{ar}=e^{a-U}$. 
 For fixed $U$, the ratio in parentheses tends
to $1$ as $r\to\infty$, proving~\eqref{eq:global-radial-tail}.

For~\eqref{eq:global-angular-tail}, fix $T\ge1$.
Starting with~\eqref{eq:coefficient-integral-zero-bound},  expanding $W$ by~\eqref{eq:weight-function}, we have
\begin{equation}\label{eq:global-angular-integral-reduction}
 \begin{aligned}
 \frac{n(E_{r,T},1/f)}{D(r)}
 &\le\frac{C}{D(r)}\int_{\overline B(E_{r,T},1)}
       W(w)\,\mathrm dA(w)\\
 &\le\frac{C}{D(r)}\sum_{\nu\in\Omega_*}
       \int_{\substack{w\in\overline B(E_{r,T},1)\\|w|\ge r/2}}
       (1+|w|)^{\beta_*}e^{\Rea(\nu w)}\,\mathrm dA(w)
       +o_T(1).
 \end{aligned}
\end{equation}
Since $\overline B(E_{r,T},1)\subset\{|w|\le r+1\}$,
\eqref{eq:exponential-disk-integral} shows that all terms with
$|\nu|<a$, or $|\nu|=a$ and $b<\beta_*$, contribute $o(D(r))$.
For each remaining term, the integral over $|w|\le r/2$ is
$O(r^{\beta_*+1/2}e^{ar/2})=o(D(r))$ by the same estimate.

To bound the integrals remaining in~\eqref{eq:global-angular-integral-reduction},
fix $\nu\in\Omega_*$ and use
the rotated polar coordinates $w=e^{-i\arg\nu}\rho e^{i\theta}$,
where $\rho=|w|$ and $\theta=\arg(\nu w)\in[-\pi,\pi]$;
thus $\Rea(\nu w)=a\rho\cos\theta$.
If $w\in\overline B(E_{r,T},1)$ and $\rho\ge r/2$,
then $\rho\le r+1$. The unit enlargement changes the angular
distance from the ray $\arg(\nu w)=0$ by at most $C/r$.
Put $\theta^*=\sqrt{a\rho}\,\theta$. The definition of $E_{r,T}$ then gives
\[
 |\theta|\ge\frac{T}{\sqrt{ar}}-\frac Cr,
 \qquad
 |\theta^*|
 \ge\frac{T}{\sqrt2}-\frac C{\sqrt r}
 \ge\frac{T}{2\sqrt2}
\]
for sufficiently large $r$, with $T$ fixed.
Using $1-\cos\theta\ge c\theta^2$ on $[-\pi,\pi]$
and substituting $\theta^*=\sqrt{a\rho}\,\theta$ gives
\begin{equation}\label{eq:global-angular-gaussian-bound}
 \begin{aligned}
 \int_{\substack{|\theta|\le\pi\\
        |\theta|\ge T/\sqrt{ar}-C/r}}
       e^{a\rho\cos\theta}\,\dd\theta
 &\le\frac{e^{a\rho}}{\sqrt{a\rho}}
       \int_{|\theta^*|\ge T/(2\sqrt2)}e^{-c(\theta^*)^2}\,\dd\theta^*\\
 &\le C(1+\rho)^{-1/2}e^{a\rho}e^{-cT^2}.
 \end{aligned}
\end{equation}
For the last inequality, use
$e^{-c(\theta^*)^2}\le e^{-cT^2/16}e^{-c(\theta^*)^2/2}$ when
$|\theta^*|\ge T/(2\sqrt2)$ and integrate the second factor over $\R$.

In polar coordinates $\mathrm dA=\rho\,\dd\rho\,\dd\theta$.
Substituting~\eqref{eq:global-angular-gaussian-bound} into
\eqref{eq:global-angular-integral-reduction}, and using the radial
estimate in~\eqref{eq:exponential-disk-integral}, yields
\[
 \begin{aligned}
 \frac{n(E_{r,T},1/f)}{D(r)}
 &\le\frac{Ce^{-cT^2}}{D(r)}\int_{r/2}^{r+1}
       \rho(1+\rho)^{\beta_*-1/2}e^{a\rho}\,\dd\rho+o_T(1)\\
 &\le Ce^{-cT^2}+o_T(1),
 \end{aligned}
\]
where the finite number of frequencies in $\Omega_*$ is absorbed
into $C$. Taking the upper limit with $T$ fixed gives
\[
 \limsup_{r\to\infty}\frac{n(E_{r,T},1/f)}{D(r)}
 \le Ce^{-cT^2}.
\]
This proves~\eqref{eq:global-angular-tail}, with $C,c>0$
independent of $T$.

Let $K_{U,T}=[0,U]\times[-T,T]$ be a compact rectangle in the
$(u,\tau)$-plane. Let $N_r(U,T)$ count, with multiplicity, the
nonzero zeros in $|z|\le r$ outside the union, over
$\omega\in\Omega_*$, of the regions $0\le u_r\le U$,
$|\tau_r|\le T$ defined by~\eqref{eq:global-radial-angular-coordinates}.
By~\eqref{eq:global-radial-tail} and~\eqref{eq:global-angular-tail},
\[
 \begin{aligned}
 N_r(U,T)&\le n(r-U/a,1/f)+n(E_{r,T},1/f),\\
 \limsup_{r\to\infty}\frac{N_r(U,T)}{D(r)}
 &\le Ce^{-U}+Ce^{-cT^2}.
 \end{aligned}
\]
Thus increasing $U,T$ makes the normalized count outside these
rectangles arbitrarily small. For fixed $\varepsilon>0$,
points outside all sectors $|\arg(\omega z)|<\varepsilon$ belong
to $E_{r,T}$ once $T/\sqrt{ar}<\varepsilon$.
Letting first $r\to\infty$ and then $T\to\infty$
in~\eqref{eq:global-angular-tail} gives an $o(D(r))$ zero count
outside these fixed sectors.
For $g\in C_b(X)$, Urysohn's lemma
\cite[Theorem~2.12, p.~39]{Rudin1987} gives
$\chi_{U,T}\in C_c(X)$ with $0\le\chi_{U,T}\le1$ and
$\chi_{U,T}=1$ on $K_{U,T}$.
The convergence already proved applies to $g\chi_{U,T}$.
Since $1-\chi_{U,T}$ vanishes on $K_{U,T}$, the error in the
normalized zero sum has upper limit at most
$C\|g\|_\infty(e^{-U}+e^{-cT^2})$ by~\eqref{eq:global-radial-tail}
and~\eqref{eq:global-angular-tail}.
For the limiting integral, the complement of $K_{U,T}$ is
$\{u>U\}\cup\{|\tau|>T\}$ in $X$, so
\[
 \begin{aligned}
 &C_\omega\left|\int_X g(1-\chi_{U,T})
       e^{-u}\frac{e^{-\tau^2/2}}{\sqrt{2\pi}}\dd u\dd\tau\right|\\
 &\qquad\le C_\omega\|g\|_\infty
       \left(e^{-U}+\frac1{\sqrt{2\pi}}
             \int_{|\tau|>T}e^{-\tau^2/2}\dd\tau\right)\\
 &\qquad\le C\|g\|_\infty(e^{-U}+e^{-cT^2}).
 \end{aligned}
\]
Here $\int_U^\infty e^{-u}\dd u=e^{-U}$, and the angular
integral is bounded by the same Gaussian tail estimate used
in~\eqref{eq:global-angular-gaussian-bound}.
Letting first $r\to\infty$ and then $U,T\to\infty$ therefore proves weak convergence for
all bounded continuous test functions, as stated in the theorem.
Finally,
\[
 \frac{P_\omega}{4\pi^2}
 \int_0^\infty e^{-u}\dd u
 \int_{\R}e^{-\tau^2/2}\dd\tau
 =\frac{P_\omega}{2\pi\sqrt{2\pi}}=C_\omega.
\]
Summing over $\Omega_*$ proves~\eqref{eq:exact-global-count}.
\end{proof}

\begin{example}\label{ex:repeated-root-global-count}
The equation
\[
 f^{(4)}-2z^2e^{2z}f''+(z^4e^{4z}+e^{(1+i)z})f=0,\qquad z\in\C,
\]
has $a=\beta_*=1$, $\Omega_*=\{1\}$, and $p_1(\mu)=(\mu^2-1)^2$.
The two double roots give the segment $[-1,1]$, so $P_1=4$ and,
by~\eqref{eq:global-perimeter-constant},
$C_1=4/(2\pi\sqrt{2\pi})=\sqrt2/\pi^{3/2}$.
By~\eqref{eq:exact-global-count}, every initial vector outside the
exceptional cone gives
$n(r,1/f_{\mathbf c})\sim \sqrt2\,r^{3/2}e^r/\pi^{3/2}$.
\end{example}

\subsection{Proof of \hyperref[thm:perturbation-stability]{Theorem~\ref*{thm:perturbation-stability}}}
\label{subsec:proof-perturbation-stability}

\begin{proof}
To reuse the proof of Theorem~\ref{thm:exact-global}, we verify
that the perturbations preserve the limiting characteristic
polynomials, the uniform growth estimates for the fundamental row,
and the global zero bounds and tail estimates.

We first verify the coefficient estimate~\eqref{eq:moving-coefficient-limit}
of Lemma~\ref{lem:moving-inverse}, with the same limits $c_k$.
Choose $b>0$ with
$\max_j\tau_1(R_j)/(n-j)<b<a$. By the definition of exponential type,
\begin{equation}\label{eq:perturbation-global-growth}
 |R_j(z)|\le C_b e^{(n-j)b|z|},\qquad z\in\C,
 \quad 0\le j\le n-2.
\end{equation}
Constants may depend on the fixed coefficients and perturbations,
but not on the initial vector or branch index.

Fix $\omega\in\Omega_*$ and put $t=\omega z$, $\beta=\beta_*$.
Choose $\kappa>0$ so that
$\gamma=1-(b/a)\sqrt{1+\kappa^2}>0$.
On $t=x+iy$, $|y|\le\kappa x$, the added coefficient
$S_j(t)=\omega^{j-n}R_j(t/\omega)$ satisfies
\begin{equation}\label{eq:perturbation-sector-decay}
 \left|\frac{S_j(t)}{(\Phi'(t))^{n-j}}\right|
 \le C|t|^{-\beta(n-j)}
 e^{(n-j)((b/a)|t|-x)}
 \le Ce^{-(n-j)\gamma x}
\end{equation}
for large $x$, since $\Phi'(t)=t^\beta e^t$ and $\beta\ge0$.
Increase $L$ in Lemma~\ref{lem:moving-inverse} so that its inverse
branches lie in this sector. With $X=\log|s|$ and
$W_\ell=\log_Js+2\pi i\ell$, its proof gives, uniformly for
$X>L(1+|\ell|)$,
\[
 \Rea\psi_\ell(s)=X-\beta\log|W_\ell|+O(\log X/X),
 \qquad |W_\ell|\asymp X.
\]
Let $\widetilde b_{k,\ell}$ denote the transformed coefficients.
The added terms in~\eqref{eq:phase-coefficients} are products of
$S_j/(\Phi')^{n-j}$ and $C_{j,k}/(\Phi')^j$.
For $k\ge1$, the recurrence~\eqref{eq:chain-coefficients} gives
$C_{j,k}=(\Phi')^kP_{j,k}(1/t)$ with $P_{j,k}$ polynomial, so the
second factor is bounded for $k\le j$; for $k=0$, only $j=0$
contributes. Thus~\eqref{eq:perturbation-sector-decay} yields
\[
 |\widetilde b_{k,\ell}(s)-b_{k,\ell}(s)|\le CX^N e^{-\delta X},
 \qquad 0\le k\le n-2,
\]
for some $N,\delta>0$, while
$\widetilde b_{n-1,\ell}=b_{n-1,\ell}$.
Consequently~\eqref{eq:moving-coefficient-limit} holds with the
same limits $c_k$, and $p_\omega,H_\omega,P_\omega,C_\omega$ are unchanged.

We next recover the growth estimate~\eqref{eq:quantified-row-growth}
of Lemma~\ref{lem:uniform-fundamental-row} and the perimeter
average~\eqref{eq:physical-potential-average} of
Lemma~\ref{lem:physical-window} for the perturbed fundamental row.
Let $\mathbf G$ be this row,
normalized at $z_0$, and let $\widetilde{\mathbf X}_\ell$ be the
matrix of its derivatives in $s$ on branch $\ell$, as
in~\eqref{eq:branch-fundamental-matrix}.
By~\eqref{eq:perturbation-global-growth}, the companion matrix in
the variable $z$ still has norm at most $Ce^{C|z|}$, as
in~\eqref{eq:original-companion-growth}.
Put $\rho_{\ell,\varepsilon}=e^{L_\varepsilon(1+|\ell|)}$.
Gronwall's inequality along the segment from $z_0$ to
$\psi_\ell(\rho_{\ell,\varepsilon}e^{i\theta})/\omega$ gives
the same bounds as~\eqref{eq:original-matrix-gronwall}.
Expressing the derivatives in $s$ in terms of those in $z$
using~\eqref{eq:derivative-change-matrix} and
\eqref{eq:derivative-change-norm-bounds} then yields
\begin{equation}\label{eq:perturbation-entrance}
 \log^+\bigl\|\widetilde{\mathbf X}_\ell
 (\rho_{\ell,\varepsilon}e^{i\theta})^{\pm1}\bigr\|
 \le C_\varepsilon e^{C_\varepsilon(1+|\ell|)},
\end{equation}
uniformly for $\theta$ in compact subintervals of $J$.
This is the perturbed version of~\eqref{eq:global-entrance}.
Together,~\eqref{eq:moving-coefficient-limit} for the perturbed
coefficients and~\eqref{eq:perturbation-entrance} provide the bounds
used in the proof of Lemma~\ref{lem:uniform-fundamental-row}; hence its conclusions,
including~\eqref{eq:quantified-row-growth}, hold for
$\mathbf G_\ell=\mathbf G\circ(\psi_\ell/\omega)$.
Put $R=ar$ and $Q=R^{\beta_*}e^R$. A fixed compact rectangle
in the $(\xi,\eta)$-plane corresponds to a region meeting
$O(\sqrt R)$ inverse branches, and $\sqrt R=o(\log Q)$.
Thus the unchanged coordinates and the proof of
Lemma~\ref{lem:physical-window} give
\begin{equation}\label{eq:perturbation-window}
 Q^{-1}\log\|\mathbf G(\omega^{-1}(R+\xi+i\sqrt R\eta))\|
 \ \rightharpoonup\ \frac{P_\omega}{2\pi}e^\xi,
\end{equation}
with the local bounds of that lemma, against continuous test functions
with compact support.

Finally, we verify the integral zero bound
\eqref{eq:coefficient-integral-zero-bound} of
Lemma~\ref{lem:coefficient-integral-zero-bound} and the tail estimates
\eqref{eq:global-radial-tail}--\eqref{eq:global-angular-tail}
for the perturbed equation.
Put $\widetilde W(z)=W(z)+e^{b|z|}$,
where $W$ is the unperturbed weight~\eqref{eq:weight-function}.
The local coefficient bound~\eqref{eq:weight-local-coefficient-bound}
and~\eqref{eq:perturbation-global-growth} give
$|A_j(w)+R_j(w)|\le C\widetilde W(z)^{n-j}$ for $|w-z|\le5$.
The rescaling and integration argument of
Lemma~\ref{lem:coefficient-integral-zero-bound} therefore yields
\[
 n(E,1/g)\le C\int_{\overline B(E,1)}\widetilde W(z)\,\mathrm dA(z)
\]
for every bounded Borel set $E$ and every nonzero perturbed solution $g$.
For each fixed $d\ge0$,
\begin{equation}\label{eq:perturbation-negligible-mass}
 \int_{|z|\le r+d}e^{b|z|}\,\mathrm dA(z)
 =O_d(re^{br})=o(D(r)).
\end{equation}
The disk estimates for $W$ thus give $n(r,1/g)\le CD(r)$ uniformly
in $g$. Moreover, the regions of integration used to prove
\eqref{eq:global-radial-tail}--\eqref{eq:global-angular-tail}
lie in $|z|\le r+1$. The integrals of the added term $e^{b|z|}$
are therefore $o(D(r))$ by~\eqref{eq:perturbation-negligible-mass},
so those two tail bounds
hold with $g$ in place of $f$, with constants independent of $g$
and the fixed tail parameters $U,T\ge1$.

These three estimates allow us to apply the argument for the limiting zero measures
from the proof of Theorem~\ref{thm:exact-global}.
The components of $\mathbf G$ are linearly independent and have
no common zero. Apply Lemma~\ref{lem:pluripolar-projection} with
$r_m=\sqrt m$, $R_m=ar_m$, $Q_m=R_m^{\beta_*}e^{R_m}$ and the
same affine maps; $\sum_mQ_m^{-1}<\infty$ gives one complete
pluripolar cone $\widetilde{\mathcal E}$ for all $\omega\in\Omega_*$.
Outside it, the loss estimate~\eqref{eq:pluripolar-projection-loss}
and~\eqref{eq:perturbation-window} give the convergence of the
measures counting zeros in~\eqref{eq:cartesian-global-zeros} along $r_m$.
The uniform disk bound permits passage to all radii and to
continuous test functions with compact support. Restricting to
$|z|\le r$ and changing to $(u,\tau)$ yield~\eqref{eq:perturbed-joint-law}
on compact sets. The tail bounds just proved extend it to bounded
continuous tests and give $o(D(r))$ zeros outside the fixed sectors.
The limiting mass in each sector is the unchanged $C_\omega$;
summing proves~\eqref{eq:perturbed-global-count}.
\end{proof}

\section{Further consequences and questions}\label{sec:questions}

\subsection{Zeros near critical rays}

Put $E=\{0\}\cup\mathcal F$, with $\mathcal F$ as in~\eqref{eq:weighted},
and let $h(\theta)=\max_{q\in E}\Rea(qe^{i\theta})$.
Define the set $\Theta_{\mathrm c}\subset[0,2\pi)$ of \emph{critical angles} by
\[
 \theta\in\Theta_{\mathrm c}
 \iff \exists\,q_1\ne q_2\in E:\quad
 \Rea(q_1e^{i\theta})=\Rea(q_2e^{i\theta})=h(\theta).
\]
Geometrically, $\Theta_{\mathrm c}$ consists of outward normal directions
to the edges of the conjugate hull $\conv\{\overline q:q\in E\}$
(the two opposite normal directions when this hull is a segment).
This set is finite and nonempty. Define
\[
 \begin{gathered}
 \mathcal R=\bigcup_{\theta\in\Theta_{\mathrm c}}
       \{\rho e^{i\theta}:\rho\ge0\},\qquad
 a_{\mathrm c}=\max_{\theta\in\Theta_{\mathrm c}}h(\theta),\\
 \mathcal T_C=\{z:\operatorname{dist}(z,\mathcal R)
                          \le C\log(2+|z|)\}.
 \end{gathered}
\]

\begin{proposition}\label{cor:critical-ray-zero-bound}
We have $0\le a_{\mathrm c}<a$. For each $C>0$ there are
constants $K,N>0$, depending only on the equation and $C$, such
that every nonzero solution of~\eqref{eq:reduced} satisfies
\[
 n(\mathcal T_C\cap\{|z|\le r\},1/f)
       \le K r^N e^{a_{\mathrm c}r},\qquad r\ge1.
\]
In particular, this count is $o(r^\eta e^{ar})$ for every fixed real $\eta$.
\end{proposition}

\begin{proof}
Since $0\in E$, we have $h\ge0$. If $h(\theta)=a$, each maximizer
satisfies $a=\Rea(qe^{i\theta})\le|q|\le a$, forcing
$q=ae^{-i\theta}$. Thus $\theta\notin\Theta_{\mathrm c}$, and finiteness
of $\Theta_{\mathrm c}$ gives $a_{\mathrm c}<a$.

Set $F_r=\mathcal T_C\cap\{|z|\le r\}$ and
$d_*=\max_{j,\lambda}\deg P_{j,\lambda}/(n-j)$.
For $w\in\overline B(F_r,1)$, a nearest point $\rho e^{i\theta}$
on a suitable critical ray satisfies $0\le\rho\le r+1$ and
$|w-\rho e^{i\theta}|\le C\log(2+r)+1$. Hence, for $q\in E$,
\[
 \Rea(qw)=\rho\Re(qe^{i\theta})
 +\Re\!\bigl(q(w-\rho e^{i\theta})\bigr)\le a_{\mathrm c}(r+1)+a\{C\log(2+r)+1\}.
\]
The weight~\eqref{eq:weight-function} therefore satisfies
$W(w)\le K_1r^{d_*+Ca}e^{a_{\mathrm c}r}$ on
$\overline B(F_r,1)\subset\{|w|\le r+1\}$. By
Lemma~\ref{lem:coefficient-integral-zero-bound},
\[
 n(F_r,1/f)\le K_2\int_{\overline B(F_r,1)}W(w)\,\mathrm dA(w)
       \le K_3r^{d_*+Ca+2}e^{a_{\mathrm c}r}.
\]
The final assertion follows from $a_{\mathrm c}<a$.
\end{proof}

If Borel sets $U,\widetilde U\subset\C$ satisfy
\[
 (U\setminus\widetilde U)\cup(\widetilde U\setminus U)
 \subseteq\mathcal T_C,
\]
Proposition~\ref{cor:critical-ray-zero-bound} shows that their zero counts
in $|z|\le r$ differ by $o(r^\eta e^{ar})$ for every fixed real $\eta$.
If a zero-counting formula with a leading term on the scale
$r^\eta e^{ar}$ has already been proved in $U$, then the count in
$\widetilde U$ has the same leading term.

\subsection{Periodic coefficients and exceptional initial values}

Periodicity gives a more specific description of the exceptional
initial values. In the subclass below, all polynomial multipliers
are constant, so~\eqref{eq:beta} and
\eqref{eq:dominant-weighted-frequencies} give
\[
 \beta_*=0,\qquad \Omega_*=\Omega_a.
\]
For finite-order solutions in this class, see~\cite{LWW}.

\begin{proposition}\label{prop:periodic-finite-exception}
Suppose that the coefficients in~\eqref{eq:reduced} have the form
\[
 A_j(z)=\sum_{k\in K_j}a_{j,k}e^{k\lambda z},\qquad
 z\in\C,\quad 0\le j\le n-2,
\]
where $\lambda\ne0$, each $K_j\subset\Z$ is finite, the $a_{j,k}$
are constants, and at least one nonzero frequency occurs.
Then the exceptional cone $\mathcal E$ in
Theorem~\ref{thm:exact-global} can be chosen to be a finite union
of proper complex linear subspaces of $\C^n$.
A nonzero initial vector $\mathbf c$ satisfies the global counting
formula~\eqref{eq:exact-global-count} and the joint
limit~\eqref{eq:joint-global-law} for every $\omega\in\Omega_a$
if and only if $\mathbf c\notin\mathcal E$.
\end{proposition}

\begin{proof}
We choose $C=e$ in~\eqref{eq:phase}, so that
$s=\Phi(\omega z)=e^{\omega z}$.
Since every $\omega/\lambda$ is rational, a common period $P$ of the
coefficients can be chosen so that
$q_\omega=\omega P/(2\pi i)\in\Z\setminus\{0\}$ for every
$\omega\in\Omega_a$. For the fundamental matrix $\mathbf H$
in~\eqref{eq:original-fundamental-matrix}, normalized as
in~\eqref{eq:fundamental-initial-matrix}, periodicity gives
\[
 \mathbf H(z+P)=\mathbf H(z)\mathbf M,\qquad z\in\C,\quad \mathbf M\in\operatorname{GL}(n,\C).
\]
Fix $\omega\in\Omega_a$ and choose a branch of $\log s$ on an
exterior sector. The inverse branches of $s=e^{\omega z}$ are
\[
 z_\ell(s)=\frac{\log s+2\pi i\ell}{\omega},\qquad \ell\in\Z.
\]
Write $\ell=h+m q_\omega$, where $h\in\{0,\ldots,|q_\omega|-1\}$ and $m\in\Z$.
Substituting $\ell=h+m q_\omega$ into the formula above and using
$\omega P=2\pi i q_\omega$ gives
$
 z_\ell(s)
       =z_h(s)+mP.
$
Every inverse branch is therefore a period translate of one of the
$|q_\omega|$ branches $z_h$ for $h\in\{0,\ldots,|q_\omega|-1\}$.
Iterating the periodicity identity gives
$\mathbf H(z+mP)=\mathbf H(z)\mathbf M^m$ for $m\in\Z$.
Its first row is $\mathbf F(z+mP)=\mathbf F(z)\mathbf M^m$.
\phantomsection\label{par:periodic-branch-translation}
Multiplying on the right by $\mathbf c$ and using
$f_{\mathbf c}=\mathbf F\mathbf c$ yields
\[
 f_{\mathbf c}(z_\ell(s))
 =\mathbf F(z_h(s)+mP)\mathbf c
 =\mathbf F(z_h(s))\mathbf M^m\mathbf c
 =f_{\mathbf M^m\mathbf c}(z_h(s)).
\]
Namely, the value of the original solution at $z_h(s)+mP$ equals
the value at $z_h(s)$ of the solution with initial vector
$\mathbf M^m\mathbf c$.
For $K=\max\{1,\|\mathbf M\|,\|\mathbf M^{-1}\|\}$, submultiplicativity gives
\begin{equation}\label{eq:periodic-translation-norm-bounds}
 K^{-|m|}\|\mathbf c\|\le\|\mathbf M^m\mathbf c\|
 \le K^{|m|}\|\mathbf c\|,
\end{equation}
so $\log\|\mathbf M^m\mathbf c\|=O_{\mathbf c}(1+|m|)$ for
$\mathbf c\ne\mathbf0$.

On each representative branch, use the derivative vector $\mathbf Y$
in~\eqref{eq:phase-derivative-vector} and the system
\eqref{eq:phase-companion-system}. By Lemma~\ref{lem:phase-coefficient-limit},
$\mathbf A(s)\to\mathbf A_0$, whose characteristic polynomial is
$p_\omega$ in~\eqref{eq:characteristic}.
Divide the $s$-plane at the rays $\Rea((\mu-\nu)s)=0$ for distinct
roots, subdividing further if necessary so that each sector $S$
has width less than $\pi$. Write $\alpha=(\omega,h,S)$ for the
frequency, representative branch $z_h$ and chosen sector.
The $j$th column of $\mathbf X_\alpha(s)$ consists of $f_j(z_h(s))$
and its $s$-derivatives of orders $1,\ldots,n-1$, so that
$\mathbf Y(s)=\mathbf X_\alpha(s)\mathbf d$ for constant $\mathbf d\in\C^n$.

Fix a ray $s=\rho e^{i\theta_0}$ in $S$ and a point
$\rho_0e^{i\theta_0}\in S$.
For $\gamma$ strictly between two consecutive values
$\Rea(e^{i\theta_0}\mu)$, the substitution
$\mathbf w(\rho)=e^{-\gamma(\rho-\rho_0)}\mathbf Y(\rho e^{i\theta_0})$
gives a system with limiting matrix
$e^{i\theta_0}\mathbf A_0-\gamma\mathbf I$.
Its eigenvalues are $e^{i\theta_0}\mu-\gamma$.
Roots with $\Rea(e^{i\theta_0}\mu)<\gamma$ give negative real parts
(the stable part), whereas those with
$\Rea(e^{i\theta_0}\mu)>\gamma$ give positive real parts
(the unstable part).
The integral construction in the proof of
Lemma~\ref{lem:slower-growth-subspace} applies to this matrix as well:
the solutions $\mathbf w$ satisfying
$\sup_{\rho\ge\rho_0}\|\mathbf w(\rho)\|<\infty$ have initial values
in a subspace of dimension
$\sum_{\Rea(e^{i\theta_0}\mu)<\gamma}m_\mu$, where $m_\mu$ is the
algebraic multiplicity of $\mu$ as a root of $p_\omega$.
For every nonzero $\mathbf Y$, \eqref{eq:jet-ray-limit} and the
choice of $\gamma$ give
\[
 \lim_{\rho\to\infty}\frac{\log\|\mathbf w(\rho)\|}{\rho}
 =\ell_v(\theta_0)-\gamma\ne0.
\]
Thus $\mathbf w$ is bounded if and only if the growth rate
$\ell_v(\theta_0)$ of $\mathbf Y(\rho e^{i\theta_0})$ is less than $\gamma$.

Order the distinct roots $\mu_1,\ldots,\mu_q$ so that
$\sigma_j:=\Rea(e^{i\theta_0}\mu_j)$ satisfy
$\sigma_1<\cdots<\sigma_q$.
For $1\le j<q$, choose $\gamma_j\in(\sigma_j,\sigma_{j+1})$.
Let $E_j$ consist of $\mathbf0$ and the coefficient vectors $\mathbf d$
whose solutions have growth rate below $\gamma_j$.
This is the image of the initial-value subspace above under
$\mathbf X_\alpha(\rho_0e^{i\theta_0})^{-1}$.
Since $\gamma_j<\gamma_{j+1}$ for $1\le j<q-1$, any growth rate
below $\gamma_j$ is also below $\gamma_{j+1}$, and hence
$E_j\subset E_{j+1}$. Therefore
\[
 \{\mathbf0\}=E_0\subset E_1\subset\cdots\subset E_q=\C^n,
 \qquad \dim E_j=\sum_{r=1}^j m_{\mu_r}.
\]
Choose a basis of $E_1$, extend it to a basis of $E_2$, and continue
until $E_q=\C^n$. At step $j$, exactly $m_{\mu_j}$ new vectors
are added. Taking these vectors as columns gives the invertible matrix
\[
 \mathbf V_\alpha=
 \bigl(\mathbf v_{1,1},\ldots,\mathbf v_{1,m_{\mu_1}},
 \ \ldots,\quad
 \mathbf v_{q,1},\ldots,\mathbf v_{q,m_{\mu_q}}\bigr).
\]
Each new vector $\mathbf v_{j,k}$ lies in $E_j\setminus E_{j-1}$,
so the corresponding solution
$\mathbf G_{\alpha,j,k}(s):=\mathbf X_\alpha(s)\mathbf v_{j,k}$
has growth rate exactly $\sigma_j$ at $\theta_0$.
By Lemma~\ref{lem:indicator-limit}, its radial growth rate
\[
 \ell_{j,k}(\theta):=
 \lim_{\rho\to\infty}
 \frac{\log\|\mathbf G_{\alpha,j,k}(\rho e^{i\theta})\|}{\rho}
 \in\{\Rea(\mu e^{i\theta}):p_\omega(\mu)=0\}
\]
is continuous in $\theta$. A change from $\mu_j$ to another root
would therefore require their real projections to agree at an
intermediate angle, which is excluded in $S$.
Hence $\ell_{j,k}(\theta)=\Rea(\mu_j e^{i\theta})$ throughout $\arg S$.
Taking $R=|s|$ and $\zeta=s/|s|$ in
\eqref{eq:scaled-scalar-and-jet-potentials},
\eqref{eq:jet-potential-local-limit} gives, for the columns of
$\mathbf G_\alpha=\mathbf X_\alpha\mathbf V_\alpha$,
\[
 \log\|\mathbf G_{\alpha,j,k}(s)\|
 =\Rea(\mu_j s)+o(|s|),\qquad 1\le k\le m_{\mu_j},
\]
locally uniformly in angle. In particular, $\mu_j$ corresponds
to exactly $m_{\mu_j}$ columns.
Recall that $s=e^{\omega z}$ on each branch, so the chain rule gives
$\mathrm d/\mathrm dz=\omega s\,\mathrm d/\mathrm ds$.
Thus the lower triangular matrix $\mathbf T(\omega z)$
in~\eqref{eq:derivative-change-matrix} has diagonal entries
$(\omega s)^k$, $0\le k<n$, and
\[
 \det\mathbf T(\omega z)
 =\prod_{k=0}^{n-1}(\omega s)^k
 =(\omega s)^{n(n-1)/2}.
\]
Since $\det\mathbf H\equiv1$ and
$\mathbf G_\alpha=\mathbf X_\alpha\mathbf V_\alpha$, this gives
\[
 \det\mathbf G_\alpha(s)
 =\det\mathbf V_\alpha\,(\omega s)^{-n(n-1)/2}.
\]
Cramer's rule and Hadamard's inequality
\cite[\S0.8.3, p.~24; Problem~2.1.P23, p.~93]{HornJohnson}
bound each inverse row using the other columns.
For example, when $n=3$, write
$\mathbf G_\alpha=(\mathbf g_1,\mathbf g_2,\mathbf g_3)$,
with corresponding roots $\nu_1,\nu_2,\nu_3$, counted with multiplicity.
If $\mathbf r_1$ is the first row of $\mathbf G_\alpha^{-1}$, then
\[
 \begin{aligned}
 \|\mathbf r_1(s)\|
 &\le\frac{\|\mathbf g_2(s)\|\,\|\mathbf g_3(s)\|}
               {|\det\mathbf G_\alpha(s)|}\\
 &=\exp\bigl(\Rea((\nu_2+\nu_3)s)+o(|s|)\bigr)
 =\exp\bigl(-\Rea(\nu_1s)+o(|s|)\bigr).
 \end{aligned}
\]
Here $\nu_1+\nu_2+\nu_3=0$, and the polynomial factor from the
reciprocal determinant contributes only $O(\log|s|)=o(|s|)$
to the logarithm.
In general dimension, the same estimate uses the product of the norms
of all columns except $\mathbf G_{\alpha,j,k}$.
Since $\sum_{j=1}^q m_{\mu_j}\mu_j=0$, the corresponding inverse row
has norm at most $\exp(-\Rea(\mu_js)+o(|s|))$.
Partition the rows of $\mathbf V_\alpha^{-1}$ according to the root groups:
\[
 \mathbf V_\alpha^{-1}
 =\bigl(\mathbf L_{\alpha,\mu_1}^{\mathsf T},\ldots,
        \mathbf L_{\alpha,\mu_q}^{\mathsf T}\bigr)^{\mathsf T},
 \qquad
 \mathbf L_{\alpha,\mu_j}\in\C^{m_{\mu_j}\times n}.
\]
Since the rows of $\mathbf V_\alpha^{-1}$ are linearly independent,
each $\mathbf L_{\alpha,\mu_j}$ has rank $m_{\mu_j}$.
These choices are fixed independently of $\mathbf c$.
Write $\mathbf Y=\mathbf X_\alpha\mathbf d=\mathbf G_\alpha\mathbf a$,
where $\mathbf a=\mathbf V_\alpha^{-1}\mathbf d$ and
$\mathbf a_\mu=\mathbf L_{\alpha,\mu}\mathbf d$ is its group of coordinates
corresponding to $\mu$. Put
$M(s)=\max_\mu e^{\Rea(\mu s)}\|\mathbf a_\mu\|$.
On each closed subsector, for any $\varepsilon>0$ and all sufficiently
large $|s|$, the column estimates and the triangle inequality give
the first bound below, while $\mathbf a=\mathbf G_\alpha^{-1}\mathbf Y$
and the inverse row estimates give the second:
\[
 \begin{aligned}
 \|\mathbf Y(s)\|&\le C e^{\varepsilon|s|}M(s),\\
 e^{\Rea(\mu s)}\|\mathbf a_\mu\|
 &\le C e^{\varepsilon|s|}\|\mathbf Y(s)\|\qquad\text{for every }\mu.
 \end{aligned}
\]
Taking the maximum over $\mu$ in the second inequality gives the
lower bound $\|\mathbf Y(s)\|\ge C^{-1}e^{-\varepsilon|s|}M(s)$.
Taking logarithms and using the arbitrariness of $\varepsilon$ yields
\begin{equation}\label{eq:periodic-block-growth}
 \log\|\mathbf X_\alpha(s)\mathbf d\|
 =\max_\mu\{\Rea(\mu s)+\log\|\mathbf L_{\alpha,\mu}\mathbf d\|\}
   +o(|s|),\qquad \mathbf d\ne\mathbf0,
\end{equation}
with error uniform in $\mathbf d$ on closed subsectors; here
$\log0=-\infty$.

\bigskip
The bounds~\eqref{eq:periodic-translation-norm-bounds} control
$\|\mathbf M^m\mathbf c\|$, but do not prevent individual groups
$\mathbf L_{\alpha,\mu}\mathbf M^m\mathbf c$ from vanishing or being
much smaller. To apply~\eqref{eq:periodic-block-growth} with
$\mathbf d=\mathbf M^m\mathbf c$, we estimate these groups, writing
$\mathbf L$ for a fixed $\mathbf L_{\alpha,\mu}$.
\phantomsection\label{par:periodic-sequence-claim}
We claim that there is an integer $N\ge1$, depending only on $\mathbf M$, such that,
for each $j\in\{0,\ldots,N-1\}$,  the sequence with $m=j+Nk$, $k\in\Z$,
is either identically zero or admits a constant $C_{\mathbf c}>0$
such that the failure sets
\[
 B_{N_0}:=\left\{k\in\Z:|k|\le N_0,\quad
 \|\mathbf L\mathbf M^{j+Nk}\mathbf c\|
 <e^{-C_{\mathbf c}(1+|k|)}\right\},\qquad N_0\in\mathbb N,
\]
satisfy
\begin{equation}\label{eq:periodic-failure-count}
 \#B_{N_0}=o(N_0)\qquad(N_0\to\infty).
\end{equation}
In this case, the proportion of failures among $-N_0,\ldots,N_0$
tends to zero.

To see this, let $p_M$ be the number of distinct eigenvalues of the
period matrix $\mathbf M$, and denote them by
$\lambda_1,\ldots,\lambda_{p_M}$, which are nonzero. Choosing each argument in $[0,2\pi)$, put
\[
 \eta_r:=\frac{\arg\lambda_r}{2\pi}\in[0,1),\qquad
 \lambda_r=|\lambda_r|e^{2\pi i\eta_r},\qquad 1\le r\le p_M.
\]
These eigenvalues and their algebraic multiplicities need not coincide
with the roots $\mu_j$ and multiplicities $m_{\mu_j}$ of $\mathbf A_0$.
Choose a $\mathbb Q$-basis $\{1,\vartheta_1,\ldots,\vartheta_d\}$ of
$\operatorname{span}_{\mathbb Q}\{1,\eta_1,\ldots,\eta_{p_M}\}$.
This space has $\mathbb Q$-dimension $d+1$.
The basis elements are linearly independent over $\mathbb Q$,
and each $\eta_r$ has a unique expression
\[
 \eta_r=b_{r0}+\sum_{\ell=1}^d b_{r\ell}\vartheta_\ell,
 \qquad b_{r\ell}\in\mathbb Q,\quad 1\le r\le p_M.
\]
For example, since $\sqrt8=2\sqrt2$, the space
$\operatorname{span}_{\mathbb Q}\{1,\sqrt2,\sqrt3,\sqrt8\}$ has basis
$\{1,\sqrt2,\sqrt3\}$, hence dimension $3$ and $d=2$.
Choose $N\ge1$ so that $Nb_{r\ell}\in\Z$ for all $r$ and $\ell$,
and write $\boldsymbol\vartheta=(\vartheta_1,\ldots,\vartheta_d)$.
For fixed $j$, the Jordan block power formula
\cite[\S3.2.5, p.~180]{HornJohnson} gives, for all sufficiently large
positive $k$,
\[
 \mathbf M^{j+Nk}
 =\sum_{r=1}^{p_M}(\lambda_r^N)^k
   \sum_{b=0}^{n-1}k^b\mathbf B_{r,b},
\]
with constant matrices $\mathbf B_{r,b}$.
Multiplying on the left by $\mathbf L$ and on the right by $\mathbf c$
replaces these matrices by the constant vectors
$\mathbf u_{r,b}:=\mathbf L\mathbf B_{r,b}\mathbf c$.
Put $\varrho_r=|\lambda_r|^N$ and $a_{r\ell}=Nb_{r\ell}\in\Z$. Then
\[
 (\lambda_r^N)^k
 =\varrho_r^k\exp\!\left(2\pi i
      \sum_{\ell=1}^d a_{r\ell}(k\vartheta_\ell)\right),
\]
since $e^{2\pi ik a_{r0}}=1$.
For each distinct modulus $\varrho$ and $0\le b<n$, define, for
$\mathbf x=(x_1,\ldots,x_d)\in\R^d$,
\[
 \mathbf Q_{\varrho,b}(\mathbf x)
 :=\sum_{\{r:\,\varrho_r=\varrho\}}\mathbf u_{r,b}
       \exp\!\left(2\pi i\sum_{\ell=1}^d a_{r\ell}x_\ell\right).
\]
Since $\bigl|\exp(2\pi i\sum_{\ell=1}^d a_{r\ell}x_\ell)\bigr|=1$
for real $\mathbf x$, the triangle inequality gives
\[
 \|\mathbf Q_{\varrho,b}(\mathbf x)\|
 \le\sum_{\{r:\,\varrho_r=\varrho\}}\|\mathbf u_{r,b}\|<\infty,
\]
uniformly in $\mathbf x$.
Combining terms with the same modulus and degree gives the full sum
\[
 \mathbf L\mathbf M^{j+Nk}\mathbf c
 =\sum_{\varrho}\sum_{b=0}^{n-1}
   \varrho^k k^b\mathbf Q_{\varrho,b}(k\boldsymbol\vartheta).
\]
If any $\mathbf Q_{\varrho,b}$ is nonzero, choose the largest such
modulus $\varrho$, then the largest degree $b$ at that modulus,
and put $\mathbf Q=\mathbf Q_{\varrho,b}$.
After division by $\varrho^k k^b$, the other groups have lower degree
at the same modulus, giving $O(k^{-1})$, or smaller moduli, such that the ratio of exponential factors is the \(k\)th power of a fixed number in \((0,1)\), and hence decays exponentially. Hence
\begin{equation}\label{eq:periodic-leading-group}
 \mathbf L\mathbf M^{j+Nk}\mathbf c
 =\varrho^k k^b\bigl(\mathbf Q(k\boldsymbol\vartheta)+O(k^{-1})\bigr),
 \qquad \varrho>0,\quad 0\le b<n,
\end{equation}
Since $1,\vartheta_1,\ldots,\vartheta_d$ are linearly independent
over $\mathbb Q$, the sequence $(k\boldsymbol\vartheta)_{k\ge1}$ is
uniformly distributed modulo $1$; see~\cite[Chapter~1, Theorem~6.1
and Example~6.1, p.~48]{KuipersNiederreiter}.
Hence, for every continuous $\varphi:\R^d\to\C$ that is $1$-periodic
in each coordinate,
\begin{equation}\label{eq:periodic-equidistribution}
 \frac1{N_0}\sum_{k=1}^{N_0}\varphi(k\boldsymbol\vartheta)
 \longrightarrow\int_{[0,1]^d}\varphi(\mathbf x)\,\dd\mathbf x.
\end{equation}
The zero set of $\mathbf Q$ in $[0,1]^d$ has Lebesgue measure zero.
For $d=0$ this is immediate, since $\mathbf Q$ is a nonzero constant.
For $d\ge1$, choose a scalar component $q\not\equiv0$.
Since $\{\mathbf Q=\mathbf0\}\subseteq\{q=0\}$, it suffices to prove
the assertion for $q$.
In one variable, write $q(x)=\sum_{m=-M}^M c_m e^{2\pi imx}$.
With $z=e^{2\pi ix}$,
\[
 z^Mq(x)=\sum_{m=-M}^M c_m z^{m+M}=:P(z)\not\equiv0.
\]
The polynomial $P$ has finitely many roots, so $q$ has finitely
many zeros on $[0,1]$.

For $d\ge2$, assume the assertion is known for $d-1$ variables.
Writing $\mathbf x=(\mathbf x',t)$ and collecting the frequencies
in the last variable gives
\[
 q(\mathbf x',t)=\sum_{m=-M}^M q_m(\mathbf x')e^{2\pi imt},
 \qquad \mathbf x'\in[0,1]^{d-1},\quad t\in[0,1].
\]
Each $q_m$ is a trigonometric polynomial in $d-1$ variables, and
some $q_{m_0}\not\equiv0$.
Thus $Z_0:=\{\mathbf x'\in[0,1]^{d-1}:q_{m_0}(\mathbf x')=0\}$
has $(d-1)$-dimensional measure zero by induction.
For $\mathbf x'\notin Z_0$, the function $t\mapsto q(\mathbf x',t)$
is not identically zero and has finitely many zeros on $[0,1]$;
for $\mathbf x'\in Z_0$, its zero set has length at most $1$.
Writing $Z:=\{\mathbf x\in[0,1]^d:q(\mathbf x)=0\}$, put
$F(\mathbf x'):=\int_0^1\mathbf1_Z(\mathbf x',t)\,\dd t$.
Thus $F=0$ outside $Z_0$, while $0\le F\le1$ on $Z_0$.
Writing $\operatorname{meas}_r$ for $r$-dimensional Lebesgue measure,
Fubini's theorem gives
\[
 \begin{aligned}
 \operatorname{meas}_d(Z)
 &=\int_{[0,1]^{d-1}\setminus Z_0}F(\mathbf x')\,\dd\mathbf x'
   +\int_{Z_0}F(\mathbf x')\,\dd\mathbf x'\\
 &\le0+\int_{Z_0}1\,\dd\mathbf x'
 =\operatorname{meas}_{d-1}(Z_0)=0.
 \end{aligned}
\]
For $\varepsilon>0$, put
$A_\varepsilon:=\{\mathbf x\in[0,1]^d:\|\mathbf Q(\mathbf x)\|\le\varepsilon\}$.
Define the continuous function
\[
 \varphi_\varepsilon(\mathbf x)=
 \begin{cases}
 1,&\|\mathbf Q(\mathbf x)\|\le\varepsilon,\\
 2-\|\mathbf Q(\mathbf x)\|/\varepsilon,
   &\varepsilon<\|\mathbf Q(\mathbf x)\|\le2\varepsilon,\\
 0,&\|\mathbf Q(\mathbf x)\|>2\varepsilon.
 \end{cases}
\]
It is $1$-periodic in each coordinate and satisfies
$\mathbf1_{A_\varepsilon}\le\varphi_\varepsilon\le\mathbf1_{A_{2\varepsilon}}$
on $[0,1]^d$.
The sets $A_{2\varepsilon}$ decrease to the zero set of $\mathbf Q$
as $\varepsilon\to0$ and have measure at most $1$.
Continuity from above of Lebesgue measure therefore gives
\[
 \begin{aligned}
 0\le\int_{[0,1]^d}\varphi_\varepsilon\,\dd\mathbf x
 =\int_{A_{2\varepsilon}}\varphi_\varepsilon\,\dd\mathbf x
 \le\operatorname{meas}_d(A_{2\varepsilon})
 \longrightarrow0\qquad(\varepsilon\to0).
 \end{aligned}
\]
Fix $\varepsilon>0$. Taking $\varphi=\varphi_\varepsilon$
in~\eqref{eq:periodic-equidistribution} gives
\[
 \frac1{N_0}\sum_{k=1}^{N_0}\varphi_\varepsilon(k\boldsymbol\vartheta)
 \longrightarrow\int_{[0,1]^d}\varphi_\varepsilon(\mathbf x)\,\dd\mathbf x.
\]
Let $\mathcal I_{N_0}:=\{1\le k\le N_0:
\|\mathbf Q(k\boldsymbol\vartheta)\|\le k^{-1/2}\}$ and put
$K_\varepsilon:=\lceil\varepsilon^{-2}\rceil$, independently of $N_0$.
There are at most $K_\varepsilon-1$ positive integers $k<K_\varepsilon$;
for $k\ge K_\varepsilon$, we have $k^{-1/2}\le\varepsilon$.
Whenever $\|\mathbf Q(k\boldsymbol\vartheta)\|\le\varepsilon$,
the value $\varphi_\varepsilon(k\boldsymbol\vartheta)$ is $1$;
all other values are nonnegative.
Thus Markov's inequality for counting measure
\cite[Exercise~1.4.36(vi), p.~100]{TaoMeasure} gives
\[
 \begin{aligned}
 \#\mathcal I_{N_0}
 \le K_\varepsilon+
 \#\{K_\varepsilon\le k\le N_0:
       \|\mathbf Q(k\boldsymbol\vartheta)\|\le\varepsilon\}\le K_\varepsilon+
       \sum_{k=1}^{N_0}\varphi_\varepsilon(k\boldsymbol\vartheta).
 \end{aligned}
\]
Dividing by $N_0$ and letting $N_0\to\infty$ with $\varepsilon$ fixed,
so that $K_\varepsilon/N_0\to0$, yields
\[
 0\le\limsup_{N_0\to\infty}\frac{\#\mathcal I_{N_0}}{N_0}
 \le\int_{[0,1]^d}\varphi_\varepsilon(\mathbf x)\,\dd\mathbf x.
\]
The integral tends to zero as $\varepsilon\to0$, hence
$\#\mathcal I_{N_0}/N_0\to0$, or equivalently
$\#\mathcal I_{N_0}=o(N_0)$.
At the other indices, where \(\|\mathbf Q(k\boldsymbol\vartheta)\|>k^{-1/2}\), the \(O(k^{-1})\) term in (7.4) has norm at most \(k^{-1/2}/2\) for sufficiently large \(k\), yielding the required lower bound as \(k\to+\infty\).
Apply $\mathbf M^{-1}$ for $k\to-\infty$.
If every group vanishes, the sequence is identically zero.
This proves the claim, with~\eqref{eq:periodic-failure-count} in the
nonzero case.

\bigskip
Apply the Cayley--Hamilton theorem and the inverse polynomial formula
\cite[Theorem~2.4.3.2 and Corollary~2.4.3.4, pp.~109--110]{HornJohnson}
to the invertible matrix $\mathbf M^N$. They give scalar coefficients
$d_{k,r}\in\C$ such that
\[
 (\mathbf M^N)^k
 =\sum_{r=0}^{n-1}d_{k,r}(\mathbf M^N)^r,
 \qquad k\in\Z.
\]
For fixed $\alpha,\mu,j,\mathbf c$, define the vector sequence
$\mathbf y_k:=\mathbf L_{\alpha,\mu}\mathbf M^{j+Nk}\mathbf c$,
$k\in\Z$. Multiplying the preceding identity on the left by
$\mathbf L_{\alpha,\mu}\mathbf M^j$ and on the right by $\mathbf c$ gives
\[
 \begin{aligned}
 \mathbf y_k
 =\mathbf L_{\alpha,\mu}\mathbf M^j(\mathbf M^N)^k\mathbf c=\sum_{r=0}^{n-1}d_{k,r}
   \mathbf L_{\alpha,\mu}\mathbf M^{j+Nr}\mathbf c
 =\sum_{r=0}^{n-1}d_{k,r}\mathbf y_r.
 \end{aligned}
\]
Hence $\mathbf y_0=\cdots=\mathbf y_{n-1}=\mathbf0$
if and only if $\mathbf y_k=\mathbf0$ for all $k\in\Z$. Define
\begin{equation}\label{eq:periodic-vanishing-spaces}
 W_{\alpha,\mu,j}
 =\bigcap_{k=0}^{n-1}
       \ker(\mathbf L_{\alpha,\mu}\mathbf M^{j+Nk}),
 \qquad 0\le j<N.
\end{equation}
These  $W_{\alpha,\mu,j}$ are proper complex linear subspaces because
$\mathbf L_{\alpha,\mu}$ has positive rank.

For fixed $\mathbf c\ne\mathbf0$, define
\[
 \Lambda_{\alpha,j}(\mathbf c)
 =\{\mu:p_\omega(\mu)=0,\ \mathbf c\notin W_{\alpha,\mu,j}\}.
\]
This set is nonempty because $\mathbf V_\alpha^{-1}\mathbf M^j$
is invertible. Put $r:=\#\Lambda_{\alpha,j}(\mathbf c)\le n$.
For the $i$th root, the
\hyperref[par:periodic-sequence-claim]{preceding claim} gives a constant
$C_i>0$ and an exceptional set $B_{i,N_0}$ of integers
$m=j+Nk$, $|m|\le N_0$.
Since $r$ is fixed, \eqref{eq:periodic-failure-count} gives
\[
\#\left(\bigcup_{i=1}^{r}B_{i,N_0}\right)
\le\sum_{i=1}^{r}\#B_{i,N_0}
=o(N_0).
\]
Set $C_{\mathbf c}:=\max_{1\le i\le r}C_i$. For $m=j+Nk$, $|m|\le N_0$, outside this union,
\begin{equation}\label{eq:periodic-root-coordinate-lower-bound}
 \|\mathbf L_{\alpha,\mu}\mathbf M^m\mathbf c\|
 \ge e^{-C_{\mathbf c}(1+|m|)},
 \qquad \mu\in\Lambda_{\alpha,j}(\mathbf c).
\end{equation}
If a root $\mu\notin\Lambda_{\alpha,j}(\mathbf c)$, then
\eqref{eq:periodic-vanishing-spaces} and the Cayley--Hamilton argument give
$\mathbf L_{\alpha,\mu}\mathbf M^{j+Nk}\mathbf c=\mathbf0$
for all $k\in\Z$. For the remaining roots,
\eqref{eq:periodic-translation-norm-bounds} gives
$\|\mathbf L_{\alpha,\mu}\mathbf M^m\mathbf c\|
\le\|\mathbf L_{\alpha,\mu}\|K^{|m|}\|\mathbf c\|$.
Combining this with~\eqref{eq:periodic-root-coordinate-lower-bound},
for retained $m$ with $|m|\le C\sqrt{\log T}$ we obtain
\begin{equation}\label{eq:periodic-root-log-negligible}
 \max_{\mu\in\Lambda_{\alpha,j}(\mathbf c)}
 \frac{\bigl|\log\|\mathbf L_{\alpha,\mu}\mathbf M^m\mathbf c\|\bigr|}{T}
 \le\frac{C'_{\mathbf c}(1+C\sqrt{\log T})}{T}
 \longrightarrow0.
\end{equation}
Substitute $s=T\zeta$, $\mathbf d=\mathbf M^m\mathbf c$
in~\eqref{eq:periodic-block-growth}, omit the zero groups
($\log0=-\infty$), and divide by $T$.
By~\eqref{eq:periodic-root-log-negligible},
\begin{equation}\label{eq:periodic-selected-growth}
 \frac1T\log\|\mathbf X_\alpha(T\zeta)\mathbf M^m\mathbf c\|
 \longrightarrow H_{\alpha,j,\mathbf c}(\zeta)
 :=\max_{\mu\in\Lambda_{\alpha,j}(\mathbf c)}\Rea(\mu\zeta).
\end{equation}
The convergence is locally uniform in $\zeta$ and uniform in retained $m$, by the
uniformity in $\mathbf d$ of~\eqref{eq:periodic-block-growth}.

To pass to the scalar logarithms, let
$v_m(s):=f_{\mathbf M^m\mathbf c}(z_h(s))$, the first component of
$\mathbf X_\alpha(s)\mathbf M^m\mathbf c$, and put
\[
 U_{T,m}(\zeta):=\frac1T\log|v_m(T\zeta)|,\qquad
 V_{T,m}(\zeta):=\frac1T
   \log\|\mathbf X_\alpha(T\zeta)\mathbf M^m\mathbf c\|.
\]
The bound for the first component and the Cauchy comparison
\eqref{eq:cauchy-log-comparison} give, on each closed disk contained
in the current sector of the $\zeta$-plane,
\[
 U_{T,m}(\zeta)\le V_{T,m}(\zeta)
 \le\sup_{|\eta-\zeta|\le\delta}U_{T,m}(\eta)
       +\frac{\log C_\delta}{T},
\]
where $C_\delta$ is independent of $T$ and $m$.
By~\eqref{eq:periodic-selected-growth}, $U_{T,m}$ is locally uniformly
bounded above, and the second inequality excludes locally uniform
convergence to $-\infty$.
By the subharmonic compactness theorem
\cite[Theorem~2.4.8, p.~28]{DinhSibonyCurrents}, applied as in the
\hyperref[par:indicator-subharmonic-compactness]{proof of Lemma~\ref*{lem:indicator-limit}},
every sequence $U_{T_\nu,m_\nu}$ with $T_\nu\to\infty$ and retained
$m_\nu$ has a subsequence satisfying
\begin{equation}\label{eq:periodic-scalar-subsequence}
 U_{T_{\nu_\ell},m_{\nu_\ell}}
 \longrightarrow\widetilde H
 \quad\text{in }L^1_{\mathrm{loc}}.
\end{equation}
As in that proof, $U_{T,m}\le V_{T,m}$ gives
$\widetilde H\le H_{\alpha,j,\mathbf c}$ by the submean inequality;
the Cauchy comparison and Hartogs' lemma give the reverse inequality.
Thus every limit in~\eqref{eq:periodic-scalar-subsequence} equals
$H_{\alpha,j,\mathbf c}$, and hence, for every compact set $K$ in this sector,
\begin{equation}\label{eq:periodic-selected-scalar-limit}
 \int_K|U_{T,m}(\zeta)-H_{\alpha,j,\mathbf c}(\zeta)|
 \,\mathrm dA(\zeta)\longrightarrow0,
\end{equation}
uniformly over the retained $m$.

We next control two contributions: the omitted integers counted
in~\eqref{eq:periodic-failure-count}, and small angular neighbourhoods
in the $\zeta$-plane. For all integers $m$ in the same range, including the omitted ones,
Gronwall's inequality on a
finite sector cover, including the boundary rays, gives, on each
compact set $K$ away from $0$,
\[
 \begin{aligned}
 e^{-C_KT}\|\mathbf M^m\mathbf c\|
 \le\|\mathbf X_\alpha(T\zeta)\mathbf M^m\mathbf c\|\le e^{C_KT}\|\mathbf M^m\mathbf c\|,\qquad \zeta\in K.
 \end{aligned}
\]
Taking logarithms and dividing by $T$, we obtain
$-C_K+o(1)\le V_{T,m}\le C_K+o(1)$, since
$\log\|\mathbf M^m\mathbf c\|=o(T)$.
The Cauchy comparison~\eqref{eq:cauchy-log-comparison} then bounds
$U_{T,m}$ locally from above and excludes locally uniform convergence
to $-\infty$. Hence every sequence $U_{T_\nu,m_\nu}$ with
$T_\nu\to\infty$ has an $L^1_{\mathrm{loc}}$-convergent subsequence,
by~\cite[Theorem~2.4.8, p.~28]{DinhSibonyCurrents}, also for omitted $m_\nu$.
In particular, $\int_K|U_{T,m}|\,\mathrm dA\le C_{K,\mathbf c}$
for large $T$; otherwise a sequence with unbounded integrals could
not have an $L^1(K)$-convergent subsequence.
Moreover, for every $\varepsilon>0$ there are $\delta>0$ and $T_0$
such that, for $T\ge T_0$ and all $|m|\le C\sqrt{\log T}$,
\begin{equation}\label{eq:periodic-small-area-control}
 E\subset K,\quad \operatorname{area}(E)<\delta
 \quad\Longrightarrow\quad
 \int_E|U_{T,m}|\,\mathrm dA<\varepsilon
 \qquad(E\text{ measurable}).
\end{equation}
Indeed, if this failed, we could choose $T_\nu\to\infty$, $m_\nu$
and $E_\nu\subset K$ with $\operatorname{area}(E_\nu)\to0$ but
$\int_{E_\nu}|U_{T_\nu,m_\nu}|\,\mathrm dA\ge\varepsilon$.
After taking and relabeling a subsequence, put
$u_\nu:=U_{T_\nu,m_\nu}\to u$ in $L^1(K)$. Then
\[
 \varepsilon\le\int_{E_\nu}|u_\nu|\,\mathrm dA
 \le\int_K|u_\nu-u|\,\mathrm dA
       +\int_{E_\nu}|u|\,\mathrm dA\longrightarrow0,
\]
a contradiction; the last term tends to zero because $u\in L^1(K)$.

For $y=\theta+2\pi\ell=\Ima(\omega z)$, where
$\ell=h+q_\omega(j+Nk)$ and $\theta=\arg s$, put
$b_\omega:=2\pi|q_\omega|N$ and define
\[
 h_{\omega,\mathbf c}(y)
 :=\max_{\mu\in\Lambda_{\alpha,j}(\mathbf c)}\Rea(\mu e^{i\theta}).
\]
Replacing $y$ by $y+b_\omega$ changes $m$ by
$\operatorname{sgn}(q_\omega)N$, so this function is $b_\omega$-periodic.
Since $h_{\omega,\mathbf c}(y)\le h_\omega(y)$, Cauchy's perimeter formula gives
\[
 \begin{aligned}
 A_\omega(\mathbf c):=\frac1{b_\omega}
 \int_0^{b_\omega}h_{\omega,\mathbf c}(y)\,\dd y	
 \le\frac1{2\pi}\int_0^{2\pi}h_\omega(\theta)\,\dd\theta
 =\frac{P_\omega}{2\pi}.
 \end{aligned}
\]
The periodic averaging in the proof of Lemma~\ref{lem:physical-window}
applies with this period and mean.
With $R=ar$, $T=e^R$ and $D(r)=e^R\sqrt R$, take
$\operatorname{supp}\phi\subset[-L,L]^2$, which meets $O(\sqrt R)$ branches.
For each fixed $\alpha$, let $\mathcal O_R$ be the omitted integers $m$;
then~\eqref{eq:periodic-failure-count} gives $\#\mathcal O_R=o(\sqrt R)$.
By the \hyperref[par:periodic-branch-translation]{branch translation identity},
the normalized scalar logarithm on the branch $\ell=h+mq_\omega$
is $U_{T,m}(\zeta)$.
As in the proof of Lemma~\ref{lem:physical-window},
$\theta=\sqrt R\eta-2\pi\ell$ and $\zeta=e^{\xi+i\theta}$ give
$\dd\xi\dd\eta=R^{-1/2}|\zeta|^{-2}\,\mathrm dA_\zeta$.
The images lie in fixed compact annular sectors $K$, where
$|\zeta|^{-2}\le e^{2L}$.
Thus the absolute contribution of the omitted branches is at most
\[
 \frac{e^{2L}\|\phi\|_\infty}{\sqrt R}
 \sum_{m\in\mathcal O_R}\int_K|U_{T,m}|\,\mathrm dA
 \le C_{K,\mathbf c,\phi}\frac{\#\mathcal O_R}{\sqrt R}=o(1),
\]
using the uniform integral bound above; summing over the finitely many
$\alpha$ preserves this limit.
Separately, for the angular neighbourhoods, let
$E_\gamma:=K\cap\bigcup_{\theta_0}
\{\zeta:|\arg\zeta-\theta_0|<\gamma\}$, where $\theta_0$ ranges over
the finitely many sector boundary angles.
Each angular interval has length $2\gamma$, so integration of
$\rho\,\dd\rho\,\dd\theta$ over a fixed annulus gives
$\operatorname{area}(E_\gamma)\le C_K\gamma$.
Thus~\eqref{eq:periodic-small-area-control} gives
\begin{equation}\label{eq:periodic-omitted-angle-integral}
 \lim_{\gamma\to0}\limsup_{T\to\infty}
 \sup_{|m|\le C\sqrt{\log T}}
 \int_{E_\gamma}|U_{T,m}|\,\mathrm dA=0.
\end{equation}
Each branch contributes a factor $R^{-1/2}$ through $\dd\eta$;
summing over $O(\sqrt R)$ branches preserves this vanishing limit.
On the retained branches outside these neighbourhoods,
\eqref{eq:periodic-selected-scalar-limit} and
$H_{\alpha,j,\mathbf c}(e^{\xi+i\theta})=e^\xi h_{\omega,\mathbf c}(y)$
give the approximation by $e^\xi h_{\omega,\mathbf c}(\sqrt R\eta)$.
Since $|h_{\omega,\mathbf c}|\le\max_{p_\omega(\mu)=0}|\mu|$,
the comparison function is bounded on $[-L,L]^2$ and its integrals
over the omitted pieces also vanish.
The preceding estimates therefore give this approximation in
$L^1([-L,L]^2)$ on the whole rectangle.
Periodic averaging with mean $A_\omega(\mathbf c)$ therefore gives
\[
 \begin{aligned}
 \int_{\R^2} e^{-R}
 \log\left|f_{\mathbf c}\left((R+\xi+i\sqrt R\eta)/\omega\right)\right|
 \phi(\xi,\eta)\,\dd\xi\dd\eta\longrightarrow
 A_\omega(\mathbf c)\int_{\R^2}e^\xi\phi(\xi,\eta)\,\dd\xi\dd\eta
 \end{aligned}
\]
for every continuous compactly supported $\phi$.
The zero-measure identity~\eqref{eq:planar-zero-measure} gives limiting
normalized zero density $A_\omega(\mathbf c)e^\xi/(2\pi)$ in $(\xi,\eta)$;
positivity implies $A_\omega(\mathbf c)\ge0$.
The coordinate change in the proof of Theorem~\ref{thm:exact-global} and the tail bounds
\eqref{eq:global-radial-tail}--\eqref{eq:global-angular-tail} therefore give
\eqref{eq:joint-global-law} with
$C_\omega(\mathbf c)=A_\omega(\mathbf c)/\sqrt{2\pi}$, and
\[
 0\le C_\omega(\mathbf c)\le C_\omega,\qquad
 \frac{n(r,1/f_{\mathbf c})}{D(r)}
 \longrightarrow\sum_\omega C_\omega(\mathbf c).
\]

We now identify the exceptional initial vectors.
For each $\alpha=(\omega,h,S)$, let $\mu_\alpha$ be the unique root
such that $\Rea(\mu_\alpha s)>\Rea(\nu s)$ for $s\in S$ and every
distinct root $\nu$. Define
\[
 \mathcal E:=\bigcup_{\alpha=(\omega,h,S)}\ \bigcup_{j=0}^{N-1}
 W_{\alpha,\mu_\alpha,j}.
\]
There are finitely many frequencies $\omega$, representative branches $h$
and sectors $S$, and $j$ has only $N$ possible values.
Each space in this union is defined in~\eqref{eq:periodic-vanishing-spaces}
and satisfies
$W_{\alpha,\mu_\alpha,j}\subset
\ker(\mathbf L_{\alpha,\mu_\alpha}\mathbf M^j)\ne\C^n$,
since $\mathbf M^j$ is invertible and
$\operatorname{rank}\mathbf L_{\alpha,\mu_\alpha}=m_{\mu_\alpha}>0$.
Thus $\mathcal E$ is a finite union of proper complex linear subspaces.
It is a complex cone, and it is complete pluripolar because the function
\[
 \mathbf c\longmapsto\sum_{\alpha,j}
 \log\left(\sum_{k=0}^{n-1}
 \|\mathbf L_{\alpha,\mu_\alpha}\mathbf M^{j+Nk}\mathbf c\|^2\right)
\]
is plurisubharmonic and has $-\infty$ set exactly $\mathcal E$.
It is not identically $-\infty$, since finitely many proper linear
subspaces cannot cover $\C^n$.

If $\mathbf c\notin\mathcal E$, then
$\mu_\alpha\in\Lambda_{\alpha,j}(\mathbf c)$ for every $\alpha,j$.
The selected maximum therefore equals the maximum over all roots,
so $h_{\omega,\mathbf c}=h_\omega$ except at the finitely many
boundary angles in each period. Hence
$A_\omega(\mathbf c)=P_\omega/(2\pi)$ and
$C_\omega(\mathbf c)=C_\omega$ for every $\omega\in\Omega_a$.
The counting formula and all joint limits of
Theorem~\ref{thm:exact-global} follow.

Conversely, if $\mathbf c\in\mathcal E\setminus\{\mathbf0\}$,
choose $\alpha=(\omega,h,S)$ and $j$ such that
$\mathbf c\in W_{\alpha,\mu_\alpha,j}$.
Then $\mu_\alpha\notin\Lambda_{\alpha,j}(\mathbf c)$.
Choose a closed interval $I\subset\arg S$ of positive length and put
\[
 \delta:=\min_{\substack{\theta\in I,\ p_\omega(\nu)=0\\
                         \nu\ne\mu_\alpha}}
 \Rea((\mu_\alpha-\nu)e^{i\theta})>0.
\]
For $y\in I+2\pi(h+q_\omega j)$, the selected maximum omits
$\mu_\alpha$, so $h_\omega(y)-h_{\omega,\mathbf c}(y)\ge\delta$.
By periodicity,
\[
 \frac{P_\omega}{2\pi}-A_\omega(\mathbf c)
 =\frac1{b_\omega}\int_0^{b_\omega}
   \bigl(h_\omega(y)-h_{\omega,\mathbf c}(y)\bigr)\,\dd y
 \ge\frac{\delta|I|}{b_\omega}>0.
\]
Thus $C_\omega(\mathbf c)<C_\omega$ for this frequency, whereas
the other constants do not exceed their full values.
The total leading coefficient is therefore strictly smaller than
in~\eqref{eq:exact-global-count}, so that formula fails.
This proves both implications and the proposition.
\end{proof}

\subsection{Examples of entire perturbations}
\label{subsec:continuous-frequency}

For $0\le j\le n-2$, let $K_j\subset\C$ be nonempty and compact,
and let $\mu_j$ be a finite complex Borel measure supported on $K_j$.
Set
\[
 R_j(z)=\int_{K_j}e^{\lambda z}\dd\mu_j(\lambda),\qquad
 \sigma_j=\max_{\lambda\in K_j}|\lambda|,\qquad z\in\C.
\]
Differentiation under the integral shows that $R_j$ is entire, and
$|R_j(z)|\le|\mu_j|(K_j)e^{\sigma_j|z|}$, where $|\mu_j|$ is the
total variation measure. Thus
\[
 \max_{0\le j\le n-2}\frac{\sigma_j}{n-j}<a
\]
implies~\eqref{eq:perturbation-type-gap}, so all conclusions of
Theorem~\ref{thm:perturbation-stability} apply with the unperturbed
$D(r)$ and $C_\omega$. Finite sums of point masses recover finite
exponential sums; the measures may also have continuous parts.

For example, on $\C$, perturb~\eqref{eq:model-original} to
\[
 g''+\left(-\frac{1+e^{2z}}4+H(z)\right)g=0,
 \qquad
 H(z)=\begin{cases}\sinh z/z,&z\ne0,\\1,&z=0.\end{cases}
\]
The identity
\[
 H(z)=\frac12\int_{-1}^1e^{tz}\dd t,\qquad z\in\C,
\]
gives $\tau_1(H)\le1$, while $H(x)\sim e^x/(2x)$ as
$x\to+\infty$ gives equality. Hence its weighted type is
$1/2<a=1$. The unperturbed equation has $\beta_*=0$,
$\Omega_*=\{1\}$ and $p_1(\mu)=\mu^2-1/4$, whose conjugate-root
hull has perimeter $P_1=2$. Therefore, outside a complete
pluripolar cone in $\C^2$,
\[
 n(r,1/g_{\mathbf c})\sim\frac{\sqrt r\,e^r}{\pi\sqrt{2\pi}}.
\]
The joint limit~\eqref{eq:perturbed-joint-law} holds with
$u=r-|\zeta|$, $\tau=\sqrt r\arg\zeta$ and
$C_1=1/(\pi\sqrt{2\pi})$.

The strict type gap in~\eqref{eq:perturbation-type-gap} cannot be replaced
by a non-strict one without additional assumptions.
For example, on $\C$ the equation $f''-e^{2z}f=0$ has $a=1$,
but adding $R_0=e^{2z}$, with $\tau_1(R_0)/2=1$, gives $g''=0$.
Its nonzero solutions have at most one zero.

\needspace{10\baselineskip}
\begin{example}\label{ex:subunit-order-perturbations}
Let
\[
 H(z)=\prod_{k=1}^{\infty}\left(1+\frac{z}{k^2}\right),\qquad
 K(z)=\sum_{k=0}^{\infty}\frac{z^k}{(k!)^2}={}_0F_1(;1;z).
\]
The canonical product $H$ and the hypergeometric function $K$
satisfy $M(r,H)\le e^{\pi\sqrt r}$ and $M(r,K)\le e^{2\sqrt r}$;
see~\cite[Eqs.~\href{https://dlmf.nist.gov/4.36.E1}{(4.36.1)}
and~\href{https://dlmf.nist.gov/16.2.E1}{(16.2.1)}]{DLMF}.
Thus both have order at most $1/2$.
For arbitrary fixed $L_0,L_1\in\C$, consider
\[
 f'''+L_1H(z)f'+\bigl(e^{3z}+L_0K(z)\bigr)f=0,
 \qquad z\in\C.
\]
Here $a=1$, $\beta_*=0$ and
$p_1(\mu)=\mu^3+1$, whose root hull has perimeter $3\sqrt3$.
Corollary~\ref{cor:subunit-order-perturbations}
and~\eqref{eq:global-perimeter-constant} therefore give, outside
a complete pluripolar cone in $\C^3$,
\[
 n(r,1/f_{\mathbf c})\sim
 \frac{3\sqrt3}{2\pi\sqrt{2\pi}}\sqrt r\,e^r.
\]
\end{example}

\subsection{Osculating hyperplanes}

We return to equations of the form~\eqref{eq:reduced} to give a
geometric application of Theorem~\ref{thm:intro}.
Consider the third-order equation
\[
 f'''+A_1(z)f'+A_0(z)f=0.
\]
\par\pagebreak[3]
Let $V$ be its three-dimensional solution space, and choose a basis
$u_1,u_2,u_3$. The associated projective curve is
\[
 F(z)=[u_1(z):u_2(z):u_3(z)]\in\mathbb P^2.
\]
Here $[X_1:X_2:X_3]$ denotes the \emph{homogeneous coordinates}
of a projective point: a nonzero triple and any nonzero scalar
multiple represent the same point.
Set $\mathbf J(z)=\bigl(u_j^{(k-1)}(z)\bigr)_{k,j=1}^{3}$.
The full Wronskian $\det\mathbf J(z)$ never vanishes, so the $u_j$
have no common zero and $F$ is defined at every $z$.
For a nonzero solution $f=\sum_{j=1}^3c_ju_j$, let $H_f$ be the
projective line $c_1X_1+c_2X_2+c_3X_3=0$.
Substituting $X_j=u_j(z)$ gives
\[
 F(z)\in H_f\iff\sum_{j=1}^3c_ju_j(z)=0\iff f(z)=0.
\]
The \emph{intersection multiplicity} at a parameter $z_0$ is
$\operatorname{ord}_{z_0}f$.
The tangent line at $z_0$ is determined by
$f(z_0)=f'(z_0)=0$. Put $\mathbf c=(c_1,c_2,c_3)^{\mathsf T}$
and $\alpha=f''(z_0)$. The initial-value relation is
\[
 \mathbf J(z_0)\mathbf c
 =\bigl(f(z_0),f'(z_0),f''(z_0)\bigr)^{\mathsf T}
 =(0,0,\alpha)^{\mathsf T}.
\]
Since $\mathbf J(z_0)$ is invertible and $\mathbf c\ne\mathbf0$,
we have $\alpha\ne0$.
Solving for $\mathbf c$ and expanding $f$ at $z_0$ gives
\[
 \mathbf c=\alpha\mathbf J(z_0)^{-1}(0,0,1)^{\mathsf T},\qquad
 f(z)=\frac{\alpha}{2}(z-z_0)^2+O\bigl((z-z_0)^3\bigr).
\]
Taking $\alpha=1$ defines a solution $f_{z_0}$.
The solutions satisfying $f(z_0)=f'(z_0)=0$ form the
one-dimensional space $\operatorname{span}_{\C}\{f_{z_0}\}\subset V$.
Now fix a two-dimensional subspace $E\subset V$, such as
$E=\operatorname{span}_{\C}\{u_1,u_2\}$.
Keeping $E$ fixed, we vary $z_0$ and seek the values of \(z_0\) for which \(f_{z_0}\in E\), or equivalently
$\operatorname{span}_{\C}\{f_{z_0}\}\subset E$.

For a concrete example, take $f'''+e^z f=0$ and
$E=\{f\in V:f'(0)=0\}$, where $V$ is its solution space.
A basis of $E$ consists of the solutions with initial vectors
$(1,0,0)^{\mathsf T}$ and $(0,0,1)^{\mathsf T}$ at $0$.
Taking \(z_0=0\) in the normalized family \(f_{z_0}\) gives \(f_0\), whose initial vector is \((0,0,1)^{\mathsf T}\); hence \(f_0\in E\). The formula with $\mathbf J(z_0)^{-1}$ gives
holomorphic dependence on $z_0$, so Taylor's formula yields, as
$z_0\to0$,
\[
 f_{z_0}'(0)
 =f_{z_0}'(z_0)-z_0f_{z_0}''(z_0)+O(z_0^2)
 =-z_0+O(z_0^2).
\]
Thus $f_{z_0}\notin E$ for all sufficiently small nonzero $z_0$.

More generally, let $n\ge3$, let $V$ be the solution space
of~\eqref{eq:reduced}, and choose a basis $u_1,\ldots,u_n$ of $V$.
For $F(z)=[u_1(z):\cdots:u_n(z)]\in\mathbb P^{n-1}$, a nonzero
solution $f=\sum_{j=1}^nc_ju_j$ defines the hyperplane
$H_f:\sum_{j=1}^nc_jX_j=0$, of complex dimension $n-2$.
Its intersection multiplicity at $z_0$ is again
$\operatorname{ord}_{z_0}f$.
The \emph{osculating hyperplane} at the parameter $z_0$ is the
unique hyperplane $H_f$ satisfying
\[
 f(z_0)=f'(z_0)=\cdots=f^{(n-2)}(z_0)=0,
 \qquad f^{(n-1)}(z_0)\ne0.
\]
Indeed, the matrix
$\bigl(u_j^{(k)}(z_0)\bigr)_{0\le k<n,\,1\le j\le n}$ is invertible,
so these vanishing conditions determine a one-dimensional solution
subspace, and every nonzero member has a zero of multiplicity $n-1$ at $z_0$.
The normalization $f^{(n-1)}(z_0)=1$ defines $f_{z_0}$; its
hyperplane $H_{f_{z_0}}$ is the osculating hyperplane at $z_0$.
For $n=3$ this is a tangent line, and for $n=4$ an osculating plane;
see~\cite[\S2.1]{HSZ} for the general correspondence.

For the following result, fix an arbitrary subspace $E\subset V$ of dimension
$n-1$, choose a basis $u_1,\ldots,u_{n-1}$ of $E$, and extend it
to a basis of $V$. Use this basis to define $F$, and put
\[
 W_E(z)=\det\bigl(u_j^{(k)}(z)\bigr)_
 {0\le k\le n-2,\;1\le j\le n-1}.
\]
For instance, replacing $u_1,u_2$ by $2u_1,u_1+u_2$
in the chosen basis doubles $W_E$.
In general, a change of basis multiplies $W_E$ by the nonzero
determinant of the constant matrix relating the two bases.
If $W_E(z)=(z-z_0)^m h(z)$ near $z_0$, with $h$ analytic and
$h(z_0)\ne0$, then the Wronskian in the new basis satisfies
\[
 \widetilde W_E(z)=cW_E(z)=c(z-z_0)^m\, h(z),\qquad
  h(z_0)\ne0,
\]
where $c\ne0$ is this determinant ($c=2$ in the preceding example).
Thus $z_0$ remains a zero of multiplicity $m$.
The zeros of $W_E$ and their multiplicities therefore depend only on $E$.
Set $p=[0:\cdots:0:1]$. For a nonzero solution
$f=\sum_{j=1}^n c_ju_j$, we have
\[
 f\in E\quad\Longleftrightarrow\quad c_n=0
 \quad\Longleftrightarrow\quad p\in H_f.
\]
For this fixed $E$, we seek the points $z_0\in\C$ at which some
nonzero $f\in E$ has a zero of multiplicity $n-1$.
The determinant \(W_E(z_0)\) vanishes exactly when a nonzero
$f\in E$ satisfies $f(z_0)=\cdots=f^{(n-2)}(z_0)=0$.
Invertibility of the full initial-value matrix makes such a solution
a nonzero multiple of $f_{z_0}$, with zero multiplicity exactly $n-1$.
Thus
\[
 W_E(z_0)=0\quad\Longleftrightarrow\quad
 \exists f\in E\setminus\{0\}:\ \operatorname{ord}_{z_0}f=n-1.
\]
Equivalently, the osculating hyperplane $H_{f_{z_0}}$ passes through $p$.
The next proposition applies Theorem~\ref{thm:intro} to the adjoint
equation satisfied by $W_E$ to count these points.

\begin{proposition}
\label{prop:osculating-hyperplanes}
$W_E\not\equiv0$.
Fix $\omega\in\Omega_a$, write $\beta=\beta(\omega)$, and
choose $J$ and an inverse branch as in Theorem~\ref{thm:intro}.
There is a nondecreasing step function $M_E$ on $J$ such that
\[
 n_{U(\alpha_1,\alpha_2)}(r,1/W_E)
 =r^\beta e^{ar}
 \bigl(M_E(\alpha_2)-M_E(\alpha_1)\bigr)
 +o(r^\beta e^{ar}),
\]
whenever $\alpha_1<\alpha_2$ lie in $J$ and neither endpoint
is a jump point of $M_E$.
\end{proposition}

The count is taken in the $z$-plane, with weight
$\operatorname{ord}_{z_0}W_E$ at $z_0$.
This weight need not equal the intersection multiplicity
$n-1$.

\begin{proof}
Linear independence of the analytic functions $u_j$ gives
$W_E\not\equiv0$.

To apply Theorem~\ref{thm:intro} to $W_E$, we check its
differential equation and the corresponding weighted frequencies.
Since the full Wronskian is a nonzero constant, the cofactor
equations give the normalized formal adjoint equation on $\C$
\[
 W_E^{(n)}+
 \sum_{j=0}^{n-2}(-1)^{n+j}(A_jW_E)^{(j)}=0;
\]
compare~\cite[Proposition~2.6]{HSZ} for the relation between
the dual curve and the adjoint operator.
The coefficient of $W_E^{(k)}$ is
\[
 B_k=\sum_{j=k}^{n-2}
 (-1)^{n+j}\binom jk A_j^{(j-k)},
 \qquad 0\le k\le n-2.
\]
These coefficients remain exponential sums, and the
$W_E^{(n-1)}$ term is absent.
A nonzero frequency $\lambda$ contributed by $A_j$ to $B_k$,
where $k<j$, satisfies
\[
 \frac{|\lambda|}{n-k}
 <\frac{|\lambda|}{n-j}\le a.
\]
Thus terms with $j>k$ have weighted frequency of modulus
strictly less than $a$. At the maximal weighted frequencies,
only the term $j=k$ remains, and the polynomial multiplier
changes by the sign $(-1)^{n+k}$.
The quantities $a$, $\Omega_a$, and $\beta(\omega)$ for
$\omega\in\Omega_a$ are unchanged, and the same change of variable and
fixed inverse branch can be used, after increasing the cutoff
if necessary. Applying Theorem~\ref{thm:intro} to the adjoint
equation proves the formula.
\end{proof}

\subsection{Further questions}

Suppose a nonzero solution has infinitely many zeros in a fixed
strip $U$, but the coefficient of $r^\beta e^{ar}$
in~\eqref{eq:intro-count} vanishes.
When can further edges of the degree polygon determine a smaller
asymptotic scale for $n_U(r,1/f)$?
Compare~\cite{GSW} for growth orders with polynomial coefficients.

For general coefficients in~\eqref{eq:reduced}, are the initial
vectors for which the conclusions of Theorem~\ref{thm:exact-global}
fail contained in a countable union of proper complex linear
subspaces? Proposition~\ref{prop:periodic-finite-exception} gives
a finite union in the periodic subclass.

\section*{Acknowledgements}
The author thanks Xiong-Feng Liu for his comments on the asymptotic
analysis in an early version of the manuscript. The author acknowledges support from the China Scholarship
Council (CSC No.~202508440485).

\section*{Partial Lean formalization}
The author thanks the research group at the University of Manchester
for providing a platform for learning Lean formalization and autoformalization.
The supplementary Lean~4 project contains machine-checked counterparts of
Proposition~\ref{prop:sectorial-counting},
Corollary~\ref{cor:exposed-frequency-vertices}, and
Lemmas~\ref{lem:phase-coefficient-limit},
\ref{lem:radial-system-growth},
\ref{lem:indicator-limit},
\ref{lem:local-root-selection},
\ref{lem:radial-transfer},
\ref{lem:slower-growth-subspace},
\ref{lem:coefficient-integral-zero-bound}, and
\ref{lem:pluripolar-projection}.
The correspondence between the manuscript statements, their hypotheses,
and the Lean declarations is documented in the coverage table.
Theorems~\ref{thm:intro}--\ref{thm:perturbation-stability}
have not been fully formalized.

The source code and coverage table are maintained in the
\href{https://github.com/xyli9596/zero-density-lean-supplement}{supplementary repository},
which will be made public after this manuscript appears on arXiv.
The scope reported here refers specifically to
\href{https://github.com/xyli9596/zero-density-lean-supplement/tree/c76593e434064df5bd62c37432587aacddb0cb27}{version~0.1.0 (commit \texttt{c76593e43406})};
later versions may extend the coverage.

\section*{Declaration of AI assistance}
In an earlier version of this manuscript, the author used conformal
pullbacks to study the zeros of infinite-order solutions in logarithmic
strips. GPT-6 Astra (OpenAI) helped refine this construction during
revision. Building on the author's study of linearly independent
solutions and the initial-value space, GPT-6.1 Sol (OpenAI) suggested
describing exceptional initial values by complex cones and exploring
pluripotential theory. This suggestion drew on discussions of classical
potential theory in a separate research project of the author.
ChatGPT and Codex also helped develop and check arguments, explore
examples, prepare figures and edit the English, and suggested changes
to the order of sections. The author is grateful to the mathematical
community for its generosity in sharing knowledge and takes
responsibility for the accuracy of the final statements, proofs and
references.

\end{document}